\documentclass[11pt,reqno]{amsart}
\usepackage{amssymb,amscd,amsbsy,mathrsfs}

\usepackage[russian]{babel}
\Russian
\usepackage[cp1251]{inputenc}
\usepackage[OT2,T1]{fontenc}

\newcommand{\lb}{\linebreak}

\renewcommand{\a}{\alpha}
\renewcommand{\b}{\beta}

\renewcommand{\d}{\delta}
\newcommand{\e}{\varepsilon}
\newcommand{\vk}{\varkappa}
\newcommand{\z}{\zeta}

\renewcommand{\l}{\lambda}

\newcommand{\s}{\sigma}
\renewcommand{\t}{\tau}

\renewcommand{\f}{\varphi}
\renewcommand{\o}{\omega}

\newcommand{\G}{\Gamma}
\newcommand{\D}{\Delta}
\renewcommand{\L}{\Lambda}

\renewcommand{\O}{\Omega}

\newcommand{\n}{\nabla}

\newcommand{\A}{{\mathscr A}}

\newcommand{\E}{{\mathscr E}}
\newcommand{\cd}{{\mathscr D}}
\newcommand{\F}{{\mathscr F}}

\newcommand{\h}{{\mathscr H}}

\newcommand{\K}{{\mathscr K}}

\newcommand{\M}{{\mathscr M}}

\newcommand{\X}{{\mathscr X}}

\newcommand{\1}{{\bf 1}}

\renewcommand{\C}{{\Bbb C}}
\newcommand{\T}{{\Bbb T}}
\newcommand{\pp}{{\Bbb P}}
\newcommand{\dd}{{\Bbb D}}
\newcommand{\R}{{\Bbb R}}
\newcommand{\Z}{{\Bbb Z}}

\newcommand{\0}{{\boldsymbol{0}}}

\newcommand{\bs}{\boldsymbol}

\newcommand{\m}{{\boldsymbol m}}
\newcommand{\bS}{{\boldsymbol S}}

\newcommand{\rf}[1]{(\ref{#1})}

\newcommand{\df}{\stackrel{\mathrm{def}}{=}}
\newcommand{\dist}{\operatorname{dist}}

\newcommand{\re}{\operatorname{Re}}
\newcommand{\spn}{\operatorname{span}}
\newcommand{\supp}{\operatorname{supp}}
\newcommand{\clos}{\operatorname{clos}}

\newcommand{\trace}{\operatorname{trace}}
\newcommand{\rank}{\operatorname{rank}}
\newcommand{\const}{\operatorname{const}}

\newcommand{\eeq}{\end{equation}}
\newcommand{\beq}{\begin{equation}}
\newcommand{\bay}{\begin{eqnarray}}
\newcommand{\ba}{\begin{align*}}
\newcommand{\ea}{\end{align*}}
\newcommand{\ey}{\end{eqnarray}}
\newcommand{\bey}{\begin{eqnarray*}}
\newcommand{\eey}{\end{eqnarray*}}

\newcommand{\eq}{\Leftrightarrow}
\newcommand{\imp}{\Rightarrow}
\newcommand{\be}{\infty}

\newcommand{\bl}{\blacksquare}
\newcommand{\ess}{\operatorname{ess}}

\newcommand{\Pf}{{\bf Доказательство. }}
\newcommand{\im}{\operatorname{Im}}
\renewcommand{\re}{\operatorname{Re}}
\newcommand{\ov}{\overline}

\newtheorem{thm}{\hspace{\parindent}Теорема}[section]

\newtheorem{cor}[thm]{\hspace{\parindent}Следствие}
\newtheorem{lem}[thm]{\hspace{\parindent}Лемма}

\newcommand{\qm}{\quad\mbox{и}\quad}

\newcommand\Li{{\rm Lip}}
\newcommand\OL{{\rm OL}}
\newcommand\CM{{\rm CM}}
\newcommand\fM{\frak M}
\newcommand\fF{\frak F}

\newcommand\dg{\frak D}

\newcommand\CL{{\rm CL}}
\newcommand\OD{{\rm OD}}
\newcommand\sgn{\operatorname{sgn}}
\newcommand{\nn}{{\Bbb N}}
\newcommand\mB{\mathcal{B}}
\newcommand\mD{\mathcal{D}}
\newcommand\mM{\mathcal{M}}
\newcommand\mS{\mathcal{S}}
\newcommand\mT{\mathcal{T}}
\newcommand\mE{\mathcal{E}}
\newcommand\mP{\mathcal{P}}

\newcommand{\dlR}{{\operatorname{Aut}}(\widehat{\Bbb R})}
\newcommand{\lC}{{\operatorname{Aut}}(\Bbb C)}
\newcommand{\lR}{{\operatorname{Aut}}(\Bbb R)}
\newcommand{\dlC}{{\operatorname{Aut}}(\widehat{\Bbb C})}

\newcommand\zh{|\!\!|\!|}

\begin{document}

\

\newcommand{\vse}{\vspace{.2in}}
\numberwithin{equation}{section}

\title{Операторно липшицевы функции}
\author{А.Б. Александров и В.В.Пеллер}
\thanks{Исследования первого автора частично поддержаны грантом РФФИ 14-01-00198;
исследования второго автора частично поддержаны грантом NSF DMS 130092}

\begin{abstract}
Целью этого обзора является подробное изучения операторно липшицевых функций. Непрерывная функция $f$ на вещественной прямой $\R$ называется операторно липшицевой, если $\|f(A)-f(B)\|\le\const\|A-B\|$ для любых самосопряжённых операторов $A$ и $B$. Приводятся достаточные условия и необходимые условия для операторной липшицевости. Изучается также класс операторно дифференцируемых функций на $\R$. Далее, рассматривается класс операторно липшицевых функций на замкнутых подмножествах плоскости, а также вводится класс коммутаторно липшицевых функций на таких подмножествах. Для изучения этих классов функций важную роль играют двойные операторные интегралы и мультипликаторы Шура.
\end{abstract}

\maketitle

{\bf
\footnotesize
\tableofcontents
\normalsize
}

\setcounter{section}{0}
\section{\bf Введение}
\setcounter{equation}{0}
\label{intr}

\

Одна из важнейших задач теории возмущений состоит в исследовании, насколько изменятся  функции $f(A)$ от оператора $A$ при малых возмущениях оператора. В частности, естественным образом возникает задача описать класс непрерывных функций $f$ на вещественной прямой $\R$ таких, что справедливо неравенство
\bay
\label{OLf}
\|f(A)-f(B)\|\le\const\|A-B\|
\ey
для произвольных (ограниченных) самосопряжённых операторов $A$ и $B$ в гильбертовом пространстве. Такие функции $f$ называются {\it операторно липшицевыми}. Напомним, что функции от самосопряжённых (нормальных) операторов определяются, как интегралы этих функций по спектральным мерам операторов, см. \cite{R}.

Класс операторно липшицевых функций на $\R$ мы будем обозначать символом $\OL(\R)$. Отметим здесь, что если $f$ -- операторно липшицева функция, то неравенство \rf{OLf} справедливо и для неограниченных самосопряжённых операторов $A$ и $B$ с ограниченной разностью, см. теорему \ref{anbnrn} ниже; причём с той же самой константой в правой части неравенства \rf{OLf}. Минимальное значение этой константы является, по определению, нормой $\|f\|_{\OL}=\|f\|_{\OL(\R)}$ функции $f$ 
в пространстве $\OL(\R)$
(строго говоря, полунормой, а нормой она станет после отождествления функций отличающихся на постоянную функцию). 

Легко видеть, что если $f$ -- операторно липшицева функция, то она {\it липшицева}, т.е.
$$
|f(x)-f(y)|\le\const|x-y|,
$$
при всех вещественных $x$ и $y$ (мы используем символ $\Li(\R)$ для обозначения класса липшицевых функций на $\R$). Обратное утверждение неверно. Фарфоровская в \cite{F1} построила пример  липшицевой функции, не являющейся операторно липшицевой. Позже в работах \cite{Mc} и \cite{Kat} было установлено, что липшицева функция $x\mapsto|x|$ не операторно липшицева.

Операторно липшицевы функции играют важную роль во многих вопросах теории операторов и математической физике. В частности, они возникают при изучении задачи применимости формулы следов Лифшица--Крейна:
\bay
\label{LiKr}
\trace\big((f(A)-f(B)\big)=\int_\R f'(t)\xi(t)\,dt
\ey
(см. \cite{Kr}).
Здесь $A$ и $B$ -- самосопряжённые операторы в гильбертовом пространстве такие, что оператор $A-B$ ядерный (т.е. $A-B\in\bS_1$), а $\xi$ -- функция класса $L^1(\R)$ ({\it функция спектрального сдвига}), которая определяется только операторами $A$ и $B$. Ясно, что правая часть равенства \rf{LiKr} имеет смысл для любой липшицевой функции $f$. Что касается левой части равенства, то, как показывает пример Фарфоровской \cite{F2}, условия $A-B\in\bS_1$ и $f\in\Li(\R)$ не гарантируют того, что $f(A)-f(B)\in\bS_1$, и, стало быть, для применимости формулы следов \rf{LiKr} для всех пар самосопряжённых операторов с ядерной разностью необходимо наложить более сильное условие на $f$. По крайней мере, функция $f$ должна обладать следующим свойством:
\bay
\label{yad}
A-B\in\bS_1\quad\Longrightarrow\quad f(A)-f(B)\in\bS_1
\ey
для самосопряжённых операторов $A$ и $B$. Для функций $f$ 
на $\R$ свойство \rf{yad}
имеет место для произвольных (не обязательно ограниченных) самосопряжённых операторов в том и только в том случае, если функция $f$ операторно липшицева (см. теорему \ref{tsentrez} ниже).
Оказывается (см. недавнюю работу \cite{Pe7}), что операторная липшицевость функции $f$ 
не только необходима для справедливости формулы следов \rf{yad} для произвольных (необязательно ограниченных) самосопряжённых операторов $A$ и $B$ с ядерной разностью, но и достаточна.

Класс операторно липшицевых функций обладает некоторыми специфическими свойствами. Так, например, операторно липшицевы функции всюду дифференцируемы, но не обязательно непрерывно дифференцируемы (см. теорему \ref{pro} и пример 7 в \S\:\ref{prim}).

Оказывается, что операторная липшицевость функции может быть охарактеризована в терминах мультипликаторов Шура (см. \S\:\ref{razdrazn}).
Мы увидим, что непрерывная функция $f$ на $\R$ является операторно липшицевой тогда и только тогда, когда она всюду дифференцируема и разделённая разность $\dg f$,
$$
(\dg f)(x,y)\df\frac{f(x)-f(y)}{x-y},\quad x,\:y\in\R,
$$
является мультипликатором Шура.

Аналогичным образом можно рассмотреть задачу для функций на окружности и унитарных операторов. Непрерывная функция $f$ на единичной окружности $\T$ называется {\it операторно липшицевой}, если 
$$
\|f(U)-f(V)\|\le\const\|U-V\|
$$
для произвольных унитарных операторов $U$ и $V$.

В главе I этого обзора мы обсудим необходимые условия и достаточные условия для операторной липшицевости функций на прямой $\R$ и окружности $\T$.
Отметим здесь, что в случае самосопряжённых операторов ключевую роль играет неравенство
\bay
\label{sig}
\|f(A)-f(B)\|\le\const\s\|f\|_{L^\be}\|A-B\|
\ey
для любых самосопряжённых операторов $A$ и $B$ с ограниченной разностью и для любой ограниченной функции $f$ на $\R$ такой, что её преобразование Фурье сосредоточено на промежутке $[-\s,\s]$, $\s>0$. Это неравенство было получено в работах \cite{Pe1} и \cite{Pe3}.
Позже в работе \cite{AP4} было установлено, что неравенство \rf{sig} справедливо с константой 1. 

По аналогии с операторно липшицевыми функциями естественно рассмотреть операторно гёльдеровы функции. Пусть $0<\a<1$. Говорят, что функция $f$ на $\R$ {\it является операторно гёльдеровой порядка} $\a$, если имеет место неравенство
$$
\|f(A)-f(B)\|\le\const\|A-B\|^\a
$$
для произвольных самосопряжённых операторов $A$ и $B$ в гильбертовом пространстве. Однако оказалось (см \S\:\ref{oHold}), что здесь ситуация сильно отличается от случая операторно липшицевых оценок: функция $f$ является операторно гёльдеровой порядка $\a$ в том и только в том случае, когда она входит в класс $\L_\a(\R)$ гёльдеровых функций порядка $\a$, т.е.
$$
|f(x)-f(y)|\le\const|x-y|^\a,\quad x,\;y\in\R.
$$

В главе II мы обсуждаем двойные операторные интегралы, т.е. выражения вида
$$
\iint\Phi(x,y)\,dE_1(x)T\,dE_2(y).
$$
Здесь $\Phi$ -- ограниченная измеримая функция, $T$ - ограниченный линейный оператор в гильбертовом пространстве, а $E_1$ и $E_2$ -- спектральные меры.
Двойные операторные интегралы появились в работе Ю.Л. Далецкого и С.Г. Крейна \cite{DK} и подверглись систематическому изучению в работах М.Ш. Бирмана и М.З. Соломяка \cite{BS1}--\cite{BS3}. Уже в этих работах обнаружилось, какую важную роль двойные операторные интегралы играют в теории возмущений. Двойные операторные интегралы для произвольных ограниченных операторов $T$ определены в том случае, когда функция $\Phi$ является {\it мультипликатором Шура} по отношению к $E_1$ и $E_2$. В главе II мы изучаем пространство таких мультипликаторов Шура. При этом мы сначала изучаем так называемые {\it дискретные мультипликаторы Шура}, а затем с их помощью изучаем мультипликаторы Шура по отношению к спектральным мерам.

Далее, в главе III мы рассмотрим операторно липшицевы функции на замкнутых подмножествах комплексной плоскости $\C$. Для замкнутого подмножества $\fF$ комплексной плоскости $\C$ класс $\OL(\fF)$ {\it операторно липшицевых функций на} $\fF$ состоит из непрерывных функций $f$ на $\fF$ таких, что
\bay
\label{neognor}
\|f(N_1)-f(N_2)\|\le\const\|N_1-N_2\|
\ey
для произвольных нормальных операторов $N_1$ and $N_2$, спектры которых содержатся в
$\fF$. Мы также подробно изучаем класс {\it коммутаторно липшицевых функций на}
$\fF$, т.е. непрерывных функций $f$ на $\fF$ таких, что
$$
\|f(N_1)R-Rf(N_2)\|\le\const\|N_1R-RN_2\|
$$
для любого ограниченного оператора $R$ и произвольных нормальных операторов $N_1$ и $N_2$ со спектрами в $\fF$.
Для изучения этих классов функций мы используем результаты главы II.

При исследовании свойств операторно липшицевых функций на всей плоскости
так же, как и в случае самосопряжённых операторов, ключевую роль играет следующее обобщение неравенства \rf{sig}:
\bay
\label{nors}
\|f(N_1)-f(N_2)\|\le\const\s\|f\|_{L^\be}\|N_1-N_2\|
\ey
для произвольных нормальных операторов $N_1$ и $N_2$ с ограниченной разностью и для всякой ограниченной функции $f$ на $\R^2$ такой, что её преобразование Фурье сосредоточено в $[-\s,\s]\times[-\s,\s]$. Заметим, однако, что доказательство неравенства \rf{sig}, полученное в \cite{Pe1} и \cite{Pe3}, не распространяется на случай нормальных операторов, и в \cite{APPS} был найден новый метод получения таких оценок. 

Мы также получим достаточное условие для коммутаторной липшицевости функций, заданных на собственных замкнутых подмножествах плоскости, в терминах интегралов Коши мер на дополнение этого множества. Это условие было найдено в работе \cite{A2}. С помощью этого условия мы выводим достаточное условия Арази--Бартмана--Фридмана \cite{ABF} для коммутаторной липшицевости функций аналитических в круге, а также его аналог для функций, аналитических в полуплоскости.

Наконец, мы приведём в третьей главе результаты, группирующиеся вокруг результатов Э. Киссина и В.С. Шульмана \cite{KS3} о свойствах коммутаторно липшицевых функций на единичной окружности 
$\T$, допускающих аналитическое продолжение в единичных круг $\dd$.

В заключение обзора мы поместили раздел ``Заключительные замечания'', в котором кратко обсуждаются результаты, не вошедшие в этот обзор.

Авторы выражают искреннюю благодарность В.С. Шульману за полезные замечания.

\

\section{\bf Предварительные сведения}
\setcounter{equation}{0}
\label{Prel}

\

{\bf1. Классы Бесова.} Пусть $w$ -- бесконечно дифференцируемая функция на  $\R$ такая, что
\bay
\label{w}
w\ge0,\quad\supp w\subset\left[1/2,2\right],\quad\mbox{и} \quad w(s)=1-w\left(\frac s2\right)\quad\mbox{при}\quad s\in[1,2].
\ey

Определим функции $W_n$, $n\in\Z$, на $\R^d$ равенством 
$$
\big(\F W_n\big)(x)=w\left(\frac{\|x\|}{2^n}\right),\quad n\in\Z, \quad x=(x_1,\cdots,x_d),
\quad\|x\|\df\Big(\sum_{j=1}^dx_j^2\Big)^{1/2},
$$
где $\F$ -- {\it преобразование Фурье}, определённое на $L^1\big(\R^d\big)$ равенством
$$
\big(\F f\big)(t)=\!\int\limits_{\R^d} f(x)e^{-{\rm i}(x,t)}\,dx,\!\quad 
x=(x_1,\cdots,x_d),
\quad t=(t_1,\cdots,t_d), \!\quad(x,t)\df \sum_{j=1}^dx_jt_j.
$$
Ясно, что
$$
\sum_{n\in\Z}(\F W_n)(t)=1,\quad t\in\R^d\setminus\{0\}.
$$

Каждому медленно растущему распределению $f$ из ${\mathscr S}^\prime\big(\R^d\big)$, сопоставим последовательность $\{f_n\}_{n\in\Z}$,
\bay
\label{fn}
f_n\df f*W_n.
\ey
Формальный ряд $\sum_{n\in\Z}f_n$, являющийся разложением типа Винера--Пэли функции $f$, не обязательно сходится к $f$. 
Сначала мы определим (однородный) класс Бесова $\dot B^s_{p,q}\big(\R^d\big)$,
$s\in\R$, $0<p,\,q\le\be$, как пространство распределений 
$f$ класса ${\mathscr S}^\prime(\R^n)$ таких, что
\bay
\label{Wn}
\{2^{ns}\|f_n\|_{L^p}\}_{n\in\Z}\in\ell^q(\Z),\quad
\|f\|_{B^s_{p,q}}\df\big\|\{2^{ns}\|f_n\|_{L^p}\}_{n\in\Z}\big\|_{\ell^q(\Z)}.
\ey
В соответствии с этим определением $\dot B^s_{p,q}(\R^d)$ содержит все полиномы; при этом $\|f\|_{B^s_{p,q}}=0$ для всякого полинома $f$. Кроме того, распределение $f$ определено последовательностью $\{f_n\}_{n\in\Z}$
единственным образом с точностью до полинома. Легко видеть, что ряд 
$\sum_{n\ge0}f_n$ сходится в\footnote{Здесь и далее мы считаем, что 
в пространстве ${\mathscr S}^\prime(\R^d)$ введена слабая топология
$\s\big({\mathscr S}^\prime(\R^d),{\mathscr S}(\R^d)\big)$} ${\mathscr S}^\prime(\R^d)$. 
Однако ряд $\sum_{n<0}f_n$, вообще говоря, может расходиться. Можно доказать, тем не менее, что ряды
\bay
\label{ryad}
\sum_{n<0}\frac{\partial^r f_n}{\partial x_1^{r_1}\cdots\partial x_d^{r_d}}\qquad \mbox{при}\quad r_j\ge0,\quad
1\le j\le d,\quad\sum_{j=1}^dr_j=r,
\ey
равномерно сходятся в $\R^d$, если $r\in\Z_+$ и 
$r>s-d/p$. Отметим, что при $q\le1$ ряды в \rf{ryad}
сходятся равномерно при более слабом условии $r\ge s-d/p$.

Теперь мы можем определить модифицированное (однородное) пространство Бесова $B^s_{p,q}\big(\R^d\big)$. Скажем, что распределение $f$
принадлежит классу $B^s_{p,q}(\R^d)$, если выполнено условие \rf{Wn} и
$$
\frac{\partial^r f}{\partial x_1^{r_1}\cdots\partial x_d^{r_d}}
=\sum_{n\in\Z}\frac{\partial^r f_n}{\partial x_1^{r_1}\cdots\partial x_d^{r_d}}\quad
\mbox{при}\quad 
r_j\ge0,\quad\mbox
1\le j\le d,\quad\sum_{j=1}^dr_j=r,
$$
в пространстве ${\mathscr S}^\prime\big(\R^d\big)$, где $r$ --
минимальное неотрицательное целое число, при котором
$r>s-d/p$ ($r\ge s-d/p$, если $q\le1$). Теперь функция $f$ однозначно определяется последовательностью $\{f_n\}_{n\in\Z}$ с точностью до полинома степени
меньше, чем $r$. При этом полином $g$ входит в 
$B^s_{p,q}\big(\R^d\big)$ тогда и только тогда, когда $\deg g<r$. 

В случае $p=q$ мы используем обозначение $B_p^s(\R^d)$ для $B_{p,p}^s(\R^d)$.

Рассмотрим теперь шкалу $\L_\a(\R^d)$, $\a>0$, {\it классов Гёльдера--Зигмунда}. Их можно определить равенством $\L_\a(\R^d)\df B_\be^\a(\R^d)$. 

%
%
%

Классы Бесова допускают много других описаний.
Приведём описание в терминах конечных разностей.
Для $h$ из $\R^d$, определим разностный оператор $\D_h$ равенством
$(\D_hf)(x)=f(x+h)-f(x)$, $x\in\R^d$.



Пусть $s>0$,  $m\in\Z$ и $m-1\le s<m$.  
Пусть $p,q\in[1,+\be]$.
Класс Бесова $B_{p,q}^s\big(\R^d\big)$ может быть определён, как множество
функций $f$ из $L^1_{\rm loc}\big(\R^d\big)$ таких, что
$$
\int_{\R^d}|h|^{-d-sq}\|\D^m_h f\|_{L^p}^q\,dh<\be,~\; q<\be;\qquad
\sup_{h\not=0}\frac{\|\D^m_h f\|_{L^p}}{|h|^s}<\be,~\; q=\be.
$$
Однако при этом определении классы Бесова могут содержать полиномы более высоких степеней, чем в случае определения в терминах свёрток с функциями $W_n$.

Пространства $B_{pq}^s$ может быть определено в терминах интеграла Пуассона.
Пусть $P_d(x,t)$ ядро Пуассона, решающее задачу Дирихле для полупространства
$\R^{d+1}_+\df\{(x,t): x\in\R^d, t>0\}$, т. е. $P_d(x,t)=c_dt(|x|^2+t^2)^{-\frac{d+1}2}$,
$c_d=\pi^{-\frac{d+1}2}\G(\frac{d+1}2)$. С каждой функцией 
$f\in L^1\big(\R^d,(\|x\|+1)^{-(d+1)}\,dx\big)$
можно связать интеграл Пуассона $\mP f$
$$
(\mP f )(x,t)=\int_{\R^d}P_d(x-y,t)f(y)\,dy.
$$
Тогда для любого целого положительного $m$ имеет место равенство
$$
\frac{\partial^m (\mP f)}{\partial^m t}(x,t)=\int_{\R^d}\frac{\partial^m P_d(x-y,t)}{\partial^m t}f(y)\,dy.
$$
Заметим, что второй интеграл имеет смысл при всех  
$f\in L^1\big(\R^d,(\|x\|+1)^{-(d+m+1)}\,dx\big)$.
Это обстоятельство позволяет корректно определить $\frac{\partial^m}{\partial t^m}\mP$ для 
функций $f$ класса $L^1\big(\R^d,(\|x\|+1)^{-(d+m+1)}\,dx\big)$.

Пусть $m\in\Z$, \mbox{$m-1\le s<m$}, $1\le p,\:q\le+\be$. Тогда мы можем определить пространство
$B_{pq}^s$ как множество всех функций 
$f\in L^1\big(\R^d,(\|x\|+1)^{-(d+m+1)}\,dx\big)$
таких, что 
$$
\left(\int_0^{\be}t^{(m-s)q-1}\left\|\Big(\frac{\partial^m}{\partial t^m}\mP f\Big)(\cdot,t)\right\|_{L^p(\R^d)}^q\,dt\right)^{\frac1q}
<+\be,\quad q<+\be,
$$
$$
\sup_{t>0}t^{m-s}\left\|\Big(\frac{\partial^m}{\partial t^m}\mP f\Big)(\cdot,t)\right\|_{L^p(\R^d)}<+\be,\quad
q=+\be.
$$
При этом определении классы Бесова тоже могут содержать полиномы более высоких степеней, 
чем в случае определения в терминах свёрток с функциями $W_n$.
Отметим ещё, что это определения в терминах интеграла Пуассона при определённых условиях
и оговорках работает и в том случае, когда $p<1$  или $q<1$.

\medskip

%

Перейдём теперь к {\it классам Бесова функций на единичной окружности $\T$}. 
Пусть $w$ -- функция, удовлетворяющая условиям \rf{w}.
Определим тригонометрические полиномы $W_n$, $n\ge0$, равенством
$$
W_n(\z)\df\sum_{j\in\Z}w\left(\frac{|j|}{2^n}\right)\z^j,\quad n\ge1,
\quad W_0(\z)\df\sum_{\{j:\,|j|\le1\}}\z^j,\quad\z\in\T.
$$
Если $f$ -- распределение на $\T$ положим
$$
f_n=f*W_n,\quad n\ge0,
$$
и скажем, что $f$ входит в класс Бесова $B_{p,q}^s(\T)$, $s\in\R$, 
$0<p,\,q\le\be$, если
\bay
\label{Bperf}
\big\{2^{ns}\|f_n\|_{L^p}\big\}_{n\ge0}\in\ell^q.
\ey

Пусть $s\in\R$, $s>\max\{0,1/p-1\}$, $m$ -- целое положительное число,
такое, что $m>\max\{s,s+1/p-1\}$. Тогда обобщённая функция $f$ на окружности $\T$ принадлежит пространству
$B_{p,q}^s(\T)$  в том и только в том случае, когда 
$$
\int_0^{1}r(1-r^2)^{(m-s)q-1}\left\|\frac{\partial^m}{\partial r^m}\Big((\mP f)(r\z)\Big)\right\|_{L^p(\T)}^q\,dr
<+\be,\quad q<+\be,
$$
$$
\sup_{r\in[0,1)}(1-r^2)^{m-s}\left\|\frac{\partial^m}{\partial r^m}\Big((\mP f)(r\z)\Big)\right\|_{L^p(\T)}<+\be,\quad q=+\be,
$$
где $\mP f$  обозначает интеграл Пуассона обобщённой функции $f$.

В определениях классов Бесова в терминах интеграла Пуассона мы рассматривали производную 
порядка $m$ по переменной $t$ в первом случае и по переменной $r$ --  во втором.
Хорошо известно, что в обоих случаях мы получили бы равносильное определение,
если бы потребовали конечность аналогичных выражение для всевозможных (в том числе
и смешанных) частных
производных порядка $m$.

Мы отсылаем читателя к \cite{Pee} и \cite{T} для более подробной информации о классах Бесова.

\medskip

{\bf2. Классы Шаттена--фон Неймана.} Для ограниченного линейного оператора $T$ в 
гильбертовом пространстве его {\it сингулярные числа} $s_j(T)$, $j\ge0$, определяются равенством
$$
s_j(T)\df\inf\big\{\|T-R\|:~\rank R\le j\big\}.
$$
{\it Класс Шаттена--фон Неймана} $\bS_p$, $0<p<\be$, состоит по определению из
операторов $T$, для которых
$$
\|T\|_{\bS_p}\df\Big(\sum_{j\ge0}\big(s_j(T)\big)^p\Big)^{1/p}<\be.
$$
При $p\ge1$ это -- нормированный идеал операторов в гильбертовом пространстве. 
Класс $\bS_1$ называется {\it классом ядерных операторов}. Если $T$ -- ядерный оператор в гильбертовом пространстве $\h$, то его {\it след} $\trace T$ определяется равенством
$$
\trace T\df\sum_{j\ge0}(Te_j,e_j),
$$
где $\{e_j\}_{j\ge0}$ -- ортонормированный базис в $\h$. Правая часть не зависит от выбора базиса.

Класс $\bS_2$ называется {\it классом Гильберта--Шмидта}. Он образует гильбертово пространство со скалярным произведением
$$
(T,R)_{\bS_2}\df\trace(TR^*).
$$

При $p\in(1,\be)$ двойственное пространство $(\bS_p)^*$ можно изометрическим образом отождествить с пространством $\bS_{p'}$, $1/p+1/p'=1$, с помощью билинейной формы
$$
\langle T,R\rangle\df\trace(TR).
$$

Пространство, сопряжённое к $\bS_1$, можно отождествить с пространством всех ограниченных линейных операторов с помощью этой же билинейной формы, в то время как пространство, сопряжённое к пространству компактных операторов, отождествляется с $\bS_1$.

Отсылаем читателя к \cite{GK} за более подробной информацией.

\medskip

{\bf3. Операторы Ганкеля.} Для функции $\f$ класса $L^2$ на единичной окружности 
$\T$ {\it оператор Ганкеля} $H_\f$ определяется на плотном в пространстве Харди $H^2$ множестве полиномов равенством $H_\f f\df\pp_-\f f$, где $\pp_-$ ортогональный проектор из $L^2$ на $H^2$. По теореме Нехари оператор $H_\f$ продолжается до ограниченного оператора из $H^2$ в $H^2_-$ тогда и только тогда, когда существует функция $\psi$ класса $L^\be$ на $\T$ такая, что её коэффициенты Фурье $\hat\psi(n)$ удовлетворяют равенству $\hat\psi(n)=\hat\f(n)$ при $n<0$. В свою очередь, по теореме Ч. Феффермана последнее свойство эквивалентно тому, что 
$\pp_-\f$ входит в класс ${\rm BMO}$. 

Оператор Ганкеля $H_\f$ входит в класс Шаттена--фон Неймана $\bS_p$ в том и только в том случае, когда функция $\pp_-\f$ входит в класс Бесова $B_p^{1/p}(\T)$. При 
$p\ge1$ это было доказано в работе \cite{Pe0}, а при $p\in(0,1)$ в \cite{Pe<1}, см. также \cite{Pek} и \cite{Se}, где даны другие доказательства при $p<1$.

Легко видеть, что оператор $H_\f$ имеет в базисах
$\{z^j\}_{j\ge0}$ и $\{\ov z^k\}_{k\ge1}$ матрицу 
\lb$\{\hat\f(-j-k-1)\}_{j\ge0,\,k>1}$. Такие матрицы, то есть матрицы вида
$\{\a_{j+k}\}_{j,k\ge0}$ называются {\it ганкелевыми матрицами}. Критерий принадлежности классу $\bS_p$ операторов Ганкеля может быть переформулирован следующим образом: {\it оператор в пространстве $\ell^2$ с ганкелевой матрицей $\{\a_{j+k}\}_{j,k\ge0}$ входит в класс $\bS_p$, $p>0$, в том и только в том случае, когда функция $\sum_{j\ge0}\a_jz^j$ входит в класс Бесова
$B_p^{1/p}(\T)$}.

Аналогичным образом определяются операторы Ганкеля на классе Харди $H^2(\C_+)$ функций в верхней полуплоскости. При этом имеют место аналоги перечисленных выше результатов для операторов Ганкеля на классе $H^2$ функций в единичном круге.

Мы отсылаем читателя к монографии \cite{Pe} за доказательствами перечисленных выше результатов и более подробной информацией об операторах Ганкеля.

\

\begin{center}
\bf\large Обозначения
\end{center}
\label{obozna}

\setcounter{section}{0}

\addtocontents{toc}{\vspace*{.03cm}\hspace*{.55cm}\textbf{{\bf Обозначения}}\hfill\pageref{obozna}\\[.25cm]}

\

Здесь приводится список некоторых обозначений, используемых в этом обзоре.

\medskip

$\OL(\fF)$ -- пространство операторно липшицевых функций на замкнутом подмножестве $\fF$ комплексной плоскости $\C$;

$\CL(\fF)$ -- пространство коммутаторно липшицевых функций на замкнутом подмножестве $\fF$ комплексной плоскости $\C$;

$\OD(\R)$ -- пространство операторно дифференцируемых функций на $\R$;

$\mB(\h)$ -- пространство ограниченных линейных операторов в гильбертовом пространстве $\h$;

$\mB(\h_1,\h_2))$ -- пространство ограниченных линейных операторов из гильбертова пространства $\h_1$ в гильбертово пространство $\h_2$;

$\mB_{\rm sa}(\h)$ -- пространство ограниченных самосопряжённых операторов в гильбертовом пространстве $\h$;

$\m$ -- нормированная мера Лебега на единичной окружности $\T$;

$\m_2$ -- мера Лебега на плоскости.

\

\begin{center}
\bf\Large Глава I
\end{center}

\medskip

\begin{center}
\bf\Large 
Операторно липшицевы функции на прямой и на окружности. По первому кругу
\end{center}
\label{1per}
\renewcommand{\thesection}{1.\arabic{section}}
\setcounter{section}{0}

\addtocontents{toc}{\vspace*{-.33cm}\hspace*{.15cm}\textbf{{\bf Глава I}. Операторно липшицевы функции на прямой и на окружности.\\[.3cm] \hspace*{1.48cm}По первому кругу}\hfill\pageref{1per}}

\

В этой вводной главе мы рассмотрим операторно липшицевы функции на вещественной прямой 
$\R$ и на единичной окружности $\T$. 

Позже, в главе III, мы подвергнем класс операторно липшицевых функций более подробному изучению, а также перейдём к изучению
операторно липшицевых функций на замкнутых подмножествах комплексной плоскости $\C$.

Мы используем обозначение $\OL(\R)$ для класса операторно липшицевых функций на $\R$ и для функции $f$ из $\OL(\R)$ и положим
$$
\|f\|_{\OL(\R)}\df\sup\left\{\frac{\|f(A)-f(B)\|}{\|A-B\|}:~
A~\mbox{и}~B~\mbox{--~самосопряжённые операторы},~A\ne B\right\}.
$$
Аналогичным образом вводится пространство $\OL(\T)$ операторно липшицевых функций на 
окружности $\T$.

Оказывается, что класс $\OL(\R)$ имеет несколько необычные свойства. В частности, функции из этого класса должны быть дифференцируемы всюду на $\R$, а также должны иметь производную на бесконечности, то есть предел
$\lim_{|t|\to\be}\frac{f(t)}{t}$
должен существовать (см. теорему \ref{pro} ниже). Заметим, что отсюда вытекает упомянутый во введении результат Макинтоша--Като: функция $x\mapsto|x|$ не является операторно липшицевой. С другой стороны, функции класса $\OL(\R)$ не обязаны быть непрерывно дифференцируемыми. В частности, функция $x\mapsto x^2\sin(1/x)$, не будучи непрерывно дифференцируемой, является операторно липшицевой, см. теорему \ref{teorox2f1x} ниже. 

Эту главу мы начнём с элементарных примеров операторно липшицевых функций (см. \S\:\ref{prim}).

Мы введём в \S\:\ref{odif} класс операторно дифференцируемых функций 
и класс локально операторно дифференцируемых функций. Причём оказывается, что при определении этих классов неважно, рассматривать ли дифференцируемость в смысле Гато или в смысле Фреше. Оказывается, что (локально) операторно дифференцируемые функции должны быть непрерывно дифференцируемы. Мы увидим, что операторно дифференцируемые функции обязательно должны быть операторно липшицевыми. Однако не всякая операторно липшицева функция является операторно дифференцируемой.

Помимо операторно липшицевых функций мы рассмотрим в \S\:\ref{CoL} {\it коммутаторно липшицевы функции}, то есть такие функции $f$ на $\R$, при которых имеет место неравенство
$$
\|f(A)R-Rf(B)\|\le\const\|AR-RB\|
$$
для любых самосопряжённых операторов $A$ и $B$ (опять, неважно, ограниченных или необязательно ограниченных) и для любого ограниченного оператора $R$. 
{\it Коммутаторно липшицева норма}
$\|f\|_{\CL(\R)}$ функции $f$ определяется, как минимальная константа, при которой это неравенство справедливо. 
Аналогичным образом можно определить коммутаторно липшицевы функции на единичной окружности, если мы заменим самосопряжённые операторы на унитарные.

Оказывается, что для функций {\it на прямой (равно как и для функций на окружности) класс коммутаторно липшицевых функции совпадает с классом операторно липшицевых функций}. Позже, в главе III мы увидим, что для функций на произвольном замкнутом подмножестве плоскости $\R^2$ это уже не так.

Мы получим в этой главе достаточное условие для операторной липшицевости на прямой и на окружности (см. \S\:\ref{Dost}), равно как и необходимое (см. \S\:\ref{Neob}) и сравним их друг с другом.

Наряду с операторно липшицевыми функциями было бы естественно рассмотреть класс {\it операторно гёльдеровых функций порядка} $\a$, $0<\a<1$, то есть класс функций $f$, для которых справедливо неравенство
$$
\|f(A)-f(B)\|\le\const\|A-B\|^\a,
$$
для самосопряжённых операторов $A$ и $B$ в гильбертовом пространстве. Однако  этот термин оказывается коротко живущим, ибо, {\it в отличие от случая операторно липшицевых функций, всякая функция $f$ на $\R$ класса Гёльдера порядка $\a$ с необходимостью является операторно гёльдеровой порядка $\a$,} см. 
\S\:\ref{oHold} 
(напомним, что $f$ принадлежит {\it классу Гёльдера} $\L_\a(\R)$, если
$|f(x)-f(y)|\le\const|x-y|^\a$ при всех $x$ и $y$ из $\R$).

\

\section{\bf Элементарные примеры операторно липшицевых функций}
\setcounter{equation}{0}
\label{prim}

\

В этом параграфе мы приведём примеры операторно липшицевых функция на прямой и простые достаточные условия для операторной липшицевости. В дальнейшем мы получим и другие достаточные условия.

\medskip

{\bf Пример 1.}  {\it При всех $\l$ из $\C\setminus\R$ функция $(\l-x)^{-1}$ является операторно липшицевой на $\R$
и $\|(\l-x)^{-1}\|_{\OL(\R)}=|\im\l|^{-2}$}.

\medskip

\Pf Из резольвентного тождества Гильберта
$$
(\l I-A)^{-1}-(\l I-B)^{-1}=(\l I-A)^{-1}(A-B)(\l I-B)^{-1}
$$
мгновенно вытекает, что $\|(\l-x)^{-1}\|_{\OL(\R)}\le|\im\l|^{-2}$. Остаётся заметить, что \lb
$\|(\l-x)^{-1}\|_{\OL(\R)}\ge\|(\l-x)^{-1}\|_{\Li(\R)}=|\im\l|^{-2}$. $\bl$

%

\medskip

{\bf Пример 1$'$.}  {\it При всех $\l$ из $\C\setminus\T$ функция $(\l-z)^{-1}$ является операторно липшицевой на $\T$
и $\|(\l-z)^{-1}\|_{\OL(\T)}=(|\l|-1)^{-2}$}.

\medskip

{\bf Пример 2.}  {\it Функция $x\mapsto\log(1+{\rm i}x)$ является операторно липшицевой на $\R$ и
$\|\log(1+{\rm i}x)\|_{\OL(\R)}=1$. Здесь $\log$ обозначает главную ветвь логарифма.}

\medskip

\Pf
Ясно, что
$$
\log(1+{\rm i}x)=\int_0^{+\be}\left(\frac1{1+t}-\frac1{1+t+{\rm i}x}\right)\,dt.
$$
Отсюда следует, что
\begin{align*}
\|\log(1+{\rm i}x)\|_{\OL(\R)}
&\le\int_0^{+\be}\left\|\frac1{1+t}-\frac1{1+t+{\rm i}x}\right\|_{(\OL)(\R)}dt\\[.2cm]
&=\int_0^{+\be}\left\|\frac1{1+t+{\rm i}x}\right\|_{(\OL)(\R)}dt
=\int_0^{+\be}\frac{dt}{(1+t)^2}=1.
\end{align*}
С другой стороны, неравенство $\|\log(1+{\rm i}x)\|_{\OL(\R)}\ge1$ очевидно, ибо 
$$
\|\log(1+{\rm i}x)\|_{\OL(\R)}\ge\|\log(1+{\rm i}x)\|_{\Li(\R)}=1.\quad\bl
$$

\medskip

Аналогичным образом можно доказать, что
при всех $\l$ из $\C\setminus\R$ имеет место равенство
$\|\log(\l-x)\|_{\OL(\R)}=|\im\l|^{-1}$, где $\log(\l-x)$ обозначает любую
регулярную на $\R$ ветвь функции $\log(\l-z)$.

%
%

\medskip

{\bf Пример 3.}  {\it Функция $\arctg$ является операторно липшицевой функцией на $\R$,
и $\|\arctg\|_{\OL(\R)}=1$}.

\medskip

\Pf Как и в предыдущем примере, 
достаточно проверить, что \lb$\|\arctg\|_{\OL(\R)}\le1$. Для этого достаточно
заметить, что $\arctg x=\im\log(1+{\rm i}x)$ при всех $x\in\R$. $\bl$

\medskip

{\bf Пример 4.} {\it При любом натуральном числе $n$ имеет место следующее равенство:
$$
\|(\l -x)^{-n}\|_{\OL(\R)}=n|\im\l|^{-n-1}\quad\mbox{при всех}\quad\l\in\C\setminus\R.
$$}

\medskip

\Pf Подставляя в элементарное тождество
\bay
\label{xny}
X^n-Y^n=\sum_{k=1}^nX^{n-k}(X-Y)Y^{k-1}
\ey
$X=(\l I-A)^{-1}$ и $Y=(\l I-B)^{-1}$, получаем
$$
(\l I-A)^{-n}-(\l I-B)^{-n}=\sum_{k=1}^n(\l I-A)^{k-n}\big((\l I-A)^{-1}-(\l I-B)^{-1}\big)(\l I-B)^{1-k}.
$$
Следовательно, для любых самосопряжённых операторов $A$ и $B$
имеем:
\begin{align*}
\|(\l I-A)^{-n}&-(\l I-B)^{-n}\|\\[.2cm]
&\le\sum_{k=1}^n\|(\l I-A)^{k-n}\|\cdot\big\|\big((\l I-A)^{-1}-(\l I-B)^{-1}\big)
\big\|\cdot\|(\l I-B)^{1-k}\|\\[.2cm]
&\le\sum_{k=1}^n|\im\l|^{k-n}|\im\l|^{-2}\|A-B\|\cdot|\im\l|^{1-k}=n|\im\l|^{-n-1}\|A-B\|.
\end{align*}
Таким образом, мы доказали, что
$\|(\l I-x)^{-n}\|_{\OL(\R)}\le n|\im\l|^{-n-1}$.
Остаётся заметить, что $\|(\l -x)^{-n}\|_{\OL(\R)}\ge\|(\l -x)^{-n}\|_{\Li(\R)}=n|\im\l|^{-n-1}$. $\bl$

%
%

\medskip

{\bf Пример 5.} {\it Для каждого вещественного числа  $a$  функция 
$x\mapsto e^{{\rm i}ax}$ на $\R$ является операторно липшицевой
 и $\|e^{{\rm i}ax}\|_{\OL(\R)}=|a|$}.

\medskip
 
\Pf Снова достаточно установить только оценку сверху. Кроме того, можно
считать, что $a=1$. Пусть $A$ и $B$ -- самосопряжённые операторы.  Тогда
$$
\Big(e^{{\rm i}tA}e^{-{\rm i}tB}\Big)'={\rm i}Ae^{{\rm i}tA}e^{-{\rm i}tB}-{\rm i}e^{{\rm i}tA}e^{-{\rm i}tB}B
={\rm i}e^{{\rm i}tA}(A-B)e^{-{\rm i}tB},
$$
откуда
\begin{align*}
\|e^{{\rm i}A}-e^{{\rm i}B}\|&=\|e^{{\rm i}A}e^{-{\rm i}B}-I\|
=\left\|{\rm i}\int_0^1e^{{\rm i}tA}(A-B)e^{-{\rm i}tB}\,dt\right\|\\[.2cm]
&\le\int_0^1\big\|e^{{\rm i}tA}(A-B)e^{-{\rm i}tB}\big\|\,dt
=\int_0^1\|A-B\|\,dt=\|A-B\|.\quad\bl
\end{align*}

Во всех приведённых выше примерах имеет место
равенство $\|f\|_{\OL(\R)}=\|f'\|_{L^\be(\R)}$, что скорее является
исключением, чем правилом.

Из примера 5 мгновенно вытекает следующее утверждение.

\begin{thm}
\label{char}
Пусть $f$ -- первообразная преобразования Фурье $\F\mu$
комплексной борелевской меры $\mu$  на $\R$. Тогда   $f\in\OL(\R)$
и $\|f\|_{\OL(\R)}\le\|\mu\|$.
\end{thm}

\Pf  Можно считать, что $f(0)=0$. Тогда
\begin{align*}
f(x)&=\int_0^x(\F\mu)(t)\,dt=
\int_0^x\left(\int_\R e^{{-\rm i}st}\,d\mu(s)\right)\,dt\\[.2cm]
&=\int_0^1\left(\int_\R xe^{{-\rm i}stx}\,d\mu(s)\right)\,dt
={\rm i}\int_\R\frac{e^{{-\rm i}sx}-1}{s}\,d\mu(s).
\end{align*}
Следовательно,
$$
\|f\|_{\OL(\R)}\le\int_\R\left\|\frac{e^{{-\rm i}sx}-1}{s}\right\|_{\OL(\R)}\,d|\mu|(s)
\le\int_\R d|\mu|(s)=\|\mu\|.\quad\bl
$$

\begin{cor}
\label{corchar}
Пусть $f\in C^1(\R)$. Предположим, что функция $f'$ является  положительно определённой.
Тогда $\|f\|_{\OL(\R)}=\|f\|_{\Li(\R)}=f'(0)$.
\end{cor}

\Pf В силу классической теоремы Бохнера, см., например, \cite{I}, положительно определённую 
функцию $f'$  можно представить в виде  $f'=\F\mu$, где $\mu$ -- конечная 
борелевская положительная мера на $\R$. Остаётся заметить, что
$$
\|\mu\|=f'(0)=|f'(0)|\le\|f\|_{\Li(\R)}\le\|f\|_{\OL(\R)}\le \|\mu\|,
$$
где последнее неравенство вытекает из теоремы \ref{char}. $\bl$

Следует отметить, что в этом параграфе по существу все приведённые выше примеры
явного вычисления операторной полунормы в пространстве $\OL(\R)$ так или иначе
основаны на следствии \ref{corchar}.

Тем не менее, можно построить пример функции $f\in\OL(\R)$ такой,
что $\|f\|_{\OL(\R)}=\|f\|_{\Li(\R)}=f'(0)$ и $f$ не удовлетворяет условиям следствия  \ref{corchar}.

С другой стороны, если $\|f\|_{\OL(\R)}=\|f\|_{\Li(\R)}=\frac{f(a)-f(0)}{a}=1$ при  $a\in\R$, $a\ne0$,
то $f(x)=x+f(0)$ при всех $x\in\R$.

%
%

Пример 5 допускает ещё одно обобщение, так называемое 
операторное неравенство Бернштейна. Об этом будет идти речь в \S\:\ref{Bern}.
В частности, в \S\:\ref{Bern} будет показано, что $L^\be(\R)\cap \E_\s\subset\OL(\R)$, где символом $\E_\s$ обозначается пространство целых функций экспоненциального типа не выше, чем $\s$. 

Рассмотрим теперь примеры операторно липшицевых функций на
единичной \lb окружности $\T$.

\medskip

{\bf Пример 6.}  {\it Пусть $n\in\Z$.
Тогда $\|z^n\|_{\OL(\T)}=|n|$ для всех $n\in\Z$}.

\medskip

\Pf  Достаточно рассмотреть случай, когда $n>0$. В этом случае
всё сводится к проверке следующего неравенства:
$$
\|U^n-V^n\|\le n\|U-V\|
$$
для любых унитарных операторов $U$ и $V$.  Чтобы его доказать,
достаточно подставить  $X=U$ и $Y=V$ в тождество \rf{xny}. $\bl$

\medskip

Из этого примера мгновенно вытекает аналог теоремы \ref{char} для окружности.

\begin{thm} 
Пусть $f$ -- непрерывная функция на единичной окружности $\T$
такая, что $\sum\limits_{n\in\Z}|n|\cdot|\hat f(n)|<\be$. Тогда $f\in\OL(\T)$
и $\|f\|_{\OL(\T)}\le\sum\limits_{n\in\Z}|n|\cdot|\hat f(n)|$.
\end{thm}

Отметим, что более сильные результаты будут вскоре приведены в 
\S\:\ref{Dost}.

\medskip 

{\bf Пример 7.} Функция $x\mapsto x^2\sin\frac1x$ является операторно липшицевой.
Чтобы убедиться в этом, докажем следующую теорему:

\begin{thm}
\label{teorox2f1x}
Пусть $f\in\OL(\R)$  и $f(0)=0$.
Положим
$$
g(x)\df\left\{\begin{array}{ll}x^{2}f(x^{-1}),&\text {если}\,\,\,x\not=0,\\[.2cm]
\quad0,&\text {если}\,\,\,\,x=0.
\end{array}\right.
$$
Тогда $g\in\OL(\R)$ и
\bay
\label{nervox21x}
\frac13\|f\|_{\OL(\R)}\le\|g\|_{\OL(\R)}\le3\|f\|_{\OL(\R)}.
\ey
\end{thm}

\medskip

\Pf Достаточно доказать только второе неравенство, поскольку первое неравенство
сводится ко второму. Можно считать, что $\|f\|_{\OL(\R)}=1$. 
Как мы отмечали во введении к этой главе, для функций на прямой операторная липшицевость эквивалентна коммутаторной липшицевости, см. \S\:\ref{CoL}
и \S\:\ref{oplip} ниже. Более того, соответствующие нормы совпадают: 
$\|f\|_{\OL(\R)}=\|f\|_{\CL(\R)}$. Поэтому достаточно доказать, что неравенство
\bay
\label{71}
\|f(A)R-Rf(A)\|\le\|AR-RA\|
\ey
для любого оператора $R$ и любого самосопряжённого оператора $A$
влечёт, что
$$
\|g(A)R-Rg(A)\|\le3\|AR-RA\|
$$
для любого оператора $R$  и любого самосопряжённого оператора $A$.
Предположим сначала, что оператор $A$ обратим. Этот случай сводится к следующему утверждению:
\bay
\label{72}
\|A^2f(A^{-1})R-RA^2f(A^{-1})\|\le3\|AR-RA\|
\ey
для любого оператора $R$  и любого самосопряжённого оператора $A$.
Имеем:
\begin{align*}
f(A^{-1})A^2R&-RA^2f(A^{-1})=f(A^{-1})A(AR-RA)\\[.2cm]
&+f(A^{-1})ARA-ARAf(A^{-1})
+(AR-RA)Af(A^{-1}). 
\end{align*}
Ясно, что
$$
\|Af(A^{-1})\|\le\sup_{t\ne0}|t^{-1} f(t)|\le\|f\|_{\Li(\R)}\le\|f\|_{\OL(\R)}=1.
$$
Следовательно, 
$$
\|f(A^{-1})A(AR-RA)\|\le\|AR-RA\|\quad\mbox{и}\quad\|(AR-RA)Af(A^{-1})\|\le\|AR-RA\|.
$$
Наконец, подставляя в \rf{71} вместо операторов $R$  и $A$ операторы $ARA$  и $A^{-1}$, получаем
$$
\|f(A^{-1})ARA-ARAf(A^{-1})\|\le\|A^{-1}ARA-ARAA^{-1}\|=\|AR-RA\|,
$$
откуда мгновенно следует \rf{72}.
Дабы рассмотреть общий случай, достаточно заметить, что для любого положительного
числа $\d$ найдётся обратимый самосопряжённый оператор $A_\d$ такой, что
$AA_\d=A_\d A$ и $\|A-A_\d\|<\d$. Тогда получаем:
\begin{align*}
\|g(A)R&-Rg(A)\|\le\|g(A)-g(A_\d)\|\cdot\|R\|+\|g(A_\d)R-Rg(A_\d)\|\\[.2cm]
&+\|g(A_\d)-g(A)\|\cdot\|R\|
\le2\d\|R\|\cdot\|g\|_{\Li(\R)}+3\|A_\d R-RA_\d\|\\[.2cm]
&\le6\d\|R\|\cdot\|f\|_{\Li(\R)}+3\|AR-RA\|+6\d\|R\|
\le3\|AR-RA\|+12\d\|R\|
\end{align*}
при всех $\d>0$. $\bl$

\medskip

{\bf Замечание.} Теперь ясно, ввиду примера 5, что
функция $g$, определённая равенством
$g(x)=x^2\sin\frac1x$, является операторно липшицева.  {\it Функция $g$ даёт нам пример операторно липшицевой функции, не являющейся непрерывно дифференцируемой}. Напомним (см. теорему \ref{pro} ниже), что всякая операторно липшицевая функция на $\R$ 
должна быть дифференцируемой в каждой точке.

Отметим ещё, что в \cite{A2} доказано, что подмножество вещественной прямой является множеством точек разрыва производной операторно липшицевой функции в том и только в том случае, когда оно является множеством первой категории
типа  $F_\s$. Другими словами, множества точек разрыва 
производных операторно липшицевых функций устроены так же, как множества точек разрыва функций первого класса Бэра.

\medskip

В \S\:\ref{dlp} получен результат, более общий, чем теорема \ref{teorox2f1x}, в котором вместо дробно-линейного преобразования $x\mapsto x^{-1}$ рассматриваются произвольные дробно-линейные преобразования.

\

\section{\bf Операторная липшицевость в сравнении с операторной \\
дифференцируемостью}
\setcounter{equation}{0}
\label{odif}

\

Пусть $H$ -- функция со значениями в банаховом пространстве $X$, заданная на подмножестве $\L$ вещественной прямой $\R$. Функция $H$ называется {\it липшицевой}, если существует неотрицательное число
$c$ такое, что
\bay
\label{Xlip+}
\|H(s)-H(t)\|_X\le c|s-t|,\quad s,\:t\in\L.
\ey
Множество всех таких функций обозначим
 символом $\Li(\L,X)$. Пусть $\|H\|_{\Li(\L,X)}$  обозначает наименьшую
 константу $c$, удовлетворяющую условию \rf{Xlip+}. 
Как обычно, мы полагаем $\|H\|_{\Li(\L,X)}\df\be$, если $H\not\in\Li(\L,X)$.
 
Пусть $f$ -- непрерывная функция на вещественной прямой $\R$.
С каждым самосопряжённым оператором $A$ и ограниченным самосопряжённым
оператором $K$ мы связываем функцию $H_{f,A,K}$, $H_{f,A,K}(t)=f(A+tK)-f(A)$,
заданную при всех $t$ из $\R$ таких, что $f(A+tK)-f(A)\in\mB(\h)$.

Заметим, что если $f\in\OL(\R)$, то $H_{f,A,K}\in\Li(\R,\mB(\h))$ и
$\|H_{f,A,K}\|_{\Li(\R,\mB(\h))}\le\|K\|\cdot\|f\|_{\OL(\R)}$.
Легко видеть, что имеет место следующее утверждение.

\begin{lem} 
\label{HfAK+}
Пусть $f$ --  непрерывная функция на $\R$. Тогда
\begin{align*}
\|f\|_{\OL(\R)}&=\sup\big\{\|H_{f,A,K}\|_{\Li(\R,\mB(\h))}:~A,~K\in\mB_{\rm sa}(\h),\|K\|=1\big\}\nonumber\\[.2cm]
&=\sup\big\{\|H_{f,A,K}\|_{\Li(\R,\mB(\h))}:~K\in\mB_{\rm sa}(\h),~
\|K\|=1,~A^*=A \big\}.\quad\bl
\end{align*} 
\end{lem}

 
Нам понадобится следующее хорошо известное элементарное утверждение.
Для удобства читателя мы приведём здесь одно из возможных доказательств.
 
\begin{lem} 
\label{verh+}
Пусть $H$ --  функция со значениями в банаховом пространстве $X$, заданная на 
невырожденном промежутке $\L$, $\L\subset\R$. Тогда
$$
\|H\|_{\Li(\L,X)}=\sup_{t\in\L}\varlimsup_{h\to0}\frac{\|H(t+h)-H(t)\|_X}{|h|}.
$$
\end{lem}

\Pf Неравенство $\ge$ очевидно.
Чтобы доказать противоположное неравенство, достаточно убедиться в том, что
неравенство \rf{Xlip+} справедливо, как только $c$ удовлетворяет условию
\bay
\label{cXlip+}
c>\sup_{t\in\L}\varlimsup_{h\to0}\frac{\|H(t+h)-H(t)\|_X}{|h|}.
\ey
Зафиксируем такое число $c$ и произвольную точку $t$ из $\L$.  Пусть $\L_t$ --множество
всех точек $s$ из $\L$, удовлетворяющих неравенству \rf{Xlip+}. Из \rf{cXlip+} мгновенно
вытекает, что множество $\L_t$
является одновременно открытым и замкнутым в $\L$. 
Кроме того, $\L_t\ne\varnothing$,
ибо $t\in\L$. Следовательно, $\L_t=\L$ в силу связности промежутка $\L$. $\bl$


\begin{thm} 
\label{lippoint+}
Пусть $f$ -- непрерывная функция на вещественной прямой $\R$.
Предположим, что
$$
\varlimsup_{t\to0}\frac{\|f(A+tK)-f(A)\|}{|t|}<+\be
$$
для любого (не обязательно ограниченного) самосопряжённого оператора $A$ и любого ограниченного
самосопряжённого оператора $K$. Тогда $f\in\OL(\R)$.
\end{thm}

\Pf Предположим, что $f\not\in\OL(\R)$. Тогда для любого неотрицательного числа $c$ существуют операторы $A,\:K\in\mB_{\rm sa}(\h)$ и число $a>0$ такие, что  $\|K\|=1$  и 
$\|H_{f,A,K}\|_{\Li(\R,\mB(\h))}>c$. Следовательно, в силу леммы \ref{verh+},
$$
\varlimsup\limits_{t\to0}\frac{\|H_{f,A,K}(t)-H_{f,A,K}(t_0)\|}{|t-t_0|}>c
$$
для некоторого числа $t_0$ из $[0,1]$.

Таким образом, при всех $n$ из $\Z_+$ найдутся операторы 
$A_n,\:K_n$ из $\mB_{\rm sa}(\h)$ такие, что $\|K_n\|=1$ и 
$$
\varlimsup_{t\to0}\frac{\|f(A_n+tK_n)-f(A_n)\|}{|t|}>n.
$$
Рассмотрим в гильбертовом пространстве $\ell^2(\h)$ самосопряжённые операторы $\A$ и ${\mathcal K}$, определённые равенствами
$$
\A(x_0,x_1,x_2,\cdots)=(A_0x_0,A_1x_1,A_2x_2,\cdots),\quad
(x_0,x_1,x_2,\cdots)\in\ell^2(\h)б
$$
и
$$
{\mathcal K}(x_0,x_1,x_2,\cdots)=(K_0x_0,K_1x_1,K_2x_2,\cdots),\quad
(x_0,x_1,x_2,\cdots)\in\ell^2(\h).
$$
Тогда 
$$
\varlimsup_{t\to0}\frac{\|f(\A+t\mathcal K)-f(\A)\|}{|t|}\ge
\varlimsup_{t\to0}\frac{\|f(A_n+tK_n)-f(A_n)\|}{|t|}>n
$$
для любого целого неотрицательного числа $n$, и мы приходим к противоречию. $\bl$


\medskip

{\bf Замечание.} Из доказательства теоремы \ref{lippoint+} можно увидеть, что имеют место 
следующие равенства
\begin{align*}
\|f\|_{\OL(\R)}&=\sup
\left\{\varlimsup_{t\to0}\frac{\|f(A+tK)-f(A)\|}{|t|}:~
A,~K\in\mB_{\rm sa}(\h),~\|K\|_{\mB_{\rm sa}(\h)}=1\right\}\\[.2cm]
&=\sup\big\{\|H_{f,A,K}\|_{\Li(\R)}:
~A,~K\in\mB_{\rm sa}(\h), \|K\|_{\mB_{\rm sa}(\h)}=1\big\}.
\end{align*}

\medskip

Прежде,
чем сформулировать следующую теорему, заметим, что функция $H_{f,A,K}$ дифференцируема
при всех самосопряжённых операторах $A$ и  $K$ тогда и только тогда, когда
она дифференцируема в нуле для любых самосопряжённых операторов $A$ и  $K$ (как обычно,
оператор $K$ предполагается ещё и ограниченным).

В доказательстве следующей теоремы мы используем один результат, который будет доказан в главе 2 (теорема \ref{sildif}). 

\begin{thm} 
\label{4i+}
Пусть $f$ -- непрерывная функция на $\R$. Тогда следующие утверждения 
равносильны:

{\rm (а)} $f\in\OL(\R)$;

{\rm (б)} $H_{f,A,K}\in\Li(\R,\mB(\h))$ для любого самосопряжённого оператора $A$ и любого ограниченного самосопряжённого оператора $K$;

{\rm (в)} функция $H_{f,A,K}$ дифференцируема, как функция из $\R$
в пространство $\mB(\h)$, наделённое слабой операторной топологией,
для любого самосопряжённого оператора $A$ и любого ограниченного самосопряжённого оператора $K$;

{\rm (г)}  функция $H_{f,A,K}$ дифференцируема, как функция из $\R$
в пространство $\mB(\h)$, наделённое сильной операторной топологией,
для любого самосопряжённого оператора $A$ и любого ограниченного самосопряжённого оператора $K$.
\end{thm}

\Pf Импликации (а)$\Longrightarrow$(б) и (г)$\Longrightarrow$(в) очевидны.
Импликация \lb(а)$\Longrightarrow$(г) вытекает из теоремы \ref{sildif} ниже.
Наконец, импликации (в)$\Longrightarrow$(а) и (б)$\Longrightarrow$(а) 
мгновенно следуют из теоремы \ref{lippoint+}. $\bl$

\medskip

Обозначим символом $\OL_{\rm loc}(\R)$ пространство непрерывных на $\R$ функций $f$
таких, что $f\big|[-a,a]\in\OL([-a,a])$ для всех $a>0$ и символом
$\Li_{\rm loc}(\R,\mB(\h))$ пространство непрерывных на $\R$ функций $f$
таких, что $f\big|[-a,a]\in\Li([-a,a],\mB(\h))$ для всех $a>0$.
Все результаты этого параграфа имеют естественный аналоги и для этих пространств.

\begin{thm} 
\label{loclippoint+}
Пусть $f$ -- непрерывная функция на вещественной прямой $\R$.
Предположим, что
$$
\varlimsup_{t\to0}\frac{\|f(A+tK)-f(A)\|}{|t|}<\be
$$
для любых $A,~K\in\mB_{\rm sa}(\h)$. Тогда    $f\in\OL_{\rm loc}(\R)$.
\end{thm}

\Pf Предположим, что $f\not\in\OL_{\rm loc}(\R)$. Тогда $f\not\in\OL([-a,a])$
для некоторого числа $a>0$. Следовательно, для любого числа $c\ge0$
существуют операторы $A,\:K$ из $\mB_{\rm sa}(\h)$ такие, что $\|A\|\le a$, $\|A+K\|\le a$
и $\|f(A+K)-f(A)\|>c\|K\|$. Повторяя аргументы доказательства теоремы \ref{lippoint+},
мы придём к противоречию, построив самосопряжённые операторы $\A$ и $\A+\mathcal K$  такие, что
$$
\|\A\|\le a,\quad\|\A+\mathcal K\|\le a\quad\mbox{и}\quad
\varlimsup_{t\to0}\frac{\|f(\A+t\mathcal K)-f(\A)\|}{|t|}=\be.\quad\bl
$$

\begin{thm} 
\label{loc4i+}
Пусть $f$ -- непрерывная функция на $\R$. Тогда следующие утверждения 
равносильны:

{\rm (а)} $f\in\OL_{\rm loc}(\R)$;

{\rm (б)} $H_{f,A,K}\in\Li_{\rm loc}(\R,\mB(\h))$ для любых $A,~K\in\mB_{\rm sa}(\h)$;

{\rm (в)}  для любых $A,\:K$ из $\mB_{\rm sa}(\h)$ функция $H_{f,A,K}$ дифференцируема как функция из $\R$
в пространство $\mB(\h)$,  наделённое слабой операторной топологией;

{\rm (г)} для любых $A,\:K$ из $\mB_{\rm sa}(\h)$ функция $H_{f,A,K}$ дифференцируема как функция из $\R$ в пространство $\mB(\h)$,  наделённое сильной операторной топологией.
\end{thm}

Эта теорема может быть доказана аналогично теореме \ref{4i+},
только вместо теоремы \ref{lippoint+} нужно использовать теорему \ref{loclippoint+}.

Отметим здесь, что в работе \cite{KS} показано, что условие (а) в теореме 
\ref{loc4i+} эквивалентно дифференцируемости по норме по всем компактным направлениям для ограниченных самосопряжённых операторов.

Из теоремы  \ref{4i+} следует, что если $f$ -- непрерывная функция на $\R$, то  
$f\in\OL(\R)$ в том и только в том случае, когда  для любого самосопряжённого 
оператора  $A$ и любого ограниченного самосопряжённого оператора $K$ существует
предел
\bay
\label{silndif+}
\lim_{t\to0}\frac1t(f(A+tK)-f(A))\df \bs{d}_{f,A}K
\ey
в сильной операторной топологии. Из теоремы \ref{sildif} ниже следует также, что $\bs{d}_{f,A}$ --
ограниченный линейный оператор, действующий из $\mB_{\rm sa}(\h)$ в $\mB(\h)$.

Аналогичные утверждения имеют место и для функций $f\in\OL_{\rm loc}(\R)$
с той лишь разницей, что теперь оператор $A$ должен быть ещё и ограниченным.

Из теоремы \ref{loc4i+} вытекает, что если $f$ -- непрерывная функция на $\R$, то 
$f\in\OL_{\rm loc}(\R)$ в том и только в том случае, когда  
для любых операторов
$A,\:K$ из $\mB_{\rm sa}(\h)$ предел \rf{silndif+} существует в сильной операторной топологии, при
этом $\bs{d}_{f,A}$ --
ограниченный линейный оператор из $\mB_{\rm sa}(\h)$ в $\mB(\h)$.

\begin{thm} 
Пусть $f\in\OL_{\rm loc}(\R)$. Тогда 
\begin{align*}
\|f\|_{\OL(\R)}&=\sup_{A\in\mB_{\rm sa}(\h)}\|\bs{d}_{f,A}\|\nonumber\\[.2cm]
&=\sup\big\{\|\bs{d}_{f,A}\|:~A~\mbox{ -- самосопряжённый оператор}\big\}.
\end{align*} 
\end{thm}

Как обычно, равенство $\|f\|_{\OL(\R)}=\be$ означает, что  
$f\not\in\OL(\R)$.

\medskip

\Pf  Достаточно воспользоваться леммой \ref{HfAK+}. $\bl$

\begin{thm} 
\label{operdif+}
Пусть $f$ -- непрерывная функция на $\R$.
Предположим, что для любого самосопряжённого 
оператора  $A$ и любого ограниченного самосопряжённого оператора $K$
предел  {\em\rf{silndif+}} существует по операторной норме.
Тогда \lb$f\in\OL(\R)\cap C^1(\R)$,
отображение $K\mapsto f(A+K)-f(A)$ ($K\in\mB_{\rm sa}(\h)$) дифференцируемо
по Фреше в точке $\0$ при каждом самосопряжённом операторе $A$,
и его дифференциал в  точке $\0$  равен $\bs{d}_{f,A}$.
\end{thm}

\Pf Включение $f\in\OL(\R)$ вытекает из теоремы \ref{4i+}.
Из теоремы \ref{sildif}
следует, что $\bs{d}_{f,A}$ -- ограниченный линейный оператор, действующий
из $\mB_{\rm sa}(\h)$ в $\mB(\h)$.
Проверим {\it дифференцируемость по Фреше}, то есть, что $\bs{d}_{f,A}$ -- линейный ограниченный оператор (уже доказано) и
$$
\lim_{t\to0}\frac1t\|f(A+tK)-f(A)-t\bs{d}_{f,A}K\|=0
$$
равномерно по всем $K$ из единичной сферы пространства
$\mB_{\rm sa}(\h)$.  

Для этого достаточно убедиться в том, что 
\bay
\label{lodF+}
\lim_{n\to\be}\frac1{t_n}\|f(A+t_nK_n)-f(A)-t_n\bs{d}_{f,A}K_n\|=0
\ey
для любой последовательности ненулевых вещественных чисел $\{t_n\}_{n\ge0}$, стремящейся к нулю, и
для любой последовательности 
самосопряжённых операторов $\{K_n\}_{n\ge0}$ такой, что
$\|K_n\|=1$ при всех $n$.



Рассмотрим в гильбертовом пространстве $\ell^2(\h)$ самосопряжённый оператор $\A$, определённый равенством
$$
 \A(x_0,x_1,x_2,\cdots)=(Ax_0,Ax_1,Ax_2,\cdots),\quad
 (x_0,x_1,x_2,\cdots)\in\ell^2(\h).
$$
Определим также ограниченный самосопряжённый оператор ${\mathcal K}$ равенством
$$
 {\mathcal K}(x_0,x_1,x_2,\cdots)=(K_0x_0,K_1x_1,K_2x_2,\cdots),\quad
 (x_0,x_1,x_2,\cdots)\in\ell^2(\h).
$$
Применяя условие теоремы к операторам $\A$ и $\mathcal K$,  получаем:
\bay
\label{lodG+}
\lim_{n\to\be}\frac1{t_n}\|f(\A+t_n{\mathcal K})-f(\A)-t_n\bs{d}_{f,\A}{\mathcal K}\|=0.
\ey
Ясно, что $\bs{d}_{f,\A}{\mathcal K}$ -- это ортогональная сумма операторов 
$\bs{d}_{f,A}K_n$, $n\ge0$, и, таким образом, \rf{lodF+} является следствием 
равенства \rf{lodG+}. 

Докажем, наконец, что $f\in C^1(\R)$.
Проверим непрерывность функции $f'$ в произвольной точке $t_0\in\R$.
Пусть $A$ -- оператор умножения на $x$  в пространстве $L^2([x_0-1,x_0+1])$.
Положим $K\df I$.
Тогда по условию теоремы существует предел по операторной норме 
$$
\lim_{t\to0} t^{-1}(f(A+tI)-f(A)).
$$
Следовательно, существует предел в $L^\be([x_0-1,x_0+1])$
$$
\lim_{t\to0} t^{-1}(f(x+t)-f(x))=f'(x),
$$
откуда $f\in C^1(t_0-1,t_0+1)$.
$\bl$

\medskip

{\bf Определение.}
Функция  $f$, удовлетворяющая условиям теоремы \ref{operdif+},
называется {\it  операторно дифференцируемой}. Обозначим символом $\OD(\R)$  множество всех операторно дифференцируемых функций
на $\R$.

\medskip 

Напомним, что для функций, определённых на банаховых пространствах, есть разные понятия дифференцируемости: существование слабой производной по Гато, существование дифференциала Гато, дифференцируемость по Фреше. Однако, как видно из теоремы \ref{operdif+}, в случае операторной дифференцируемости функций на прямой все эти определения эквивалентны. Отметим, что эквивалентность операторной дифференцируемости по Фреше и существования дифференциала Гато, являющегося линейным ограниченным оператором, доказана в работе \cite{KS}. 

Нетрудно видеть, что следующее утверждение можно доказать примерно так же, как теорему \ref{operdif+}.

\begin{thm} 
\label{locoperdif+}
Пусть $f$ -- непрерывная функция на $\R$.
Предположим, что для любых $A,~K\in\mB_{\rm sa}$ 
предел  {\em\rf{silndif+}} существует по операторной норме.
Тогда $f\in\OL_{\rm loc}(\R)\cap C^1(\R)$,
отображение $K\mapsto f(A+K)-f(A)$, $K\in\mB_{\rm sa}$, дифференцируемо
по Фреше в точке $\0$ при всех $A$ из $\mB_{\rm sa}$,
и его дифференциал в  точке $\0$  равен $\bs{d}_{f,A}$.
\end{thm}

Если функция  $f$ удовлетворяет условиям теоремы \ref{locoperdif+}, то говорят, что она {\it локально операторно дифференцируема}, и пишут 
$f\in{\rm OD}_{\rm loc}(\R)$. 

\medskip

Обратим внимание на то, что теоремы \ref{operdif+} и \ref{locoperdif+}, в частности, утверждают, что 
{\it если $f\in{\rm OD}_{\rm loc}(\R)$, то функция $f$ непрерывно дифференцируема и входит в класс $\OL_{\rm loc}(\R)$, а, если $f\in{\rm OD}(\R)$, то
$f\in\OL(\R)$}.

\medskip

{\bf Замечание.} Отметим, что не будучи непрерывно дифференцируемой, функция $g$, $g(x)=x^2\sin\frac1x$, не может быть операторно дифференцируемой. Таким образом, в примере 7 в \S\:\ref{prim} нельзя заменить класс операторно липшицевых функций классом операторно дифференцируемых функций.
Действительно, нетрудно проверить, что функция $x\mapsto\sin x=\im e^{{\rm i}x}$ является операторно дифференцируемой.

\medskip

Нашей ближайшей целью является доказательство непрерывной зависимости (по операторной норме) дифференциала $\bs{d}_{f,A}$ от оператора $A$ для (локально) операторно дифференцируемых функций $f$. Следующий результат был получен в работе \cite{KS}.

\begin{thm}
\label{ndlod+}
Пусть $f$ -- локально операторно дифференцируемая функция, и пусть $c>0$. Тогда
для любого положительного числа $\e$ существует положительное число $\d$ такое, что 
$$
\|\bs{d}_{f,A}-\bs{d}_{f,B}\|\le\e,
$$
как только
$A$ и $B$ - самосопряжённые операторы такие, что
$\|A\|\le c$, $\|B\|\le c$ и $\|A-B\|\le\d$.
\end{thm}

Докажем сначала следующую лемму, полученную в работе \cite{KS}.

\begin{lem}
\label{01+}
Пусть функция $f$ -- локально операторно дифференцируемая функция.  
Тогда для любых положительных чисел $c$ и $\e$ 
найдётся число $\d>0$
такое, что 
$$
\|f(A+K)-f(A)-\bs{d}_{f,A}K\|\le\e\|K\|,
$$
как только $A$ и $K$ - самосопряжённые операторы такие, что $\|K\|\le\d$ и $\|A\|\le c$.
\end{lem}

\Pf Предположим противное. Тогда для некоторых положительных чисел $c$ и $\e$ найдутся
последовательности самосопряжённых операторов $\{A_n\}_{n\ge0}$ и $\{K_n\}_{n\ge0}$
такие, что $\|K_n\|\to0$, $\|A_n\|\le c$ и 
\bay
\label{>e+}
\|f(A_n+K_n)-f(A_n)-\bs{d}_{f,A_n}K_n\|>\e\|K_n\|,\quad n\ge0.
\ey
Определим ограниченный самосопряжённый оператор $\A$ в $\ell^2(\h)$ равенством
\bay
\label{oprA+}
 \A(x_0,x_1,x_2,\cdots)=(A_0x_0,A_1x_1,A_2x_2,\cdots),\quad
 (x_0,x_1,x_2,\cdots)\in\ell^2(\h).
\ey
Тогда $\|A\|\le c$. Поскольку
функция $f$ дифференцируема по Фреше в точке $A$, существует 
число $\d>0$ такое, что
\bay
\label{A_Ke+}
\|f(A+K)-f(A)-\bs{d}_{f,A}K\|\le\e\|K\|
\ey
для любого самосопряжённого оператора $K$, удовлетворяющего условию $\|K\|\le\d$. 
Определим теперь оператор ${\mathcal K}_n$ в $\ell^2(\h)$ равенством
\bay
\label{oprKn+}
{\mathcal K}_n(x_0,x_1,x_2,\cdots)=(\0,\,\cdots,\0,K_nx_n,\0,\0,\cdots),\quad
 (x_0,x_1,x_2,\cdots)\in\ell^2(\h).
\ey

Применяя неравенство \rf{A_Ke+} при достаточно больших $n$, получаем
$$
\|f(A_n+K_n)-f(A_n)-\bs{d}_{f,A_n}K_n\|=
\|f(\A+{\mathcal K}_n)-f(\A)-\bs{d}_{f,\A}{\mathcal K}_n\|\le\e\|{\mathcal K}_n\|
=\e\|K_n\|,
$$
что противоречит неравенству \rf{>e+}. $\bl$

\medskip

{\bf Доказательство теоремы \ref{ndlod+}.} Пусть $c$, $\e$ и $\d$ обозначают то же, что в лемме \ref{01+}. Рассмотрим  самосопряжённые
операторы $A$ и $B$ такие, что $\|A\|\le c$, $\|B\|\le c/2$ и $\|B-A\|\le\min\{\d/2,c/2\}$. Пусть $K$ --  самосопряжённый оператор
такой, что $\|K\|=\d/2$. Тогда $\|B+K\|\le c$, $\|B-A\|\le\|K\|$ и $\|B-A+K\|\le2\|K\|$. Следовательно,
$$
\|f(B+K)-f(B)-\bs{d}_{f,B}K\|\le\e\|K\|,
$$
$$
\|f(B)-f(A)-\bs{d}_{f,A}(B-A)\|\le\e\|B-A\|\le\e\|K\|
$$
и
$$
\|f(B+K)-f(A)-\bs{d}_{f,A}(B-A+K)\|\le\e\|B-A+K\|\le2\e\|K\|.
$$
Используя равенство 
$$
\bs{d}_{f,A}(B-A+K)=\bs{d}_{f,A}(B-A)+\bs{d}_{f,A}K,
$$
получаем
\begin{align*}
\|\bs{d}_{f,B}K-\bs{d}_{f,A}K\|&\le\|\bs{d}_{f,B}K-f(B+K)+f(B)\|\\[.2cm]
&+\|f(B+K)-f(A)-\bs{d}_{f,A}(B-A+K)\|\\[.2cm]
&+\|\bs{d}_{f,A}(B-A)-f(B)+f(A)\|\le4\e\|K\|,
\end{align*}
откуда вытекает, что $\|\bs{d}_{f,B}-\bs{d}_{f,A}\|\le4\e$. $\bl$

\medskip

Перейдём теперь к случаю операторно дифференцируемых функций. Скажем, что {\it не обязательно ограниченные самосопряжённые операторы $A$ и $B$ эквивалентны}, если  существует оператор $K$ из $\mB_{\rm sa}(\h)$ такой, что $B=A+K$. Для операторов из одного класса можно ввести метрику $\dist(A,B)\df\|B-A\|$.


\begin{thm} 
\label{odineo+}
Пусть $f$ -- операторно дифференцируемая функция на $\R$. Тогда
отображение
$$
A\mapsto\bs{d}_{f,A}
$$
на каждом классе эквивалентности непрерывно в операторной норме.
\end{thm}

\begin{lem} 
\label{ravnots+}
Пусть $f$ -- операторно дифференцируемая функция на $\R$. Тогда для любого положительного числа $\e$  
найдётся положительное число $\d$
такое, что 
$$
\|f(A+K)-f(A)-\bs{d}_{f,A}K\|\le\e\|K\|
$$
для любого (не обязательно ограниченного) самосопряжённого оператора $A$ и для любого самосопряжённого оператора $K$, норма которого не больше, чем $\d$.
\end{lem}

\Pf Предположим противное. Тогда для некоторого числа $\e>0$ найдутся
две последовательности самосопряжённых операторов $\{A_n\}_{n=1}^\be$ и $\{K_n\}_{n=1}^\be$
такие, что $\|K_n\|\to0$ и 
\bay
\label{verotsep+}
\|f(A_n+K_n)-f(A_n)-\bs{d}_{f,A_n}K_n\|>\e\|K_n\|
\ey
при всех $n\ge1$. Пусть $\A$ и ${\mathcal K}_n$ -- операторы в пространстве $\ell^2(\h)$,
определённые формулами \rf{oprA+} и \rf{oprKn+}. Поскольку
функция $f$ дифференцируема по Фреше в точке $\A$, существует положительное
число $\d$ такое, что
$$
\|f(\A+K)-f(\A)-\bs{d}_{f,\A}K\|\le\e\|K\|
$$
для всех самосопряжённых операторов $K$ нормы не выше, чем $\d$. 
Применяя это неравенство к оператору ${\mathcal K}_n$
при достаточно больших $n$, получаем:
$$
\|f(A_n+K_n)-f(A_n)-\bs{d}_{f,A_n}K_n\|=
\|f(\A+{\mathcal K}_n)-f(\A)-\bs{d}_{f,\A}{\mathcal K}_n\|\le\e\|{\mathcal K}_n\|=\e\|K_n\|,
$$
что противоречит неравенству \rf{verotsep+}. $\bl$

\medskip

{\bf Доказательство теоремы \ref{odineo+}.}
Пусть $\e$ и $\d$ обозначают то же, что в лемме \ref{ravnots+}. Рассмотрим самосопряжённые
операторы $A$ и $B$ такие, что $\|B-A\|\le\d/2$. Пусть $K$ --  самосопряжённый оператор
такой, что $\|K\|=\d/2$. Тогда  $\|B-A\|\le\|K\|$ и $\|B-A+K\|\le2\|K\|$. Следовательно,
$$
\|f(B+K)-f(B)-\bs{d}_{f,B}K\|\le\e\|K\|,
$$
$$
\|f(B)-f(A)-\bs{d}_{f,A}(B-A)\|\le\e\|B-A\|\le\e\|K\|
$$
и
$$
\|f(B+K)-f(A)-\bs{d}_{f,A}(B-A+K)\|\le\e\|B-A+K\|\le2\e\|K\|.
$$
Используя равенство 
$$
\bs{d}_{f,A}(B-A+K)=\bs{d}_{f,A}(B-A)+\bs{d}_{f,A}K,
$$
получаем:
\begin{align*}
\|(d_Bf)(K)-(d_Af)(K)\|&\le\|\bs{d}_{f,B}K-f(B+K)+f(B)\|\\[.2cm]
&+\|f(B+K)-f(A)-\bs{d}_{f,A}(B-A+K)\|\\[.2cm]
&+\|\bs{d}_{f,A}(B-A)-f(B)+f(A)\|\le4\e\|K\|
\end{align*}
для всех самосопряжённых операторов $K$ таких, что $\|K\|=\d/2$,
откуда вытекает, что $\|\bs{d}_{f,B}-\bs{d}_{f,A}\|\le4\e$
 всякий раз, когда $\|B-A\|\le\d/2$. $\bl$

\begin{thm}
Пусть  $f\in\OL_{\rm loc}(\R)$. Тогда функция $f$ локально операторно дифференцируема в том и только в том случае,
если отображение $A\mapsto \bs{d}_{f,A}$ действует непрерывно,
как отображение 
из банахова пространства $\mB_{\rm sa}(\h)$
в банахово пространство ограниченных операторов из $\mB_{\rm sa}(\h)$ в 
$\mB(\h)$.
\end{thm}

\Pf Из теоремы \ref{ndlod+} следует, что достаточно проверить только, что непрерывность
отображения $A\mapsto \bs{d}_{f,A}$ влечёт операторную дифференцируемость.
Заметим, что $H_{f,A,K}'(s)= \bs{d}_{f,A+sK}K$ (производная берётся в сильной
операторной топологии). Следовательно,
\bay
\label{prirost+}
f(A+K)-f(A)=\int_0^1(\bs{d}_{f,A+sK}K)\,ds,
\ey
где интеграл понимается в следующем смысле:
$$
(f(A+K)-f(A))u=\int_0^1\big((\bs{d}_{f,A+sK}K)u\big)\,ds
$$
для любого $u\in\h$. Применяя тождество \rf{prirost+} к оператору $tK$  вместо $K$, получаем:
$$
t^{-1}(f(A+K)-f(A))-\bs{d}_{f,A}K=\int_0^1\big((\bs{d}_{f,A+stK}-\bs{d}_{f,A})K\big)\,ds.
$$
Будем считать, что $\|K\|=1$. Тогда из последнего тождества следует, что 
$$
\|t^{-1}(f(A+K)-f(A))-\bs{d}_{f,A}K\|\le\int_0^1\|\bs{d}_{f,A+stK}-\bs{d}_{f,A}\|\,ds
$$
Остаётся заметить, что
$$
\lim_{t\to0}\int_0^1\|\bs{d}_{f,A+stK}-\bs{d}_{f,A}\|\,ds=0
$$
равномерно по всем самосопряжённым операторам $K$ единичной нормы
в силу непрерывности отображения $A\mapsto \bs{d}_{f,A}$ в точке $A$. $\bl$

\medskip

Аналогичным образом можно доказать следующее утверждение.

\begin{thm}
\label{neprclass+}
Пусть  $f\in\OL_{\rm loc}(\R)$. Тогда функция $f$ операторно дифференцируема в том
и только в том случае, если отображение $A\mapsto \bs{d}_{f,A}$ непрерывно в операторной норме на каждом классе эквивалентности.
\end{thm}

\begin{thm} Множество ${\rm OD}(\R)$ является замкнутым подпространством
пространства  $\OL(\R)$. 
\end{thm}

 \Pf Нам нужно доказать, что если $\lim\limits_{n\to\be}f_n=f$  в пространстве  $\OL(\R)$
и $f_n\in{\rm OD}(\R)$  при всех $n$, то $f\in{\rm OD}(\R)$.
Из теорем  \ref{sildif} и \ref{muld} следует, что
$$
\|\bs{d}_{f_n,A}-\bs{d}_{f,A}\|=\|\bs{d}_{f_n-f,A}\|
\le\|\dg(f_n-f)\|_{\fM(\R\times\R)}=\|f_n-f\|_{\OL(\R)}\to0
$$
при $n\to\be$. Таким образом, $\lim\limits_{n\to\be}\bs{d}_{f_n,A}=\bs{d}_{f,A}$ по норме равномерно по всем самосопряжённым операторам $A$.
Теперь можно сослаться на теорему \rf{neprclass+}, поскольку непрерывность сохраняется
при равномерной сходимости. $\bl$

Отметим здесь, что в случае функций на конечных интервалах замкнутость множества операторно дифференцируемых функций в пространстве операторно липшицевых функций была установлена в \cite{KS}.


\

\section{\bf Коммутаторная липшицевость}
\setcounter{equation}{0}
\label{CoL}

\

Напомним, что непрерывная функция $f$ на $\R$ называется {\it коммутаторно липшицевой}, если 
\bay
\label{fAR-RfA}
\|f(A)R-Rf(A)\|\le\const\|AR-RA\|
\ey
для любого ограниченного самосопряжённого оператора $A$ и для любого ограниченного оператора $R$. Как и в определении операторно липшицевых функций, если функция $f$ коммутаторно липшицева, то неравенство \rf{fAR-RfA} справедливо для любого (не обязательно ограниченного) оператора $A$ и для любого ограниченного оператора $R$, см. теорему \ref{anbnrn}.

Нетрудно установить, что если функция $f$ коммутаторно липшицева, то в неравенстве 
\rf{fAR-RfA} можно заменить коммутаторы на квазикоммутаторы, т.е. имеет место более общее неравенство
$$
\|f(A)R-Rf(B)\|\le\const\|AR-RB\|
$$
для любого ограниченного оператора $R$ и любых самосопряжённых операторов $A$ и $B$.
Более того, оба эти условия эквивалентны тому, что функция $f$ операторно липшицева.
Это вытекает из следующего результата.

\begin{thm}
\label{komlipf}
Пусть $f$ -- непрерывная функция на $\R$. Тогда следующие условия эквивалентны:

{\em(а)} $\|f(A)-f(B)\|\le\|A-B\|$ для любых самосопряжённых операторов $A$ и $B$;

{\em(б)} $\|f(A)R-Rf(A)\|\le\|AR-RA\|$ для любого самосопряжённого оператора $A$ и любого ограниченного линейного оператора $R$;

{\em(в)} $\|f(A)R-Rf(B)\|\le\|AR-RB\|$ для произвольных самосопряжённых операторов $A$ и $B$ и для любого ограниченного линейного оператора $R$.
\end{thm}

Мы выведем теорему \ref{komlipf} из более общего результата для функций нормальных операторов в \S\:\ref{oplip}. Заметим, однако, что в случае функций нормальных операторов коммутаторная липшицевость отнюдь не эквивалентна операторной липшицевости.

\

\section{\bf Операторные неравенства Бернштейна}
\setcounter{equation}{0}
\label{Bern}

\

В этом параграфе мы приводим элементарное доказательство результата работы
\cite{Pe3}, который заключается в том, что функции из $L^\be(\R)$, преобразование Фурье которых имеет компактный носитель, являются операторно липшицевыми. Более того, мы получим так-называемое операторное неравенство Бернштейна с константой 1. При этом мы следуем подходу, приведённому в работе \cite{AP4}. Мы также получим аналогичные результаты для функций на окружности.

Отметим, что в \S\:\ref{Dost} мы увидим из этих результатов, что достаточным условием для операторной липшицевости является принадлежность классу Бесова $B_{\be,1}^1(\R)$.

Пусть $\s>0$. Говорят, что целая функция $f$ является функцией {\it экспоненциального
типа не выше} $\s$, если для любого положительного числа $\e$
найдётся число $c>0$ такое, что $|f(z)|\le c e^{(\s+\e)|z|}$ при всех $z\in\C$.

Обозначим через $\E_\s$ множество всех
целых функций экспоненциального типа не выше $\s$. Хорошо известно, что
$\E_\s\cap L^\be(\R)=\{f\in L^\be(\R):\supp\F f\subset[-\s,\s]\}$.

Отметим ещё, пространство $\E_\s\cap L^\be(\R)$
совпадает с множеством всех целых функций $f$ таких, что  $f\in L^\be(\R)$ и 
\bay
\label{fEs}
|f(z)|\le e^{\s|\im z|}\|f\|_{L^\infty(\R)},\quad z\in\C,
\ey
см, например, \cite{L}, стр. 97.

Неравенство Бернштейна (см. \cite{Be}) утверждает, что
$$
\sup_{x\in \R}|f^\prime(x)|\le\s\sup_{x\in \R}|f(x)|
$$
для любой функции $f$ из $\mathscr E_\s\cap L^\be(\R)$. Отсюда вытекает, что
\bay
\label{corber}
|f(x)-f(y)|\le\s\|f\|_{L^\infty(\R)}|x-y|,\quad f\in\mathscr E_\s\cap L^\be(\R),\quad x,~y\in\R,
\ey
где $\|f\|_{L^\infty(\R)}\df\sup\limits_{x\in \R}|f(x)|$.

Бернштейн также доказал в \cite{Be} следующее усиление неравенства \rf{corber}:
\bay
\label{bern2}
|f(x)-f(y)|\le\b(\s(|x-y|))\|f\|_{L^\infty(\R)},\quad f\in\mathscr E_\s\cap L^\be(\R),\quad x,~y\in\R,
\ey
где
$$
\b(t)\df\left\{\begin{array}{ll}2\sin(t/2),&\text {если}\,\,\,\,0\le t\le\pi,\\[.1cm]
2,&\text {если}\,\,\,\,t>\pi.
\end{array}\right.
$$
Заметим, что $\b(t)\le\min(t,2)$  при всех $t\ge0$.

Ясно, что неравенства \rf{bern2} превращаются в равенство для функции $f(z)=e^{{\rm i}\s z}$.

%
%
%

Пусть $X$ -- комплексное банахово пространство. Обозначим через $\E_\s(X) $  пространство
всех целых $X$-значных функций $f$ экспоненциального
типа не выше $\s$, т. е. удовлетворяющих следующему условию:
для любого положительного числа $\e$
найдётся число $c>0$ такое, что $\|f(z)\|_X\le c e^{(\s+\e)|z|}$ при всех $z\in\C$.

\medskip

{\bf Неравенство Бернштейна для векторнозначных функций.} {\em Пусть $f$ -- функция из 
$\E_\s(X)\cap L^\be(\R,X)$, где $\s>0$. 
Тогда 
\bay
\label{bn}
\|f(x)-f(y)\|_X\le\b(\s(|x-y|))\|f\|_{L^\infty(\R,X)}\le \s\|f\|_{L^\be(\R,X)}|x-y|
\ey
для всех $x,y\in\R$.}

 \medskip
 
 Векторная версия неравенства Бернштейна сводится к скалярной версии при помощи
теоремы Хана--Банаха.

\medskip

\begin{thm}
\label{opbern}
Пусть $f\in\mathscr E_\s\cap L^\be(\R)$.  Тогда
\bay
\label{opnerBer}
\|f(A)-f(B)\|\le\b(\s(\|A-B\|))\|f\|_{L^\infty}\le\s\|f\|_{L^\infty}\|A-B\|
\ey
для любых (ограниченных) самосопряжённых операторов $A$ и $B$. В частности,
$\|f\|_{\OL(\R)}\le\s\|f\|_{L^\be(\R)}$.
\end{thm}

{\bf Доказательство теоремы \ref{opbern}.} Пусть $A$ и $B$ -- самосопряжённые
операторы в гильбертовом пространстве $\h$. Нам нужно доказать, что
$$
\|f(A)-f(B)\|\le\b_\s(\|A-B\|)\|f\|_{L^\infty}.
$$
Положим $F(z)=f(A+z(B-A))$. Ясно, что $F$ -- целая функция со значениями в пространстве операторов
$\mB(\h)$ и $\|F(t)\|\le\|f\|_{L^\be(\R)}$ для всех $t\in\R$. Из неравенства фон Неймана (см. \cite{SNF}) следует, что 
$F\in\E_{\s\|B-A\|}(\mB(\h))$. Чтобы закончить доказательство, остаётся применить
неравенство Бернштейна \rf{bn}  к векторнозначной функции  $F$ при $x=0$ и  $y=1$.  $\bl$

\medskip

Ранее в работе \cite{Pe3} было доказано, что 
\bay
\label{PeR}
\|f\|_{\OL(\R)}\le\const\s\|f\|_{L^\be(\R)},\quad f\in\mathscr E_\s\cap L^\be(\R)
\ey
и, в частности, $\E_\s\cap L^\be(\R)\subset\OL(\R)$. Отсюда следует, что
для любой функции $f\in\E_\s\cap \Li(\R)$  функция $f'$ операторно липшицева.

Следующий пример показывает, что $\E_\s\cap \Li(\R)\not\subset\OL(\R)$.

\medskip

{\bf Пример.} {\it Рассмотрим функцию $f(x)\df\int_0^x{\rm Si}(t)\,dt$, где $\rm Si$  
обозначает интегральный синус,
$$
{\rm Si}(x)\df\int_0^x\frac{\sin t}t dt.
$$
Ясно, что $f\in\E_1\cap \Li(\R)$,  но функция $f$ не может быть
операторно липшицевой, см. теоремы {\rm\ref{olm}} и {\rm\ref{pro}} ниже, поскольку не существует предела 
$\lim\limits_{|x|\to\be}x^{-1}f(x)$ (на самом деле 
$\lim\limits_{x\to\be}x^{-1}f(x)=\lim\limits_{x\to\be}{\rm Si}(x)=\frac\pi2=-\lim\limits_{x\to-\be}x^{-1}f(x)$)}.

\medskip

Интересно отметить, что если функцию $f$ из этого примера чуть ``испортить'',
заменив её функцией $g(x)\df\int_0^x{\rm Si}(|t|)\,dt$, то она станет операторно лиишицевой.
Достаточно убедиться в том, что функция $g(x)-\frac\pi2x$ является операторно липшицевой.
Это следует из того (см. предложение 7.8 статьи \cite{BS3}),
что производная этой функции принадлежит пространству $L^2(\R)\cap\Li(\R)$
(это также можно вывести из теоремы \ref{Besdiffer} ниже).

\medskip

Установим теперь операторные аналоги неравенств Бернштейна для унитарных 
операторов.

\begin{lem} 
\label{IUV}
Пусть $U$ и $V$ -- унитарные операторы. Тогда существует самосопряжённый 
оператор $A$ такой, что $V=e^{{\rm i}A}U$, $\|A\|\le\pi$ и 
$\b(\|A\|)=\|U-V\|$.
\end{lem}

\Pf Определим оператор $A$ равенством $A=\arg(VU^{-1})$, где функция $\arg$
определена на окружности $\T$  равенством $\arg(e^{{\rm i}s}) =s, s\in[-\pi,\pi)$.
Очевидно, что 
$\b(\|A\|)=\|I-e^{{\rm i}A}\|=\|U-V\|$. $\bl$

\begin{thm}  
\label{uniber}
Пусть $f$ -- тригонометрический полином степени не более, чем $n$.
Тогда
$$
\|f(U)-f(V)\|\le n\|f\|_{L^\be(\T)}\|U-V\|
$$
для любых унитарных операторов $U$ и $V$.
\end{thm}

\Pf Пусть $A$ -- самосопряжённый оператор, существование которого утверждается
в лемме \ref{IUV}. Положим $\Phi(z)\df f(e^{{\rm i}zA}U)$, $z\in\C$, где тем же символом $f$ мы обозначаем
аналитическое продолжение в $\C\setminus\{0\}$ тригонометрического полинома $f$.
Ясно, что $\Phi$ -- целая функция со значениями в пространстве операторов
$\mB(\h)$ и $\|\Phi(t)\|\le\|f\|_{L^\be(\T)}$ для всех $t\in\R$. Из неравенства фон Неймана (см. \cite{SNF}) следует, что 
$\Phi\in\E_{n\|A\|}(\mB(\h))$.   Применяя теперь векторное неравенство Бернштейна,
получаем:
$$
\|f(U)-f(V)\|=\|\Phi(1)-\Phi(0)\|\le\b(n\|A\|)\|f\|_{L^\be(\T)}.
$$
Остаётся заметить, что
$$
\b(n\|A\|)\le n \b(\|A\|)=n\|U-V\|.\quad\bl
$$

\medskip

Отметим, что в работе \cite{Pe1} было доказано, что
$$
\|f(U)-f(V)\|\le\const n\|f\|_{L^\be(\T)}\|U-V\|
$$
для любого тригонометрического полинома $f$ степени $n$ и для любых унитарных операторов $U$ и $V$.


\medskip

{\bf Замечание.} 
Из доказательства теоремы \ref{uniber} видно, что 
$$
\|f(U)-f(V)\|\le\b(n\|A\|)\|f\|_{L^\be(\T)}=\b\left(2n\arcsin\frac{\|U-V\|}2\right)\|f\|_{L^\be(\T)}.
$$
Эта оценка неулучшаема для функции $f(z)=z^n$, поскольку
$$
\sup\big\{|z_1^n-z_2^n|:~z_1\in\T,~z_2\in\T,~
|z_1-z_2|<n\big\}=\b\big(2n\arcsin(\d/2)\big),\quad\d\in(0,2].
$$

\

\section{\bf Необходимые условия для операторной липшицевости}
\setcounter{equation}{0}
\label{Neob}

\

В этом параграфе для функций на прямой и на окружности мы получим необходимые условия для операторной липшицевости. Эти необходимые условия по существу были получены в работах \cite{Pe1} и \cite{Pe3}, в которых использовались другие методы. Здесь для достижения цели помимо критерия ядерности операторов Ганкеля
(см. \S\:\ref{Prel}, раздел 3), на который также опираются работы \cite{Pe1} и \cite{Pe3}, используются и результаты параграфа \ref{olriolt+} этого обзора о поведении производных операторно липшицевых функций при дробно линейных преобразованиях.

Для доказательства следующего результата мы воспользуемся результатами параграфа 
\ref{yadiyadcom} о поведении функций от операторов при ядерных возмущениях.

\begin{thm}
\label{B11}
Пусть $f$ -- операторно липшицева функция на окружности $\T$. 
Тогда $f\in B_1^1(\T)$.
\end{thm}

\Pf Ввиду замечания к теореме \ref{tsentrez}, функция $f$ обладает свойством:
$$
f(U)-f(V)\in\bS_1,
$$
как только $U$ и $V$ -- унитарные операторы такие, что $U-V\in\bS_1$.

Определим операторы $U$ и $V$ в пространстве $L^2(\T)$ равенством
$$
Uf=\bar z f\quad\mbox{and}\quad Vf=\bar zf-2(f,\1)\bar z,\quad f\in L^2.
$$
Легко видеть, что $U$ и $V$ - унитарные операторы и
$$
\rank(V-U)=1.
$$
Также нетрудно проверить, что при $n\ge0$ имеет место равенство:
$$
V^nz^j=\left\{
\begin{array}{ll}z^{j-n},&j\ge n,\\[.1cm]
-z^{j-n},&0\le j<n,\\[.1cm]
z^{j-n},&j<0.
\end{array}\right.
$$
Отсюда получаем, что для любой непрерывной функции $f$ на $\T$
справедливо соотношение:
\begin{align*}
\big((f(V)-f(U))z^j,z^k\big)&=\sum_{n>0}\hat f(n)\big((V^nz^j,z^k)-(z^{j-n},z^k)\big)\\
&+\sum_{n<0}\hat f(n)\big((V^nz^j,z^k)-(z^{j-n},z^k)\big)\\
&=-2\left\{\begin{array}{ll}\hat f(j-k),&j\ge0,~k<0,\\[.1cm]
\hat f(j-k), &j<0,~k\ge0,\\[.1cm]
0,&\mbox{otherwise}.
\end{array}\right.
\end{align*}

Таким образом, если $f(U)-f(V)\in\bS_1$, то операторы в $\ell^2$
с ганкелевыми матрицами
$$
\{\hat f(j+k)\}_{j\ge0,k\ge1}\qm\{\hat f(-j-k)\}_{j\ge0,k\ge1}
$$
входят в класс $\bS_1$. Теперь можно воспользоваться критерием ядерности операторов Ганкеля (см. \S\:\ref{Prel}, раздел 3) и заключить, что $f\in B_{1}^{1}(\T)$. $\bl$

\medskip

Отметим здесь, что конструкция в доказательстве теоремы \ref{B11} позаимствована из работы \cite{AP3}.

Нам удобно в этом параграфе ввести обозначение $\zh M\zh$ для нормы матрицы $M$.

Чтобы сформулировать следствие теоремы \ref{B11}, нам понадобится банахово пространство
$(\OL)_{\rm loc}'(\T)$, которое будет подробно рассмотрено в \S\:\ref{olriolt+}.
Здесь мы только отметим, что $(\OL)_{\rm loc}'(\T)=\{f'+c\ov z: f\in\OL(\T),\:c\in\C\}$ (см. следствие \ref{rt4}), причём,
как всегда в этой статье, производная понимается в комплексном смысле, т. е.
$f'(\z)\df\lim_{\t\to\z}(\t-\z)^{-1}(f(\t)-f(\z))$.

\begin{cor} 
\label{sumgrad}
Пусть $u$ -- интеграл Пуассона функции $f\in(\OL)_{\rm loc}'(\T)$. Тогда 
$\zh\n u\,\zh\in L^1(\dd)$ и 
$\big\|\zh\n u\,\zh\big\|_{L^{1}(\dd)}\le\const\|f\|_{(\OL)_{\rm loc}'(\T)}$. 
\end{cor}



\Pf  Пусть  $f=g'$, где $g\in\OL(\T)$. Тогда $f\in B_1^0(\T)$, поскольку $g\in B_1^1(\T)$, и достаточно
воспользоваться характеризацией класса Бесова $B_1^0(\T)$ в терминах 
гармонического продолжения,  см \S \ref{Prel}. Остаётся заметить, что следствие очевидно
для функции $f(z)=z^{-1}=\ov z$. $\bl$


\medskip

Чтобы сформулировать более сильное необходимое условие для операторной липшицевости,
нам понадобится понятие меры Карлесона. Пусть $\mu$ -- положительная борелевская мера 
в открытом единичном круге $\dd$. Хорошо известная теорема Карлесона утверждает,
что каждая функция класса Харди $H^p$ входит в пространство $L^p(\mu)$ ($0<p<+\be$)
в том и только в том случае, когда
для любой
точки $\z$ единичной окружности $\T$  и любого положительного числа $r$ 
выполняется неравенство
$$
\mu\{z\in\dd:\,|z-\z|<r\}\le\const r.
$$
Такие меры $\mu$ называются {\it мерами Карлесона} в круге $\dd$.
Заметим, что условие Карлесона не зависит от $p\in(0,+\be)$.
Подробнее о мерах Карлесона см., например, \cite{Ni1}, \cite{N}. Нам понадобится следующая
эквивалентная переформулировка условия Карлесона:
\bay
\label{Vin}
\sup_{a\in\dd}\int_\dd\frac{1-|a|^2}{|1-\ov a z|^2}\,d\mu(z)<+\be,
\ey
см., например,  \cite{Ni1}, лекция VII. 
Заметим, что условие \rf{Vin} означает, что $\|k_a\|_{L^2(\mu)}\le\const\|k_a\|_{H^2}$
для всех $a\in\dd$, где $k_a(z)\df(1-z\ov a)^{-1}$ -- воспроизводящее ядро гильбертова пространства $H^2$.

{\it Обозначим символом ${\rm CM}(\dd)$ пространство всех комплексных мер Радона $\mu$
в $\dd$ таких, что $|\mu|$  -- мера Карлесона, и символом $\|\mu\|_{{\rm CM}(\dd)}$  норму 
тождественного оператора вложения из $H^1$ в  $L^1(|\mu|)$}. Хорошо известно, что
(квази)норма тождественного оператора вложения из $H^p$ в $L^p(|\mu|)$ равна
$\|\mu\|_{{\rm CM}(\dd)}^{1/p}$  при всех $p\in(0,+\infty)$.

Всё вышесказанное о мерах Карлесона в круге $\dd$  имеет естественные аналоги для
полуплоскости $\C_+$.
В этом случае условие Карлесона для положительной борелевской меры $\mu$ в $\C_+$
переписывается следующими образом:
$$
\mu\{z\in\C_+:\,|z-t|<r\}\le\const r
$$
для всех $t\in\R$ и всех $r>0$. Аналогом условия \rf{Vin} является неравенство
$$
\sup_{a\in\C_+}\int_{\C_+}\frac{\im a}{|z-\ov a |^2}\,d\mu(z)<\be.
$$
В частности, точно также мы можем ввести пространство ${\rm CM}(\C_+)$
и норму в нём.

Пусть $f$ --  (обобщённая) функция на единичной окружности $\T$.
{\it Обозначим символом $\mP f$  интеграл Пуассона функции  $f$}.


\begin{thm} 
Пусть  $f\in(\OL)'_{\rm loc}(\T)$.
Тогда $\zh\n (\mP f)\,\zh\,d\m_2\in\CM(\dd)$.
\end{thm}

\Pf Пусть $f\in(\OL)_{\rm loc}'(\T)$. Тогда из теоремы \ref{kp} и следствия \ref{sumgrad}
вытекает, что
\bay
\label{afia}
\int_\dd|\!\!|\!|((\nabla u)\circ\f)(z)|\!\!|\!|\cdot|\f'(z)|\,d\m_2
\le\const\|f\|_{(\OL)_{\rm loc}'(\T)}
\ey
для любого дробно-линейного автоморфизма единичного круга $\dd$,
где $u=\mP f$. Возьмём в качестве $\f$
функцию $\f(z)=(1-\ov a z)^{-1}(a-z)$, где $a\in\dd$.  Делая в 
интеграле \rf{afia} замену переменной $z=\f(w)$, получаем
$$
\sup_{a\in\dd}
\int_\dd|\!\!|\!|(\nabla u)(w)|\!\!|\!|\frac{1-|a|^2}{|1-\ov a w|^2}\,d\m_2(w)\le\const\|f\|_{(\OL)'(\T)}.
$$
Таким образом, мера $|\!\!|\!|\nabla u\,|\!\!|\!|\,d\m_2=\zh\n (\mP f)\,\zh\,d\m_2$ удовлетворяет
условию \rf{Vin}. $\bl$

\begin{thm} 
Пусть  $f\in\OL(\T)$.
Тогда\footnote{Здесь и далее ${\rm Hess}$ обозначает гессиан, т. е.
матрицу из частных производных второго порядка.} 
$\zh{\rm Hess}\, \mP f\,\zh\,d\m_2\in\CM(\dd)$.
\end{thm}


Перейдём теперь к интегралам Пуассона функций на $\R$.
Если функция $f$ принадлежит пространству $L^1(\R,(1+x^2)^{-1}\,dx)$, то интеграл
Пуассона определяется стандартным образом.   
Нам понадобится интеграл Пуассона фунцции $f$ такой, что $f'\in  L^1(\R,(1+x^2)^{-1}\,dx)$.
Ясно, что достаточно рассмотреть случай вещественной функции $f$. Пусть $u$ -- интеграл 
Пуассона функции $f'$. Обозначим буквой $v$ гармоническую функцию, сопряжённую с функцией $u$.
Функция $u+{\rm i}v$  имеет первообразную $F$ такую, что граничные значения функции $\re F$
совпадают с функцией $f$ всюду на $\R$. Функция $F$ не определяется однозначно, поскольку
гармонически сопряжённая функция $v$ не определяется однозначно. Семейство $\{v+c\}_{c\in\R}$
состоит из всех гармонических функций, сопряжённых с функцией $u$. Нужная нам первообразная
функции $u+{\rm i}(v+c)$ имеет вид $F+c{\rm i}z+{\rm i}\a$, где $\a\in\R$.
Заметим, что $\re(F+c{\rm i}z+{\rm i}\a)=\re F-cy$. Таким образом, естественно определить
интеграл Пуассона функции $f$ как класс функций $\{\re F-cy\}_{c\in\R}$. Поскольку ${\rm Hess}\,y=0$,
гессиан интеграла Пуассона ${\rm Hess}\,\mP f$ функции $f$ определяется однозначно.

\begin{thm}
Пусть  $f\in(\OL)'(\R)$.
Тогда $|\!\!|\!|\n \mP f\,|\!\!|\!|\,d\m_2\in\CM(\C_+)$.
\end{thm}

\Pf 
Пусть $f\in(\OL)'(\R)$. Тогда из теоремы \ref{pere} и следствия \ref{sumgrad}
вытекает, что
\bay
\label{afia1}
\int_\dd|\!\!|\!|((\nabla u)\circ\f)(z)|\!\!|\!|\cdot|\f'(z)|\,d\m_2\le\const\|f\|_{(\OL)'(\R)},
\ey
для любого автоморфизма $\f$ из $\dlC$ такого, что $\f(\dd)=\C_+$,
где $u=\mP f$.
Возьмём в качестве $\f$
функцию $\f(z)=(1-z)^{-1}(a-\ov a z)$, где $a\in\C_+$.  Делая в 
интеграле \rf{afia} замену переменной $z=(w-\ov a)^{-1}(w-a)$, получаем
$$
\sup_{a\in\C_+}
\int_{\C_+}|\!\!|\!|(\nabla u)(w)|\!\!|\!|\frac{2\im a}{|w-\ov a|^2}\,d\m_2(w)\le\const\|f\|_{(\OL)_{\rm loc}'(\T)}.
$$
Последнее условие равносильно карлесоновости меры
$|\!\!|\!|\nabla u\,|\!\!|\!|\,d\m_2$. $\bl$

\begin{thm} 
Пусть  $f\in\OL(\R)$.
Тогда $|\!\!|\!|{\rm Hess}\, \mP f\,|\!\!|\!|\,d\m_2\in\CM(\C_+)$.
\end{thm}

Необходимые условия операторной липшицевости, приведённые выше, были получены первоначально в работах \cite{Pe1} и \cite{Pe3}. Именно, в работе \cite{Pe1} было доказано, что если $f\in\OL(\T)$, то оба оператора Ганкеля $H_f$ и $H_{\ov f}$ отображают пространство Харди $H^1$ в класс Бесова $B_1^1(\T)$
(класс таких функций $f$ обозначен в  \cite{Pe1} символом $L$). С. Семмес заметил, что $f\in L$ тогда и только тогда, когда мера 
$\zh{\rm Hess}\,\mP f\zh\,d\m_2$ является карлесоновой; см. \cite{Pe4}, где изложено доказательство этой эквивалентности. Аналогичное утверждение имеет место и для функций на $\R$, см.
\cite{Pe3}. В работе \cite{Pe1} также показано, что необходимое условие для операторной липшицевости, обсуждаемое выше, не является достаточным. Более того, оно не является достаточным даже для липшицевости.

Рассмотрим теперь пространства $\pp_+(b_\infty^{-1}(\T))$ и 
$\pp_+(b^{-1}_{1,\be}(\T))$, являющиеся замыканиями множества аналитических полиномов в пространствах Бесова $B_\infty^{-1}(\T)$ и $B^{-1}_{1,\be}(\T)$.
Хорошо известно, что эти пространства допускают следующие описания в терминах аналитического продолжения в единичный круг:
$$
\pp_+(b_\infty^{-1}(\T))=
\big\{h:~\lim_{r\to1^-}(1-r)\|h(rz)\|_{L^\be(\T)}=0\big\};
$$
$$
\pp_+(b^{-1}_{1,\be}(\T))=
\big\{h:~\lim_{r\to1^-}(1-r)\|h(rz)\|_{L^1(\T)}=0\big\}.
$$
В работе \cite{Pe1} также отмечено, что пространство $\pp_+L$ является двойственным к пространству аналитических функций $g$ в круге $\dd$, допускающих представление 
\bay
\label{predprobes}
\hspace*{-.4cm}g=\sum_n\f_n\psi_n,~\;\mbox{где}~\; \f_n\in H^1,~
\psi_n\in\pp_+(b_\infty^{-1}(\T)),~\;
\sum_n\|f_n\|_{H^1}\|\psi_n\|_{\pp_+(b_\infty^{-1}(\T))}<\be.
\ey
Очевидно, что такие функции $g$ входят в пространство $\pp_+(b^{-1}_{1,\be}(\T))$, двойственное к которому отождествляется естественным образом с пространством Бесова $\pp_+B^1_{\be,1}(\T)$ аналитических в $\dd$ функций. Тем не менее, не всякая функция из \lb$\pp_+(b^{-1}_{1,\be}(\T))$ представима в виде \rf{predprobes}. Иначе мы имели бы равенство $L=B_{\be,1}^1(\T)$, что невозможно, ибо условие $f\in L$, будучи необходимым для операторной липшицевости, не является достаточным.

\medskip

{\bf Замечание.}
Отметим здесь, что пространство 
$$
L=\{f\in {\rm BMO}(\R):~\zh{\rm Hess}\,\mP f\zh\,d\m_2\in{\rm CM}(\C_+)\}
$$
является предельным пространством шкалы пространств Трибеля--Лизоркина 
\lb$F_{p,q}^s(\R)$ и обозначается символом $F_{\be,1}^1(\R)$ 
(см. \cite {FJ}, \S\:5). Аналогичным образом определяется пространство Трибеля--Лизоркина $F_{\be,1}^1(\T)$ функций на $\T$. Полученные здесь необходимые условия операторной липшицевости можно переформулировать следующим образом:
$\OL(\R)\subset F_{\be,1}^1(\R)$ и $\OL(\T)\subset F_{\be,1}^1(\T)$.

\

\section{\bf Достаточное условие для операторной липшицевости в терминах классов Бесова}
\setcounter{equation}{0}
\label{Dost}

\

В этом параграфе мы покажем, что функции класса Бесова $B_{\be,1}^1(\R)$ (см. \S\:\ref{Prel}) являются операторно липшицевыми. Мы также получим аналогичный результат для функций на окружности. Приводимые здесь доказательства отличаются от оригинальных доказательств работ \cite{Pe1} и \cite{Pe3} и основаны на операторных неравенствах Бернштейна (см. \S\:\ref{Bern}).

\begin{thm}
\label{Besdost}
Пусть $f\in B_{\be,1}^1(\R)$. Тогда функция $f$ операторно липшицева  и
\bay
\label{BeOLsa}
\|f(A)-f(B)\|\le\const\|f\|_{B_{\be,1}^1}\|A-B\|
\ey
для любых самосопряжённых операторов $A$ и $B$ с ограниченной разностью $A-B$.
\end{thm}

\Pf Как мы отмечали во введении (см. теорему \ref{anbnrn} ниже), достаточно доказать неравенство \rf{BeOLsa}
для ограниченных самосопряжённых операторов $A$ и $B$.

Не умаляя общности, мы можем считать, что $f(0)=0$.
Рассмотрим функции $f_n=f*W_n$, определённые равенством \rf{fn}. Пусть $N\in\Z$. 
Положим $g_n\df f_n-f_n(0)$. Из определения класса $B_{\be,1}^1(\R)$ 
(см. \S\:\ref{Prel}) вытекает, что 
$\sum_{n=-\be}^\be g'_n=f'$, 
причём ряд сходится равномерно на $\R$. Отсюда вытекает, что ряд 
$\sum\limits_{n=-\be}^\be g_n$ сходится равномерно на каждом компактном
подмножестве в $\R$. Стало быть, 
$$
\sum_{n=-\be}^\be g_n(A)=f(A)\quad\mbox{и}\quad\sum_{n=-\be}^\be g_n(B)=f(B),
$$
причём ряды сходятся абсолютно по операторной норме. 

Поскольку очевидно, что $f_n\in \E_{2^n+1}\cap L^\be(\R)$, операторное неравенство Бернштейна \rf{opnerBer} позволяет заключить, что
\begin{align*}
\|f(A)-f(B)\|&\le\left\|\sum_{n=-\be}^\be\big((g_n(A)-g_n(B)\big)\right\|
=\left\|\sum_{n=-\be}^\be\big((f_n(A)-f_n(B)\big)\right\|\\[.2cm]
&\le\sum_{n=-\be}^\be2^{n+1}\|f_n\|_{L^\be}\|A-B\|
\le\const\|f\|_{B_{\be,1}^1}\|A-B\|.\quad\bl
\end{align*}

Аналогичным образом можно доказать следующий аналог теоремы \ref{Besdost} для функций на единичной окружности.

\begin{thm}
\label{Besokr}
Пусть $f\in B_{\be,1}^1(\T)$. Тогда $f$ -- операторно липшицева функция и
$$
\|f(U)-f(V)\|\le\const\|f\|_{B_{\be,1}^1}\|U-V\|,\quad f\in B_{\be,1}^1(\T),
$$
для любых унитарных операторов $U$ и $V$.
\end{thm}

Следующая формулировка позволяет нам объединить необходимые условия, полученные в \S\:\ref{Neob}, с достаточными условиями этого параграфа.

\begin{thm}
$B_{\be,1}^1(\T)\subset\OL(\T)\subset F_{\be,1}^1(\T)$ и 
$B_{\be,1}^1(\R)\subset\OL(\T)\subset F_{\be,1}^1(\R)$.
\end{thm}

Оказывается, что функции класса $B_{\be,1}^1(\R)$ не только операторно липшицевы, но и операторно дифференцируемы.

\begin{thm}
\label{Besdiffer}
Пусть $f$ - функция класса $B_{\be,1}^1(\R)$. Тогда $f$ -- операторно дифференцируемая функция.
\end{thm}

Мы отсылаем читателя к работам \cite{Pe3} и \cite{Pe5}
для доказательства теоремы \ref{Besdiffer}.

\

\section{\bf Операторно гёльдеровы функции}
\setcounter{equation}{0}
\label{oHold}

\

В этом параграфе мы обсудим другие приложения операторных неравенств Бернштейна, полученных в \S\:\ref{Bern}. Мы покажем, что класс операторно гёльдеровых функций порядка $\a$, $0<\a<1$, совпадает с классом гёльдеровых функций порядка $\a$.
Мы также кратко остановимся на случае произвольных модулей непрерывности.
Результаты этого параграфа получены в \cite{AP1} и \cite{AP2}. Другой подход к этим задачам был найден в \cite{FN}, где были получены аналогичные результаты для функций класса Гёльдера с несколько худшими константами и несколько более слабый результат для произвольных модулей непрерывности.

Хорошо известно, что (скалярное) неравенство Бернштейна играют важную роль в теории приближений, см., например, \cite{Ah}, \cite{Dz}, \cite{NI}, \cite{Ti}. Речь идёт об описании 
тех или иных свойств типа гладкости в терминах приближений хорошими функциями.
 Прямые теоремы теории приближений дают
оценки приближений функций из данного функционального пространства  $X$ (обычно в том или ином 
смысле гладких функций) хорошими 
функциями. Обратные теоремы позволяют по оценкам приближений данной функции
$f$ хорошими функциями сделать вывод о принадлежности данной функции 
тому или иному функциональному пространству.
В случае, когда для функционального пространства $X$ прямые теоремы ``стыкуются''
с обратными теоремами, получается полное описание пространства $X$ в терминах приближений.

В этом параграфе мы будем рассматривать функциональные пространства на
единичной окружности $\T$ и на вещественной прямой $\R$. В качестве хороших функций
в первом случае берутся пространства $\mathcal P_n$ тригонометрических многочленов степени
не выше $n$, а во втором -- пространства целых функций $\E_\s$ экспоненциального типа
не выше $\s$. Мы будем рассматривать только равномерные приближения.

Классические неравенства Бернштейна играют решающую роль при доказательстве обратных теорем теории приближений. Нетрудно понять, что если в таком доказательстве
воспользоваться операторной версией неравенства Бернштейна, то мы получим соответствующую
гладкость функции $f$  на множестве унитарных операторов, если речь идёт о функциях
на окружности, или на множестве самосопряжённых операторов, если речь идёт о функциях
на прямой.

Проиллюстрируем это на примере. 
Классическая теорема Джексона утверждает, что если функция $f$ 
входит в класс Гёльдера $\L_\a(\T)$, $0<\a<1$, то
$$
\dist(f,\mathcal P_n)\le\const(n+1)^{-\a}\|f\|_{\L_\a}.
$$

Бернштейн доказал, что верно и обратное, т. е. если для некоторой непрерывной функции $f$ на $C(\T)$ при $\a\in(0,1)$ выполняются неравенства
$$
\dist(f,\mathcal P_n)\le c(n+1)^{-\a},\quad n\in\Z_+,
$$
для положительного числа $c$, то  $f\in\L_\a(\T)$.

Приведём стандартное доказательство этого результата Бернштейна.
Не умаляя общности, можно считать, что $c=1$. 
Для $n\ge0$ найдётся тригонометрический полином $f_n$  
такой, что $\deg f_n<2^n$  и $\|f-f_n\|_{C(\T)}\le 2^{-\a n}$. Ясно, что 
$$
\|f_n-f_{n-1}\|_{C(\T)}\le\|f-f_n\|_{C(\T)}+\|f-f_{n-1}\|_{C(\T)}\le2^{-\a n}(1+2^\a)\le3\cdot2^{-\a n}.
$$
Следовательно, 
$$
\|f_n-f_{n-1}\|_{\Li(\T)}\le2^n\|f_n-f_{n-1}\|_{C(\T)}\le3\cdot2^{(1-\a)n}
$$
в силу неравенства Бернштейна. Принимая во внимание очевидное равенство
\lb$\|f_0\|_{\Li(\T)}=0$, получаем
$$
\|f_N\|_{\Li(\T)}\le\sum_{n=1}^N\|f_n-f_{n-1}\|_{\Li(\T)}\le3\sum_{n=1}^N2^{(1-\a)n}\le\frac3{1-2^{\a-1}}2^{(1-\a)N},\quad N\in\Z_+.
$$
Пусть $\z,\,\tau\in\T$. Выберем $N\in\Z_+$
таким образом, чтобы $2^{-N}<|\z-\tau|\le2^{1-N}$. 
Тогда
\begin{align*}
|f(\z)-f(\tau)|&\le|f(\z)-f_N(\z)|+|f_N(\z)-f_N(\tau)|+|f_N(\xi)-f(\tau)|\\[.2cm]
&\le2\|f-f_N\|_{L^\be}+\|f_N\|_{\Li}|\z-\tau|\le2\cdot 2^{-\a N}+\frac3{1-2^{\a-1}}2^{(1-\a)N}|\z-\tau|\\[.2cm]
&\le2|\z-\tau|^\a+
\frac{3\cdot2^{1-\a}}{1-2^{\a-1}}|\z-\tau|^\a\le\frac8{1-2^{\a-1}}|\z-\tau|^{\a}.\quad\bl
\end{align*}

Следующая теорема говорит о том, что всякая функция $f$ из $\L_\a(\T)$, $0<\a<1$, является 
операторно гёльдеровой порядка $\a$, что резко контрастирует со случаем липшицевых функций.

\begin{thm} 
\label{43}
Пусть $f\in\L_\a(\T)$, где $\a\in(0,1)$. 
Тогда существует константа $c$ такая, что
$$
\|f(U)-f(V)\|\le\frac{8c}{1-2^{\a-1}}\|f\|_{\L_\a}\|U-V\|^\a
$$
для любых унитарных операторов $U$ и $V$.
\end{thm}

\Pf Пусть $f\in\L_\a(\T)$. Сначала применяем соответствующую прямую теорему теории приближений (в данном случае теорему Джексона).  В силу этой теоремы существует константа $c$ такая, что
$$
\dist(f,\mathcal P_n)\le c(n+1)^{-\a}\|f\|_{\L_\a},\quad n\in\Z_+.
$$ 
Повторяя теперь почти слово в слово доказательство соответствующей обратной теоремы 
(оно было изложено перед формулировкой этой теоремы), заменяя $\z$ и $\tau$  унитарными операторами $U$ и $V$ и применяя операторное 
неравенство Бернштейна вместо скалярного неравенства Бернштейна, получаем искомый результат.  $\bl$

\medskip


Аналогичное утверждение верно и в случае прямой.

\begin{thm} 
\label{44}
Пусть $f\in\L_\a(\R)$, $\a\in(0,1)$. 
Тогда существует константа $c$ такая, что
$$
\|f(A)-f(B)\|\le\frac{8c}{1-2^{\a-1}}\|f\|_{\L_\a}\|A-B\|^\a
$$
для любых самосопряжённых операторов $A$ и $B$.
\end{thm}

Доказательство основано на аналогичном описании функций класса $\L_\a(\R)$
в терминах приближения целыми функциями экспоненциального типа.

\medskip

{\bf Прямая теорема для пространства $\L_\a(\R)$.} {\em Пусть $f\in\L_\a(\R)$, $0<\a<1$. Тогда существует 
абсолютная константа $c$ такая, что 
$$
\inf\{h\in\E_\s:\|f-h\|_{L^\be(\R)}\}\le c\s^{-\a}\|f\|_{\L_\a(\T)}
$$
для любого $\s>0$.}

\medskip

{\bf Обратная теорема для пространства $\L_\a(\R)$.} {\em Пусть $0<\a<1$. Пусть $f$ -- непрерывная функция на 
вещественной прямой $\R$ такая, что $\lim\limits_{|x|\to\be} x^{-1}f(x)=0$. Предположим, что
$$
\inf\{h\in\E_\s:\|f-h\|_{L^\be(\R)}\}\le c\s^{-\a}
$$
для некоторого положительного числа $c$ и для всех $\s>0$. Тогда $f\in\L_\a(\R)$ и 
$$
\|f\|_{\L_\a(\R)}\le\frac{5c}{1-2^{\a-1}}.
$$

}

\medskip

Более подробно доказательства теорем \ref{43} и \ref{44} изложены в статье \cite{AP2}.
В этой же работе получен и ряд других результатов, основанных, в конечном итоге, на
некоторых результатах теории приближений.

В частности, аналоги теорем \ref{43} и \ref{44} получены там для всех $\a>0$. 
Приведём ещё некоторые результаты, полученные в \cite{AP2}, которые также могут
быть доказаны методами теории приближений.

Функцию $\o:[0,+\be)\to\R$ будем называть {\it модулем непрерывности}, если она
является  неотрицательной неубывающей непрерывной функцией такой,
что $\o(0)=0$, $\o(x)>0$ при $x>0$ и
$\o(x+y)\le\o(x)+\o(y)$ для всех $x,y\in[0,+\be)$.

Обозначим через $\L_\o(\R)$ пространство всех непрерывных функций на $\R$
таких, что 
$$
\|f\|_{\L_\o(\R)}\df\sup_{x\ne y}\frac{|f(x)-f(y)|}{\o(|x-y|)}<+\be.
$$
Аналогичным образом можно определить пространство $\L_\o(\T)$.

Положим 
\bay
\label{omzv}
\o_*(x)\df x\int_x^\be \frac{\o(t)}{t^2}\,dt.
\ey

\begin{thm} 
\label{prmone}
Пусть $f\in\L_\o(\R)$, где $\o$ --  модуль непрерывности.
Тогда 
$$
\|f(A)-f(B)\|\le c\|f\|_{\L_\o(\R)}\o_*(\|A-B\|)
$$  
для любых самосопряжённых операторов
$A$ и  $B$, где $c$ -- абсолютная постоянная.
\end{thm}

Аналогичный результат имеет место и для функций $f\in\L_\o(\T)$.

\

\section{\bf Гёльдеровские функции при возмущениях операторами класса Шаттена--фон Неймана}
\setcounter{equation}{0}
\label{SchavN}

\

Мы рассмотрим в этом параграфе ещё одно приложение операторных неравенств Бернштейна, приведённых в \S\:\ref{Bern}. Пусть $f$ -- функция класса Гёльдера
$\L_\a(\R)$, \lb$0<\a<1$, и пусть $p\ge1$. Предположим, что самосопряжённый оператор $B$ является возмущением самосопряжённого оператора $A$ оператором класса
$\bS_p$, т.е. $B-A\in\bS_p$. Зададим себе вопрос, что мы можем сказать об операторе $f(A)-f(B)$. Этот вопрос подробно изучался в работе \cite{AP3}.
Здесь мы сформулируем результат работы \cite{AP3} в случае $p>1$.

\begin{thm}
\label{aHolSp}
Пусть $p>1$ и $0<\a<1$. Тогда существует положительное число $c$ такое, что
$$
\|f(A)-f(B)\|_{\bS_{p/\a}}\le c\|f\|_{\L_\a}\|A-B\|_{\bS_p}^\a
$$
для произвольных самосопряжённых операторов $A$ и $B$, разность которых входит в класс $\bS_p$.
\end{thm}

Мы опускаем здесь доказательство теоремы \ref{aHolSp} и отсылаем читателя за доказательством к работе \cite{AP3}. Случай $p=1$ также подробно рассмотрен в \cite{AP3}. Заметим здесь, что при $p=1$ заключение теоремы \ref{aHolSp} неверно. В работе 
\cite{AP3} также получен аналог теоремы \ref{aHolSp} при всех положительных $\a$,
а также рассматриваются более общие вопросы возмущений операторами из симметрично нормируемых идеалов.

\

\begin{center}
\bf\Large Глава II
\end{center}

\medskip

\begin{center}
\bf\Large 
Мультипликаторы Шура и \\ двойные операторные интегралы
\end{center}
\label{muSch}

\addtocontents{toc}{\vspace*{.3cm}\hspace*{-.35cm}\textbf{{\bf Глава II}. Мультипликаторы Шура и двойные операторные интегралы}\hfill\pageref{muSch}}

\renewcommand{\thesection}{2.\arabic{section}}
\setcounter{section}{0}

\

В этой главе мы будем изучать мультипликаторы Шура, как дискретные, так и мультипликаторы Шура по отношению к спектральным мерам. Мы используем описание дискретных мультипликаторов Шура, основанное на теореме Гротендика (см. книгу
Ж. Пизье \cite{Pi} и работу \cite{Pi2}). Этот результат уточняется в случае, когда исходная функция задана на произведении топологических пространств и является непрерывной по каждой переменной. Получено также уточнение этого результата для борелевских функций, заданных на произведении топологических пространств.

Затем мы определяем двойные операторные интегралы и мультипликаторы Шура по отношению к спектральным мерам. Изучение таких мультипликаторов Шура в случае борелевских функций на произведении топологических пространств сводится к дискретным мультипликаторам Шура.

\

\section{\bf Дискретные мультипликаторы Шура}
\setcounter{equation}{0}
\label{dmsh}

\

Обозначим через $\ell^p({\mT})$ пространство
комплексных функций $\a: t\mapsto\a_t$, заданных на не обязательно счётном или конечном
множестве $\mT$ и таких, что $\sum\limits_{t\in\mT}|\a_t|^p <\be$,
с нормой
$\|\a\|_p =\Big(\sum_{t\in\mT}|\a_t|^p\Big)^{1/p}$, где $p\in[1,+\be)$.
При $p=\be$ пространство $\ell^p(\mT)$ состоит из всех ограниченных комплексных 
функций $\a: t\mapsto\a_t$, заданных на $\mT$, и $\|\a\|_\be = \sup_{t\in\mT}|\a_t|$. 
В тех случаях, когда необходимо явно
указать множество $\mT$, на котором задано семейство $\a$, мы будем писать $\|\a\|_{\ell^p(\mT)}$
вместо $\|\a\|_p$. Обозначим через $c_0(\mT)$ подпространство
пространства $\ell^\be(\mT)$, состоящее из функций $\a$, стремящихся к нулю в бесконечности.

Пусть ${\mS}$ и $\mT$ -- произвольные непустые множества. Каждому ограниченному оператору
$A:\ell^2(\mT)\to\ell^2(\mS)$ соответствует единственная матрица $\{a(s,t)\}_{(s,t)\in\mS\times\mT}$ такая, что
$(Ax)_s=\sum\limits_{t\in T}a(s,t)x_t$ для всех $x=\{x_t\}_{t\in\mT}$ из $\ell^2(\mT)$. В этом случае мы говорим, что
{\it матрица $\{a(s,t)\}_{(s,t)\in\mS\times\mT}$ определяет ограниченный оператор} $A:\ell^2(\mT)\to\ell^2(\mS)$.  Положим
$$
\|\{a(s,t)\}_{(s,t)\in\mS\times\mT}\|\df\|A\|
\quad\mbox{и}\quad \|\{a(s,t)\}_{(s,t)\in\mS\times\mT}\|_{\bS_1}\df\|A\|_{\bS_1}.
$$
В случае, если $A\not\in\bS_1(\ell^2(\mT),\ell^2(\mS))$, будем  считать, что последняя норма равна $\be$. 
Если матрица $\{a(s,t)\}_{(s,t)\in\mS\times\mT}$ не определяет ограниченный
оператор из $\ell^2(\mT)$ в $\ell^2(\mS)$, будем считать, что её операторная норма (равно, как и ядерная) равна $\be$.
Обозначим через $\mB(\mS\times\mT)$ {\it множество всех матриц $\{a(s,t)\}_{(s,t)\in\mS\times\mT}$,
задающий ограниченный оператор  из} $\ell^2(\mT)$ в $\ell^2(\mS)$. В некоторых случаях
мы будем писать $\|\{a(s,t)\}_{(s,t)\in\mS\times\mT}\|_{\mB(\mS\times\mT)}$ вместо 
$\|\{a(s,t)\}_{(s,t)\in\mS\times\mT}\|$
и $\|\{a(s,t)\}_{(s,t)\in\mS\times\mT}\|_{\bS_1(\mS\times\mT)}$ вместо 
$\|\{a(s,t)\}_{(s,t)\in\mS\times\mT}\|_{\bS_1}$.

Матрица ${\Phi}=\{\Phi(s,t)\}_{(s,t)\in\mS\times\mT}$ называется {\it мультипликатором Шура} 
пространства
$\mB(\mS\times\mT)$,
если для любой матрицы
$A=\{a(s,t)\}_{(s,t)\in\mS\times\mT}$ из $\mB(\mS\times\mT)$ матрица
${\Phi}\star A\df\{\Phi(s,t)a(s,t)\}_{(s,t)\in\mS\times\mT}$ также
принадлежит пространству $\mB(\mS\times\mT)$.

{\it Множество всех мультипликаторов Шура пространства $\mB(\mS\times\mT)$ обозначим символом} $\fM(\mS\times\mT)$. Из теоремы о замкнутом графике легко вывести, что мультипликаторы Шура
порождают ограниченные операторы в $\mB(\mS\times\mT)$. Положим
$$
\|\Phi\|_{\fM(\mS\times\mT)}\df
\sup\{\|\Phi\star A\|:~A\in\mB(\mS\times\mT),~  \|A\|_\mB\le1\}.
$$
Отсюда по двойственности получаем, что
\bay
\label{S1}
\|\Phi\|_{\fM(\mS\times\mT)}=\sup\{\|\Phi\star A\|_{\bS_1}:~A\in\mB(\mS\times\mT),~
\|A\|_{\bS_1}\le1\}.
\ey

Легко видеть, что
$$
\|A\|_{\mB(\mS\times\mT)}=\sup\|A\|_{\mB(\mS_0\times\mT_0)}, \quad
\|A\|_{\bS_1(S\times\mT)}=\sup\|A\|_{\bS_1(\mS_0\times\mT_0)}
$$
и
$$
\|\Phi\|_{\fM(\mS\times\mT)}=\sup\|{\bs\f}\|_{\fM(\mS_0\times\mT_0)},
$$
где супремумы берутся по всем конечным подмножествам $\mS_0$ и $\mT_0$
множеств $\mS$ и $\mT$.

Отметим ещё, что 
$\|\Phi\|_{\ell^\be(\mS\times\mT)}\le\|\Phi\|_{\fM(\mS\times\mT)}$.
Легко видеть, что неравенство превращается в равенство для любой
матрицы $\{\Phi(s,t)\}_{(s,t)\in\mS\times\mT}$ ранга $1$.
Есть и другие классы матриц, для которых это неравенство превращается в равенство.
Например, если каждая строка (или каждый столбец) матрицы $\Phi$ содержит
не более одного ненулевого элемента, то  
$\|\Phi\|_{\ell^\be(\mS\times\mT)}=\|\Phi\|_{\fM(\mS\times\mT)}$.

Нам понадобится ещё одна характеристика матрицы
$\Phi$. Положим
$$
\|\Phi\|_{\fM_0(\mS\times\mT)}\df
\sup\{\|\Phi\star A\|:~A\in\mB(\mS\times\mT),~
\|A\|\le1,~a(t,t)=0~\;\mbox{при}\;~t\in\mS\cap\mT\}.
$$
Обозначим символом $\fM_0(\mS\times\mT)$  множество всех матриц 
$\Phi=\{\Phi(s,t)\}_{(s,t)\in\mS\times\mT}$ таких, что $\|\Phi\|_{\fM_0(\mS\times\mT)}<+\be$.
Очевидно, $\|\Phi\|_{\fM_0(\mS\times\mT)}\le\|\Phi\|_{\fM(\mS\times\mT)}$.
Легко видеть, что
$$
\|\Phi\|_{\fM_0(\mS\times\mT)}=\sup\|\Phi\|_{\fM_0(\mS_0\times\mT_0)},
$$
где супремум берётся по всем конечным подмножествам $\mS_0$ и $\mT_0$
множеств $\mS$ и $\mT$.


Заметим ещё, что если матрицы $\Phi=\{\Phi(s,t)\}_{(s,t)\in\mS\times\mT}$
и $\Psi=\{\Psi(s,t)\}_{(s,t)\in\mS\times\mT}$ совпадают вне ``диагонали''
$\{(t,t):t\in\mS\cap\mT\}$,  то
$$
\|\Phi-\Psi\|_{\fM_0(\mS\times\mT)}=0\quad
\text{и}\quad\|\Phi\|_{\fM_0(\mS\times\mT)}=\|\Psi\|_{\fM_0(\mS\times\mT)}.
$$

Прежде, чем сформулировать следующее утверждение, отметим, что
$\|\Phi\|_{\fM_0(\mS\times\mT)}=\|\Phi\|_{\fM(\mS\times\mT)}$, если
$\mS\cap\mT=\varnothing$.

\begin{lem} 
Пусть $\Phi\in\ell^\be(\mS\times\mT)$, где $\mS$ и $\mT$ -- произвольные
множества такие, что $\mS\cap\mT\ne\varnothing$.
Тогда
$$
\max\left\{\|\Phi\|_{\fM_0(\mS\times\mT)},\|\Phi(t,t)\|_{\ell^\be(\mS\cap\mT)}\right\}
\le\|\Phi\|_{\fM(\mS\times\mT)}
\le2\|\Phi\|_{\fM_0(\mS\times\mT)}+\|\Phi(t,t)\|_{\ell^\be(\mS\cap\mT)}.
$$
\end{lem}

\Pf  Первое неравенство очевидно. Докажем второе.
Обозначим символом $\chi$ характеристическую функцию множества
$\{(s,t)\in\mS\times\mT:s=t\}$. Легко видеть, что $\|\chi\|_{\fM(\mS\times\mT)}=1$,
 откуда $\|1-\chi\|_{\fM(\mS\times\mT)}\le\|1\|_{\fM(\mS\times\mT)}+\|\chi\|_{\fM(\mS\times\mT)}=2$.
 Пусть $A\in\mB(S\times T)$ и $\|A\|\le1$. Тогда
 $$
 \Phi\star A=\Phi\star(1-\chi)\star A+\Phi\star\chi\star A.
 $$
 Остаётся заметить, что
 $$
 \|\Phi\star(1-\chi)\star A\|\le\|\Phi\|_{\fM_0(S\times T)}
 \|(1-\chi)\star A\|\le2\|\Phi\|_{\fM_0(S\times T)}
 $$
 и
 $$
 \|\Phi\star\chi\star A\|\le\|\Phi\star\chi\|_{\fM(S\times T)}
 =\|\Phi(t,t)\|_{\ell^\be(S\cap T)}.\quad \bl
 $$

 \begin{cor}
 \label{32}
 Если $\Phi(t,t)=0$ при всех $t\in\mS\cap\mT$, то
 $$
\|\Phi\|_{\fM_0(\mS\times\mT)}\le\|\Phi\|_{\fM(\mS\times\mT)}
\le2\|\Phi\|_{\fM_0(\mS\times\mT)}.
$$
 \end{cor}

 \begin{lem}
 \label{33}
 Пусть $\mS$ -- хаусдорфово топологические пространство. Предположим, что
 никакая точка множества $\mS\cap\mT$ не является изолированной точкой
 в пространстве $\mS$.
 Тогда  $\|\Phi\|_{\fM_0(\mS\times\mT)}=\|\Phi\|_{\fM(\mS\times\mT)}$ для любой
 функции $\Phi\in\ell^\be(\mS\times\mT)$, непрерывной по переменной $s\in\mS$.
 \end{lem}

 \Pf Достаточно доказать, что
 $\|\Phi\|_{\fM(\mS\times\mT)}\le\|\Phi\|_{\fM_0(\mS\times\mT)}$
 или, что то же самое, $\|\Phi\|_{\fM(\mS_0\times\mT_0)}
 \le\|\Phi\|_{\fM_0(\mS\times\mT)}$
 для всех конечных подмножеств $\mS_0$ и $\mT_0$ множеств $\mS$ и $\mT$.
 Зафиксируем конечные подмножества $\mS_0$ и $\mT_0$ множеств $\mS$ и $\mT$.
 Ясно, что
 для любого положительного числа  $\e$ существует возмущение $\widetilde\mS_0$
 множества $\mS_0$ такое, что $\widetilde\mS_0\cap\mT_0=\varnothing$ и
 $\|\Phi\|_{\fM(\mS_0\times\mT_0)}<\e+\|\Phi\|_{\fM(\widetilde\mS_0\times\mT_0)}$.
 Следовательно,
 $$
 \|\Phi\|_{\fM(\mS_0\times\mT_0)}<\e+\|\Phi\|_{\fM(\widetilde\mS_0\times\mT_0)}=
 \e+\|\Phi\|_{\fM_0(\widetilde\mS_0\times\mT_0)}\le\e+\|\Phi\|_{\fM_0(\mS\times\mT)}
 $$
 для любого $\e>0$. $\bl$
 
 \medskip

Мы собираемся рассмотреть аналог пространства $\fM_{0}(\mS\times\mT)$,
определённого через $\bS_1$-норму вместо операторной нормы.
Для этого мы положим
$$
\|\Phi\|_{\fM_{0,\bS_1}(\mS\times\mT)}\df
\sup\{\|\Phi\star A\|_{\bS_1}:~A\in\mB(\mS\times\mT),~
\|A\|_{\bS_1}\le1,~a(t,t)=0~\mbox{для}~t\in\mS\cap\mT\}
$$
и $\fM_{0,\bS_1}(\mS\times\mT)\df\{\Phi: \|\Phi\|_{\fM_{0,\bS_1}(\mS\times\mT)}<+\be\}$.
При этом следует отметить, что соответствующий аналог пространства 
$\fM(\mS\times\mT)$ нам не нужно определять, ибо он совпадает с этим же пространством
$\fM(\mS\times\mT)$ в силу равенства \rf{S1}.

Точно так же, как для операторной нормы, можно доказать следующие утверждения.

\begin{lem}
\label{32s1l}
Пусть $\Phi\in\ell^\be(\mS\times\mT)$, где $\mS$ и $\mT$ -- произвольные
множества такие, что $\mS\cap\mT\ne\varnothing$.
Тогда
\begin{align*}
\max\left\{\|\Phi\|_{\fM_{0,\bS_1}(\mS\times\mT)},\|\Phi(t,t)\|_{\ell^\be(\mS\cap\mT)}\right\}
&\le\|\Phi\|_{\fM(\mS\times\mT)}\\[.2cm]
&\le2\|\Phi\|_{\fM_{0,\bS_1}(\mS\times\mT)}+\|\Phi(t,t)\|_{\ell^\be(\mS\cap\mT)}.
\end{align*}
Если $\Phi(t,t)=0$ при всех $t\in\mS\cap\mT$, то
 $$
\|\Phi\|_{\fM_{0,\bS_1}(\mS\times\mT)}\le\|\Phi\|_{\fM(\mS\times\mT)}
\le2\|\Phi\|_{\fM_{0,\bS_1}(\mS\times\mT)}.
$$
\end{lem}

\begin{cor}
\label{32s1}
Если $\Phi(t,t)=0$ при всех $t\in\mS\cap\mT$, то
$$
\|\Phi\|_{\fM_{0,\bS_1}(\mS\times\mT)}\le\|\Phi\|_{\fM(\mS\times\mT)}
\le2\|\Phi\|_{\fM_{0,\bS_1}(\mS\times\mT)}.
$$
\end{cor}

\begin{lem}
\label{33s1}
Пусть $\mS$ -- хаусдорфово топологические пространство. Предположим, что
никакая точка множества $\mS\cap\mT$ не является изолированной точкой
в пространстве $\mS$.
Тогда  $\|\Phi\|_{\fM_{0,\bS_1}(\mS\times\mT)}=\|\Phi\|_{\fM(\mS\times\mT)}$ для любой
функции $\Phi$ из $\ell^\be(\mS\times\mT)$, непрерывной по переменной $s\in\mS$.
\end{lem}

\

\section{\bf Описание дискретных мультипликаторы Шура}
\setcounter{equation}{0}
\label{dmsh+}

\

\begin{thm}
\label{41}
Пусть $\{u_s\}_{s\in\mS}$ и $\{v_t\}_{t\in\mT}$ -- семейства  векторов в
(не обязательно сепарабельном) гильбертовом пространстве $\h$
такие, что $\|u_s\|\cdot\| v_t\|\le 1$ для всех $s$ из $\mS$ и
 $t$ из $\mT$. Положим
$\Phi(s,t)\df(u_s, v_t)$, $s\in\mS$, $t\in\mT$.  Тогда $\Phi \in\fM(\mS\times\mT)$
и $\|\Phi\|_{\fM(\mS\times\mT)}\le 1$.
\end{thm}

\Pf В силу \rf{S1} достаточно доказать, что
$$
\|\{a(s,t)(u_s,v_t)\}\|_{\bS_1} \le\|\{a(s,t)\}\|_{\bS_1}
$$
для любой матрицы $\{a(s,t)\}$, задающей ядерный оператор. Ясно, что достаточно
ограничиться случаем, когда $\rank\{a(s,t)\}=1$. Кроме того, можно считать, что
$\|a(s,t)\|_{\bS_1}=1$. Тогда $a(s,t)=\a_s\b_t$ для некоторых $\a\in\ell^2(\mS)$
и $\b\in\ell^2(\mT)$ таких, что $\|\a\|_{\ell^2(\mS)}=\|\b\|_{\ell^2(\mT)}=1$.
Пусть $\{e_j\}_{j\in J}$ -- ортонормированный базис в гильбертовом пространстве
$\h$. Положим $\hat x(j)\df(x,e_j)$, $j\in J$. Тогда
\begin{align*}
\|\{\a_s\b_t(u_s,v_t)\}\|_{\bS_1}&\le\sum_{j\in J}\|\{\a_s\b_t\hat u_s(j)\ov{\hat v_t(j)})\}\|_{\bS_1}\\[.2cm]
&=\sum_{j\in J}\|\{\a_s\hat u_s(j)\}\|_{\ell^2(\mS)}\|\{\b_t\ov{\hat v_t(j)}\}\|_{\ell^2(\mT)}\\[.2cm]
&\le\left(\sum_{j\in J}\|\{\a_s\hat u_s(j)\}\|_{\ell^2(\mS)}^2\right)^{\frac12}
\left(\sum_{j\in J}\|\{\b_t\ov{\hat v_t(j)}\}\|_{\ell^2(\mT)}^2\right)^{\frac12}.
\end{align*}
Ясно, что,
\begin{align*}
\sum_{j\in J}\|\{\a_s\hat u_s(j)\}\|_{\ell^2(\mS)}^2
&=\sum_{j\in J}\sum_{s\in\mS}|\a_s|^2|\hat u_s(j)|^2\\[.2cm]
&=\sum_{s\in\mS}|\a_s|^2\sum_{j\in J}|\hat u_s(j)|^2
=\sum_{s\in\mS}|\a_s|^2\|u_s\|^2
\le\sup_{s\in\mS}\|u_s\|^2.
\end{align*}
Аналогичным образом
$$
\sum_{j\in J}\|\{\b_t\ov{\hat v_t(j)}\}\|_{\ell^2(\mT)}^2\le\sup_{t\in\mT}\|v_t\|^2.
$$
Отсюда
$$
\|\{\a_s\b_t(u_s,v_t)\}\|_{\bS_1}
\le\sup_{s\in\mS}\|u_s\|\sup_{t\in\mT}\|u_t\|\le 1. \quad  \bl
$$

Весьма нетривиальным является обратное утверждение,
которое содержится в теореме 5.1 монографии \cite{Pi}, см. также \cite{Pi2}.
Этот результат мы приводим без доказательства.

\begin{thm}
\label{Pi}
Пусть $\Phi\df\{\Phi(s,t)\}$ -- мультипликатор Шура на $\mB(\mS\times\mT)$ и
$\|\Phi\|_\fM\le1$. Тогда найдутся два семейства векторов
$\{u_s\}_{s\in\mS}$ и $\{v_t\}_{t\in\mT}$  в
(не обязательно сепарабельном) гильбертовом пространстве $\h$
такие, что $\|u_s\|\le1$ при всех $s\in\mS$, $\| v_t\|\le 1$  всех
 $t\in\mT$ и
$$
\Phi(s,t)\df(u_s, v_t),\quad s\in\mS,\quad t\in\mT.
$$
\end{thm}

%

%
%
%

{\bf Замечание  к теореме \ref{Pi}.}  {\it В этой теореме можно дополнительно потребовать, чтобы линейная оболочка как семейства $\{u_s\}_{s\in\mS}$,
так и семейства $\{ v_t\}_{t\in\mT}$ была всюду плотна в гильбертовом пространстве $\h$}. Действительно,
пусть $\h_1$ -- замыкание линейной оболочки семейства $\{ v_t\}_{t\in\mT}$ и пусть $P_1$ -- ортогональный проектор на $\h_1$.
Тогда
$\{P_1u_s\}_{s\in\mS}$ и $\{ v_t\}_{t\in\mT}$ -- семейства в гильбертовом пространстве $\h_1$ такие,
что $\Phi(s,t)=(P_1u_s,v_t)$ для всех $(s,t)\in\mS\times\mT$.
Пусть теперь $\h_2$ -- замыкание линейной оболочки семейства $\{P_1u_s\}_{s\in\mS}$, а $P_2$ -- ортогональный проектор на $\h_2$.
Тогда
$\{P_1u_s\}_{s\in\mS}$
и $\{P_2 v_t\}_{t\in\mT}$ -- семейства векторов в  $\h_2$ такие, что
 $\Phi(s,t)=(P_1u_s,P_2 v_t)$ при $(s,t)\in\mS\times\mT$. Ясно, что линейные оболочки семейств
$\{P_1u_s\}_{s\in\mS}$ и $\{P_2 v_t\}_{t\in\mT}$ плотны в гильбертовом пространстве $\h_2$. 

\medskip

Напомним, что отображение $g$ из топологического пространства $\mT$ в гильбертово пространство $\h$
называется {\it слабо непрерывным}, если оно непрерывно из 
$\mT$ в гильбертово пространство $\h$, наделённое слабой топологией, т.е. функция $t\mapsto(g(t),u)$ непрерывна на $\mT$
 для любого вектора $u$ из $\h$.
Отметим, что если $\h$-значная функция $g$ ограничена, то непрерывность
функции $t\mapsto(g(t),u)$ достаточно проверить только для векторов $u$ из подмножеств, линейная оболочка которых плотна
в $\h$.

Следующая теорема содержится в результатах работ \cite{KS3} и \cite{A4}.

\begin{thm}
\label{topolx}
Пусть $\Phi\in\fM(\mS\times\mT)$, где $\mS$ и $\mT$ топологические пространства.
Предположим, что функция $\Phi$ непрерывна по каждой переменной. Тогда существуют два семейства
$\{u_s\}_{s\in\mS}$ и $\{ v_t\}_{t\in\mT}$ в (не обязательно сепарабельном) гильбертовом пространстве
$\h$ такие, что

{\rm а)} линейная оболочка семейства $\{u_s\}_{s\in\mS}$ плотна в $\h$;

{\rm б)} линейная оболочка семейства $\{ v_t\}_{t\in\mT}$ плотна в $\h$;

{\rm в)} отображение $s\mapsto u_s$ слабо непрерывно;

{\rm г)} отображение $t\mapsto v_t$ слабо непрерывно;

{\rm д)} $\|u_s\|^2\le\|\Phi\|_{\fM(\mS\times\mT)}$ для всех $s\in\mS$;

{\rm е)} $\| v_t\|^2\le\|\Phi\|_{\fM(\mS\times\mT)}$ для всех $t\in\mT$;

{\rm ж)} $\f (s,t)=(u_s, v_t)$ для всех $(s,t)\in\mS\times\mT$.
\end{thm}

\Pf В силу теоремы \ref{Pi} и замечания к этой теореме найдутся 
семейства $\{u_s\}_{s\in\mS}$ и $\{ v_t\}_{t\in\mT}$,
удовлетворяющие условиям а), б), д), е) и ж).
Ясно, что функция $s\mapsto(u_s,h)$ непрерывна для $h= v_t$, где
$t\in\mT$. Следовательно, функция $s\mapsto(u_s,h)$ непрерывна при всех $h\in\h$ в силу б).
Таким образом, отображение $s\mapsto u_s$ слабо непрерывно. Аналогично из a) выводится слабая
непрерывность отображения $t\mapsto v_t$. $\bl$

\medskip

{\bf Замечание.}  Если хотя бы одно из пространств $\mS$ и $\mT$ сепарабельно, то и пространство $\h$
тоже сепарабельно.
Действительно, достаточно заметить, что если, например, пространство $\mS$ сепарабельно, то замыкание линейной
оболочки семейства $\{u_s\}_{s\in S}$ сепарабельно.

\medskip

Это замечание позволяет получить следующее утверждение.

\begin{thm}
\label{npkp}
Пусть $\Phi\in\fM(\mS\times \mT)$, где $ \mathcal S$ и $ \mathcal T$ --
топологические пространства, причём хотя бы одно из них сепарабельно.
Предположим, что функция $\Phi$ непрерывна по каждой переменной. 
Тогда существует последовательность непрерывных функций $\{\f_n\}_{n\ge1}$ на множестве $ \mathcal S$
 и последовательность непрерывных функций $\{\psi_n\}_{n\ge1}$ на множестве $ \mathcal T$ такие, что 
$$
\sum_{n=1}^\be|\f_n(s)|^2\le\|\Phi\|_{\fM(\mathcal \mathcal S\times \mathcal T)},\quad
\sum_{n=1}^\be|\psi_n(t)|^2\le\|\Phi\|_{\fM(\mathcal \mathcal S\times \mathcal T)}
$$
и
$$
\sum_{n=1}^\be\f_n(s)\psi_n(t)=\Phi(s,t)
\quad\mbox{при всех }~s\in\mathcal S~\mbox{~и~}~t\in\mathcal T.
$$
\end{thm}

\Pf  Пусть $\{u_s\}_{s\in \mathcal S}$ и $\{ v_t\}_{t\in \mathcal T}$ обозначают  два 
семейства в гильбертовом пространстве $\h$, существование которых утверждается
в теореме \ref{topolx}. Из замечания к этой теореме следует, что пространство $\h$
сепарабельно. Пусть \lb $\{e_n\}_{n=1}^N$ -- ортонормированный базис в $\h$, где $0\le N\le\be$. Остаётся положить $\f_n(s)\df(u_s,e_n)$ и
$\psi_n(t)\df(e_n,v_t)$;  при этом если $N<\be$, то $\f_n(s)\df\psi_n(t)\df0$
при $n>N$. $\bl$

\medskip

{\bf Определение.}
Отображение $g$, действующее из топологического пространства $\mT$ в гильбертово пространство $\h$, называется {\it слабо измеримым по Борелю}, если
функция $t\mapsto(g(t),u)$ измерима по Борелю на $T$ для любого вектора $u$ из
$\h$.

\medskip

Легко видеть, что условие измеримости по Борелю функции $t\mapsto(g(t),u)$
достаточно проверить только для векторов $u$ из подмножества 
гильбертова пространства $\h$, линейная оболочка которого плотна
в $\h$.

\begin{thm}
\label{topolbor}
Пусть $\mS$ и $\mT$ -- топологические пространства, а 
$\Phi\in\fM(\mS\times\mT)$.
Предположим, что функция $\Phi$ измерима по Борелю по каждой
переменной. Тогда существуют два семейства
$\{u_s\}_{s\in\mS}$ и $\{ v_t\}_{t\in\mT}$ в (не обязательно сепарабельном) гильбертовом пространстве
$\h$ такие, что

{\rm а)} линейная оболочка семейства $\{u_s\}_{s\in\mS}$ плотна в $\h$;

{\rm б)} линейная оболочка семейства $\{ v_t\}_{t\in\mT}$ плотна в $\h$;

{\rm в)} отображение $s\mapsto u_s$ измеримо по Борелю в слабом смысле;

{\rm г)} отображение $t\mapsto v_t$ измеримо по Борелю в слабом смысле;

{\rm д)} $\|u_s\|^2\le\|\Phi\|_{\fM(S\times T)}$ для всех $s\in S$;

{\rm е)} $\| v_t\|^2\le\|\Phi\|_{\fM(S\times T)}$ для всех $t\in T$;

{\rm ж)} $\Phi(s,t)=(u_s, v_t)$ для всех $(s,t)\in S\times T$.
\end{thm}

Доказательство теоремы \ref{topolbor}, по существу, слово в слово повторяет доказательство теоремы
\ref{topolx}  с той лишь разницей, что вместо непрерывности речь должна идти
об измеримости по Борелю.

%

\begin{thm}
\label{topolborizm}
Пусть $\mS$ и $\mT$ -- топологические пространства, а 
$\Phi$ -- борелевская функция из $\fM(\mS\times\mT)$. 
Предположим, что $\mu$   и $\nu$ -- борелевские
$\s$-конечные меры на $\mS$ и $\mT$. Тогда существуют последовательности
$\{\f_k\}_{k\ge1}$ и $\{\psi_k\}_{k\ge1}$ такие, что

{\rm а)} $\f_k\in L^\be(\mu)$ и $\psi_k\in L^\be(\nu)$ при всех $k\ge1$;

{\rm б)} $\sum_{k=1}^\be|\f_k(s)|^2\le\|\Phi\|_{\fM(S\times T)}$ при $\mu$-почти всех $s\in\mS$;

{\rm в)} $\sum_{k=1}^\be|\psi_k(t)|^2\le\|\Phi\|_{\fM(S\times T)}$ при $\nu$-почти всех $t\in\mT$;

{\rm г)} $\Phi(s,t)=\sum_{k=1}^\be\f_k(s)\psi_k(t)$  при $\mu\otimes\nu$-почти всех $(s,t)\in\mS\times\mT$.
\end{thm}

\Pf Ясно, что можно считать, что $\|\Phi\|_{\fM(S\times T)}=1$. Пусть 
$\{u_s\}_{s\in\mS}$ и $\{ v_t\}_{t\in\mT}$ -- семейства в гильбертовом пространстве $\h$,
существование которых утверждается в теореме \ref{topolbor}. Пусть $\{e_j\}_{j\in J}$ --
ортонормированный базис в гильбертовом пространстве $\h$. Положим
$\f_j(s)\df(u_s,e_j)$  и $\psi_j(t)\df(e_j,v_t)$. Тогда функции $\f_j$ и $\psi_j$ измеримы по Борелю,
$$
\sum\limits_{j\in J}|\f_j(s)|^2\le1\quad\mbox{при всех }~s\in\mS,
\quad\sum\limits_{j\in J}|\psi_j(t)|^2\le1\quad\mbox{при всех }~t\in\mT
$$
и 
$$
\Phi(s,t)=\sum_{j\in J}\f_j(s)\psi_j(t)\quad\mbox{при всех }~(s,t)\in\mS\times\mT.
$$

Отсюда мгновенно вытекает доказываемая теорема в случае, когда
множество $J$  не более, чем счётно. Чтобы рассмотреть случай произвольного
множества $J$,
положим $\Psi(s,t)\df\sum\limits_{j\in J}|\f_j(s)|\cdot|\psi_j(t)|$. По неравенству Коши--Буняковского
$$
\Psi(s,t)\le\left(\sum\limits_{j\in J}|\f_j(s)|^2\right)^{\frac12}\left(\sum\limits_{j\in J}|\psi_j(t)|^2\right)^{\frac12}
\le1.
$$
Можно считать, что меры $\mu$  и $\nu$ вероятностные. 
Положим $J_s\df\{j\in J: \f_j(s)\ne0\}$, где $s\in\mS$. Заметим,
что множество $J_s$ не более, чем счётно при каждом $s\in\mS$, поскольку $\sum_{j\in J}|\f_j(s)|^2\le1$.
Легко видеть, что для всех $s$ из $\mS$ 
\begin{align*}
\sum_{j\in J}|\f_j(s)|\int_\mT|\psi_j(t)|\,d\nu(t)&=\sum_{j\in J_s}|\f_j(s)|\int_\mT|\psi_j(t)|\,d\nu(t)\\[.2cm]
&=\int_\mT\Big(\sum_{j\in J_s}|\f_j(s)|\cdot|\psi_j(t)|\Big)\,d\nu(t)=\int_\mT\Psi(s,t)\,d\nu(t).
\end{align*}
Чтобы проинтегрировать теперь по $s$, мы введём не более, чем счётное множество
$J_\flat\df\Big\{j\in J:\int_T|\psi_j(t)|\,d\nu(t)\ne0\Big\}$. Тогда
\begin{align*}
\sum_{j\in J}\int_\mS|\f_j(s)|\,d\mu(s)\int_\mT|\psi_j(t)|\,d\nu(t)&=
\sum_{j\in J_\flat}\int_\mS|\f_j(s)|\,d\mu(s)\int_\mT|\psi_j(t)|\,d\nu(t)\\[.2cm]
&=
\int_\mS\left(\int_\mT\Psi(s,t)\,d\nu(t)\right)\,d\mu(s).
\end{align*}
Теперь ясно, что 
$$
\int_\mS\left(\int_\mT\sum_{j\in J\setminus J_\flat}|\f_j(s)|\cdot|\psi_j(t)|\,d\nu(t)\right)\,d\mu(s)=0.
$$
Отсюда и из неравенства
$$
\left|\Phi(s,t)-\sum_{j\in J_\flat}\f_j(s)\psi_j(t)\right|\le\sum_{j\in J\setminus J_\flat}|\f_j(s)|\cdot|\psi_j(t)|
$$
вытекает, что
$$
\sum_{j\in J_\flat}\f_j(s)\psi_j(t)=\Phi(s,t)
$$
для $\mu\otimes\nu$-почти всех $(s,t)\in\mS\times\mT$. $\bl$

\medskip

Рассмотрим теперь некоторые примеры мультипликаторов Шура.
Пусть $\mM(\T^2)$ обозначает пространство всех комплексных борелевских мер на
двухмерном торе $\T^2$ с нормой $\|\mu\|_{\mM(\T^2)}\df|\mu|(\T)$.

\medskip

{\bf Пример.}  Пусть $\mu\in\mM(\T^2)$.
Тогда $\{\widehat\mu(m,n)\}\in\fM(\Z^2)$ и $\|\widehat\mu\|_{\fM(\Z\times\Z)}\le\|\mu\|$.

\medskip

Это утверждение является очевидным следствием теоремы \ref{41}. Ясно, что далеко
не каждый мультипликатор Шура ${\bs a}\in\fM(\Z\times\Z)$ представим в виде
${\bs a}=\widehat\mu$, где $\mu\in\mM(\T^2)$.
Рассмотрим, например, случай, когда матрица
${\bs a}=\{a_{mn}\}_{m,n\in\Z}$ состоит из одинаковых столбцов (или строк). Пусть для определённости
$a_{mn}=t_n$ при всех $m,n\in\Z$. Тогда
${\bs a}\in\fM(\Z\times\Z)$ в том и только в том случае, когда
${\bs a}\in\ell^\be(\Z\times\Z)$ и $\|{\bs a}\|_{\fM(\Z\times\Z)}=\|{\bs a}\|_{\ell^\be(\Z\times\Z)}=\|\{t_n\}\|_{\ell^\be}$.
Разумеется, отнюдь не каждая такая матрица ${\bs a}$ с ограниченными коэффициентами
представима в виде ${\bs a}=\widehat\mu$,  где $\mu\in\mM(\T^2)$.

С другой стороны, если предположить, что матрица ${\bs a}=\{a_{mn}\}_{m,n\in\Z}$
является матрицей {\it Лорана}, т. е. $a_{mn}=t_{m-n}$, то ситуация резко изменится.

\begin{thm}
\label{Boch}
Пусть $A=\{a_{mn}\}_{m,n\in\Z}$ -- матрица Лорана.
Тогда $A\in\fM(\Z^2)$ в том и только в том случае, когда  $a_{mn}=\widehat\mu(m-n)$
для  некоторой меры $\mu$ из $\mM(\T)$;
при этом $\|{\bs a}\|_{\fM(\Z^2)}=\|\mu\|_{\mM(\T)}$.
\end{thm}

Все эти результаты переносятся на локально компактные абелевы группы. При этом в случае недискретной абелевой группы $G$ в формулировке аналога
теоремы \ref{Boch} нужно потребовать непрерывность соответствующих функций.

\begin{thm}
\label{BochR}
Пусть $h$ - непрерывная функция на $\R$.
Тогда матрица $A=\{h(s-t)\}_{s,t\in\R}$ принадлежит пространству
$\fM(\R\times\R)$ в том и только в том случае, когда существует комплексная борелевская мера $\mu$ такая, что $h=\F\mu$. При этом
$\|A\|_{\fM(\R\times\R)}=\|\mu\|_{\mM(\R)}$.
\end{thm}

%

\

\section{\bf Двойные операторные интегралы}
\setcounter{equation}{0}
\label{mSid}

\

Двойные операторные интегралы -- это выражения вида
\bay
\label{doi}
\int_\mS\int_\mT\Phi(s,t)\,dE_1(s)T\,dE_2(t),
\ey
где $E_1$ и $E_2$ -- спектральные меры в сепарабельном гильбертовом пространстве $\h$, 
$\Phi$ -- ограниченная измеримая функция, а $T$ - линейный непрерывный оператор в $\h$.

Двойные операторные интегралы появились в работе \cite{DK}. В работах Бирмана и Соломяка
\cite{BS1}, \cite{BS2} и \cite{BS3} была создана стройная теория двойных операторных интегралов. Идея Бирмана и Соломяка состояла в том, чтобы сначала определить двойные операторные интегралы вида \rf{doi} для произвольных ограниченных измеримых функций $\Phi$ и для операторов $T$ класса Гильберта--Шмидта $\bS_2$. Для этого рассматривается спектральная мера $\E$, которая принимает значение в множестве ортогональных проекторов 
в гильбертовом пространстве $\bS_2$ и определяется равенством
$$
\E(\L\times\D)T=E_1(\L)TE_2(\D),\quad T\in\bS_2,
$$
где $\L$ и $\D$ --измеримые подмножества множеств $\mS$ и $\mT$. 
Ясно, что умножение слева на $E_1(\L)$ коммутирует с умножением справа на
$E_2(\D)$. В работе \cite{BS4}  было доказано, что $\E$ продолжается до спектральной меры
на $\mS\times\mT$. В этой ситуации двойной операторный интеграл \rf{doi} определяется равенством 
$$
\int_\mS\int_\mT\Phi(s,t)\,d E_1(s)T\,dE_2(t)\df
\left(\,\,\int_{\mS\times\mT}\Phi\,d\E\right)T.
$$
Из этого определения сразу же вытекает, что
$$
\left\|\int_\mS\int_\mT\Phi(s,t)\,dE_1(s)T\,dE_2(t)\right\|_{\bS_2}
\le\|\Phi\|_{L^\be}\|T\|_{\bS_2}.
$$

Если функция $\Phi$ обладает свойством
$$
T\in\bS_1\quad\Longrightarrow\quad\int_\mS\int_\mT\Phi(s,t)\,d E_1(s)T\,dE_2(t)\in\bS_1,
$$
то говорят, что $\Phi$ -- {\it мультипликатор Шура пространства $\bS_1$ относительно
спектральных мер $E_1$ и $E_2$}.

Чтобы определить двойные операторные интегралы \rf{doi} для ограниченных операторов $T$,
рассмотрим трансформатор
$$
Q\mapsto\int_{\mT}\int_{\mS}\Phi(t,s)\,d E_2(t)\,Q\,dE_1(s),\quad Q\in\bS_1,
$$
и предположим, что функция $(y,x)\mapsto\Phi(y,x)$ -- мультипликатор Шура пространства 
$\bS_1$ относительно $E_2$ и $E_1$.
В этом случае трансформатор
\bay
\label{tra}
T\mapsto\int_\mS\int_\mT\Phi(s,t)\,d E_1(s)T\,dE_2(t),\quad T\in \bS_2,
\ey
продолжается по двойственности до ограниченного линейного трансформатора в пространстве ограниченных линейных операторов в $\h$. При этом говорят, что $\Phi$ -- {\it мультипликатор Шура (относительно $E_1$ и $E_2$) пространства ограниченных линейных операторов}.
Обозначим пространство таких мультипликаторов Шура символом $\fM(E_1,E_2)$.
Норма функции $\Phi$ в $\fM(E_1,E_2)$ определяется, как норма трансформатора
\rf{tra} в пространстве ограниченных линейных операторов.

Легко видеть, что если функция $\Phi$ на $\mS\times\mT$ входит в {\it проективное тензорное произведение} $L^\be(E_1)\hat\otimes L^\be(E_2)$ пространств $L^\be(E_1)$ 
и $L^\be(E_2)$ (т.е. $\Phi$ допускает представление
$$
\Phi(s,t)=\sum_{n\ge0}\f_n(s)\psi_n(t),\quad\mbox{где}\quad
\sum_{n\ge0}\|\f_n\|_{L^\be(E_1)}\|\psi_n\|_{L^\be(E_2)}<\be),
$$
то $\Phi\in\fM(E_1,E_2)$.
Для таких функций $\Phi$ имеем
$$
\int_\mS\int_\mT\Phi(s,t)\,dE_1(s)T\,dE_2(t)=
\sum_{n\ge0}\left(\,\int_\mS\f_n\,dE_1\right)T\left(\,\int_\mT\psi_n\,dE_2\right).
$$

Более общо, $\Phi\in\fM(E_1,E_2)$, если функция $\Phi$
принадлежит {\it интегральному проективному тензорному произведению} $L^\be(E_1)\hat\otimes_{\rm i}L^\be(E_2)$ пространств $L^\be(E_1)$ и $L^\be(E_2)$, т.е. функция $\Phi$ допускает представление
\bay
\label{ipt}
\Phi(s,t)=\int_\O \f(s,w)\psi(t,w)\,d\s(w),
\ey
где $(\O,\s)$ -- $\s$-конечное пространство с мерой, $\f$ -- измеримая функция на $\mS\times \O$,\lb
$\psi$ -- измеримая функция на $\mT\times \O$ и
$$
\int_\O\|\f(\cdot,w)\|_{L^\be(E_1)}\|\psi(\cdot,w)\|_{L^\be(E_2)}\,d\s(w)<\be.
$$

Оказывается, что все мультипликаторы Шура могут быть получены таким образом (см. теорему \ref{tomSc}).

Ещё одно достаточное условие для того, чтобы функция была мультипликатором Шура, может быть сформулировано в терминах тензорного произведения Хогерупа  $L^\be(E_1)\!\otimes_{\rm h}\!L^\be(E_2)$, которое определяется, как пространство функций $\Phi$ вида
\bay
\label{FiH}
\Phi(s,t)=\sum_{n\ge0}\f_n(s)\psi_n(t),
\ey
где $\{\f_n\}_{n\ge0}\in L_{E_1}^\be(\ell^2)$ и $\{\psi_n\}_{n\ge0}\in L_{E_2}^\be(\ell^2)$.
При этом
$$
\|\Phi\|_{L^\be(E_1)\otimes_{\rm h}\!L^\be(E_2)}\df\inf
\Big\|\sum_{n\ge0}|\f_n|^2\Big\|_{L^\be(E_1)}^{1/2}
\Big\|\sum_{n\ge0}|\psi_n|^2\Big\|_{L^\be(E_2)}^{1/2},
$$
где инфимум берётся
по всем представлениям функции $\Phi$ вида \rf{FiH}. 
Нетрудно проверить, что если $\Phi\in L^\be(E_1)\!\otimes_{\rm h}\!L^\be(E_2)$, то $\Phi\in\fM(E_1,E_2)$ и
\bay
\label{Haagrazl}
\iint\Phi(s,t)\,dE_1(s)T\,dE_2(t)=
\sum_{n\ge0}\Big(\int\f_n\,dE_1\Big)T\Big(\int\psi_n\,dE_2\Big),
\ey
причём ряд в правой части сходится в слабой операторной топологии и
$$
\|\Phi\|_{\fM(E_1,E_2)}\le\|\Phi\|_{L^\be(E_1)\otimes_{\rm h}L^\be(E_2)}.
$$

Как видно из следующей теоремы, условие
$\Phi\in L^\be(E_1)\!\otimes_{\rm h}\!L^\be(E_2)$ не только достаточно, но и необходимо. 

\begin{thm}
\label{tomSc} 
Пусть $\Phi$ -- измеримая функция на
$\mS\times\mT$ и пусть $\mu$ и $\nu$ -- положительные $\s$-конечные меры на $\mS$ и $\mT$, которые взаимно абсолютно непрерывны с мерами $E_1$ и $E_2$. Следующие утверждения эквивалентны:

{\rm (а)} $\Phi\in\fM(E_1,E_2)$;

{\rm (б)} $\Phi\in L^\be(E_1)\hat\otimes_{\rm i}L^\be(E_2)$;

{\rm (в)} $\Phi\in L^\be(E_1)\!\otimes_{\rm h}\!L^\be(E_2)$;

{\rm (г)} существуют $\s$-конечная мера $\s$ на множестве $\O$, измеримые функции $\f$ на $\mS\times\O$ и $\psi$ на $\mT\times\O$ такие, что имеет место
{\em\rf{ipt}} и
\bay
\label{bs}
\left\|\left(\int_\O|\f(\cdot,w)|^2\,d\s(w)\right)^{1/2}\right\|_{L^\be(E_1)}
\left\|\left(\int_\O|\psi(\cdot,w)|^2\,d\s(w)\right)^{1/2}\right\|_{L^\be(E_2)}<\be;
\ey

{\rm (д)} если интегральный оператор
$f\mapsto\int k(x,y)f(y)\,d\nu(y)$ из $L^2(\nu)$ в $L^2(\mu)$ входит в класс $\bS_1$, то в этот же класс входит и интегральный оператор
\lb$f\mapsto\int\Psi(x,y)k(x,y)f(y)\,d\nu(y)$.
\end{thm}

Импликации (г)$\imp$(а)$\eq$(д) установлены в \cite{BS3}.
В случае матричных мультипликаторов Шура импликация (а)$\imp$(б) была доказана в работе \cite{Ben}. Мы отсылаем читателя к работе \cite{Pe1} за доказательством эквивалентности утверждений (а), (б) и (г). Эквивалентность (в) и (г) доказывается элементарно.

Легко видеть, что условия (а) - (д) эквивалентны и тому, что функция $\Phi$ является мультипликатором Шура класса $\bS_1$. 

Отметим также, что с таким же успехом можно определить двойные операторные интегралы вида \rf{doi}, в которых $E_1$ и $E_2$ -- спектральные меры в разных гильбертовых пространствах, а $T$ -- оператор из одного пространства в другое.

\medskip

{\bf Замечание.} Из теорем \ref{topolborizm} и \ref{tomSc} с лёгкостью вытекает, что если $\mS$ и $\mT$ -- топологические пространства, а $\Phi$ -- борелевская функция на $\mS\times\mT$ класса $\fM(\mS\times\mT)$ (т.е. $\Phi$ -- дискретный мультипликатор Шура), то $\Phi\in\fM(E_1,E_2)$ для любых борелевских спектральных мер $E_1$ и $E_2$ на $\mS$ и $\mT$.

\medskip

%
%
%
%

Двойные операторные интегралы также могут быть определены и по отношению к полуспектральным мерам. Напомним, что {\it полуспектральная мера} $\mE$ 
на измеримом пространстве $(\X,{\frak B})$ -- это отображение, заданное на $\s$-алгебре ${\frak B}$, которое принимает значения в множестве ограниченных операторов в гильбертовом пространстве $\h$, счётно аддитивное в сильной операторной топологии и такое, что
$$
\mE(\D)\ge\0,\quad\D\in{\frak B},\quad\mE(\varnothing)=\0
\quad\mbox{и}\quad\mE(\X)=I.
$$
По теореме Наймарка \cite{Nai} всякая полуспектральная мера $\mE$ имеет 
{\it спектральную дилатацию}, т.е. спектральную меру $E$, заданную на том же измеримом
пространстве $(\X,{\frak B})$, принимающую значения в множестве ортогональных проекторов в гильбертовом пространстве $\K$, содержащем $\h$, и такую, что
$$
\mE(\D)=P_\h E(\D)\big|\h,\quad\D\in{\frak B},
$$
где $P_\h$ -- ортогональный проектор в $\K$ на подпространство $\h$. Такую спектральную дилатацию можно выбрать {\it минимальной} в том смысле, что 
$$
\K=\clos\spn\{E(\D)\h:~\D\in{\frak B}\}.
$$ 
Отметим здесь, что в работе \cite{MM} показано, что если $E$ -- минимальная  спектральная дилатация полуспектральной меры $\E$, то $E$ и $\E$ взаимно абсолютно непрерывны.

Интегралы по полуспектральным мерам определяются следующим образом:
$$
\int_\X \f(x)\,d\mE(x)=P_\h\left(\int_\X \f(x)\,d E(x)\right)\Big|\h,\quad
\f\in L^\be(\mE)\df L^\be(E).
$$

Если $(\X_1,{\frak B}_1)$ и $(\X_2,{\frak B}_2)$ -- полуспектральные меры, $E_1$
и $E_2$ -- их спектральные дилатации в гильбертовом пространствах $\K_1$ и $\K_2$, а функция 
$\Phi$ на $\X_1\times\X_2$ удовлетворяет одному из условий теоремы \ref{tomSc}, то двойной операторный интеграл по полуспектральным мерам $\mE_1$ и $\mE_2$ определяется равенством 
$$
\int_{\!\X_1}\!\!\int_{\!\X_2}\Phi(x_1,x_2)\,d\E_1(x_1)Q\,d\E_2(x_2)=
P^{[1]}_\h\int_{\!\X_1}\!\!\int_{\!\X_2}\Phi(x_1,x_2)\,dE_1(x_1)\big(QP^{[2]}_{\h}\big)\,dE_2(x_2)\Big|\h
$$
для произвольного ограниченного оператора $Q$ в $\h$. 
Здесь $P^{[1]}_\h$ и $P^{[2]}_\h$ -- ортогональные проекторы из $\K_1$ и $\K_2$
на $\h$.
Если $\Phi\in L^\be(E_1)\!\otimes_{\rm h}\!L^\be(E_2)$, то
имеет место равенство
\bay
\label{doipolcpm}
\iint\Phi(x_1,x_2)\,d\mE_1(x_1)T\,d\mE_2(x_2)=
\sum_{n\ge0}\Big(\int\f_n\,d\mE_1\Big)T\Big(\int\psi_n\,d\mE_2\Big),
\ey
где $T\in\mB(\h)$, а $\f_n$ и $\psi_n$ -- функции из представления \rf{FiH}.

Двойные операторные интегралы по полуспектральным мерам были введены в работах \cite{Pe*}, см. также \cite{Pe9}.

\

\begin{center}
\bf\Large Глава III
\end{center}

\medskip

\begin{center}
\bf\Large 
Операторно липшицевы функции \\ на подмножествах плоскости
\end{center}
\label{OL2per}

\addtocontents{toc}{\vspace*{.3cm}\hspace*{-.35cm}\textbf{{\bf Глава III}. Операторно липшицевы функции на подмножествах плоскости}\hfill\pageref{OL2per}}

\renewcommand{\thesection}{3.\arabic{section}}
\setcounter{section}{0}

\

В этой главе мы изучаем операторно липшицевы и коммутаторно липшицевы функции на замкнутых подмножествах комплексной плоскости. При этом значительную роль играют мультипликаторы Шура. Мы предлагаем два метода разностных и коммутаторных оценок. Первый метод использует дискретные мультипликаторы Шура и аппроксимацию операторами с конечным спектром. Второй метод основан на применении двойных операторных интегралов.

\

\section{\bf Операторно липшицевы и коммутаторно липшицевы функции на замкнутых подмножествах плоскости}
\setcounter{equation}{0}
\label{oplip}

\

Мы рассмотрим в этом параграфе классы операторно липшицевых функций и коммутаторно липшицевых функций на замкнутых подмножествах плоскости. Мы увидим, что в отличие от случая функций на прямой и на окружности, эти два класса отнюдь не обязаны совпадать. При определении этих классов будем здесь рассматривать только ограниченные операторы. В следующем параграфе мы увидим, что рассмотрение не обязательно ограниченных операторов приводит к тем же классам функций.

Пусть $\fF$ -- непустое подмножество комплексной плоскости $\C$. Обозначим
символом $\Li(\fF)$  пространство функций $f:\fF\to\C$,
удовлетворяющих {\it условию Липшица}:
\bay
\label{lipsc}
|f(z)-f(w)|\le C|z-w|,\quad z,w\in\C.
\ey
Наименьшую из констант $C\ge0$, удовлетворяющих
условию \rf{lipsc}, обозначим через $\|f\|_{\Li(\fF)}=\|f\|_{\Li}$.
Положим $\|f\|_{\Li}\df\be$, если $f\not\in\Li(\fF)$.

Обычно мы будем требовать замкнутость множества $\fF$.

Из спектральной теоремы для пары {\it коммутирующих} нормальных операторов легко следует, что неравенство
\bay
\label{1comlip}
\|f(N_1)-f(N_2)\|\le \|f\|_{\Li(\fF)}\|N_1-N_2\|
\ey
выполняется для любых коммутирующих нормальных операторов $N_1$ и $N_2$, спектры которых
содержатся в $\fF$.

Комплексную непрерывную функцию $f$, заданную на непустом замкнутом множестве $\fF$,
$\fF\subset\C$,
будем называть {\it операторно липшицевой}, если существует положительное число
$C$ такое, что
\bay
\label{lipid}
\|f(N_1)-f(N_2)\|\le C\|N_1-N_2\|
\ey
для любых нормальных операторов $N_1$ и $N_2$
со спектрами, лежащими в $\fF$.
Пространство всех операторно липшицевых функций
на $\fF$ обозначим символом $\OL(\fF)$. Наименьшую из констант $C$, удовлетворяющих
условию \rf{lipid},
обозначим через $\|f\|_{\OL(\fF)}=\|f\|_{\OL}$. Положим $\|f\|_{\OL}=\be$, если $f\not\in\OL(\fF)$.

Если функция $f$ задана на более широком множестве ${\frak G}\supset\fF$, то для краткости мы будем
обычно писать $f\in\OL(\fF)$ и $\|f\|_{\OL(\fF)}$ вместо $f|\fF\in\OL(\fF)$ и $\|f|\fF\|_{\OL(\fF)}$.
Это же соглашение будет применяться и для других функциональных пространств.

Легко видеть, что $\OL(\fF)\subset\Li(\fF)$ и $\|f\|_{\Li(\fF)}\le\|f\|_{\OL(\fF)}$
для любой функции $f\in\OL(\fF)$.  Мы увидим в \S\:\ref{orav}, что равенство $\OL(\fF)=\Li(\fF)$
выполняется только для конечных множеств $\fF$.

Заметим, что если $f\in\OL(\fF)$ и $\|f\|_{\OL}\le1$,  то
\bay
\label{unit}
\|f(N_1)U-Uf(N_2)\|\le\|N_1U-UN_2\|
\ey
для всех унитарных операторов $U$ и всех
нормальных операторов $N_1$  и $N_2$ таких, что $\s(N_1),\s(N_2)\subset\fF$.
Чтобы убедиться в этом, достаточно применить неравенство \rf{lipid}
для $C=1$ к нормальным операторам $U^*N_1U$ и $N_2$. Обратно, если неравенство \rf{unit}
выполняется для всех унитарных операторов $U$ и всех
нормальных операторов $N_1$  и $N_2$ таких, что $N_1=N_2$ и $\s(N_1)=\s(N_2)\subset\fF$, то
$f\in\OL(\fF)$ и $\|f\|_{\OL(\fF)}\le1$.
Действительно, применяя неравенство \rf{unit} к операторам
$$
\mathcal N_1=\mathcal N_2= \left(\begin{array}{cc}
    N_1 & \0 \\
    \0 & N_2 \\
  \end{array}\right)
\quad\text{и}\quad
\mathcal U= \left(\begin{array}{cc}
    \0 & I \\
    I & \0 \\
  \end{array}\right),
$$
получаем:
$$
\|f(N_1)-f(N_2)\|\le\|N_1-N_2\|.
$$
Заметим, что в этом рассуждении мы обошлись только самосопряжёнными унитарными операторами
$\mathcal U$, т. е. нормальными операторами $\mathcal U$ такими, что $\mathcal U^2=I$
или, что то же самое, унитарными операторами со спектрами в  $\{-1,1\}$.

Справедлива следующая теорема:

\begin{thm}
\label{Nlip0}
Пусть $f$ -- непрерывная функция на замкнутом подмножестве $\fF$ комплексной плоскости $\C$.
Тогда следующие утверждения равносильны:

{\em(а)} $\|f(N_1)-f(N_2)\|\le\|N_1-N_2\|$ для всех нормальных операторов $N_1$ и $N_2$ таких, что
$\s(N_1),\s(N_2)\subset\fF$;

{\em(б)} $\|f(N_1)U-Uf(N_2)\|\le\|N_1U-UN_2\|$ для всех унитарных операторов $U$
и всех  нормальных операторов $N_1$ и $N_2$ таких, что $\s(N_1),\s(N_2)\subset\fF$;

{\em(в)} $\|f(N)U-Uf(N)\|\le\|NU-UN\|$ для всех самосопряжённых унитарных операторов $U$ и всех
нормальных операторов $N$ таких, что $\s(N)\subset\fF$;

{\em(г)} $\|f(N)A-Af(N)\|\le\|NA-AN\|$ для всех самосопряжённых операторов $A$ и всех
нормальных операторов $N$ таких, что $\s(N)\subset\fF$.
\end{thm}

\Pf Эквивалентность утверждений (а), (б) и (в), по существу, была доказана
перед формулировкой теоремы. Импликация (г)$\Longrightarrow$(в) очевидна.
Остаётся доказать, что (в) влечёт (г). Обозначим через $\frak X$ множество
всех операторов $R$ таких, что $\|f(N)R-Rf(N)\|\le\|NR-RN\|$ для всех
нормальных операторов $N$ со спектром в $\fF$. Ясно, что множество $\frak X$ замкнуто по норме и
$\a U+\b I\in\frak X$
для всех унитарных операторов $U$ и всех $\a,\b\in\C$. Теперь, чтобы доказать,
что любой самосопряжённый оператор $A$ принадлежит множеству $\frak X$,
достаточно заметить, что оператор \lb$(I-\e{\rm i}A)(I+\e{\rm i}A)^{-1}$ является
унитарным при всех  $\e$ из $(-\|A\|^{-1},\|A\|^{-1})$ и
$$
A=\lim_{\e\to0}\frac1{2\e{\rm i}}\left(I-(I-\e{\rm i}A)(I+\e{\rm i}A)^{-1}\right).\quad\bl
$$

\medskip

{\bf Замечание.}  Унитарный оператор $U$ является самосопряжённым
в том и только в том случае, когда он представим в виде $U=2P-I$, где $P$ -- ортогональный проектор.
Таким образом, утверждение (в) теоремы \ref{Nlip0} можно переписать
следующим образом:

{\it$\|f(N)P-Pf(N)\|\le\|NP-PN\|$ для всех ортогональных проекторов $P$ и всех
нормальных операторов $N$ таких, что $\s(N)\subset\fF$}.

\medskip

Имеет место также следующее утверждение:

\begin{thm}
\label{Nlip}
Пусть $f$ -- непрерывная функция, заданная на замкнутом подмножестве $\fF$ комплексной плоскости $\C$.
Тогда следующие утверждения равносильны:

{\em(а)} $\|f(N_1)-f(N_2)\|\le\|N_1-N_2\|$ для всех нормальных операторов $N_1$ и $N_2$ таких, что
$\s(N_1),\s(N_2)\subset\fF$;

{\em(б)} $\|f(N)R-Rf(N)\|\le\max\big\{\|NR-RN\|,\|N^*R-RN^*\|\big\}$ для всех операторов
$R\in\mB(\h)$ и всех
нормальных операторов $N$ таких, что $\s(N)\subset\fF$;

{\em(в)} $\|f(N_1)R-Rf(N_2)\|\le\max\big\{\|N_1R-RN_2\|,\|N_1^*R-RN_2^*\|\big\}$ для всех операторов
$R\in\mB(\h)$
и всех  нормальных операторов $N_1$ и $N_2$ таких, что $\s(N_1),\s(N_2)\subset\fF$.
\end{thm}

\Pf Докажем сначала импликацию (а)$\Longrightarrow$(б). Пусть выполнено
утверждение (а). Тогда из теоремы \ref{Nlip0} следует, что $\|f(N)A-Af(N)\|\le\|NA-AN\|$ для всех самосопряжённых операторов $A$ и всех
нормальных операторов $N$ таких, что $\s(N)\subset\fF$. Применяя это утверждение
к нормальному оператору
$\left(
  \begin{array}{cc}
     N & \0 \\
     \0 & N \\
  \end{array}
\right)$
и самосопряжённому оператору
$\left(
  \begin{array}{cc}
     \0 & R \\
      R^* & \0 \\
  \end{array}
\right)$, получим
$$
\max\big\{\|f(N)R-Rf(N)\|,\|f(N)R^*-R^*f(N)\|\big\}
\le\max\big\{\|NR-RN)\|,\|NR^*-R^*N\|\big\},
$$
откуда следует (б).

Применяя (б) к нормальному оператору
$\left(
  \begin{array}{cc}
     N_1 & \0 \\
     \0 & N_2 \\
  \end{array}
\right)$
и оператору
$\left(
  \begin{array}{cc}
     \0 & R \\
     \0 & \0 \\
  \end{array}
\right)$, получим (в). 
Импликация (в)$\Longrightarrow$(а) очевидна. $\bl$

В доказательствах теорем   \ref{Nlip0} и \ref{Nlip} мы использовали стандартный приём, связанный
с переходом к блочным операторам, который позволяет в некоторых случаях переходить 
от одного оператора к паре операторов.  Этот приём будет нам полезен и в дальнейшем.
Киттанех \cite{Ki} называет этот приём трюком Бербериана, имея
в виду, видимо, статью Бербериана \cite{Ber}.

Теоремы \ref{Nlip0} и \ref{Nlip} содержатся
в теореме 3.1 статьи \cite{AP6}, но по существу в некотором виде их можно извлечь из статьи \cite{KS2}, где наряду
с операторной нормой рассматриваются также произвольные симметричные нормы.

%
%
%

Заметим, что равенство
$$
\|N_1^*R-RN_2^*\|=\|N_1R-RN_2\|,
$$
а значит, и равенство
$$
\max\big(\|N_1R-RN_2\|,\|N_1^*R-RN_2^*\|\big)=\|N_1R-RN_2\|,
$$
выполняется
в каждом из следующих частных случаев:

1) операторы $N_1$ и $N_2$ являются самосопряжёнными (это так, если $\fF\subset\R$);

2) операторы $N_1$ и $N_2$ являются унитарными (это так, если $\fF\subset\T$);

3) $R$ -- самосопряжённый оператор и $N_1=N_2$;

4) $R$ -- унитарный оператор.

Комплексную функций $f$,
непрерывную на замкнутом  множестве $\fF$, $\fF\subset\C$, будем называть
{ \it коммутаторно липшицевой} функцией, если существует число $C\ge0$
такое, что
\bay
\label{cld}
\|f(N)R-Rf(N)\|\le C\|NR-RN\|
\ey
для всех $R\in\mB(\h)$ и всех нормальных операторов $N$ со спектром в $\fF$.
Множество всех коммутаторно липшицевых функций,
заданных на $\fF$, обозначим символом $\CL(\fF)$. Наименьшую из констант $C$, удовлетворяющих
условию \rf{cld},
обозначим через $\|f\|_{\CL(\fF)}=\|f\|_{\CL}$. Положим $\|f\|_{\CL(\fF)}=\be$, если $f\not\in\CL(\fF)$.

\begin{thm}
\label{Clip}
Пусть $f$ -- непрерывная функция, заданная на замкнутом подмножестве $\fF$ комплексной плоскости $\C$.
Тогда следующие три утверждения равносильны:

{\em(а)} $\|f(N)R-Rf(N)\|\le\|NR-RN\|$ для всех операторов $R\in\mB(\h)$ и всех
нормальных операторов $N$ таких, что $\s(N)\subset\fF$;

{\em(б)} $\|f(N_1)R-Rf(N_2)\|\le\|N_1R-RN_2\|$ для всех операторов $R\in\mB(\h)$ и всех нормальных операторов
$N_1$ и $N_2$ таких, что
$\s(N_1),\:\s(N_2)\subset\fF$;

{\em(в)} $\|f(N_1)A-Af(N_2)\|\le\|N_1A-AN_2\|$ для всех самосопряжённых операторов $A$
и всех  нормальных операторов $N_1$ и $N_2$ таких, что $\s(N_1),\:\s(N_2)\subset\fF$.
\end{thm}

\Pf Чтобы доказать импликацию (а)$\Longrightarrow$(б), достаточно применить (а)
к нормальному оператору
$\left(
  \begin{array}{cc}
     N_1 & \0 \\
     \0 & N_2 \\
  \end{array}
\right)$
и оператору
$\left(
  \begin{array}{cc}
     \0 & R \\
     \0 & \0 \\
  \end{array}
\right)$.
Импликация (б)$\Longrightarrow$(в) очевидна.
Остаётся доказать, что (в) влечёт (а).
 Применяя (в) к $N_1=U^{*}NU$ и $N_2=N$, где $U$ -- унитарный оператор, получаем:
 $$
 \|f(N)UA-UAf(N)\|=\|f(U^{*}NU)A-Af(N)\|\le\|NUA-UAN\|
 $$
 для любого самосопряжённого оператора $A$, унитарного оператора $U$
 и нормального оператора $N$ такого, что $\s(N)\subset\fF$.
 Заметим, что если условие (а) выполняется для оператора $R\in\mB(\h)$,
 то оно будет выполняться и для оператора $R+\l I$, где $\l\in\C$. Таким
 образом, можно считать, что оператор $R$ обратим. Тогда, применяя
 полярное разложение к обратимому оператору $R$, получаем:
 $R=UA$, где $U$ -- унитарный оператор, $A$ -- (положительный) самосопряжённый
 оператор. $\bl$

Из теорем \ref{Nlip0} и \ref{Clip} мгновенно следует, что
$\CL(\fF)\subset\OL(\fF)$ и $\|f\|_{\OL(\fF)}\le\|f\|_{\CL(\fF)}$
для всех $f$ из $\CL(\fF)$.

\medskip

{\bf Замечание.} 
%
%
В утверждениях предыдущих теорем этого параграфа, в которых
участвуют нормальные операторы $N_1, N_2$ и ещё один оператор $U$ или $R$
(другими словами, речь идёт об утверждениях  (б) теоремы \ref{Nlip0} и  (в) теоремы \ref{Nlip}, 
(б) теоремы \ref{Clip}), 
можно считать, что нормальные операторы $N_1$ и $N_2$ действуют
в разных гильбертовых пространствах (при этом там, где участвует оператор 
$U$,  следует допустить, чтобы унитарный оператор действовал из одного гильбертова
пространства в другое).   
Это видно из доказательств этих теорем.
В качестве иллюстрации мы приведём здесь соответствующую переформулировку
утверждения (б) теоремы  \ref{Clip}:

\medskip

{\em $\|f(N_1)R-Rf(N_2)\|\le\|N_1R-RN_2\|$ для всех операторов  $R\in\mB(\mathcal \h_2,\h_1)$  и
всех нормальных операторов $N_1$ и $N_2$,
действующих соответственно в гильбертовых пространствах $\h_1$
и $\h_2$, и удовлетворяющих условию: $\s(N_1),\:\s(N_2)\subset\fF$.}

\medskip

Аналоги теорем  \ref{Nlip0}, \ref{Nlip}  и \ref{Clip} с соответствующими изменениями имеют место
по существу с теми же доказательствами для симметрично-нормированных идеалов. 
Мы рассмотрим здесь только идеал ядерных операторов $\bS_1$.
С каждым замкнутым множеством $\fF$, $\fF\subset\C$, мы связываем пространство
{\it ядерно липшицевых} (или {\it$\bS_1$-липшицевых}) функций $\OL_{\bS_1}(\fF)$ и 
{\it ядерно-коммутаторно липшицевых}
(или {\it$\bS_1$-коммутаторно липшицевых}) функций 
$\CL_{\bS_1}(\fF)$. Чтобы получить определение пространств $\OL_{\bS_1}(\fF)$
и $\CL_{\bS_1}(\fF)$, 
нужно только в неравенствах \rf{lipid} и \rf{cld}  операторную норму заменить  ядерной нормой.

Соответствующий ''объединённый'' аналог теорем   \ref{Nlip0} и \ref{Nlip} для ядерной нормы
может быть 
сформулирован следующим образом.

\begin{thm}
\label{yaNlip01}
Пусть $f$ -- непрерывная функция на замкнутом подмножестве $\fF$ комплексной плоскости $\C$.
Тогда следующие утверждения равносильны:

{\em(а)} $\|f(N_1)-f(N_2)\|_{\bS_1}\le\|N_1-N_2\|_{\bS_1}$ для всех нормальных операторов $N_1$ и $N_2$ таких, что
$\s(N_1),\:\s(N_2)\subset\fF$;

{\em(б)} $\|f(N_1)U-Uf(N_2)\|_{\bS_1}\le\|N_1U-UN_2\|_{\bS_1}$ для всех унитарных операторов $U$
и всех  нормальных операторов $N_1$ и $N_2$ таких, что $\s(N_1),\:\s(N_2)\subset\fF$;

{\em(в)} $\|f(N)U-Uf(N)\|_{\bS_1}\le\|NU-UN\|_{\bS_1}$ для всех самосопряжённых унитарных операторов $U$ и всех
нормальных операторов $N$ таких, что $\s(N)\subset\fF$;

{\em(г)} $\|f(N)A-Af(N)\|_{\bS_1}\le\|NA-AN\|_{\bS_1}$ для всех самосопряжённых операторов $A$ и всех
нормальных операторов $N$ таких, что $\s(N)\subset\fF$;

{\em(д)} $\|f(N)R-Rf(N)\|_{\bS_1}+\|\ov f(N)R-R\ov f(N)\|_{\bS_1}\le\|NR-RN\|_{\bS_1}+\|N^*R-RN^*\|_{\bS_1}$ для всех операторов
$R\in\mB(\h)$ и всех
нормальных операторов $N$ таких, что $\s(N)\subset\fF$;

{\em(е)} $\|f(N_1)R-Rf(N_2)\|_{\bS_1}+\|\ov f(N_1)R-R\ov f(N_2)\|_{\bS_1}
\le\|N_1R-RN_2\|_{\bS_1}+\lb\|N_1^*R-RN_2^*\|_{\bS_1}$ для всех операторов
$R\in\mB(\h)$
и всех  нормальных операторов $N_1$ и $N_2$ таких, что $\s(N_1),\:\s(N_2)\subset\fF$.
\end{thm}

Эта теорема может быть доказана аналогично теоремам \ref{Nlip0} и \ref{Nlip}.
Приведём теперь аналог теоремы \ref{Clip}.

\begin{thm}
\label{ClipS1}
Пусть $f$ -- непрерывная функция, заданная на замкнутом подмножестве $\fF$ комплексной плоскости $\C$.
Тогда следующие три утверждения равносильны:

{\em(а)} $\|f(N)R-Rf(N)\|_{\bS_1}\le\|NR-RN\|_{\bS_1}$ для всех операторов 
$R\in\mB(\h)$ и всех
нормальных операторов $N$ таких, что $\s(N)\subset\fF$;

{\em(б)} $\|f(N_1)R-Rf(N_2)\|_{\bS_1}\le\|N_1R-RN_2\|_{\bS_1}$ для всех операторов $R\in\mB(\h)$ и всех нормальных операторов
$N_1$ и $N_2$ таких, что
$\s(N_1),\:\s(N_2)\subset\fF$;

{\em(в)} $\|f(N_1)A-Af(N_2)\|_{\bS_1}\le\|N_1A-AN_2\|_{\bS_1}$ для всех самосопряжённых операторов $A$
и всех  нормальных операторов $N_1$ и $N_2$ таких, что $\s(N_1),\:\s(N_2)\subset\fF$.
\end{thm}

Заметим, что для самосопряжённых операторов теорему \ref{yaNlip01}
можно переформулировать следующим образом

\begin{thm}
\label{yaNlip01A}
Пусть $f$ -- вещественная непрерывная функция, заданная на замкнутом подмножестве $\fF$ вещественной 
прямой $\R$.
Тогда следующие утверждения равносильны:

{\em(а)} $\|f(A)-f(B)\|_{\bS_1}\le\|A-B\|_{\bS_1}$ для всех самосопряжённых операторов $A$ и $B$ таких, что
$\s(A),\s(B)\subset\fF$;

{\em(б)} $\|f(A)U-Uf(B)\|_{\bS_1}\le\|AU-UB\|_{\bS_1}$ для всех унитарных операторов $U$
и всех  самосопряжённых операторов $A$ и $B$ таких, что $\s(A),\:\s(B)\subset\fF$;

{\em(в)} $\|f(A)U-Uf(A)\|_{\bS_1}\le\|AU-UA\|_{\bS_1}$ для всех самосопряжённых унитарных операторов $U$ и всех
самосопряжённых операторов $A$ таких, что $\s(A)\subset\fF$;

{\em(г)} $\|f(A)R-Rf(A)\|_{\bS_1}\le\|AR-RA\|_{\bS_1}$ для всех самосопряжённых операторов $A$ и  $R$ таких, что $\s(A)\subset\fF$;

{\em(д)} $\|f(A)R-Rf(A)\|_{\bS_1}\le\|AR-RA\|_{\bS_1}$ для всех операторов
$R\in\mB(\h)$ и всех
самосопряжённых операторов $A$ таких, что $\s(A)\subset\fF$;

{\em(е)} $\|f(A)R-Rf(B)\|_{\bS_1}
\le\|AR-RB\|_{\bS_1}$ для всех операторов
$R\in\mB(\h)$
и всех  самосопряжённых операторов $A$ и $B$ таких, что $\s(A),\:\s(B)\subset\fF$.
\end{thm}

\begin{cor}
\label{sledoveshch}
Если $f$ -- вещественная непрерывная функция на замкнутом
подмножестве $\fF$  вещественной прямой $\R$, то
$\|f\|_{\OL_{\bS_1}(\fF)}=\|f\|_{\CL(\fF)}$.
\end{cor}

Отсюда следует, что для комплексной непрерывной функции $f$ имеют место следующие неравенства:
$$
\|f\|_{\OL_{\bS_1}(\fF)}\le\|f\|_{\CL_{\bS_1}(\fF)}\le2\|f\|_{\OL_{\bS_1}(\fF)}.
$$

Всё то же самое можно сказать и в случае унитарных операторов $N_1$  и  $N_2$,
т. е. в том случае, когда множество $\fF$ лежит на единичной окружности $\T$.

Ясно, что $\ov z\in\OL(\fF)$ для любого замкнутого множества $\fF$ в $\C$
и 
$\|\ov z\|_{\OL(\fF)}=1$, если множество $\fF$ содержит по крайней мере
две точки.

\medskip

{\bf Определение.}
Замкнутое подмножество $\fF$ комплексной плоскости $\C$ называется
{\it множеством Фугледе}, если $\CL(\fF)=\OL(\fF)$. 

\medskip

Это понятие было введено
Киссиным и Шульманом в работе \cite{KS2}. 
Отметим также, что это понятие там введено
для любого симметрично-нормированного операторного идеала. 

Джонсон и Вильямс \cite{JW} доказали, что каждая функция $f\in\CL(\fF)$
дифференцируема в комплексном смысле в каждой неизолированной точке множества
$\fF$, см. ниже теоремы  \ref{olm}  и  \ref{pro}. Заметим, что $\ov z\in\OL(\fF)$.
Таким образом, множество Фугледе не может иметь внутренних точек и даже не может 
содержать двух пересекающихся отрезков, не лежащих на одной прямой.
Киссин и Шульман \cite{KS2} доказали, что любая компактная кривая класса 
$C^2$ является множеством Фугледе.

Следующая теорема по существу содержится в предложении 4.5 статьи \cite{KS2}.

\begin{thm} 
\label{Fug}
Замкнутое подмножество $\fF$ комплексной плоскости $\C$ является множеством
Фугледе  в том и только в том случае, когда
$\ov z\in\CL(\fF)$. Если $\ov z\in\CL(\fF)$, то
$\|f\|_{\CL(\fF)}\le\|\ov z\|_{\CL(\fF)}\|f\|_{\OL(\fF)}$ для всех $f\in\OL(\fF)$.
\end{thm}

Очевидным образом эта теорема  вытекает из теорем \ref{Nlip} и \ref{Clip}.

\begin{cor} Пусть $\fF$ -- замкнутое подмножество комплексной плоскости
$\C$. Тогда равенство $\CL(\fF)=\OL(\fF)$ имеет место вместе с равенством
полунорм $\|\cdot\|_{\CL(\fF)}=\|\cdot\|_{\OL(\fF)}$ в том и только в том случае,
когда $\|\ov z\|_{\CL(\fF)}\le1$.
\end{cor}

Заметим, что $\|\ov z\|_{\CL(\fF)}\ge\|\ov z\|_{\OL(\fF)}=1$, если только
множество $\fF$ содержит по крайней мере две точки.
Таким образом, условие $\|\ov z\|_{\CL(\fF)}\le1$ можно заменить
условием $\|\ov z\|_{\CL(\fF)}=1$, если множество $\fF$ содержит
по крайней мере две точки.

%
%

\begin{thm}
 \label{su}
Если все точки  замкнутого подмножества  $\fF$ комплексной плоскости
$\C$ лежат на прямой или окружности, то $\fF$ -- множество Фугледе
и  \lb$\|\cdot\|_{\CL(\fF)}=\|\cdot\|_{\OL(\fF)}$.
\end{thm}

\Pf  Рассмотрим сначала случай, когда множество $\fF$ лежит на прямой.
Аффинное преобразование комплексной плоскости позволяет ограничиться
случаем, когда $\fF\subset\R$. Тогда
$\|N^*R-RN^*\|=\|NR-RN\|$ для любого оператора $R$ и любого нормального оператора $N$ со спектром
в $\fF$, поскольку в этом случае $N$ -- самосопряжённый оператор. Следовательно,
$\|\ov z\|_{\CL(\fF)}\le1$.
Аналогично, в случае, когда множество  $\fF$ лежит на окружности,
достаточно рассмотреть случай, когда  $\fF\subset\T$.
Тогда $\|N^*R-RN^*\|=\|N(N^*R-RN^*)N\|=\|NR-RN\|$
для любого оператора $R$ и любого нормального оператора $N$ со спектром
в $\fF$, поскольку теперь $N$ -- унитарный оператор. Следовательно,
$\|\ov z\|_{\CL(\fF)}\le1$. $\bl$

\medskip

Легко видеть, что нормы в $\CL(\fF)$ и $\OL(\fF)$ совпадают
тогда и только тогда, когда
\bay
\label{Km}
\|N^*R-RN^*\|=\|NR-RN\|
\ey
для всех нормальных операторов $N$ со спектром в $\fF$ и всех ограниченных операторов $R$.
Камовиц \cite{Ka} доказал, что равенство \rf{Km} выполняется для всех ограниченных 
операторов $R$ при фиксированном $N$ в том и только в том случае, когда
$N$ -- нормальный оператор со спектром, лежащим на окружности или на прямой.
Из этого результата Камовица  мгновенно следует, что
теорему \ref{su} можно обратить. Другими словами, равенство 
$\|\cdot\|_{\CL(\fF)}=\|\cdot\|_{\OL(\fF)}$
имеет место в том и только в том случае, когда множество $\fF$ лежит на окружности или на прямой.
Таким образом, результат Камовица содержит в себе теорему  \ref{su}.


Пусть $\fF_1$ и $\fF_2$ -- непустые замкнутые подмножества
комплексной плоскости $\C$.  Обозначим через $\CL(\fF_1,\fF_2)$
пространство всех непрерывных  функций $f$, заданных на $\fF=\fF_1\cup\fF_2$,
для которых существует константа $C\ge0$ такая, что
\bay
\label{cld2}
\|f(N_1)R-Rf(N_2)\|\le C\|N_1R-RN_2\|
\ey
для всех $R\in\mB(\h)$ и всех нормальных операторов $N_1$ и $N_2$ со спектрами,
лежащими в $\fF_1$ и $\fF_2$.
Наименьшую из констант $C$, удовлетворяющих
условию \rf{cld2},
обозначим через $\|f\|_{\CL(\fF_1,\fF_2)}$.
Положим $\|f\|_{\CL(\fF_1,\fF_2)}=\be$, если $f\not\in\CL(\fF_1,\fF_2)$.

Переходя к сопряжённым операторам, мы видим, что условие \rf{cld2}
равносильно условию:
$$
\|R^*\ov f(N_1)-\ov f(N_2)R^*\|\le C\|R^*N_1^*-N_2^*R^*\|.
$$

Отсюда следует, что $f\in\CL(\fF_1,\fF_2)$ в том и только в том случае,
когда $\ov f(\ov z)\in\CL(\ov\fF_2,\ov\fF_1)$ и
$\|\ov f(\ov z)\|_{\CL(\ov\fF_2,\ov\fF_1)}=\|f\|_{\CL(\fF_1,\fF_2)}$,
где для подмножества $\fF$ комплексной плоскости $\C$ символом $\ov\fF$ обозначается  множество $\{\ov\z:~\z\in\fF\}$.

Если условие \rf{cld2} переписать в терминах матриц, а затем перейти
к транспонированным матрицам, то точно так же получится, что
$\CL(\fF_1,\fF_2)=\CL(\fF_2,\fF_1)$ и
$\|\cdot\|_{\CL(\fF_1,\fF_2)}=\|\cdot\|_{\CL(\fF_2,\fF_1)}$.

\begin{thm}
\label{28}
Пусть $f$ -- непрерывная функция, заданная на объединении
\lb$\fF_1\cup\fF_2$ замкнутых подмножеств $\fF_1$ и $\fF_2$
комплексной плоскости $\C$. Тогда следующие утверждения
равносильны:

{\em(а)} $\|f(N_1)R-Rf(N_2)\|\le\|N_1R-RN_2\|$ для всех операторов $R\in\mB(\h)$ и всех
нормальных операторов $N_1$ и $N_2$ таких, что $\s(N_1)\subset\fF_1$ и
$\s(N_2)\subset\fF_2$;

{\em(б)} $\|f(N_1)R-Rf(N_2)\|\le\|N_1R-RN_2\|$ для любого оператора
$R$, действующего из произвольного гильбертова пространства $\h_2$
в произвольное гильбертово пространство $\h_1$, и всех нормальных операторов
$N_1$ и $N_2$, действующих соответственно в пространствах $\h_1$
и $\h_2$, таких, что $\s(N_1)\subset\fF_1$ и $\s(N_2)\subset\fF_2$;

{\em(в)} утверждение {\em(б)} выполняется при дополнительном предположении
о простоте спектров нормальных операторов $N_1$ и $N_2$;

{\em(г)} $\|f(N_1)A-Af(N_2)\|\le\|N_1A-AN_2\|$ для всех самосопряжённых операторов
$A$
и всех  нормальных операторов $N_1$ и $N_2$ таких, что $\s(N_1)\subset\fF_1$ и
$\s(N_2)\subset\fF_2$.

\end{thm}

\Pf Докажем сначала эквивалентность первых трёх утверждений.
Импликации (б)$\Longrightarrow$(а) и  (б)$\Longrightarrow$(в) тривиальны.
Докажем, что (а)$\Longrightarrow$(б) и  (в)$\Longrightarrow$(б).
Начнём с первой импликации.
Если пространства $\h_1$ и $\h_2$
изоморфны,
то тогда существует унитарный оператор $U:\h_1\to\h_2$.
Тогда оператор $RU$ и нормальные операторы $N_1$ и $U^*N_2U$ действуют в одном
и том же гильбертовом пространстве $\h_1$, поэтому
$$
\|f(N_1)(RU)-(RU)f(U^*N_2U)\|\le\|N_1(RU)-(RU)U^*N_2U\|
$$
в силу (а), откуда мгновенно следует нужная нам оценка.
Чтобы свести общий случай к рассмотренному выше частному,
рассмотрим операторы $\mathcal R\df\bigoplus_{j\ge1}R$, 
$\mathcal N_1\df\bigoplus_{j\ge1}N_1$
и $\mathcal N_2\df\bigoplus_{j\ge1}N_2$. Легко видеть, что неравенство
$$
\|f(N_1)R-Rf(N_2)\|\le\|N_1R-RN_2\|
$$
равносильно неравенству
$$
\|f(\mathcal N_1)\mathcal R-\mathcal Rf(\mathcal N_2)\|\le
\|\mathcal N_1\mathcal R-\mathcal R\mathcal N_2\|.
$$

Докажем теперь, что (в) влечёт (б).
Предположим противное. Тогда существуют оператор $R\in\mB(\h_2,\h_1)$ и нормальные операторы $N_1$ и $N_2$,
действующие соответственно в пространствах $\h_1$ и $\h_2$, и такие, что
$\s(N_1)\subset\fF_1$, $\s(N_2)\subset\fF_2$, $\|N_1R-RN_2\|=1$ и $\|f(N_1)R-Rf(N_2)\|>1$.
Таким образом, существуют векторы $u_0\in\h_2$ и $v_0\in\h_1$ такие,
что $\|u_0\|=1$,
$\|v_0\|=1$ и
$\left|\left((f(N_1)R-Rf(N_2))u_0,v_0\right)\right|>1$. Пусть $\h_1^0$  и
$\h_2^0$ обозначают наименьшие приводящие подпространства операторов
$N_1$ и $N_2$, содержащие соответственно $v_0$ и $u_0$.
Пусть $P$ обозначает ортогональный проектор на подпространство $\h_1^0$,
а $Q$ -- на подпространство $\h_2^0$.
Заметим, что $\|f(N_1)PRQ-PRQf(N_2)\|>1$, поскольку $\left((f(N_1)PRQ-PRQf(N_2))u_0,v_0\right)=
\left(f(N_1)R-Rf(N_2))u_0,v_0\right)$. Кроме того,
$\|N_1PRQ-PRQN_2\|=\|P(N_1R-RN_2)Q\|\le1$.
Положим $N_1^0\df N|{\h_1^0}$ и $N_2^0\df N|{\h_2^0}$.
Тогда операторы $N_1^0$ и $N_2^0$ можно рассматривать как нормальные
операторы соответственно в пространствах $\h_1^0$ и $\h_2^0$.
Ясно, что $N_1^0$ и $N_2^0$ -- нормальные операторы с простыми спектрами.
Теперь, чтобы получить противоречие, достаточно заметить, что
$\|f(N_1^0)PRQ-PRQf(N_2^0)\|>1$ и $\|N_1^0PRQ-PRQN_2^0\|\le1$.

Заметим, что импликация (а)$\Longrightarrow$(г)  тривиальна.
Остаётся доказать, что (г) влечёт (а).
Применяя (г) к нормальным операторам $U^{*}N_1U$ и $N_2$, где $U$ --
унитарный оператор, получаем:
 $$
 \|f(N_1)UA-UAf(N_2)\|\le\|N_1UA-UAN_2\|
 $$
 для любого самосопряжённого оператора $A$, унитарного оператора $U$
 и нормальных операторов $N_1$ и  $N_2$ таких, что $\s(N_1)\subset\fF_1$ и
 $\s(N_2)\subset\fF_2$.  Отсюда при помощи полярного разложения следует
 (а) для обратимых операторов $R$. Чтобы получить утверждение (а) в полном объёме,
 достаточно заметить, что множество всех обратимых операторов
 всюду плотно в банаховом пространстве $\mB(\h)$. $\bl$

Нам понадобится следующее хорошо известное элементарное утверждение, позволяющее
приближать произвольный нормальный
оператор нормальными операторами с конечным спектром.

\begin{lem}
\label{apr}
Пусть $N$ -- ограниченный нормальный оператор. Предположим, что множество $\L$, $\L\subset\C$,
является $\e$-сетью спектра $\s(N)$ оператора $N$, т.~е. для любой точки $\z\in\s(N)$
существует точка $\l\in\L$ такая, что $|\l-\z|<\e$. Тогда существует нормальный оператор
$N_0$ такой, что $NN_0=N_0N$, $\|N-N_0\|<\e$ и $\s(N_0)$ -- конечное подмножество множества
$\L$.
\end{lem}

\Pf В силу компактности спектра оператора $N$ существует конечная $\e$-сеть $\L_0$
множества $\s(N)$ такая, что $\L_0\subset \L$. Тогда мы можем найти борелевскую
функцию $\eta:\s(N)\to\L_0$ такую, что $\sup\{|z-\eta(z)|:z\in\s(N)\}<\e$.
Остаётся положить $N_0\df\eta(N)$. $\bl$

\medskip

Из этой леммы и из неравенства \rf{1comlip}  легко вытекает, что если неравенство
\rf{lipid} выполняется для всех нормальных операторов $N_1$ и  $N_2$ с конечными спектрами,
лежащими в  $\fF$, то $f\in\OL(\fF)$ и $\|f\|_{\OL(\fF)}\le C$.

Другими словами, для любой непрерывной функции $f$, заданной на замкнутом множестве  $\fF$,
$\fF\subset\C$, выполняется следующее равенство:
\bay
\label{finite}
\|f\|_{\OL(\frak F)}=\sup\big\{\|f\|_{\OL(\L)}:\L\subset\frak F,\,\,\,
\L\,\,\text{-- конечно}\big\}.
\ey
Более того, для любой непрерывной функции $f:\fF\to\C$ имеет место следующее равенство:
\bay
\label{finite1}
\|f\|_{\OL(\frak F)}=\sup\big\{\|f\|_{\OL(\L)}:\L\subset\fF_0,\,\,\,
\L\,\,\text{-- конечно}\big\},
\ey
где $\fF_0$ -- всюду плотное подмножество множества $\fF$.

 Аналогичные равенства имеют место и для коммутаторно липшицевых полунормы.

Отсюда видно, что, по существу, мы не получили бы ничего существенно нового, если бы
попытались определить пространства $\OL(\fF)$ и $\CL(\fF)$ для произвольного подмножества
$\fF$  комплексной плоскости $\C$.

Остановимся для определённости на случае пространства  $\OL(\fF)$ (случай 
пространства $\CL(\fF)$  может быть рассмотрен аналогично).
Будем говорить, что произвольная функция $f:\fF\to\C$ принадлежит
пространству  $\OL(\fF)$, если существует константа  $C\ge0$ такая, что 
неравенство \ref{lipid} имеет место для всех нормальных операторов 
$N_1$ и $N_2$  с конечными спектрами, лежащими в $\fF$.
Отметим, что требование конечности спектров позволяет определить
$f(N_1)$ и $f(N_2)$ для произвольной функции $f$.
Ясно, что $\OL(\fF)\subset\Li(\fF)$. Таким образом, каждая функция
$f$ допускает липшицево продолжение на замыкание $\clos\fF$ множества $\fF$.
Из равенства \rf{finite1} видно, что эта продолженная функция попадёт
в пространство $\OL(\clos\fF)$  и её \lb$\OL$-полунорма не изменится.
Таким образом, пространство $\OL(\fF)$ отождествляется естественным образом
с пространством $\OL(\clos\fF)$. 

Принимая во внимание это замечание, мы будем для краткости писать
$\OL(\dd)$, $\OL(\C_+)$, $\CL(\dd)$ и $\CL(\C_+)$ вместо
$\OL(\clos\dd)$, $\OL(\clos\C_+)$, $\CL(\clos\dd)$ и $\CL(\clos\C_+)$.

\

\section{\bf Ограниченные и неограниченные нормальные операторы}
\setcounter{equation}{0}
\label{ogrneogr}

\

Мы докажем в этом параграфе некоторые вспомогательные утверждения, из которых вытекает, что в определении операторно липшицевых и коммутаторно липшицевых (равно как и операторно гёльдеровых) функций мы можем рассматривать лишь ограниченные нормальные операторы или допустить к рассмотрению также неограниченные нормальные операторы. В обоих случаях мы получим те же самые классы функций с теми же нормами.

Пусть $N_1$ и $N_2$ -- не обязательно ограниченные нормальные операторы в гильбертовых 
пространствах $\h_1$ и $\h_2$ с областями определения $\mD_{N_1}$ и $\mD_{N_2}$.
Пусть $R$ -- ограниченный оператор из $\h_2$ в $\h_1$.
Мы говорим, что $N_1R-RN_2$ -- {\it ограниченный оператор}, если
$R(\mD_{N_2})\subset \mD_{N_1}$ и $\|N_1Ru-RN_2u\|\le C\|u\|$ для всех $u\in \mD_{N_2}$.
Тогда существует единственный ограниченный оператор
$K$ такой, что
$Ku=N_1Ru-RN_2u$ для всех $u\in \cd_{N_2}$. В этом случае мы пишем $K=N_1R-RN_2$.
Таким образом, $N_1R-RN_2$ -- ограниченный оператор в том и только в том случае, когда 
\bay
\label{MN}
\big|(Ru,{N_1}^*v)-(N_2u,R^*v)\big|\le C\|u\|\cdot\|v\|
\ey
для всех $u\in \mD_{N_2}$ и $v\in \mD_{N_1^*}=\mD_{N_1}$. Легко видеть, что
$N_1R-RN_2$ -- ограниченный оператор, в том и только в том случае, если $N_2^*R^*-R^*N_1^*$ -- ограниченный оператор. При этом
$(N_1R-RN_2)^*=-(N_2^*R^*-R^*N_1^*)$.
В частности, мы пишем $N_1R=RN_2$, если $R(\mD_{N_2})\subset\mD_{N_1}$ и $N_1Ru=RN_2u$ для всех $u\in\mD_{N_2}$.
Скажем, что $\|N_1R-RN_2\|=\be$, если $N_1R-RN_2$ не является ограниченным оператором.

\medskip

{\bf Замечание.} Пусть $N_1$ и $N_2$ нормальные операторы. Предположим, что оператор $N_1^*$ является
замыканием оператора ${N_1}_\flat$, а оператор $N_2$ -- оператора ${N_2}_\sharp$.
Тогда, если неравенство \rf{MN} выполняется для всех $u\in \mD_{{N_2}_\sharp}$ и $v\in\mD_{{N_1}_\flat}$,
то оно будет выполняться для всех $u\in \mD_{N_2}$ и $v\in \mD_{N_1}$.

\begin{thm}
\label{anbnrn}
Пусть $N_1$ и $N_2$ -- нормальные операторы в гильбертовых пространствах $\mathscr H_1$
и $\mathscr H_2$, а $R$ -- ограниченный оператор из $\h_2$ в $\h_1$. Тогда существуют последовательности ограниченных нормальных операторов 
$\{N_{1,n}\}_{n\ge1}$ и $\{N_{2,n}\}_{n\ge1}$, действующих в гильбертовых пространствах $\h_{1,n}$
и $\h_{2,n}$, и последовательность ограниченных операторов $\{R_n\}_{n\ge1}$,
действующих из гильбертова пространства $\h_{2,n}$ в гильбертово
пространство $\h_{1,n}$, такие, что 

{\em(а)} последовательность $\{\|R_n\|\}_{n\ge1}$ является неубывающей и
$$
\lim_{n\to\be}\|R_n\|=\|R\|;
$$

{\em(б)} $\s(N_{1,n})\subset\s(N_1)$ и $\s(N_{2,n})\subset\s(N_2)$при всех $n\ge1$;

{\em(в)} для любой непрерывной на множестве $\s(N_1)\cup\s(N_2)$ функции $f$
последовательность 
$\big\{\big\|f(N_{1,n})R_n-R_nf(N_{2,n})\big\|\big\}_{n\ge1}$ является неубывающей
и
$$
\lim_{n\to\be}\big\|f(N_{1,n})R_n-R_nf(N_{2,n})\big\|=\|f(N_1)R-Rf(N_2)\|;
$$

{\em(г)} для любой непрерывной на множестве $\s(N_1)\cup\s(N_2)$ функции $f$ такой, что
$\|f(N_1)R-Rf(N_2)\|<\be$, и любого натурального числа   $j$ последовательность сингулярных чисел
$\big\{s_j\big(f(N_{1,n})R_n-R_nf(N_{2,n})\big)\big\}_{n\ge0}$ является неубывающей и
$$
\lim_{n\to\be}s_j\big(f(N_{1,n})R_n-R_nf(N_{2,n})\big)=s_j\big(f(N_1)R-Rf(N_2)\big).
$$
\end{thm}

\Pf Не умаляя общности, можно считать, что $0\in\s(N_1)\cup\s(N_2)$.
Положим $P_{1,n}\df E_{N_1}\big(\{|\l|\le n\})$ и $P_{2,n}\df E_{N_2}\big(\{|\l|\le n\})$,
где $E_{N_1}$ и $E_{N_2}$ -- спектральные меры нормальных операторов $N_1$ и $N_2$. 
Положим 
$$
\widetilde N_{1,n}\df P_{1,n}N_1=N_1P_{1,n}=P_{1,n}N_1P_{1,n},\quad
\widetilde N_{2,n}\df P_{2,n}N_2=N_2P_{2,n}=P_{2,n}N_2P_{2,n},
$$
$$
\mathscr H_{1,n}\df P_{1,n}\mathscr H_1\quad\mbox{и}\quad\mathscr H_{2,n}\df P_{2,n}\mathscr H_2.
$$
Ясно, что $\widetilde N_{1,n}$ и $\widetilde N_{2,n}$ -- ограниченные нормальные операторы
в гильбертовых пространствах $\mathscr H_1$ и $\mathscr H_2$, а $\mathscr H_{1,n}$
и $\mathscr H_{2,n}$ -- приводящие подпространства операторов $\widetilde N_{1,n}$ и $\widetilde N_{2,n}$
соответственно. 

Положим $N_{1,n}\df\widetilde N_{1,n}\big|\mathscr H_{1,n}$ и $N_{2,n}\df\widetilde N_{2,n}\big|\mathscr H_{2,n}$.
Тогда $N_{1,n}$ и $N_{2,n}$ -- нормальные операторы, если их рассматривать как операторы,
действующие в гильбертовых пространствах $\mathscr H_{1,n}$ и $\mathscr H_{2,n}$. Оператор $R_n$ из $\mB(\h_{2,n},\h_{1,n})$  определяется равенством $R_nu\df P_{1,n}Ru=P_{1,n}R P_{2,n}u$  для $u\in\mathscr H_{2,n}$.
Утверждения (а) и (б) очевидны. Чтобы доказать остальные утверждения,
достаточно заметить, что
$$
P_{1,n}\big(f(N_1)R-Rf(N_2)\big)P_{2,n}u=\Big(f(N_{1,n})R_n-R_nf(N_{2,n})\Big)u
$$
для всех  $u$ из $\mathscr H_{2,n}$. $\bl$

\medskip

Из теоремы \ref{anbnrn} видно, что в определениях классов операторно липшицевых, коммутаторно липшицевых функций, а также в определении операторно гёльдеровых функций мы можем рассматривать не только ограниченные нормальные (самосопряжённые) операторы, но и неограниченные. При этом определяемые классы функций не изменятся, равно как и не изменятся введённые на них нормы.

\

\section{\bf Разделённые разности и коммутаторная липшицевость}
\setcounter{equation}{0}
\label{razdrazn}

\

С каждой функцией $f$, заданной на замкнутом множестве  $\fF$, $\fF\subset\C$,
мы связываем функцию
$\dg_0 f:\fF\times\fF\to\C$,
\bay
\label{razdra0}
(\dg_0 f)(z,w)\df\left\{\begin{array}{ll}\dfrac{f(z)-f(w)}{z-w},&\text {если}\,\,\,\,z\ne w,\\[.2cm]
0,&\text{если}\,\,\,\,z=w.
\end{array}\right.
\ey
Если множество  $\fF$ не имеет изолированных точек, и в каждой точке
$z$ множества $\fF$ существует конечная производная $f'(z)$
в комплексном смысле, то мы можем определить {\it разделённую разность}
$\dg f:\fF\times\fF\to\C$ равенством
$$
(\dg f)(z,w)\df\left\{\begin{array}{ll}\dfrac{f(z)-f(w)}{z-w},&\text {если}\,\,\,\,z\ne w,\\[.2cm]
f'(z),&\text{если}\,\,\,\,z=w.
\end{array}\right.
$$

\begin{thm}
 \label{f1f2}
Пусть $f$ -- непрерывная функция на объединении
$\fF_1\cup\fF_2$ замкнутых подмножеств $\fF_1$ и $\fF_2$
комплексной плоскости $\C$. Тогда $f\in\CL(\fF_1,\fF_2)$ в том и только
в том случае, если $\dg_0 f\in\fM(\fF_1\times\fF_2)$. При этом
$$
\|f\|_{\CL(\fF_1,\fF_2)}=\|\dg_0 f\|_{\fM_0(\fF_1\times\fF_2)}\le\|\dg_0 f\|_{\fM(\fF_1\times\fF_2)}\le2\|f\|_{\CL(\fF_1,\fF_2)}.
$$
\end{thm}

\Pf Мы докажем только равенство $\|f\|_{\CL(\fF_1,\fF_2)}=\|\dg_0 f\|_{\fM_0(\fF_1\times\fF_2)}$, 
поскольку всё остальное будет тогда следовать из леммы \ref{32}.
Рассмотрим сначала случай конечных множеств $\fF_1$ и $\fF_2$.
Пусть $N_1$ и $N_2$ -- нормальные операторы такие, что $\s(N_1)\subset\fF_1$
и $\s(N_2)\subset\fF_2$.
В силу теоремы  \ref{28} можно считать, что операторы $N_1$ и $N_2$
имеют простой спектр. Тогда в пространствах $\h_1$ и $\h_2$
существуют ортонормированные базисы $\{u_\l\}_{\l\in\s(N_1)}$
и $\{v_\mu\}_{\mu\in\s(N_2)}$ такие, что $N_1u_\l=\l u_\l$ при всех 
$\l\in\s(N_1)$
и $N_2v_\mu=\mu v_\mu$ при всех $\mu\in\s(N_2)$.  С каждым оператором
$X:\h_2\to\h_1$ мы связываем матрицу
$\{(Xv_\mu,u_\l)\}_{(\l,\mu)\in\s(N_1)\times\s(N_2)}$. Имеем
$$
\left((N_1R-RN_2)v_\mu,u_\l\right)=
\big(Rv_\mu,N_1^*u_\l)-(RN_2v_\mu,u_\l)=
(\l-\mu)(Rv_\mu,u_\l).
$$
Аналогично,
$$
\left((f(N_1)R-Rf(N_2))v_\mu,u_\l\right)=(f(\l)-f(\mu))(Rv_\mu,u_\l).
$$
Ясно, что
$$
\{(f(\l)-f(\mu))(Rv_\mu,u_\l)\}
=\{(\dg_0 f)(\l,\mu)\}\star
\{(\l-\mu)(Rv_\mu,u_\l)\}.
$$
Заметим, что матрица $\{a_{\l\mu}\}_{(\l,\mu)\in\s(N_1)\times\s(N_2)}$ представима в виде
$$
\{a_{\l\mu}\}_{(\l,\mu)\in\s(N_1)\times\s(N_2)}=
\{(\l-\mu)(Rv_\mu,u_\l)\}_{(\l,\mu)\in\s(N_1)\times\s(N_2)},
$$
где $R$ -- оператор, действующий из $\h_2$  в $\h_1$, в
том и только в том случае, когда \lb$a_{\l\mu}=0$ при $\l=\mu$.
Теперь равенство $\|f\|_{\CL(\fF_1,\fF_2)}=\|\dg_0 f\|_{\fM_0(\fF_1\times\fF_2)}$
в случае конечных множеств $\fF_1$ и $\fF_2$ очевидно.

Отсюда сразу получаем, что в общем случае имеет место
следующее неравенство $\|f\|_{\CL(\fF_1,\fF_2)}\ge\|\dg_0 f\|_{\fM_0(\fF_1\times\fF_2)}$.

Остаётся доказать, что $\|f\|_{\CL(\fF_1,\fF_2)}\le\|\dg_0 f\|_{\fM_0(\fF_1\times\fF_2)}$.
Ясно, что это достаточно доказать в случае компактных множеств $\fF_1$ и $\fF_2$.
Если множества $\fF_1$ и $\fF_2$ компактны, то существует последовательность
борелевских функций $h_n$, заданных на множестве $\fF_1\cup\fF_2$,
такая, что при всех $n\ge1$ множество $h_n(\fF_1\cup\fF_2)$  конечно, 
$h_n(\fF_1)\subset\fF_1$, $h_n(\fF_2)\subset\fF_2$ и
$|h_n(z)-z|<n^{-1}$ при всех $z\in\fF_1\cup\fF_2$. Тогда
\begin{align*}
\|f(N_1)R-Rf(N_2)\|&=\lim_{n\to\be}\|f(h_n(N_1))R-Rf(h_n(N_2))\|\\[.2cm]
&\le\varliminf_{n\to\be}\|\dg_0 f\|_{\fM_0(h_n(\fF_1)\times h_n(\fF_2))}
\|h_n(N_1)R-Rh_n(N_2)\|\\[.2cm]
&\le\|\dg_0 f\|_{\fM_0(\fF_1\times\fF_2)}\lim_{n\to\be}\|h_n(N_1)R-Rh_n(N_2)\|
\\[.2cm]
&=\|\dg_0 f\|_{\fM_0(\fF_1\times\fF_2)}\|N_1R-RN_2\|.\quad\bl
\end{align*}


%

\medskip

{\bf Замечание.} Неравенство $\|f\|_{\CL(\fF_1,\fF_2)}\le\|\dg_0 f\|_{\fM(\fF_1\times\fF_2)}$ можно доказать и другим способом 
с помощью двойных операторных интегралов, см. замечание к теореме \ref{doikomlip}.

\medskip

Теорема \ref{f1f2} в случае, когда $\fF_1=\fF_2$, приводит к следующему результату.

\begin{thm}
\label{olm}
Пусть $f$ -- функция, заданная на непустом замкнутом
подмножестве $\fF$ комплексной плоскости $\C$. Тогда $f\in\CL(\fF)$
в том и только в том случае, когда $\dg_0 f\in\fM(\fF\times\fF)$. При
этом 
$$
\|f\|_{\CL(\fF)}=\|\dg_0 f\|_{\fM_0(\fF\times\fF)}\le\|\dg_0 f\|_{\fM(\fF\times\fF)}\le2\|f\|_{\CL(\fF)}.
$$
\end{thm}

Заметим, что если $\dg_0 f\in\fM(\fF\times\fF)$ для функции $f$, заданной
на $\F$,  то она непрерывна и даже удовлетворяет условию Липшица.
Действительно, если $\z,\t\in\fF$, то 
$|(\dg_0 f)(\t,\z)|\le\|\dg_0 f\|_{\fM_0(\fF\times\fF)}$,
откуда $|f(\z)-f(\t)|\le\|\dg_0 f\|_{\fM_0(\fF\times\fF)}|\z-\t|$.

Следующее утверждение получено в работе \cite{JW}. 

\begin{thm}
\label{pro}
Пусть $f$ -- функция на замкнутом подмножестве
 $\fF$ комплексной плоскости $\C$ такая, что $\dg_0 f\in\fM(\fF\times\fF)$.
Тогда функция $f$ дифференцируема в комплексном смысле в каждой
неизолированной точке множества $\fF$. Кроме того, если
множество $\F$ неограничено, то существует конечный предел
$\lim\limits_{|z|\to\be}\dfrac{f(z)}{z}$.
\end{thm}

Нам понадобится элементарная лемма, которую мы приводим
без доказательства.

\begin{lem}
\label{37}
Пусть $S$ и $T$ -- произвольные множества.
Предположим, что последовательность функций 
$\{{\bs\f}_n\}$ на $S\times T$
сходится поточечно к функции $\bs\f$. Тогда
$\|\bs\f\|_{\fM(S\times T)}\le\varliminf\limits_{n\to\be}\|{\bs\f}_n\|_{\fM(S\times T)}$.
\end{lem}

{\bf Доказательство теоремы \ref{pro}.} Докажем сначала дифференцируемость
функции в каждой неизолированной точке $a$ множества $\fF$.
Не умаляя общности, можно считать, что $a=0$ и $f(0)=0$. Нам нужно доказать,
что функция $z^{-1}f(z)$ имеет конечный предел при $z\to0$. Предположим противное.
Тогда  эта функция должна иметь по крайней мере две конечных (поскольку
функция $f$ липшицева) предельных точки. Ясно, что можно считать, что
этими предельными точками являются $1$ и $-1$. Таким образом, существуют
две стремящиеся к нулю последовательности $\{\l_n\}_{n\ge1}$ и $\{\mu_n\}_{n\ge1}$  точек множества
$\fF\setminus\{0\}$ такие, что $\lim\limits_{n\to\be}\l_n^{-1}f(\l_n)=1$  и
$\lim\limits_{n\to\be}\mu_n^{-1}f(\mu_n)=-1$. Переходя, если нужно, к
подпоследовательностям, мы можем добиться выполнения следующих
условий:

a) $|\l_n|>|\mu_n|>|\l_{n+1}|$ при всех $n\ge1$;

б) $\lim\limits_{n\to\be}\mu_n^{-1}\l_n=0$ и $\lim\limits_{n\to\be}\l_{n+1}^{-1}\mu_n=0$.

Ясно, что
$\|\{(\dg_0 f)(\l_m,\mu_n)\}\|_{\fM(\nn\times\nn)}\le\|\dg_0 f\|_{\fM(\fF\times\fF)}$.
Заметим, что последовательность $\left\{\|\{(\dg_0 f)(\l_{m+k},\mu_{n+k})\}\|
_{\fM(\nn\times\nn)}\right\}_{k\ge1}$  является невозрастающей и
$$
\lim_{k\to\be}(\dg_0 f)(\l_{m+k},\mu_{n+k})=\sgn(m-n+1/2).
$$
Теперь из леммы \ref{37} следует, что  
$\|\{\sgn(m-n+1/2)\}\|_{\fM(\nn\times\nn)}<+\be$,
что противоречит теореме \ref{Boch}.

Существование конечного предела $\lim\limits_{|z|\to\be}z^{-1}f(z)$
в случае неограниченного множества $\fF$ доказывается аналогичным образом
с той лишь разницей, что теперь нужно выбирать 
последовательности
$\{\l_n\}_{n\ge1}$ и $\{\mu_n\}_{n\ge1}$, стремящиеся к бесконечности. $\bl$

\begin{cor} 
\label{tollin}
Пространство $\CL(\C)$  совпадает с множеством
всех функций вида $az+b$, где $a,b\in\C$.
\end{cor}

\Pf Ясно, что любая функция вида $az+b$, где $a,b\in\C$, принадлежит
пространству $\CL(\C)$. Обратно, из теоремы \ref{pro} следует, что
$f$ -- целая функция. Ясно, что функция $f'$ ограничена,
поскольку  $\CL(\C)\subset\OL(\C)\subset\Li(\C)$. Остаётся воспользоваться
теоремой Лиувилля. $\bl$

\begin{thm}
\label{muld}
Пусть $f$ -- функция, заданная на совершенном подмножестве $\fF$
комплексной плоскости $\C$. Тогда $f\in\CL(\fF)$ в том и только в том случае,
когда $f$ дифференцируема в комплексном смысле в каждой точке множества
$\fF$ и \lb$\dg f\in\fM(\fF\times\fF)$. При этом $\|f\|_{\CL(\fF)}=\|\dg f\|_{\fM(\fF\times\fF)}$.
\end{thm}

\Pf Если $f\in\CL(\fF)$,  то $\dg_0 f\in\fM(\fF\times\fF)$ в силу теоремы \ref{olm}
и следствия \ref{32}. Дифференцируемость функции $f$ следует из
теоремы \ref{pro}. Обратно, если $\dg f\in\fM(\fF\times\fF)$, то
$\dg_0 f\in\fM_0(\fF\times\fF)$ и можно применить теорему \ref{olm}.
Равенство  $\|f\|_{\CL(\fF)}=\|\dg f\|_{\fM(\fF\times\fF)}$ следует из теоремы \ref{olm},
леммы \ref{33} и из очевидного равенства
$\|\dg f\|_{\fM_0(\fF\times\fF)}=\|\dg_0 f\|_{\fM_0(\fF\times\fF)}$. $\bl$

\

Следующая теорема показывает, что нет необходимости рассматривать все нормальные операторы
$N_1$ и $N_2$, а достаточно ограничиться только одной парой нормальных
операторов $N_1$ и $N_2$ таких, что $\s(N_1)=\fF_1$ и $\s(N_2)=\fF_2$.
В частности, когда речь идёт о пространстве $\CL(\fF)$,  можно считать, что $N_1=N_2$,
т.е. ограничиться одним нормальным оператором  $N=N_1=N_2$ таким, что $\s(N)=\fF$.

\begin{thm} 
\label{n1n2}
Пусть $N_1$ и  $N_2$ --  нормальные операторы,
действующих в гильбертовых пространствах $\h_1$
и $\h_2$. Предположим, что непрерывная функция $f$, заданная по крайней мере
на объединении спектров $\s(N_1)\cup\s(N_2)$ этих операторов, обладает следующим 
свойством:
\bay
\label{nerkomlip1}
\|f(N_1)R-Rf(N_2)\|\le\|N_1R-RN_2\|
\ey
для всех $R\in\mB(\h_2,\h_1)$. Тогда
$f\in\CL(\s(N_1),\s(N_2))$  и $\|f\|_{\CL(\s(N_1),\s(N_2))}\le1$. 
\end{thm}

Пусть $f$ -- непрерывная функция на подмножестве комплексной плоскости. Предположим, что 
$N_1$ и $N_2$ -- нормальные операторы, действующие в 
гильбертовых пространствах $\h_1$ и $\h_2$, объединение спектров которых содержится в области определения функции $f$. Будем говорить, что пара $(N_1,N_2)$ является {\it$f$-регулярной}, если неравенство \ref{nerkomlip1} справедливо при всех
$R\in\mB(\h_2,\h_1)$.

Теорему \ref{n1n2} можно переформулировать следующим образом:

\medskip

{\em Если упорядоченная пара нормальных операторов $(N_1,N_2)$ является 
$f$-регулярной, то $f$-регулярной будет и любая пара нормальных
операторов $(M_1,M_2)$ такая, что $\s(M_1)\subset\s(N_1)$ и
$\s(M_2)\subset\s(N_2)$}.

\medskip

Докажем сначала лемму.

\begin{lem} 
\label{n1n2lem}
Пусть $(N_1,N_2)$ -- $f$-регулярная пара ограниченных нормальных операторов. 
Предположим, что оператор $M_1$  унитарно эквивалентен
сужению оператора $N_1$  на приводящее подпространство оператора $N_1$
и оператор $M_2$  унитарно эквивалентен
сужению оператора $N_2$  на приводящее подпространство оператора $N_2$.
Тогда пара $(M_1,M_2)$ является $f$-регулярной.
\end{lem}

\Pf  Пусть операторы $N_1$ и $N_2$ действуют в гильбетовых пространствах
$\h_1$ и  $\h_2$ и пусть $\K_1$ и $\K_2$ --
приводящие подпространства операторов $N_1$ и $N_2$ такие, что оператор $M_1$
унитарно эквивалентен оператору $N_1\big|\K_1$  и оператор $M_2$
унитарно эквивалентен оператору $N_2\big|\K_2$.
Пусть $M_1\in\mB(\widetilde{\h_1})$ и 
$M_2\in\mB(\widetilde{\h_2})$.

Достаточно рассмотреть следующие два частных случая:

1. $\K_1=\h_1$ и $\K_2=\h_2$.
Тогда $M_1=U_1^*N_1U_1$ и $M_2=U_2^*N_2U_2$
для некоторых унитарных операторов (мы допускаем, чтобы
унитарный оператор действовал из одного гильбертова пространства в другое).
Имеем:
\begin{align*}
\|f(M_1)R-Rf(M_2)\|&=\|U_1^*f(N_1)U_1R-RU_2^*f(N_2)U_2\|\\[.2cm]
&=\|f(N_1)U_1RU_2^*-U_1RU_2^*f(N_2)\|\\[.2cm]
&\le\|N_1U_1RU_2^*-U_1RU_2^*N_2\|
=\|M_1R-RM_2\|
\end{align*}
для любого оператора $R$ из $\mB(\widetilde{\h_2},\widetilde{\h_1})$.

2. $M_1=N_1\big|\K_1$ и $M_2=N_2\big|\K_2$. Пусть
$P_1$ -- ортогональный проектор из $\h_1$ на $\K_1$
и  $P_2$ -- ортогональный проектор из $\h_2$ на $\K_2$.
Если $R\in\mB(\K_2,\K_1)$, то
Тогда
\begin{align*}
\|f(M_1)R-Rf(M_2)\|&=\|P_1(f(M_1)R-Rf(M_2))P_2\|\\[.2cm]
&=\|P_1(f(N_1)R-Rf(N_2))P_2\|
=\|f(N_1)P_1RP_2-P_1RP_2f(N_2)\|\\[.2cm]
&\le\|N_1P_1RP_2-P_1RP_2N_2\|
=\|M_1R-RM_2\|.\quad\bl
\end{align*}

\medskip

{\bf Доказательство теоремы \ref{n1n2}.} Прежде всего, заметим, что в силу теоремы  
\ref{anbnrn}, достаточно рассмотреть случай ограниченных операторов $N_1$ и $N_2$. 
Легко видеть, что из леммы \ref{n1n2lem} и теоремы \ref{28} 
мгновенно вытекает теорема \ref{n1n2} в случае, когда спектры операторов $N_1$  и $N_2$ 
конечны. Таким образом, из леммы \ref{apr}  следует, что остаётся доказать, что
для любых конечных подмножеств $\D_1$ и $\D_2$ множеств $\s(N_1)$ и $\s(N_2)$
найдутся нормальные операторы $M_1\in\mB(\K_1)$ и $M_2\in\mB(\K_2)$
такие, что $\s(M_1)=\D_1$, $\s(D_2)=\L_2$ и
$$
\|f(M_1)R-Rf(M_2)\|\le\|M_1R-RM_2\|
$$
для всех $R$ из $\mB(\K_2,\K_1)$.
С каждым нормальным оператором $N$ мы свяжем функцию $\a_N$ такую, что
$\a_N(\z)$ есть кратность спектра оператора $N$ в изолированной точке $\z$ его спектра $\s(N)$ 
и $\a_N(\z)=\be$ в каждой неизолированной точке его спектра.

В качестве операторов $M_1$ и $M_2$ мы возьмём какие-либо нормальные
операторы в гильбертовых пространствах $\K_1$,
$\K_2$ и обладающие следующими свойствами:

1)  $\s(M_1)=\L_1$ и $\s(M_2)=\L_2$;

2)  функции $\a_{M_1}$ и $\a_{M_2}$ являются сужениями функций $\a_{N_1}$ и $\a_{N_2}$.

Пусть $\D_1^{(\e)}$ и $\D_2^{(\e)}$ обозначают замкнутые $\e$-окрестности множеств
$\D_1$ и $\D_2$. Пусть $N_1^{(\e)}$ обозначает сужение оператора $N_1$  на подпространство
$E_{N_1}\big(\D_1^{(\e)}\cap\s(N_1)\big)$, а $N_2^{(\e)}$ -- сужение оператора $N_2$  на подпространство
$E_{N_2}\big(\D_2^{(\e)}\cap\s(N_2)\big)$, где $E_{N_1}$ и $E_{N_2}$ --  спектральные меры 
операторов $N_1$ и $N_2$. Ясно, что существуют операторы $M_1^{(\e)}$
из $\mB(\K_1)$
и $M_2^{(\e)}$ из $\mB(\K_2)$ такие, что оператор $M_1^{(\e)}$ унитарно эквивалентен
оператору  $N_1^{(\e)}$, оператор $M_2^{(\e)}$ унитарно эквивалентен
оператору  $N_2^{(\e)}$, $\|M_1-M_1^{(\e)}\|\le\e$ и $\|M_2-M_2^{(\e)}\|\le\e$.
Тогда для любого оператора $R$ из $\mB(\K_2,\K_1)$ имеем
\begin{align*}
\|f(&M_1)R-Rf(M_2)\|
\le\|R\|\cdot\big\|f(M_1)-f\big(M_1^{(\e)}\big)\big\|+
\|R\|\cdot\big\|f(M_2)-f\big(M_2^{(\e)}\big)\big\|\\[.2cm]
+&\big\|f\big(M_1^{(\e)}\big)R-Rf\big(M_2^{(\e)}\big)\big\|
\le\|R\|\cdot\big\|f(M_1)-f\big(M_1^{(\e)}\big)\big\|
+\|R\|\cdot\big\|f(M_2)-f\big(M_2^{(\e)}\big)\big\|\\[.2cm]
+&\big\|M_1^{(\e)}R-RM_2^{(\e)}\big\|
\le\|R\|\cdot\big\|f(M_1)-f\big(M_1^{(\e)}\big)\big\|
+\|R\|\cdot\big\|f(M_2)-f\big(M_2^{(\e)}\big)\big\|\\[.2cm]
+&2\e\|R\|
+\|M_1R-RM_2\|.
\end{align*}
Остаётся перейти к пределу при $\e\to0$. $\bl$

\medskip

 Следующая теорема содержится в работе \cite{JW}.

\begin{thm} 
Пусть  $M$ и $N$ -- операторы в гильбертовом пространстве $\h$,
причём оператор $N$ нормальный. Тогда следующие утверждения равносильны:

{\em(а)} $M=f(N)$ для некоторой функция $f\in\CL(\s(N))$;

{\em(б)} существует константа $c$ такая, что $\|MR-RM\|\le c\|NR-RN\|$ для всех ограниченных операторов $R$;

{\em(в)} существует константа $c$ такая, что $\|MR-RM\|_{\bS_1}\le c\|NR-RN\|_{\bS_1}$ для всех ограниченных операторов $R$;

{\em(г)} для любого ограниченного оператора $T$ найдётся ограниченный оператор  $S$
такой, что $SN-NS=TM-MT$;

{\em(д)} для любого компактного оператора $T$ найдётся ограниченный оператор $S$
такой, что $SN-NS=TM-MT$;

{\em(е)} для любого оператора $T$ из $\bS_1(\h)$ найдётся оператор $S$
из $\bS_1(\h)$
такой, что $SN-NS=TM-MT$.
\end{thm}

На самом деле работа  \cite{JW}  содержит и целый ряд дополнений к этой теореме.
Например, чтобы доказать эквивалентности (б)$\Longleftrightarrow$(е) и (в)$\Longleftrightarrow$(г)$\Longleftrightarrow$(д), можно не требовать, чтобы оператор $N$ был нормальным.

\

\section{\bf Мультипликаторы Шура и операторная липшицевость}
\setcounter{equation}{0}
\label{dmshol}

\

Заметим, что если замкнутое подмножество $\fF$ является множеством Фугледе,
то $\OL(\fF)=\CL(\fF)$ в силу теоремы \ref{Fug}. Таким образом, в этом случае теорема \ref{olm}
даёт полное описание пространства $\OL(\fF)$ в терминах мультипликаторов Шура.

В частности, для множеств $\fF$, лежащих на прямой или окружности, мы имеем
полное описание пространства  $\OL(\fF)$ в терминах мультипликаторов Шура.
Причём в последнем случае и полунорма операторно липшицевой функции
выражается через норму мультипликатора Шура.

В случае, когда множество $\fF$ не является множеством  Фугледе, нам неизвестно 
полное описание операторно липшицевых функций на $\fF$  в терминах мультипликаторов Шура.

В этом случае мы можем предложить следующее достаточное условие операторной липшицевости.

\begin{thm}
\label{teoroboperlipots} 
Пусть $f$ -- функция, непрерывная на замкнутом подмножестве $\fF$
комплексной плоскости $\C$. Предположим, что существуют 
мультипликаторы Шура $\Phi_1,\Phi_2\in\fM(\fF\times\fF)$ такие, что 
$$
f(z)-f(w)=(z-w)\Phi_1(z,w)+(\ov z-\ov w)\Phi_2(z,w).
$$
Тогда $f\in\OL(\fF)$ и 
$$
\|f\|_{\OL(\fF)}\le\|\Phi_1\|_{\fM(\fF\times\fF)}+\|\Phi_2\|_{\fM(\fF\times\fF)}.
$$
\end{thm}

Эту теорему можно доказать с помощью аппроксимации операторами с конечными спектрами, как это сделано в доказательстве теоремы \ref{f1f2}. Мы опустим это доказательство, приведя взамен доказательство, основанное на двойных операторных интегралах, см. теорему \ref{olivsdoi} и замечание к ней.

\medskip

{\bf Замечание.} Теорему \ref{teoroboperlipots}  иногда удобнее применять, привлекая  вещественные
переменные: $z=x_1+{\rm i}y_1$, $w=x_2+{\rm i}y_2$.  Предположим, что существуют
мультипликаторы Шура $F_1,F_2\in\fM(\fF\times\fF)$ такие, что 
$$
f(z)-f(w)=(x_1-x_2)F_1(z,w)+(y_1-y_2)F_2(z,w).
$$
Тогда $f\in\OL(\fF)$ и 
$$
\|f\|_{\OL(\fF)}\le\frac12\|F_1+{\rm i}F_2\|_{\fM(\fF\times\fF)}+\frac12\|F_1-{\rm i}F_2\|_{\fM(\fF\times\fF)}
\le\|F_1\|_{\fM(\fF\times\fF)}+\|F_2\|_{\fM(\fF\times\fF)}.
$$

\medskip

Следующая теорема показывает, что операторно липшицева функция имеет производную
в каждой точке области определения вдоль любого невырожденного направления.

\begin{thm} 
Пусть $f\in\OL(\fF)$, где $\fF$ -- замкнутое множество в
$\C$. Тогда для любой прямой $l$ сужение $f\big|l\cap\fF$
дифференцируемо в каждой неизолированной точке множества $l\cap\fF$
и в $\be$, если множество $l\cap\fF$ не ограничено.
\end{thm}

\Pf Ясно, что $f\big|l\cap\fF\in\OL(l\cap\fF)$. Остаётся заметить, что
$\CL(l\cap\fF)=\OL(l\cap\fF)$ в силу теоремы \ref{su}, и применить
теорему \ref{muld} к функции $f\big|l\cap\fF$. $\bl$

\begin{cor}  Пусть $f\in\OL(\fF)$, где $\fF$ -- замкнутое подмножество
комплексной плоскости. Тогда функция $f$ дифференцируема по любому направлению
в каждой внутренней точке множества $\fF$.
\end{cor}

{\bf Замечание.}  Функция $f\in\OL(\fF)$ не обязана быть
дифференцируемой как функция двух вещественных переменных в каждой внутренней 
точке множества $\fF$. Например, нетрудно проверить, что функция $f$, записываемая
в полярных координатах следующим образом: $f(r,\theta)=re^{3{\rm i}\theta}$,
принадлежит пространству $\OL(\C)$, но не дифференцируема в нуле как функция
двух вещественных переменных $x$ и $y$.

\medskip

Это было отмечено в \cite{AP6},  см. также \cite{A1}.

\

\section{\bf Роль двойных операторных интегралов}
\setcounter{equation}{0}
\label{doinastene}

\

В этом параграфе мы продемонстрируем роль двойных операторных интегралов в оценках операторных разностей и (квази)коммутаторов. Мы начнём с оценок операторных разностей при возмущении самосопряжённого оператора оператором класса Гильберта--Шмидта и обсудим формулу Бирмана--Соломяка.

Далее, мы вернёмся к результатам предыдущих двух параграфов, где мы получили условия для коммутаторной липшицевости и операторной липшицевости в терминах принадлежности некоторых функций пространству дискретных мультипликаторов Шура. В этом параграфе мы приведём другое доказательство достаточности этих условий с помощью двойных операторных интегралов.
При этом мы получим полезные формулы, выражающие операторные разности и коммутаторы в терминах двойных операторных интегралов. 

Наконец, мы получим формулы для операторных производных в терминах двойных операторных интегралов.

Следующая теоремы была получена М.С. Бирманом и М.З. Соломяком в \cite{BS3}.

\begin{thm}
\label{BiSoS2}
Пусть $f$ -- липшицева функция на $\R$, а $A$ и $B$ -- самосопряжённые операторы в гильбертовом пространстве, разность которых $A-B$ входит в класс Гильберта--Шмидта $\bS_2$. Тогда имеет место формула
\bay
\label{infoABS2}
f(A)-f(B)=\int_\R\!\int_\R\big(\dg_0 f\big)(x,y)\,dE_A(x)(A-B)\,dE_B(y).
\ey
\end{thm}

Отметим, что из формулы \ref{infoABS2} непосредственно вытекает неравенство
$$
\|f(A)-f(B)\|_{\bS_2}\le\|f\|_{\Li}\|A-B\|_{\bS_2}.
$$
Иными словами, липшицевы функции являются $\bS_2$-липшицевыми. Оказывается, что 
липшицевы функции также являются $\bS_p$-липшицевыми при $p\in(1,\be)$. Это было недавно доказано в работе \cite{PS}. Напомним, что при $p=1$ соответствующее утверждение неверно. Впервые это было доказано в работе \cite{F2}. Более того, класс \lb$\bS_1$-липшицевых функций совпадает с классом операторно липшицевых функций, см. теорему \ref{tsentrez}.

Мы опускаем здесь доказательство теоремы \ref{BiSoS2} и отсылаем читателя за доказательством к работе \cite{BS3}.

Займёмся теперь коммутаторной липшицевостью. 

\begin{thm}
\label{doikomlip}
Пусть $\fF_1$ и $\fF_2$ -- замкнутые подмножества  плоскости $\C$. Предположим, что $f$ -- непрерывная функция на $\fF_1\cup\fF_2$ такая, что 
функция $\dg_0 f$, определённая равенством {\em\rf{razdra0}}, принадлежит классу
мультипликаторов Шура $\fM(\fF_1\times\fF_2)$. Тогда, если $N_1$ и $N_2$ -- нормальные операторы такие, что $\s(N_j)\subset\fF_j$, $j=1,\,2$, а $R$ -- ограниченный линейный оператор, то справедлива формула
\bay
\label{fkommudoi}
f(N_1)R-Rf(N_2)=
\int_{\!\fF_1}\!\!\int_{\!\fF_2}
\big(\dg_0 f\big)(\z_1,\z_2)\,dE_1(\z_1)(N_1R-RN_2)\,dE_2(\z_2),
\ey
где $E_j$ -- спектральная мера оператора $N_j$.
\end{thm}

{\bf Замечание.}
Из формулы \rf{fkommudoi} сразу же вытекает, что
$$
\|f(N_1)R-Rf(N_2)\|\le\|\dg_0 f\|_{\fM(E_1,E_2)}\|N_1R-RN_2\|
\le\|\dg_0 f\|_{\fM(\fF_1\times\fF_2)}\|N_1R-RN_2\|,
$$
и, в частности, $f$ -- коммутаторно липшицева функция.

\medskip 

В частном случае, когда $R$ -- единичный оператор, получаем следующий результат:

\begin{thm}
\label{doioplip}
Пусть $\fF$ -- замкнутое подмножество плоскости $\C$, а  $f$ -- непрерывная функция на $\fF$ такая, что 
 $\dg_0 f\in\fM(\fF\times\fF)$. Тогда, если $N_1$ и $N_2$ -- нормальные операторы такие, что $\s(N_j)\subset\fF$, то справедлива формула
\bay
\label{frazndoi}
f(N_1)-f(N_2)=
\int_\fF\int_\fF
\big(\dg_0 f\big)(\z_1,\z_2)\,dE_1(\z_1)(N_1-N_2)\,dE_2(\z_2).
\ey
\end{thm}

\medskip

{\bf Доказательство теоремы \ref{doikomlip}.} Предположим сначала, что операторы $N_1$ и $N_2$ ограничены. Имеем
\begin{align*}
\int_{\fF_1}\!\!\int_{\fF_2}\big(\dg_0 f\big)(\z_1,\z_2)\,dE_1(\z_1)&(N_1R-RN_2)\,dE_2(\z_2)\\[.2cm]
&=\int_{\fF_1}\!\!\int_{\fF_2}
\big(\dg_0 f\big)(\z_1,\z_2)\,dE_1(\z_1)N_1R\,dE_2(\z_2)\\[.2cm]
&-\int_{\fF_1}\!\!\int_{\fF_2}
\big(\dg_0 f\big)(\z_1,\z_2)\,dE_1(\z_1)RN_2\,dE_2(\z_2).
\end{align*}

Из определения двойных операторных интегралов следует, что
$$
\int_{\fF_1}\!\!\int_{\fF_2}
\big(\dg_0 f\big)(\z_1,\z_2)\,dE_1(\z_1)N_1R\,dE_2(\z_2)
=
\int_{\fF_1}\!\!\int_{\fF_2}
\z_1\big(\dg_0 f\big)(\z_1,\z_2)\,dE_1(\z_1)R\,dE_2(\z_2)
$$
и
$$
\int_{\fF_1}\!\!\int_{\fF_2}
\big(\dg_0 f\big)(\z_1,\z_2)\,dE_1(\z_1)RТ_2\,dE_2(\z_2)
=
\int_{\fF_1}\!\!\int_{\fF_2}
\z_2\big(\dg_0 f\big)(\z_1,\z_2)\,dE_1(\z_1)R\,dE_2(\z_2).
$$
Поскольку $(\z_1-\z_2)(\dg_0 f\big)(\z_1,\z_2)=f(\z_1)-f(\z_2)$, $\z_1\in\fF_1$, $\z_2\in\fF_2$, получаем
\begin{align*}
\int_{\fF_1}\!\!\int_{\fF_2}
\big(\dg_0 f\big)(\z_1,\z_2)&\,dE_1(\z_1)(N_1R-RN_2)\,dE_2(\z_2)\\
=&\int_{\fF_1}\!\!\int_{\fF_2}
f(\z_1)\,dE_1(\z_1)R\,dE_2(\z_2)-
\int_{\fF_1}\!\!\int_{\fF_2}
f(\z_2)\,dE_1(\z_1)R\,dE_2(\z_2).
\end{align*}
Опять, легко видеть из определения двойных операторных интегралов, что
$$
\int_{\fF_1}\!\!\int_{\fF_2}f(\z_1)\,dE_1(\z_1)R\,dE_2(\z_2)=
\left(\,\,\int_{\fF_1}f(\z_1)\,dE_1(\z_1)\right)R=f(N_1)R
$$
и
$$
\int_{\fF_1}\!\!\int_{\fF_2}f(\z_2)\,dE_1(\z_1)R\,dE_2(\z_2)=
R\int\limits_{\fF_1}f(\z_1)\,dE_1(\z_1)=Rf(N_2),
$$
откуда вытекает равенство \rf{fkommudoi}.

Пусть теперь $N_1$ и $N_2$ -- необязательно ограниченные нормальные операторы.
Отметим прежде всего, что из доказанной части теоремы \ref{doikomlip} и из теоремы \ref{anbnrn} вытекает коммутаторно липшицева оценка и, стало быть, оператор
$f(N_1)R-Rf(N_2)$ ограничен.

Положим
$$
P_k\df E_1\big(\{\z\in\C:|\z|\le k\}\big)\quad\mbox{и}\quad
Q_k\df E_2\big(\{\z\in\C:|\z|\le k\}\big),\quad k>0.
$$
Тогда
$$
N_{1,k}\df P_kN_1\quad\mbox{и}\quad N_{2,k}\df Q_kN_2~\mbox{ --}
$$
ограниченные нормальные операторы. Пусть $E_{j,k}$ -- спектральная мера оператора
$N_{j,k}$, $j=1,\,2$. 
Имеем
\begin{align*}
P_k&\left(\int_{\fF_1}\!\!\int_{\fF_2}\big(\dg_0 f\big)(\z_1,\z_2)\,
dE_1(\z_1)(N_1R-RN_2)\,dE_2(\z_2)\right)Q_k\\[.2cm]
&=P_k\left(\int_{\fF_1}\!\!\int_{\fF_2}\big(\dg_0 f\big)(\z_1,\z_2)\,
dE_{1,k}(\z_1)(P_kf(N_1)R-Rf(N_2)Q_k)\,dE_{2,k}(\z_2)\right)Q_k.
\end{align*}

Применяя \rf{fkommudoi} к ограниченным нормальным операторам $N_{1,k}$ и $N_{2,k}$, получаем
\begin{align*}
P_k\big(f(N_{1,k})R&-Rf(N_{2,k})\big)Q_k=\\[.2cm]
=&P_k\left(\int_{\fF_1}\!\!\int_{\fF_2}\big(\dg_0 f\big)(\z_1,\z_2)\,
dE_{1,k}(\z_1)(P_kN_1R-RN_2Q_k)\,dE_{2,k}(\z_2)\right)Q_k.
\end{align*}
Поскольку
$$
P_k\big(f(N_{1,k})R-Rf(N_{2,k})\big)Q_k=P_k\big(f(N_{1})R-Rf(N_{2})\big)Q_k,
$$
имеем
\begin{align*}
P_k\big(f(N_{1})R&-Rf(N_{2})\big)Q_k=\\[.2cm]
=&P_k\left(\int_{\fF_1}\!\!\int_{\fF_2}\big(\dg_0 f\big)(\z_1,\z_2)\,
dE_1(\z_1)(N_1R-RN_2)\,dE_2(\z_2)\right)Q_k.
\end{align*}
Остаётся перейти к пределу в сильной операторной топологии. $\bl$

\medskip

В случае, когда  $\fF_1=\fF_2$ и множество $\fF_1$ совершенно, можно, ввиду теоремы \ref{muld}, в формуле \ref{fkommudoi} заменить $\dg_0f$ разделённой разностью $\dg f$.

\begin{thm}
\label{infsovmn}
Пусть $\fF$  -- замкнутое совершенное подмножество плоскости $\C$, а 
$f\in\CL(\fF)$. Тогда, если $N_1$ и $N_2$ -- нормальные операторы такие, что $\s(N_j)\subset\fF$, $j=1,\,2$, а $R$ -- ограниченный линейный оператор, то справедлива формула
\bay
\label{fkommudoisov}
f(N_1)R-Rf(N_2)=
\int_{\fF}\!\int_{\fF}
\big(\dg f\big)(\z_1,\z_2)\,dE_1(\z_1)(N_1R-RN_2)\,dE_2(\z_2),
\ey
где $E_j$ -- спектральная мера оператора $N_j$.
\end{thm}

Перейдём теперь к осмыслению результатов \S\:\ref{dmshol} в свете двойных операторных интегралов. Имеет место следующее утверждение:

\begin{thm}
\label{olivsdoi}
Пусть $f$ -- функция, непрерывная на замкнутом подмножестве $\fF$
комплексной плоскости $\C$ такая, что существуют 
мультипликаторы Шура $\Phi_1,\Phi_2\in\fM(\fF\times\fF)$ такие, что 
$$
f(\z_1)-f(\z_2)=(\z_1-\z_2)\Phi_1(\z_1,\z_2)+(\ov\z_1-\ov\z_2)\Phi_2(\z_1,\z_2),
\quad\z_1,~\z_2\in\fF.
$$
Предположим, что $N_1$ и $N_2$ -- нормальные операторы, спектры которых содержатся в $\fF$.
Тогда справедлива формула
\begin{align}
\label{intforadoi}
f(N_1)-f(N_2)&=\int_\fF\int_\fF
\Phi_1(\z_1,\z_2)\,
dE_1(\z_1)(N_1-N_2)\,dE_2(\z_2)\nonumber\\[.2cm]
&+\int_\fF\int_\fF
\Phi_2(\z_1,\z_2)(\z_1,\z_2)\,
dE_1(\z_1)(N^*_1-N^*_2)\,dE_2(\z_2).
\end{align}
\end{thm}

{\bf Замечание.}
Из формулы \rf{intforadoi} легко следует, что 
$$
\|f(N_1)-f(N_2)\|\le\big(\|\Phi_1\|_{\fM(\fF\times\fF)}
+\|\Phi_2\|_{\fM(\fF\times\fF)}\big)\|N_1-N_2\|
$$
и, в частности, $f$ -- {\it операторно липшицева функция}.

\medskip

\Pf Как и в доказательстве теоремы \ref{doikomlip}, предположим сначала, что операторы $N_1$ и $N_2$ ограничены. Тогда
\begin{align*}
\int_\fF\int_\fF
\Phi_1(\z_1,\z_2)\,
&dE_1(\z_1)(N_1-N_2)\,dE_2(\z_2)\\[.2cm]
&=\int_\fF\int_\fF
\Phi_1(\z_1,\z_2)
dE_1(\z_1)N_1dE_2(\z_2)-
\int_\fF\int_\fF
\Phi_1(\z_1,\z_2)
dE_1(\z_1)N_2dE_2(\z_2)\\[.2cm]
&=\int_\fF\int_\fF
\z_1\Phi_1(\z_1,\z_2)
dE_1(\z_1)dE_2(\z_2)-
\int_\fF\int_\fF
\z_2\Phi_1(\z_1,\z_2)
dE_1(\z_1)dE_2(\z_2)\\[.2cm]
&=\int_\fF\int_\fF
(\z_1-\z_2)\Phi_1(\z_1,\z_2)\,dE_1(\z_1)\,dE_2(\z_2).
\end{align*}

Аналогичным образом,
$$
\int_\fF\int_\fF
\Phi_2(\z_1,\z_2)(\z_1,\z_2)
dE_1(\z_1)(N^*_1-N^*_2)dE_2(\z_2)=
\int_\fF\int_\fF(\ov\z_1-\ov\z_2)\Phi_2(\z_1,\z_2)
\,dE_1(\z_1)\,dE_2(\z_2).
$$
Таким образом, правая часть равенства \rf{intforadoi} равна 
\begin{align*}
\int_\fF\int_\fF(f(\z_1)&-f(\z_2))\,dE_1(\z_1)\,dE_2(\z_2)\\[.2cm]
&=\int_\fF f(\z_1)\,dE_1(\z_1)
-\int_\fF f(\z_2)\,dE_2(\z_2)=f(N_1)-f(N_2).
\end{align*}
Переход от ограниченных к неограниченным операторам осуществляется так же, как и в доказательстве теоремы \ref{doikomlip}. $\bl$

\medskip

Перейдём теперь к применению двойных операторных интегралов в задачах
операторной дифференцируемости.

\begin{thm} 
\label{sildif}
Пусть $f$ -- операторно липшицева функция на $\R$, а
$A$ и $K$ -- самосопряжённые операторы, причём оператор $K$ ограничен.
Тогда
$$
\lim_{t\to0}\frac1t(f(A+tK)-f(A))=\int_\R\int_\R(\dg f)(x,y)\,dE_A(x)K\,dE_A(y),
$$
где предел берётся в сильной операторной топологии.
\end{thm}

Нам понадобится несколько вспомогательных утверждений.
Пусть $\widehat\R\df\R\cup\{\be\}$  обозначает {\it одноточечную компактификацию
вещественной прямой} $\R$. Напомним, что  любая функция $f\in\OL(\R)$  дифференцируема всюду на $\widehat\R$, см. теорему \ref{pro}.

\begin{lem}
\label{abchag} 
Пусть $f\in\OL(\R)$. Тогда существуют две последовательности непрерывных на 
$\widehat\R$
функций $\{\f_n\}_{n\ge0}$ и $\{\psi_n\}_{n\ge0}$ такие, что 

{\rm a)} $\sum\limits_{n\ge0}|\f_n|^2\le\|f\|_{\OL(\R)}$ всюду на $\widehat\R$,

{\rm б)} $\sum\limits_{n\ge0}|\psi_n|^2\le\|f\|_{\OL(\R)}$ всюду на $\widehat\R$,

{\rm в)} $(\dg f)(x,y)=\sum\limits_{n\ge0}\f_n(x)\psi_n(y)$ при всех $x,\,y\in\R$.
\end{lem}

\Pf В силу теоремы \ref{muld} мы имеем: $\dg f\in\fM(\R\times\R)$ и $\|\dg f\|_{\fM(\R\times\R)}=\|f\|_{\OL(\R)}$.
Продолжим функцию $\dg f$ на множество $\widehat\R\times\widehat\R$, положив
$(\dg f)(x,y)=f'(\be)=\lim\limits_{t\to\be}\dfrac{f(t)}{t}$ в случае, когда $|x|+|y|=\be$. Ясно что эта продолженная на $\widehat\R\times\widehat\R$
функция $\dg f$
непрерывна по каждой переменной. Отсюда получаем:
\begin{align*}
\|\dg f\|_{\fM(\widehat\R\times\widehat\R)}&=
\sup\{\|\dg f\|_{\fM(\L_1\times\L_2)}:~\L_1,\,\L_2\subset\widehat\R, \,\,
\L_1\,\,\text {и} \,\,\L_2\,\,\,\text{конечны}\}\\[.2cm]
&=\sup\{\|\dg f\|_{\fM(\L_1\times\L_2)}:~\L_1,\,\L_2\subset\R, \,\,
\L_1\,\,\text {и} \,\,\L_2\,\,\,\text{конечны}\}=\|\dg f\|_{\fM(\R\times\R)}.
\end{align*}
Остаётся применить теорему \ref{npkp} к функции $\dg f:\widehat\R\times\widehat\R\to\C$. $\bl$

\begin{lem} 
\label{nepratk-}
Пусть $A$ и $K$ -- самосопряжённые операторы, причём оператор $K$ ограничен.
Тогда для любой функции $f$ из $C(\widehat\R)$ функция $H$, $H(t)=f(A+tK)$, действует
непрерывно из $\R$ в пространство $\mB(\h)$ с нормированной топологией.
\end{lem}

Отметим, что в работе \cite{AP2} получено существенно более сильное утверждение.

\medskip

\Pf Можно считать, что $f(\be)=0$. Тогда мы можем построить последовательность
$\{f_n\}_{n\ge0}$
функций класса $C^\be$ с компактным носителем  такую, что $f_n\to f$ равномерно.
Каждая функция $H_n$, $H_n(t)=f_n(A+tK)$, непрерывна, поскольку $f_n\in\OL(\R)$
при $n\ge0$. Остаётся заметить, что $H_n\to H$ равномерно. $\bl$

\begin{lem} 
\label{Xnun-}
Пусть $\{X_n\}_{n\ge0}$ -- последовательность в пространстве операторов $\mB(\h)$,
а $\{u_n\}_{n\ge0}$ -- в гильбертовом пространстве $\h$. Предположим, что 
\lb$\sum_{n\ge0}X_nX_n^*\le a^2I$  и $\sum_{n\ge0}\|u_n\|^2\le b^2$
для некоторых неотрицательных чисел $a$ и $b$. Тогда ряд $\sum_{n\ge0}X_nu_n$
слабо сходится и 
$$
\Big\|\sum_{n\ge0}X_nu_n\Big\|\le ab.
$$
\end{lem}

\Pf Пусть $v\in\h$ и $\|v\|=1$. Тогда
$$
\sum_{n\ge0}|(X_nu_n,v)|=\sum_{n\ge0}|(u_n,X_n^*v)|\le
\Big(\sum_{n\ge0}\|u_n\|^2\Big)^{1/2}\Big(\sum_{n\ge0}\|X_n^*v\|^2\Big)^{1/2}\le ab,
$$
откуда следует доказываемое утверждение. $\bl$

\medskip

{\bf Доказательство теоремы \ref{sildif}.} 
В силу
формул \rf{frazndoi} и \rf{Haagrazl}  теорему \ref{sildif}  можно переформулировать следующим образом:
$$
\lim_{t\to0}\sum_{n\ge0}\f_n(A+tK)K\psi_n(A)=\sum_{n\ge0}\f_n(A)K\psi_n(A)
$$
в сильной операторной топологии, где $\f_n$ и $\psi_n$ обозначают функции, существование которых утверждается в лемме \ref{abchag}.
Другими словами, нам нужно доказать, что для любого вектора $u\in\h$ мы имеем
$$
\lim_{t\to0}\sum_{n\ge0}(\f_n(A+tK)-\f_n(A))K\psi_n(A)u=\0,
$$
где ряд суммируется в слабой топологии пространства $\h$, а предел берётся в пространстве $\h$ по норме.
Будем считать, что $\|u\|=1$ и $\|f\|_{\OL(\R)}=1$.  Тогда $\sum_{n\ge0}|\f_n|^2\le1$ и
$\sum_{n\ge0}|\psi_n|^2\le1$ всюду на $\R$.

Положим $u_n\df K\psi_n(A)u$. Имеем:
$$
\sum_{n\ge0}\|u_n\|^2\le\|K\|^2\sum_{n\ge0}\|\psi_n(A)u\|^2=
\|K\|^2\sum_{n\ge0}(|\psi_n|^2(A)u,u)\le\|K\|^2<+\be.
$$
Пусть $\e>0$. Выберем  натуральное число $N$ так, чтобы 
$\sum_{n>N}\|u_n\|^2<\e^2$. Тогда из леммы  \ref{Xnun-}
следует, что
$$
\Big\|\sum_{n>N}(\f_n(A+tK)-\f_n(A))u_n\Big\|\le2\e
$$
при всех $t\in\R$. Заметим, что
$$
\Big\|\sum_{n=0}^N(\f_n(A+tK)-\f_n(A))u_n\Big\|\le\|K\|\sum_{n=0}^N\|\f_n(A+tK)-\f_n(A)\|<\e
$$
при всех достаточно близких к нулю $t$ в силу леммы \ref{nepratk-}.
Таким образом, 
$$
\Big\|\sum_{n\ge0}(\f_n(A+tK)-\f_n(A))u_n\Big\|<3\e
$$
при всех достаточно близких к нулю $t$. $\bl$

\medskip

Аналогично теореме \ref{sildif} можно  доказать следующую теорему.

\begin{thm} 
\label{sildifloc-}
Пусть $A$ и $K$ -- ограниченные самосопряжённые операторы.
Тогда
$$
\lim_{t\to0}\frac1t(f(A+tK)-f(A))=\int_\R\int_\R(\dg f)(x,y)\,dE_A(x)K\,dE_A(y),
$$
для любой функции $f\in\OL_{\rm loc}(\R)$,
где предел берётся в сильной операторной топологии.
\end{thm}

Из теоремы \ref{sildif} вытекает следующее утверждение:

\begin{thm}
Пусть $f$ -- операторно дифференцируемая функция на $\R$, а $A$ и $K$ -- самосопряжённые операторы, причём оператор $K$ ограничен. Тогда для производной 
функции $t\mapsto f(A+tK)-f(A)$ по операторной норме справедлива формула
\bay
\label{propoopno}
\frac{d}{dt}\big(f(A+tK)-f(A)\big)\Big|_{t=0}=
\int_\R\int_\R\frac{f(x)-f(y)}{x-y}\,dE_A(x)K\,dE_A(y).
\ey
\end{thm} 

В частности, формула \rf{propoopno} справедлива для любой функции $f$ класса Бесова $B_{\be,1}^1(\R)$, см. теорему \ref{Besdiffer}.

Аналогичные утверждения справедливы и для функций на окружности.

\

\section{\bf Ядерная липшицевость и ядерно-коммутаторная липшицевость}
\setcounter{equation}{0}
\label{yadiyadcom}

\

Цель этого параграфа состоит в том, чтобы доказать, что для произвольного замкнутого множества $\fF$ на плоскости классы $\CL(\fF)$ и  $\CL_{\bS_1}(\fF)$ совпадают. В частности, если $\fF\subset\R$, то совпадают и классы $\OL(\fF)$ и  
$\OL_{\bS_1}(\fF)$ (см. \S\:\ref{oplip}, где определяются классы 
$\CL_{\bS_1}(\fF)$ и $\OL_{\bS_1}(\fF)$). 

Заметим, что  определение класса $\CL_{\bS_1}(\fF)$ естественным образом обобщается до определения класса $\CL_{\bS_1}(\fF_1,\fF_2)$,
где $\fF_1$ и $\fF_2$ -- непустые замкнутые подмножества
комплексной плоскости $\C$.


\begin{lem}
\label{f1f2s1}
Пусть $f$ -- непрерывная функция на объединении
$\fF_1\cup\fF_2$ замкнутых подмножеств $\fF_1$ и $\fF_2$
комплексной плоскости $\C$. Тогда
$$
\|f\|_{\CL_{\bS_1}(\fF_1,\fF_2)}\ge\|\dg_0 f\|_{\fM_{0,\bS_1}(\fF_1\times\fF_2)}
\ge\frac12\|\dg_0 f\|_{\fM(\fF_1\times\fF_2)}.
$$
\end{lem}

\Pf Второе неравенство вытекает из следствия \ref{32s1}. Докажем первое неравенство.
Достаточно рассмотреть случай конечных множеств $\fF_1$ и $\fF_2$.
Пусть $N_1$ и $N_2$ -- нормальные операторы с простыми спектрами такими, что $\s(N_1)=\fF_1$
и $\s(N_2)=\fF_2$.
Тогда в пространствах $\h_1$ и $\h_2$
существуют ортонормированные базисы $\{u_\l\}_{\l\in\s(N_1)}$
и $\{v_\mu\}_{\mu\in\s(N_2)}$ такие, что $N_1u_\l=\l u_\l$ при всех 
$\l\in\s(N_1)$
и $N_2v_\mu=\mu v_\mu$ при всех $\mu\in\s(N_2)$.  С каждым оператором
$X:\h_2\to\h_1$ мы связываем матрицу
$\{(Xv_\mu,u_\l)\}_{(\l,\mu)\in\s(N_1)\times\s(N_2)}$. Имеем
$$
\left((N_1R-RN_2)v_\mu,u_\l\right)=
\big(Rv_\mu,N_1^*u_\l)-(RN_2v_\mu,u_\l)=
(\l-\mu)(Rv_\mu,u_\l).
$$
Аналогично,
$$
\left((f(N_1)R-Rf(N_2))v_\mu,u_\l\right)=(f(\l)-f(\mu))(Rv_\mu,u_\l).
$$
Ясно, что
$$
\{(f(\l)-f(\mu))(Rv_\mu,u_\l)\}
=\{(\dg_0 f)(\l,\mu)\}\star
\{(\l-\mu)(Rv_\mu,u_\l)\}.
$$
Заметим, что матрица $\{a_{\l\mu}\}_{(\l,\mu)\in\s(N_1)\times\s(N_2)}$ представима в виде
$$
\{a_{\l\mu}\}_{(\l,\mu)\in\s(N_1)\times\s(N_2)}=
\{(\l-\mu)(Rv_\mu,u_\l)\}_{(\l,\mu)\in\s(N_1)\times\s(N_2)},
$$
где $R$ -- оператор, действующий из $\h_2$  в $\h_1$, в
том и только в том случае, когда \lb$a_{\l\mu}=0$ при $\l=\mu$.
Теперь неравенство $\|f\|_{\CL(\fF_1,\fF_2)}\ge\|\dg_0 f\|_{\fM_{0,\bS_1}(\fF_1\times\fF_2)}$
очевидно. $\bl$

\begin{cor}
\label{napryamoi}
Пусть $f$ -- вещественная непрерывная функция на замкнутом подмножестве $\fF$  
вещественной прямой $\R$. Тогда
\bay
\label{tsepravnerav}
\|f\|_{\OL_{\bS_1}(\fF)}=\|f\|_{\CL_{\bS_1}(\fF)}\ge\|\dg_0 f\|_{\fM_{0,\bS_1}(\fF\times\fF)}
\ge\frac12\|\dg_0 f\|_{\fM(\fF\times\fF)}.
\ey
Если отказаться от вещественности функции $f$, то 
\bay
\label{vtortsepner}
\|f\|_{\OL_{\bS_1}(\fF)}\ge\frac12\|f\|_{\CL_{\bS_1}(\fF)}\ge\frac12\|\dg_0 f\|_{\fM_{0,\bS_1}(\fF\times\fF)}
\ge\frac14\|\dg_0 f\|_{\fM(\fF\times\fF)}.
\ey
\end{cor}

\Pf Равенство в \rf{tsepravnerav} вытекает из следствия \ref{sledoveshch}.
Все неравенства в \rf{tsepravnerav} уже доказаны выше. Очевидно, что 
\rf{vtortsepner} следует из \rf{tsepravnerav} $\bl$

\medskip

Перейдём теперь к основным результатам этого параграфа.

\begin{thm}
\label{sovkliyadkl}
Пусть $\fF$ -- замкнутое множество в $\C$. Тогда $\CL(\fF)\!=\!\CL_{\bS_1}(\fF)$ и
$$
\frac12\|f\|_{\CL_{\bS_1}(\fF)}\le\|f\|_{\CL(\fF)}\le2\|f\|_{\CL_{\bS_1}(\fF)},
\quad f\in\CL_(\fF).
$$
\end{thm}

\Pf Пусть $f\in\CL(\fF)$ и пусть $N_1$ и $N_2$ -- нормальные операторы со спектрами в $\fF$ такие, что $N_1R-RN_2\in\bS_1$, а $R$ - ограниченный оператор. Тогда, ввиду замечания к теореме \ref{doikomlip},
\begin{align*}
\|f(N_1)R-Rf(N_2)\|_{\bS_1}
&\le\big\|\dg_0f\|_{\fM(E_1,E_2)}\|N_1R-RN_2\|_{\bS_1}\\
&\le\big\|\dg_0f\|_{\fM(\fF\times\fF)}\|N_1R-RN_2\|_{\bS_1}
\le2\|f\|_{\CL(\fF)}\|N_1R-RN_2\|_{\bS_1}
\end{align*}
в силу теоремы \ref{olm}.
Отсюда получаем неравенство $\|f\|_{\CL_{\bS_1}(\fF)}\le2\|f\|_{\CL(\fF)}$.
С другой стороны из леммы \ref{f1f2s1} и из теоремы \ref{olm} получаем
$$
\|f\|_{\CL(\fF)}\le\big\|\dg_0f\big\|_{\fM(\fF\times\fF)}\le2\|f\|_{\CL_{\bS_1}(\fF)}.\quad\bl
$$

\medskip

Если $\fF$ -- совершенное множество, можно получить более точный результат.

\begin{thm}
\label{sovershenstvo}
Пусть $\fF$ -- замкнутое совершенное множество в $\C$. Тогда 
\lb$\|f\|_{\CL_{\bS_1}(\fF)}=\|f\|_{\CL(\fF)}$ для всех 
$f$ из $\CL(\fF)=\CL_{\bS_1}(\fF)$.
\end{thm}

\Pf В силу формулы \rf{fkommudoisov} при $f\in\CL(\fF)$ имеем
\begin{align*}
\|f(N_1)R-Rf(N_2)\|_{\bS_1}
&\le\big\|\dg f\|_{\fM(E_1,E_2)}\|N_1R-RN_2\|_{\bS_1}\\
&\le\big\|\dg f\|_{\fM(\fF\times\fF)}\|N_1R-RN_2\|_{\bS_1}
=\|f\|_{\CL(\fF)}\|N_1R-RN_2\|_{\bS_1}.
\end{align*}
Последнее равенство гарантируется теоремой \ref{muld}. Тем самым
мы доказали, что $\|f\|_{\CL_{\bS_1}(\fF)}\le\|f\|_{\CL(\fF)}$. 
Применяя теперь леммы \ref{f1f2s1}  и \ref{33s1}, а также теорему \ref{muld}, получаем:
$$
\|f\|_{\CL_{\bS_1}(\fF)}\ge\|\dg_0 f\|_{\fM_{0,\bS_1}(\fF\times\fF)}
=\|\dg f\|_{\fM(\fF\times\fF)}=\|f\|_{\CL(\fF)}.\quad\bl
$$

%

\medskip

Настало время перейти к центральному результату этого параграфа.

\begin{thm}
\label{tsentrez}
Пусть $f$ -- непрерывная функция на $\R$. Следующие утверждения эквивалентны:

{\em(а)} функция $f$ операторно липшицева;

{\em(б)} функция $f$ ядерно липшицева;

{\em(в)} $f(A)-f(B)\in\bS_1$, как только $A$ и $B$ -- самосопряжённые операторы с разностью $A-B$ из $\bS_1$.
\end{thm} 

Отметим, что в утверждении (в) {\it необходимо рассматривать не только ограниченные операторы} $A$ и $B$.

\medskip

\Pf Эквивалентность условий (а) и (б) установлена в следствии \ref{napryamoi}. Импликация (б)$\Rightarrow$(в) тривиальна. Покажем, что (в)$\Rightarrow$(б).
Предположим, что $f\not\in\CL_{\bS_1}(\R)$. Тогда можно найти последовательности самосопряжённых операторов $A_n$ и $B_n$ такие, что $A_n-B_n\in\bS_1$ и
$\|A_n-B_n\|_{\bS_1}^{-1}\|f(A_n)-f(B_n)\|_{\bS_1}\to\be$, когда $n\to\be$.
Не умаляя общности, можно считать, что $\lim_{n\to\be}\|A_n-B_n\|_{\bS_1}=0$.
Действительно, рассмотрим приращение $A_n\mapsto A_n+K_n$, где $K_n\df B_n-A_n$. 
Рассмотрим теперь следующие приращения: 
$A_n+(j/M_n)K_n\mapsto A_n+((j+1)/M_n)K_n$,
$0\le j\le M_n-1$, где $\{M_n\}$ -- последовательность натуральных чисел такая, что 
$\lim_{n\to\be}\|A_n-B_n\|_{\bS_1}/M_n=0$. Выберем теперь такое число $j$, при котором число 
$$
\|f(A_n+((j+1)/M_n)K_n)-f(A_n+(j/M_n)K_n)\|_{\bS_1}
$$
максимально. Заменим теперь пару $(A_n,B_n)$ парой 
$A_n+(j/M_n)K_n, A_n+((j+1)/M_n)K_n$. Ясно, что теперь
$$
\lim_{n\to\be}\|A_n-B_n\|_{\bS_1}^{-1}\|f(A_n)-f(B_n)\|_{\bS_1}=\be\quad\mbox{и}
\quad\lim_{n\to\be}\|A_n-B_n\|_{\bS_1}=0.
$$
Теперь достаточно, в случае необходимости, выбрать подпоследовательность
последовательности $(A_n,B_n)$ или повторить некоторые члены этой последовательности, чтобы добиться условия
$$
\sum_n\|B_n-A_n\|_{\bS_1}<\be\quad\mbox{но}\quad
\sum_n\|f(B_n)-f(A_n)\|_{\bS_1}=\be.
$$
Пусть теперь $A$ -- ортогональная сумма операторов $A_n$, а
$B$ -- ортогональная сумма операторов $B_n$. Тогда $B-A\in\bS_1$, а 
$f(B)-f(A)\not\in\bS_1$. $\bl$

\medskip

{\bf Замечание.} Аналогичное утверждение справедливо для функций на единичной окружности и унитарных операторов.

\

\section{\bf Операторно липшицевы функции на плоскости. Достаточное условие}
\setcounter{equation}{0}
\label{OLnaplBes}

\

В этом параграфе мы получим достаточное условие для того, чтобы функция $f$ на плоскости была операторно липшицевой. Это достаточное условие даётся в терминах класса Бесова $B_{\be,1}^1(\R^2)$ и аналогично теореме \ref{Besdost} для функций на
вещественной прямой. Результаты этого параграфа получены в работе \cite{APPS}.

Напомним (см. формулу \rf{frazndoi}), что в случае функций на прямой операторная липшицевость функции $f$ может быть получена из формулы
$$
f(A)-f(B)=\int_\R\int_\R\frac{f(s)-f(t)}{s-t}\,dE(s)(A-B)\,dE(t).
$$
Здесь $A$ и $B$ - самосопряжённые операторы в гильбертовом пространстве. Именно таким образом было получено первоначальное доказательство операторной липшицевости
функций класса $B_{\be,1}^1(\R)$ в \cite{Pe1} и \cite{Pe3}.

Естественно было бы попробовать поступить так же и с функциями на плоскости. 
Однако (см. лемму \ref{tollin}), если разделённая разность является мультипликатором Шура для любых спектральных мер на $\C$, то функция должна быть линейной.

В работе \cite{APPS} был применён другой метод: для нормальных операторов $N_1$ и $N_2$ разность $f(N_1)-f(N_2)$ представляется в виде суммы двойных операторных интегралов. При этом интегрируются разделённые разности по каждой переменной.

Введём следующие обозначения. Пусть $N_1$ и $N_2$ -- нормальные операторы в гильбертовом пространстве. Положим
$$
A_j\df\re N_j,\quad B_j\df\im N_j,\quad E_j\quad\mbox{--~ спектральная мера оператора}\quad N_j,\quad j=1,\,2.
$$
Другими словами, $N_j=A_j+{\rm i}B_j$, $j=1,\,2$, где $A_j$ и $B_j$ -- коммутирующие самосопряжённые операторы. 

Если $f$ -- функция на $\R^2$, имеющая всюду частные производные по обеим переменным, рассмотрим разделённые разности по каждой переменной
$$
\big(\dg_xf\big)(z_1,z_2)\df\frac{f(x_1,y_2)-f(x_2,y_2)}{x_1-x_2},
\quad z_1,\,z_2\in\C,
$$
и
$$
\big(\dg_yf\big)(z_1,z_2)\df\frac{f(x_1,y_1)-f(x_1,y_2)}{y_1-y_2},
\quad z_1,\,z_2\in\C,
$$
где
$$
x_j\df\re z_j,\quad y_j\df\im z_j,\quad  j=1,\,2.
$$
Заметим, что на множествах $\{(z_1,z_2):~x_1=x_2\}$ и
$\{(z_1,z_2):~y_1=y_2\}$
разделённые разности $\dg_xf$ и $\dg_yf$
будут пониматься, как соответствующие частные производные функции $f$.

Следующий результат даёт ключевую оценку.

\begin{thm}
\label{SMdd}
Пусть $f$ -- непрерывная ограниченная функция на $\R^2$, преобразование Фурье 
$\F f$ которой имеет компактный носитель. Тогда разделённые разности $\dg_xf$ и 
$\dg_yf$ являются мультипликаторами Шура класса $\fM(\C\times\C)$.

Более того, если
$$
\supp\F f\subset\{\z\in\C:~|\z|\le\s\},\quad\s>0,
$$
то
\bay
\label{nerxy}
\|\dg_xf\|_{\fM(\C\times\C)}\le\const\s\|f\|_{L^\be}
\quad\mbox{и}\quad\|\dg_yf\|_{\fM(\C\times\C)}\le \const\s\|f\|_{L^\be}.
\ey
\end{thm}

Из определения класса Бесова $B_{\be,1}^1(\R^2)$ и из теоремы \ref{SMdd} следует, что для любой функции $f$ класса $B_{\be,1}^1(\R^2)$ разделённые разности $\dg_xf$ и $\dg_yf$ являются мультипликаторами Шура и
\bay
\label{nerdlyaBes}
\|\dg_xf\|_{\fM(\C\times\C)}\le\const\|f\|_{B_{\be,1}^1}
\quad\mbox{и}\quad\|\dg_yf\|_{\fM(\C\times\C)}\le \const\|f\|_{B_{\be,1}^1}.
\ey

Неравенства \rf{nerdlyaBes} вместе с теоремой \ref{olivsdoi} влекут следующий результат, полученный в работе \cite{APPS} и являющийся центральным в этом параграфе.

\begin{thm}
\label{doin}
Пусть $f$ - функция класса Бесова $B_{\be,1}^1(\R^2)$. Предположим, что $N_1$ и $N_2$ -- нормальные операторы такие, что оператор $N_1-N_2$ ограничен. Тогда
\begin{align*}
f(N_1)-f(N_2)&=\iint\limits_{\C^2}\big(\dg_yf\big)(z_1,z_2)\,
dE_1(z_1)(B_1-B_2)\,dE_2(z_2)\nonumber\\[.2cm]
&+\iint\limits_{\C^2}\big(\dg_xf\big)(z_1,z_2)\,
dE_1(z_1)(A_1-A_2)\,dE_2(z_2),
\end{align*}
и имеет место неравенство
$$
\|f(N_1)-f(N_2)\|\le\const\|f\|_{B_{\be,1}^1}\|N_1-N_2\|,
$$
то есть $f$ -- операторно липшицева функция на $\C$.
\end{thm}

Чтобы доказать теорему \ref{SMdd}, мы воспользуемся одной формулой для представления разделённой разности в виде элемента тензорного произведения Хогерупа. Напомним, что символом $\mathscr E_\s$ обозначается множество целых функций (одной комплексной переменной)
экспоненциального типа не выше, чем $\s$.

\begin{lem}
\label{predosnnaformKotSh}
Пусть $\f\in\mathscr E_\s\cap L^\be(\R)$. Тогда
\begin{align}
\label{haa1}
\frac{\f(x)-\f(y)}{x-y}&=\sum_{n\in\Z}\s\cdot\frac{\f(x)-\f\big(\pi n\s^{-1}\big)}{\s x-\pi n}
\cdot\frac{\sin(\s y-\pi n)}{\s y-\pi n}.
\end{align}
Кроме того,
\bay
\label{2vy}
\sum_{n\in\Z}\frac{\big|\f(x)-\f\big(\pi n\s^{-1}\big)\big|^2}{(\s x-\pi n)^2}
\le 3\|\f\|_{L^\be(\R)}^2,
\quad x\in\R,
\ey
и
\bay
\label{hiz}
\sum_{n\in\Z}\frac{\sin^2(\s y-\pi n)}{(\s y-\pi n)^2}=1,\quad y\in\R.
\ey
\end{lem}

Мы отсылаем читателя к работе \cite{APPS}, в параграфе 5 которой даётся доказательство леммы \ref{predosnnaformKotSh}, основанное на формуле Котельникова -- Шеннона, которая, в свою очередь, основана на том, что семейство функций 
$\{(z-\pi n)^{-1}\sin(z-\pi n)\}_{n\in\Z}$
образует ортонормированный базис в $\mathscr E_1\cap L^2(\R)$,
см. \cite{L}, лекция 20, пункт 3.

\medskip

{\bf Доказательство теоремы \ref{SMdd}.} Ясно, что $f$ -- сужение на $\R^2$ целой функции двух комплексных переменных.
Более того, $f(\cdot,a),\,f(a,\cdot)\in \mathscr E_\s\cap L^\be(\R)$ при всех $a\in\R$.
Не умаляя общности, мы можем предположить, что $\s=1$. По лемме \ref{predosnnaformKotSh} 
\bey
\big(\dg_xf\big)(z_1,z_2)=\frac{f(x_1,y_2)-f(x_2,y_2)}{x_1-x_2}=
\sum_{n\in\Z}(-1)^n\frac{f(\pi n,y_2)-f(x_2,y_2)}{\pi n-x_2}\cdot\frac{\sin(x_1-\pi n)}{x_1-\pi n}
\eey
и
\bey
\big(\dg_yf\big)(z_1,z_2)=\frac{f(x_1,y_1)-f(x_1,y_2)}{y_1-y_2}=\sum_{n\in\Z}(-1)^n\frac{f(x_1,y_1)-f(x_1,\pi n)}{y_1-\pi n}\cdot\frac{\sin(y_2-\pi n)}{y_2-\pi n}.
\eey

Заметим, что выражения $\dfrac{\sin(x_1-n\pi)}{x_1-\pi n}$ и $\dfrac{f(x_1,y_1)-f(x_1,\pi n)}{y_1-\pi n}$ 
зависят от $z_1=(x_1,y_1)$ и не зависят от $z_2=(x_2,y_2)$, в то время как
выражения $\dfrac{f(\pi n,y_2)-f(x_2,y_2)}{\pi n-x_2}$ и $\dfrac{\sin(y_2-\pi n)}{y_2-\pi n}$
зависят от $z_2=(x_2,y_2)$ и не зависят от $z_1=(x_1,y_1)$. Кроме того, по лемме \ref{predosnnaformKotSh}
\bey
\sum_{n\in\Z}\frac{|f(x_1,y_1)-f(x_1,\pi n)|^2}{(y_1-\pi n)^2}\le3\|f(x_1,\cdot)\|_{L^\be(\R)}^2
\le3\|f\|_{L^\be(\C)}^2,\\
\sum_{n\in\Z}\frac{|f(\pi n,y_2)-f(x_2,y_2)|^2}{(\pi n-x_2)^2}\le3\|f(\cdot,y_2)\|_{L^\be(\R)}^2
\le3\|f\|_{L^\be(\C)}^2\\
\eey
и
$$
\sum_{n\in\Z}\frac{\sin^2(x_1-\pi n)}{(x_1-\pi n)^2}=
\sum_{n\in\Z}\frac{\sin^2(y_2-\pi n)}{(y_2-\pi n)^2}=1,
$$
что доказывает \rf{nerxy}. $\bl$

\medskip

Отметим, что неравенства \rf{nerxy} играют роль операторных неравенств Бернштейна (см. \S\:\ref{Bern}). Так же, как и в случае функций от самосопряжённых операторов, можно доказать следующие утверждения.

\begin{thm}
\label{Hanormal}
Пусть $0<\a<1$ и пусть $f$ -- функция класса Гёльдера $\L_\a(\R^2)$.
Тогда 
$$
\|f(N_1)-f(N_2)\|\le\frac{8c}{1-2^{\a-1}}\|f\|_{\L_\a}\|N_1-N_2\|^\a
$$
для некоторой константы $c>0$ и
для любых нормальных операторов $N_1$ и $N_2$ с ограниченной разностью $N_1-N_2$.
\end{thm}

Теорему \ref{Hanormal} можно обобщить на случай произвольных модулей непрерывности.

\begin{thm}
\label{Spnormal}
Пусть $0<\a<1$, $p>1$ и пусть $f$ -- функция класса Гёльдера $\L_\a(\R^2)$. Тогда существует положительное число $c$ такое, что
$$
\|f(N_1)-f(N_2)\|_{S_{p/\a}}\le c\,\|f\|_{\L_\a}\|N_1-N_2\|_{\bS_p}^\a
$$
для произвольных нормальных операторов $N_1$ и $N_2$, разность которых входит в класс Шаттена--фон Неймана $\bS_p$.
\end{thm}

Мы отсылаем читателя к работе \cite{APPS}, где можно найти доказательства этих утверждений, равно как и других результатов на эту тему.

\

\section{\bf Достаточное условие коммутаторной липшицевости в терминах
интегралов Коши}
\setcounter{equation}{0}
\label{dostol}

\

В этом параграфе мы приводим достаточное условие коммутаторной липшицевости,
полученное в работе \cite{A2}.

Пусть $\fF$ -- непустое замкнутое подмножество комплексной плоскости $\C$, причём 
$\fF\ne\C$.
Обозначим через $\M(\C\setminus \fF)$
множество всех комплексных мер Радона на $\C\setminus\fF$ таких, что
\bay
\label{uslname}
\|\mu\|_{\M(\C\setminus \fF)}\df
\sup_{z\in\fF}\int_{\C\setminus\fF}\frac{d|\mu|(\z)}{|\z-z|^2}<+\be.
\ey
Для $\mu\in\M(\C\setminus \fF)$ интеграл Коши
$$
\widehat\mu(z)=\int_{\C\setminus\fF}\frac{d\mu(\z)}{\z-z}\,\,,
$$
вообще говоря, не определён  даже при $z\in\fF$, поскольку функция $\z\mapsto(\z-z)^{-1}$
не обязана быть суммируемой по мере $|\mu|$.
С каждой фиксированной точкой $z_0\in\fF$ мы связываем модифицированный интеграл Коши следующим образом:
$$
\widehat\mu_{z_0}(z)\df\int_{\C\setminus\fF}\Big(\frac1{\z-z}-\frac1{\z-z_0}\Big)\,d\mu(\z).
$$
Из неравенства Коши--Буняковского следует, что
величина $\widehat\mu_{z_0}(z)$ корректно определена при $z\in\fF$ и
\mbox{$|\widehat\mu_{z_0}(z)|\le\|\mu\|_{\M(\C\setminus \fF)}|z-z_0|$}.
Кроме того, $\widehat\mu_{z_0}(z_1)=-\widehat\mu_{z_1}(z_0)$ и
$$
|\widehat\mu_{z_0}(z_1)-\widehat\mu_{z_0}(z_2)|=|\widehat\mu_{z_1}(z_1)-\widehat\mu_{z_1}(z_2)|
=|\widehat\mu_{z_1}(z_2)|
\le\|\mu\|_{\M(\C\setminus \fF))}|z_1-z_2|
$$
для любых $z_1,z_2\in\fF$.

Заметим, что отображение $z\mapsto(\z-z)^{-1}$
действует непрерывно из множества $\fF$ в гильбертово пространство $L^2(|\mu|)$,
наделённое слабой топологией.
Это обстоятельство позволяет легко проверить, что функция $\widehat\mu_{z_0}(z)$
дифференцируема, как функция комплексной переменной, в каждой неизолированной точке
множества $\fF$. В частности, $\widehat\mu_{z_0}(z)$ аналитична во внутренности
множества $\fF$.

Обозначим символом $\widehat\M(\fF)$
множество всех функций
$f$ на $\fF$, представимых в виде $f=c+\widehat\mu_{z_0}$, где $c$ -- константа. Положим
$$
\|f\|_{\widehat\M(\fF)}\df
\inf\{\|\mu\|_{\M(\C\setminus \fF)}:~\mu\in\M(\C\setminus \fF),\,\,\,
f-\widehat\mu_{z_0}=\const\,\,\, \text{на}\,\,\,\fF\}.
$$
Легко видеть, что определение пространства $\widehat\M(\fF)$ и
полунормы $\|f\|_{\widehat\M(\fF)}$
не зависит от выбора точки $z_0\in\fF$.


\begin{thm}
\label{coint7}
Пусть $\fF$ -- собственное замкнутое подмножество комплексной плоскости $\C$.
Тогда $\widehat\M(\fF)\subset\CL(\fF)$ и $\|f\|_{\CL(\fF)}\le\|f\|_{\widehat\M(\fF)}$
для всех $f$ из $\widehat\M(\fF)$.
\end{thm}

\Pf Пусть $\mu\in\M(\C\setminus \fF)$ и $f=\widehat\mu_{z_0}$. 
Рассмотрим разделённую разность
\begin{align*}
\frac{f(z)-f(w)}{z-w}&=
\frac1{z-w}\int_{\C\setminus\fF}\left(\frac1{\z-z}-\frac1{\z-w}\right)\,d\mu(\z)
\\[.2cm]
&=\int_{\C\setminus\fF}\frac{d\mu(\z)}{(\z-z)(\z-w)}.
\end{align*}
Неравенство \rf{uslname} означает, что эта разделённая разность удовлетворяет условию (г) теоремы \ref{tomSc}. Стало быть, она является мультипликатором Шура для любых спектральных борелевских мер на $\fF$. Причём её мультипликаторная норма не превосходит $\|\mu\|_{{\mathscr M}(\C\setminus\fF)}$.

Пусть теперь $N$ -- нормальный оператор со спектром в $\fF$, а $R$ -- ограниченный линейный оператор. Из теоремы \ref{doikomlip} немедленно вытекает, что 
\begin{align*}
f(N)R-Rf(N)=\int_{\C\setminus\fF}\int_{\C\setminus\fF}
\frac{f(z)-f(w)}{z-w}\,dE_N(z)(NR-RN)\,dE_N(w),
\end{align*}
где $E_N$ -- спектральная мера оператора $N$. Стало быть,
\begin{align*}
\|f(N)R-Rf(N)\|&\le\left\|\frac{f(z)-f(w)}{z-w}\right\|_{\fM_{E_N,E_N}}
\|NR-RN\|\\[.2cm]
&\le\|\mu\|_{{\mathscr M}(\C\setminus\fF)}\|NR-RN\|.\quad\bl
\end{align*}

%

\medskip

Пусть $\widehat\M_\be(\fF)$ обозначает пространство функций вида
$f+az$, где $f\in\widehat\M(\fF)$, $a\in\C$. Легко видеть, что линейная функция
$az$ принадлежит пространству $\widehat\M(\fF)$, если множество $\fF$ компактно. Таким образом,
$\widehat\M_\be(\fF)=\widehat\M(\fF)$ для компактных $\fF$.
В случае неограниченного множества $\fF$  нетрудно проверить, что $f'(\be)=0$ для любой
функции $f\in\widehat\M(\fF)$. Таким образом, $\widehat\M_\be(\fF)\ne\widehat\M(\fF)$
для некомпактных множеств $\fF$. Из теоремы \ref{coint7} следует, что
$\widehat\M_\be(\fF)\subset\CL(\fF)$.

Авторам неизвестно, имеет ли место равенство $\widehat\M_\be(\fF)=\CL(\fF)$
даже для таких простых множество $\fF$, как окружность или прямая.


\

\section{\bf Коммутаторно липшицевы функции в круге и в полуплоскости}
\setcounter{equation}{0}
\label{full}

\

Мы рассмотрим здесь пространства коммутаторно липшицевых функций в единичном круге $\dd$ и в верхней полуплоскости $\C_+\df\{\z\in\C:~\re\z>0\}$. В частности, мы изложим результаты работы Э. Киссина и В.С. Шульмана \cite{KS3} и их аналоги для верхней полуплоскости.

Пусть $C_A$ обозначает диск-алгебру, т. е. пространство всех функций $f$, аналитических в открытом круге $\dd$ и непрерывных в его замыкании. В работе \cite{KS3} доказано, что 
$\CL(\dd)=\{f\in C_A:~f\in\OL(\T)\}$. Следующая теорема показывает, что при этом имеет место равенство норм.

\begin{thm}
\label{KST}
Пусть $f\in\CL(\dd)$. Тогда $f\in C_A$ и 
 $\|f\|_{\CL(\dd)}=\|f\|_{\CL(\T)}=\|f\|_{\OL(\T)}$. Если  $f\in C_A$, то $f\in\CL(\dd)$
в том и только в том случае, когда \lb$f\in\OL(\T)$.
\end{thm}

\Pf Равенство $\|f\|_{\CL(\T)}=\|f\|_{\OL(\T)}$ следует из теоремы \ref{su}.
Неравенство $\|f\|_{\CL(\T)}\le\|f\|_{\CL(\dd)}$  очевидно. Остаётся доказать, что
$\|f\|_{\CL(\dd)}\le\|f\|_{\CL(\T)}$. Можно считать, что $\|f\|_{\CL(\T)}=1$.
Тогда $\|\dg f\|_{\fM(\T\times\T)}=1$ в силу теоремы \ref{muld}. Применим теперь теорему \ref{topolx}.
Тогда мы получим два семейства $\{u_\z\}_{\z\in\T}$ и $\{v_\t\}_{\t\in\T}$
в некотором гильбертовом пространстве $\h$, непрерывно зависящие от параметров, и такие, что $\|u_\z\|\le1$, 
$\|v_\t\|\le1$  и
$(\dg f)(\z,\t)=(u_\z,v_\t)$ при всех $\z,\,\t\in\T$. Рассмотрим гармонические продолжения функций
$\z\mapsto u_\z$ и $\t\mapsto v_\t$ в единичный круг, положив
$$
u_z\df\int_{\T}\frac{1-|z|^2}{|z-\z|^2}u_\z\,d\m(\z)\quad\text{и}\quad
v_w\df\int_{\T}\frac{1-|w|^2}{|w-\t|^2}v_\t\,d\m(\t)
$$
при $z,w\in\dd$. Отметим, что интегралы понимаются, как интегралы $\h$-значных функций, непрерывных в слабой топологии.
Тогда, применяя интеграл Пуассона по переменной $\z$ к обеим частям
равенства $(\dg f)(\z,\xi)=(u_\z,v_\t)$, получаем $(\dg f)(z,\t)=(u_z,v_\t)$ при всех $z\in\clos\dd$
и  $\t\in\T$. Применяя теперь интеграл Пуассона к последнему равенству, получаем:
$(\dg f)(z,w)=(u_z,v_w)$ при всех $z\in\clos\dd$
и $w\in\clos\dd$. Теперь ясно, что 
$$
\|f\|_{\CL(\dd)}=\|\dg f\|_{\fM(\clos\dd\times\clos\dd)}\le\sup_{z\in\clos\dd}\|u_z\|\sup_{w\in\clos\dd}\|v_w\|
=\sup_{\z\in\T}\|u_\z\|\sup_{\t\in\T}\|v_\t\|=1.\quad \bl
$$

\medskip

Приведём теперь аналог теоремы \ref{npkp} для функций в единичном круге.

\begin{thm} 
\label{KSTc}
Пусть $f\in\CL(\dd)$. Тогда  в диск-алгебре $C_A$ существуют 
последовательности $\{\f_n\}_{n\ge1}$ и $\{\psi_n\}_{n\ge1}$ такие, что
$$
\left(\sup_{z\in\dd}\sum_{n=1}^\be|\f_n(z)|^2\right)\left(\sup_{w\in\dd}\sum_{n=1}^\be|\psi_n(w)|^2\right)
=\|f\|_{\CL(\dd)}^2\quad\text{и}\quad
(\dg f)(z,w)=\sum_{n=1}^\be\f_n(z)\psi_n(w),
$$
при этом все ряды сходятся равномерно, когда переменные $z$ и $w$ пробегают компактное подмножество
открытого единичного круга.
\end{thm}

\Pf Можно считать, что $\|f\|_{\CL(\dd)}=1$. Вместо первого равенства достаточно
доказать неравенство $\le$, поскольку неравенство $\ge$ вытекает из теорем \ref{muld} и 
\ref{41} Пусть $\h$, $u_z$ и $v_w$ обозначают
то же, что в доказательстве теоремы \ref{KST}. Рассмотрим ортонормированный базис 
$\{e_n\}_{n=1}^\be$ в пространстве $\h$.
Положим $\f_n(z)\df(u_z,e_n)$ и $\psi_n(w)\df(e_n,v_w)$ для всех $n\ge1$. Докажем, что 
$\f_n\in C_A$ и $\psi_n\in C_A$. Обозначим символом $X$ множество векторов 
$e\in\h$ таких, что $(u_z,e)\in C_A$. Ясно, что $X$ -- замкнутое подпространство
пространства $\h$.  Заметим, что $v_\t\in X$ для всех $\t\in\T$,
поскольку $(\dg f)(\cdot,\t)\in C_A$ при всех $\t\in\T$. Таким образом, $X=\h$,
поскольку линейная оболочка семейства $\{v_\t\}_{\t\in\T}$ плотна в $\h$.
Таким образом, $(u_z,e)\in C_A$ для всех $e\in\h$. Аналогично доказывается,
что $(e,v_w)\in C_A$ для всех $e\in\h$. Остаётся доказать утверждение, относящееся
к равномерной сходимости на компактах. Заметим, что
$$
\left|\sum_{n=N}^\be\f_n(z)\psi_n(w)\right|\le\left(\sum_{n=N}^\be|\f_n(z)|^2\right)^{\frac12}
\left(\sum_{n=N}^\be|\psi_n(w)|^2\right)^{\frac12}.
$$
Таким образом, достаточно доказать равномерную сходимость на компактных 
подмножествах круга $\dd$  рядов $\sum_{n=1}^\be|\f_n(z)|^2$ и $\sum_{n=1}^\be|\psi_n(z)|^2$.
Это вытекает из следующей элементарной леммы.

\begin{lem} 
\label{kner}
Пусть  $\{h_k\}_{k=1}^\be$ --  последовательность аналитических в единичном круге $\dd$ функций.
Предположим, что функция $\sum_{k=1}^\be|h_k(z)|$ ограничена в $\dd$. Тогда ряд
$\sum_{k=1}^\be|h_k(z)|$ сходится равномерно на компактных подмножествах
открытого единичного круга.
\end{lem}

\Pf Не умаляя общности, можно считать, что $h_k\in C_A$ при всех $k$.
Положим $\d_k\df\int_\T|h_k(\z)|\,d\m(\z)$. Ясно, что ряд $\sum_{k=1}^\be\d_k$ сходится.
Остаётся заметить, что  из формулы Коши вытекает, что $|h_k(z)|\le\frac{\d_k}{1-|z|}$. $\bl$ 

\medskip 

Обозначим через $(\OL)_+(\T)$ пространство функций $f$ из $\OL(\T)$, допускающих
аналитическое продолжение в единичный круг $\dd$, непрерывное вплоть до границы.
Из теоремы \ref{pro} следует, что каждая функция $f\in\CL(\dd)$ аналитична в
$\dd$. Таким образом, из теоремы \ref{KST} вытекает следующее
результат работы \cite{KS3}.

\begin{thm} 
\label{olclt}
Оператор сужения $f\mapsto f\big|\T$ является линейной изометрией пространства
$\CL(\dd)$ на пространство $(\OL)_+(\T)$.
\end{thm}

Аналогичные результаты имеют место и для пространства $\CL(\C_+)$. 

\begin{thm} 
\label{KSR}
Пусть $f$ -- непрерывная функция в замкнутой верхней полуплоскости $\clos\C_+$.
Предположим, что функция   $f$ аналитична в открытой полуплоскости $\C_+$. 
Тогда $\|f\|_{\CL(\C_+)}=\|f\|_{\CL(\R)}=\|f\|_{\OL(\R)}$. В частности, $f\in\CL(\C_+)$
в том и только в том случае, когда $f\in\OL(\R)$.
\end{thm}


Обозначим через $C_A(\C_+)$ множество всех аналитических в $\C_+$ функций, непрерывных 
вплоть до границы и имеющих конечный предел в бесконечности.

\begin{thm} 
\label{KSRc}
Пусть $f\in\CL(\C_+)$. Тогда  в пространстве  $C_A(\C_+)$ существуют 
последовательности $\{\f_n\}_{n=1}^\be$ и $\{\psi_n\}_{n=1}^\be$ такие, что
$$
\left(\sup_{z\in\C_+}\sum_{n=1}^\be|\f_n(z)|^2\right)\left(\sup_{w\in\C_+}\sum_{n=1}^\be|\psi_n(w)|^2\right)
=\|f\|_{\CL(\C_+)}^2
$$
и
$$
(\dg f)(z,w)=\sum_{n=1}^\be\f_n(z)\psi_n(w).
$$
При этом все ряды сходятся равномерно, когда переменные $z$ и $w$ пробегают компактное подмножество
открытой верхней полуплоскости.
\end{thm}


Мы опускаем здесь доказательства теорем \ref{KSR} и \ref{KSRc}. Они аналогичны доказательствам соответствующих результатов для функций в единичном круге.

Имеет место также следующий аналог для прямой теоремы \ref{olclt}.

\begin{thm} 
\label{olclr}
Оператор сужения $f\mapsto f\big|\R$ является линейной изометрией пространства
$\CL(\C_+)$ на пространство $(\OL)_+(\R)$.
\end{thm}

Отметим, что в статье \cite{A4} получен следующий результат, по существу содержащий в себе как теорему \ref{KST},
так и теорему \ref{KSR}.

\begin{thm}
Пусть $\fF_0$ и $\fF$ -- непустые совершенные подмножества комплексной плоскости $\C$
такие, что $\fF_0\subset\fF$ и множество $\Omega\df\fF\setminus\fF_0$ открыто.
Предположим, что функция  $f_0\in\CL(\fF_0)$  допускает непрерывное продолжение
$f$ на множество $\fF$ такое, что функция $f$ аналитична на $\Omega$ и
$$
|f(z)z^{-2}|\to0\quad\text{при}\quad z\to\be
$$
в каждой неограниченной\footnote{Последнее условие выполняется автоматически, если
множество $\Omega$ ограничено.} компоненте связности множества  $\Omega$.
Тогда
$f\in\CL(\fF)$ и $\|f\|_{\CL(\fF)}=\|f_0\|_{\CL(\fF_0)}$.
\end{thm}

Авторам неизвестен ответ на следующий вопрос.
Пусть $f$ -- непрерывная функция в замкнутом единичном круге, гармоническая 
внутри круга. Предположим, что $f\in\OL(\T)$. Вытекает ли отсюда, что функция 
$f\in\OL(\dd)$?
Иными словами, является ли оператор  сужения  $f\mapsto f\big|\T$ изоморфизмом
пространства $\OL(\dd)$ на $\OL(\T)$?
Аналогичный вопрос можно поставить для полуплоскости и для других областей.

Напомним, что если $T$ -- сжатие в гильбертовом пространстве $\h$, то по теореме
Сёкефальви--Надя (см. \cite{SNF}, глава I, \S\:5) существует его {\it унитарная дилатация}, т.е. такой унитарный оператор $U$ в гильбертовом пространстве $\K$, $\h\subset\K$, что
$T^n=P_\h U^n\big|\h$, $n\ge0$. Дилатацию всегда можно выбрать минимальной.
При этом можно определить линейное и мультипликативное функциональное исчисление:
$\f\mapsto\f(T)\df P_\h\f(U)\big|\h$, $\f\in C_A$.
{\it Полуспектральная мера $\mE_T$ сжатия} $T$ определяется равенством $\mE_T(\D)\df P_\h E_U(\D)\big|h$, где $E_U$ -- спектральная мера оператора $U$, а $\D$ -- борелевское подмножество окружности $\T$. Легко видеть, что
$$
\f(T)=\int_\T\f(\z)\,d\mE(\z),\quad\f\in C_A.
$$

\begin{thm}
\label{teorokomszha}
Пусть $f\in\CL(\dd)$ и пусть $T_1$ и $T_2$ -- сжатия в гильбертовом пространстве $\h$, а  $R\in\mB(\h)$. Тогда
\bay
\label{doifszha}
f(T_1)R-Rf(T_2)=
\int_\T\int_\T\big(\dg f)(\z,\t)\,d\mE_1(\z)(T_1R-RT_2)\,d\mE_2(\t),
\ey
где $\mE_1$ и $\mE_2$ -- полуспектральные меры сжатий $T_1$ и $T_2$,
и справедливо неравенство 
\bay
\label{KSlip}
\|f(T_1)R-Rf(T_2)\|\le\|f\|_{\CL(\dd)}\|T_1R-RT_2\|.
\ey
\end{thm}

\Pf Пусть $\{\f_n\}_{n\ge1}$ и $\{\psi_n\}_{n\ge1}$ -- последовательности функций в диск-алгебре, удовлетворяющие заключению теоремы \ref{KSTc}. Ввиду формулы
\rf{doipolcpm}

\begin{align*}
\int_\T\int_\T\big(\dg f)(\z,\t)&
\,d\mE_1(\z)(T_1R-RT_2)\,d\mE_2(\t)
=\sum_{n\ge1}\f_n(T_1)(T_1R-RT_2)\psi_n(T_2)\\[.2cm]
&=
\sum_{n\ge1}T_1\f_n(T_1)R\psi_n(T_2)-\sum_{n\ge1}\f_n(T_1)R\psi_n(T_2)T_2\\[.2cm]
&=\int_\T\int_\T\z\big(\dg f)(\z,\t)
\,d\mE_1(\z)R\,d\mE_2(\t)
-\int_\T\int_\T\t\big(\dg f)(\z,\t)
\,d\mE_1(\z)R\,d\mE_2(\t)\\[.2cm]
&=\int_\T\int_\T(f(\z)-f(\t))
\,d\mE_1(\z)R\,d\mE_2(\t)=f(T_1)R-Rf(T_2),
\end{align*}
что доказывает формулу \rf{doifszha}, которая немедленно влечёт неравенство 
\rf{KSlip}. $\bl$

\medskip

Неравенство  \rf{KSlip} было доказано Киссиным и Шульманом в \cite{KS3} другим методом. Доказательство, приведённое здесь, аналогично доказательству теоремы 4.1 в работе \cite{Pe9}, см. также \cite{Pe*}.
В случае $f\in B_{\be,1}^1(\T)\cap C_A$ и $R=I$ теорема \ref{teorokomszha} была доказана в работе \cite{Pe*}. 

Аналогичное утверждение можно доказать и для диссипативных операторов, см. \cite{AP*} по поводу возмущений функций диссипативных операторов.

\

\section{\bf Операторно липшицевы функции и дробно-линейные преобразования}
\setcounter{equation}{0}
\label{dlp}

\

Пусть $\dlC$ обозначает группу Мёбиуса дробно-линейных преобразований расширенной
комплексной плоскости $\widehat \C\df\C\cup\{\be\}$. Другими словами,
$$
\dlC=\Big\{\f:~\f(z)=\frac{az+b}{cz+d}, \quad a,\,b,\,c,\,d\in \C,\quad ad-bc\ne0\Big\}.
$$

Множество всех линейных
преобразований комплексной плоскости обозначим через $\lC$, т. е.
$$
\lC=\{\f\in\dlC:~\f(\be)=\be\}=\{\f:~\f(z)=az+b, a,b\in\C, a\ne0\}.
$$

Пусть $\widehat\R$ обозначает расширенную вещественную прямую, $\widehat\R\df\R\cup\{\be\}$.
Положим 
$$
\dlR\df\{\f\in\dlC:~\f(\widehat\R)=\widehat\R\}\quad\mbox{и} \quad
\lR\df\{\f\in\lC:~\f(\R)=\R\}.
$$

С каждым дробно-линейном преобразованием $\f$ и с каждой функцией $f$,
заданной на замкнутом множестве  $\fF,\:\fF\subset\C$, мы связываем функцию
${\mathcal Q}_\f f$, определённую на замкнутом множестве 
$\fF_\f\df\C\cap\f^{-1}(\fF\cup\{\be\})$
следующим образом:
$$
({\mathcal Q}_\f f)(z)\df\left\{\begin{array}{ll}\dfrac{f(\f(z))}{\f^{\,\prime}(z)},&\text {если}\,\,\,
z\in\C,\,\,\f(z)\in\fF\,\,\,\text{и}\,\,\, \f(z)\ne\be,\\[.3cm]
\quad 0,&\text {если}\,\,\,z\in\C\,\,\,\text{и}\,\,\, \f(z)=\be.
\end{array}\right.
$$

Легко видеть, что если $\f\in\lC$, то $\fF_\f=\f^{-1}(\fF)$, ${\mathcal Q}_\f f=(\f'(0))^{-1}(f\circ\f)$,
${\mathcal Q}_\f(\OL(\fF))=\OL(\fF_\f)$, ${\mathcal Q}_\f(\CL(\fF))=\CL(\fF_\f)$, $\|{\mathcal Q}_\f f\|_{\OL(\fF_\f)}=\|f\|_{\OL(\fF)}$
для всех $f$ из $\OL(\fF)$  и $\|{\mathcal Q}_\f f\|_{\CL(\fF_\f)}=\|f\|_{\CL(\fF)}$
для всех $f$ из $\CL(\fF)$. Поэтому в дальнейшем нас будет в основном интересовать случай,
когда $\f\notin\lC$. Заметим, что если  $\fF=\C$, то $\fF_\f=\C$  при всех $\f$ из $\dlC$.
Если же  $\fF=\R$, то $\fF_\f=\R$  при всех $\f$ из $\dlR$.

Пусть $a\in\fF$, где $\fF$ -- замкнутое подмножество  комплексной плоскости $\C$.
Положим 
$$
\OL_a(\fF)\df\{f\in\OL(\fF):f(a)=0\}\quad\mbox{и}\quad
\CL_a(\fF)\df\{f\in\CL(\fF):f(a)=0\}.
$$ 
Ясно, что $\OL_a(\fF)$ и $\CL_a(\fF)$ -- банаховы пространства.

\begin{thm}
\label{zamfone}
Пусть $\fF$ -- замкнутое подмножество в $\C$, $a\in\fF$ и
пусть $\f$ -- автоморфизм из $\dlC$ такой, что $a\df\f(\be)$.
Тогда ${\mathcal Q}_\f(\OL_a(\fF))\subset\OL(\fF_\f)$ и
$$
\|{\mathcal Q}_\f f\|_{\OL(\fF_\f)}\le3\|f\|_{\OL(\fF)}
$$
для всех $f$ из $\OL_a(\fF)$.
\end{thm}

\Pf  Рассмотрим сначала частный случай, когда $\f(z)=\phi(z)\df z^{-1}$.
Тогда $a=0$ и нам нужно доказать, что ${\mathcal Q}_\phi(\OL_0(\fF))\subset\OL(\fF_\f)$
и $\|{\mathcal Q}_\phi f\|_{\OL(\fF_\f)}\le3\|f\|_{\OL(\fF)}$.
Пусть $f\in\OL_0(\fF)$. Можно считать, что $\|f\|_{\OL(\fF)}=1$. 
Тогда всё сводится к следующему неравенству:
$$
\|({\mathcal Q}_\f f)(N)R-R({\mathcal Q}_\f f)(N)\|\le3\max(\|NR-RN\|,\|N^*R-RN^*\|)
$$
для любых ограниченных операторов $N$ и $R$ таких, что оператор $N$  нормален и $\s(N)\subset\fF_\phi$.
Рассмотрим функцию 
$$
h(z)\df\left\{\begin{array}{ll}zf(z^{-1}),&\text {если}\,\,\,
z\in\fF_\phi\setminus\{0\},\\[.2cm]
\quad0,&\text {если}\,\,\,z=0.
\end{array}\right.
$$
Легко видеть, что $\sup|h|\le\|f\|_{\Li(\fF)}\le\|f\|_{\OL(\R)}=1$,  поскольку $f(0)=0$.
Заметим, что $({\mathcal Q}_\phi f)(N)=-Nh(N)$.
Таким образом, нам нужно доказать, что
$$
\|Nh(N)R-RNh(N)\|\le3\max\{\|NR-RN\|,\|N^*R-RN^*\|\}.
$$
Мы воспользуемся следующим элементарным тождеством:
\begin{multline}
\label{42}
Nh(N)R-RNh(N)=h(N)(NR-RN)\\
+h(N)RN-NRh(N)
+(NR-RN)h(N).
\end{multline}

Заметим, что 
$$
\|h(N)(NR-RN)\|\le\|NR-RN\|\le\max\{\|NR-RN\|,\|N^*R-RN^*\|\}.
$$
Аналогично оценивается и норма оператора $(NR-RN)h(N)$.

Остаётся доказать, что
$$
\|h(N)RN-NRh(N)\|\le\max\{\|NR-RN\|,\|N^*R-RN^*\|\}.
$$

Если оператор $N$ обратим, то
\begin{multline*}
\|h(N)RN-NRh(N)\|=\|f(N^{-1})NRN-NRNf(N^{-1}\|)\\
\le
\max\{\|RN-NR\|,\|(N^*)^{-1}NRN-NRN(N^*)^{-1}\|\}\\
=\max\{\|NR-RN\|,\|(N^*)^{-1}N(RN^*-N^*R)N(N^*)^{-1}\|\}
\\
=\max\{\|NR-RN\|,\|N^*R-RN^*\|\}.
\end{multline*}

Если $0$ --  предельная точка множества $\fF_\phi$ (т. е. множество $\fF$ неограничено), то на этом доказательство можно закончить, ибо в этом случае каждый нормальный оператор $N$ со спетром  в $\fF_\phi$
приближается с любой точностью нормальным оператором $M$ таким, что 
$MN=NM$ и $\s(M)\in\fF_\f\setminus\{0\}$. Это следует, например, из леммы \ref{apr}.

Пусть теперь $0$ --  изолированная точка множества $\fF_\phi$.
Рассмотрим необратимый нормальный оператор $N$ со спектром $\fF_\phi$.
Тогда оператор $N$ представим в виде $N=\0\oplus N_0$,  где $N_0$ -- обратимый нормальный оператор.
Заметим, что ${\mathcal Q}_\phi(N)=\0\oplus N_0^2f(N_0^{-1})$. Пусть $P$ обозначает ортогональный 
проектор на подпространство, в  котором действует оператор $N_0$. Легко видеть,
что
\begin{multline*}
\|h(N)RN-NRh(N)\|=\|P(h(N)RN-NRh(N))P\|\\=\|h(N)PRPN-NPRPh(N)\|
=\|h(N_0)(PRP)N_0-N_0(PRP)h(N_0)\|\\
\le\max(\|N_0(PRP)-(PRP)N_0\|,\|N_0^*(PRP)-(PRP)N_0^*\|\\
\le\max(\|NR-RN\|,\|N^*R-RN^*\|).
\end{multline*}

Переходим к общему случаю. Положим $b=\f^{-1}(\be)$.
Ясно, что $\f(z)=a+c\phi(z-b)$, где $c\in\C\setminus\{0\}$. Таким образом,
всё сводится к случаю, когда $a=b=0$, т. е.  $\f=c\phi$, поскольку параллельные переносы 
сохраняют операторно липшицеву норму. Наконец, случай $\f=c\phi$ легко сводится
к разобранному случаю $\f=\phi$. $\bl$

\medskip

{\bf Пример.} Пусть $\f(z)=z^{-1}$, $\fF=\C$, $f=\ov z$.
Тогда $f\in\OL_{0}(\C)$ и $\|f\|_{\OL(\C)}=1$. Кроме того, $({\mathcal Q}_\f f)(z)=-\ov z^{-1}z^2$
и
$$
3=\|{\mathcal Q}_\f f\|_{\Li(\T)}\le\|{\mathcal Q}_\f f\|_{\Li(\C)}\le\|{\mathcal Q}_\f f\|_{\OL(\C)}\le3\|f\|_{\OL(\C)}=3.
$$
Этот пример показывает, что $\|{\mathcal Q}_\f f\|_{\OL(\C)}=\|{\mathcal Q}_\f f\|_{\Li(\C)}=3$ и константа
$3$ в теореме \ref{zamfone} неулучшаема.

\medskip

Из теоремы \ref{zamfone} легко вытекает следующее утверждение.

\begin{thm}
\label{zamf}
Пусть $\f\in\dlC$, $a=\f(\be)$ и $b=\f^{-1}(\be)$. Пусть $\fF$ -- замкнутое
множество в $\C$, содержащее точку $a$.
Тогда ${\mathcal Q}_\f(\OL_a(\fF))=\OL_b(\fF_\f)$ и
$$
\frac13\|f\|_{\OL(\fF)}\le\|{\mathcal Q}_\f f\|_{\OL(\fF_\f)}\le3\|f\|_{\OL(\fF)}
\quad\mbox{для всех}\quad f~\mbox{~из~}~\OL_a(\fF).
$$
\end{thm}

\Pf Заметим, что $({\mathcal Q}_\f(\OL_a(\fF)))(b)=0$. Таким образом, из теоремы \ref{zamfone}
следует, что ${\mathcal Q}_\f(\OL_a(\fF))\subset\OL_b(\fF_\f)$ и $\|{\mathcal Q}_\f f\|_{\OL(\fF_\f)}\le3\|f\|_{\OL(\fF)}$.
Чтобы доказать, что ${\mathcal Q}_\f(\OL_a(\fF))\supset\OL_b(\fF_\f)$, и получить требуемую
оценку снизу для $\|{\mathcal Q}_\f f\|_{\OL(\fF_\f)}$, достаточно применить теорему \ref{zamfone}
к замкнутому множеству $\fF_\f$  и дробно-линейному преобразованию $\f^{-1}$. 
$\bl$

\medskip

Приведём ещё один близкий результат.

\begin{thm}
\label{zamfcom}
Пусть $\f\in\dlC\setminus\lC$, $a=\f(\be)$. Пусть $\fF$ -- замкнутое
подмножество комплексной плоскости $\C$ такое, что $a\not\in\fF$. Если $z_0$ -- одна из ближайших к $a$
точек множества $\fF$,
то ${\mathcal Q}_\f(\OL_{z_0}(\fF))\subset\OL(\fF_\f)$ и
$$
\|{\mathcal Q}_{\f}f\|_{\OL(\fF_\f)}
\le 5\|f\|_{\OL(\fF)}
$$
для всех $f\in\OL_{z_0}(\fF)$.
\end{thm}

\Pf 
Так же, как в доказательстве теоремы \ref{zamfone},  достаточно ограничиться
случаем $\f(z)=\phi(z)\df z^{-1}$.  Пусть $f\in\OL_{z_0}(\fF)$ и $\|f\|_{\OL_{z_0}(\fF)}=1$.
Нам нужно доказать, что
$$
\|({\mathcal Q}_\f f)(N)R-R({\mathcal Q}_\f f)(N)\|\le5\max(\|NR-RN\|,\|N^*R-RN^*\|)
$$
для любых нормальных операторов $N_1$ и $N_2$ таких, что $\s(N_1),\;\s(N_2)\subset\fF_\phi$.
Пусть $h$ обозначает то же, что и в доказательстве теоремы \ref{zamfone}. Только 
теперь мы не можем утверждать, что $\sup|h|\le1$. Имеем
\begin{multline*}
\sup_{z\in\fF_\phi}|h(z)|\le\sup\{|zf(z^{-1})|:z\in\phi^{-1}(\fF)\}=\sup\{|z|^{-1}|f(z)-f(z_0)|:z\in\fF\}\\
\le\sup\{|z|^{-1}|z-z_0|:z\in\fF\}\le\sup\{1+|z|^{-1}|z_0|:z\in\fF\}=2.
\end{multline*}
Повторяя теперь аргументы доказательства теоремы \ref{zamfone}, получаем
\begin{align*}
\|({\mathcal Q}_\f f)(N)R-R({\mathcal Q}_\f f)(N)\|
&\le(1+2\sup|h(z)|)\max(\|NR-RN\|,\|N^*R-RN^*\|)\\[.2cm]
&\le5\max(\|NR-RN\|,\|N^*R-RN^*\|). \quad\bl
\end{align*}

\medskip

{\bf Пример.} Пусть $\f(z)=z^{-1}$, $\fF=\T$, $z_0=1$, $f=1-\ov z$.
Тогда $f\in\OL_{z_0}(\T)$ и $\|f\|_{\OL(\T)}=1$. Легко проверить, что $({\mathcal Q}_\f f)(z)=z^3-z^2$ и $\|{\mathcal Q}_\f f\|_{\Li(\T)}\ge5$. Тогда
$$
5\le\|{\mathcal Q}_\f f\|_{\Li(\T)}\le\|{\mathcal Q}_\f f\|_{\OL(\T)}\le\|z^3\|_{\OL(\T)}+\|z^2\|_{\OL(\T)}=5.
$$
Этот пример показывает, что константа
$5$ в теореме \ref{zamfcom} неулучшаема.

\medskip

{\bf Замечание 1.} Можно ввести следующее обобщение оператора ${\mathcal Q}_\f$, 
положив
$$
({\mathcal Q}_{n,\f} f)(z)\df\left\{\begin{array}{ll}\dfrac{|\f'(z)|^nf(\f(z))}{(\f^{\,\prime}(z))^{n+1}},&\text {если}\,\,\,
z\in\C,\,\,\f(z)\in\fF\,\,\,\text{и}\,\,\, \f(z)\ne\be,\\[.3cm]
\quad 0,&\text {если}\,\,\,z\in\C\,\,\,\text{и}\,\,\, \f(z)=\be,
\end{array}\right.
$$
где $n\in\Z$.
Имеют место аналоги теорем \ref{zamfone}, \ref{zamf} и \ref{zamfcom} для операторов 
${\mathcal Q}_{n,\f}$ с константами, зависящими от $n$. Аналоги теорем \ref{zamfone} и \ref{zamf}
можно найти в \cite{A1}. Аналог теоремы  \ref{zamfcom}  может быть получен таким же образом.

\medskip

{\bf Замечание 2.} Доказательства теорем \ref{zamfone}, \ref{zamf} и \ref{zamfcom} работают
и для пространств коммутаторно липшицевых функций. Случай теорем \ref{zamfone} и \ref{zamf} 
разобран в \cite{A1}. Разумеется, в случае пространств $\CL(\fF)$  речь может идти 
об обобщениях на операторы ${\mathcal Q}_{n,\f}$ (см. замечание 1) только для ``тонких'' множеств $\fF$.
Например, если множество $\fF$ содержит внутренние точки, то такие обобщения
невозможны, поскольку функции класса $\CL(\fF)$ аналитичны на внутренности
множества $\fF$.

\medskip

Далее нас в основном будут интересовать случай, когда $\fF=\R$
и $\fF=\T$. В этих случаях имеет место равенство $\CL(\fF)=\OL(\fF)$ вместе с равенством полунорм.

Отметим, что из теоремы \ref{zamfone} вытекает следующее утверждение.

\begin{thm}
\label{zamrt}
Пусть $\f\in\dlC$. Предположим, что $\f(\widehat\R)=\T$. Тогда
$$
\|{\mathcal Q}_\f f\|_{\OL(\R)}\le3\|f\|_{\OL(\T)}
$$
для всех $f\in\OL_a(\T)$, где $a=\f(\be)$.
\end{thm}

\Pf Применим теорему \ref{zamf} к $\fF=\T$. Тогда $\fF_\f=\R\cup\{\f^{-1}(\be)\}$,
и мы имеем
$$
\|{\mathcal Q}_\f f\|_{\OL(\R)}\le\|{\mathcal Q}_\f f\|_{\OL(\R\cup\{\f^{-1}(\be)\})}\le3\|f\|_{\OL(\T)}
$$
для всех $f\in\OL_a(\T)$. $\bl$

\

\section{\bf Производные операторно липшицевых функций и дробно-линейные преобразования}
\setcounter{equation}{0}
\label{difdlp}

\


Пусть $\OL'(\R)\df\{f':f\in\OL(\R)\}$ и $\|f'\|_{\OL'(\R)}\df\|f\|_{\OL(\R)}$.
Тогда $\OL'(\R)$ -- банахово пространство,
состоящее из функций, заданных на расширенной прямой $\widehat\R$.

\begin{thm}
\label{dxa}
Пусть $f\in\OL(\R)$. Тогда $(x-a)^{-1}(f(x)-f(a))\in(\OL)'(\R)$ и
$$
\left\|\dfrac{f(x)-f(a)}{x-a}\right\|_{(\OL)'(\R)}\le\|f\|_{\OL(\R)}\quad
\mbox{при всех}\quad a\in\R.
$$
\end{thm}

\Pf Достаточно рассмотреть случай $a=0$ и $f(0)=0$.  Положим
$$
F(x)=\int_0^x\frac{f(t)}t\,dt=\int_0^1\frac{f(tx)}t\,dt.
$$
Нам нужно доказать, что $F\in\OL(\R)$ и $\|F\|_{\OL(\R)}\le\|f\|_{\OL(\R)}$.
Заметим, что при всех $t$ из $(0,1]$  функция $x\mapsto t^{-1}{f(tx)}$ принадлежит
банаховому пространству $\OL_0(\R)$ (см. определение в \S\:\ref{dlp}) и 
$\|t^{-1}{f(tx)}\|_{\OL(\R)}=\|f\|_{\OL(\R)}$
при всех $t$ из $(0,1]$. Следовательно,
$$
\|F\|_{\OL(\R)}\le\int_0^1\|t^{-1}f(tx)\|_{\OL(\R)}\,dt=\|f\|_{\OL(\R)}. \quad\bl
$$

\medskip

{\bf Замечание.} Аналогично можно доказать, что для любого замкнутого невырожденного
промежутка $J$ и для любой функции $f$ из $\OL(J)$ справедливо неравенство
$$
\left\|\dfrac{f(x)-f(a)}{x-a}\right\|_{(\OL)'(J)}\le\|f\|_{\OL(J)}
\quad\mbox{при всех}\quad a\in J,
$$
где $\OL'(J)\df\{g':g\in\OL(J)\}$ и $\|g'\|_{\OL'(J)}\df\|g\|_{\OL(J)}$.

\begin{thm}
\label{delab}
Пусть $f\in\OL(\R)$. Тогда $(x-a-b{\rm i})(f(x)-f(a))\in(\OL)'(\R)$ и
$$
\left\|\dfrac{f(x)-f(a)}{x-a-b{\rm i}}\right\|_{(\OL)'(\R)}\le2\|f\|_{\OL(\R)} 
\quad\mbox{при всех}\quad a,b\in\R.
$$
\end{thm}

\Pf Достаточно рассмотреть случай, когда $a=0$, $b=1$, $f(0)=0$ и $\|f\|_{\OL(\R)}=1$.  
Тогда из теоремы \ref{dxa} следует, что 
$$
\left\|\dfrac{f(x)}{x-{\rm i}}\right\|_{(\OL)'(\R)}\le\left\|\dfrac{xf(x)}{x-{\rm i}}\right\|_{\OL(\R)}.
$$
Остаётся доказать, что $\left\|\dfrac{xf(x)}{x-{\rm i}}\right\|_{\OL(\R)}\le2$. Пусть
$A$ и $B$ -- самосопряжённые операторы. Имеем
\begin{align*}
A(A-{\rm i}I)^{-1}f(A)-B(B-{\rm i}I)^{-1}f(B)
&=A(A-{\rm i}I)^{-1}(f(A)-f(B))\\
&+
{\rm i}(A-{\rm i}I)^{-1}(B-A)(B-{\rm i}I)^{-1}f(B),
\end{align*}
откуда
$$
\|A(A-{\rm i}I)^{-1}f(A)-B(B-{\rm i}I)^{-1}f(B)\|\\
\le\|f(A)-f(B)\|+\|A-B\|\cdot\|g(B)\|\le2\|A-B\|,
$$
где $g(t)=(t-{\rm i})^{-1}f(t)$. $\bl$

\begin{cor}
\label{cordelab}
Пусть $f\in\OL(\R)$. Тогда $(x-a-b{\rm i})^{-1}f(x)\in(\OL)'(\R)$ и
$$
\left\|\dfrac{f(x)}{x-a-b{\rm i}}\right\|_{(\OL)'(\R)}
\le\left(2+\frac{|f(a)|}{|b|}\right)\|f\|_{\OL(\R)} 
\quad\mbox{при всех}\quad a,b\in\R,~\;b\ne0.
$$
\end{cor}

\Pf  Можно считать, что $a=0$, $b=1$ и $\|f\|_{\OL(\R)}=1$. 
Используя теорему \ref{delab} и пример 2 в  \S\:\ref{prim}, получаем:
$$
\left\|\dfrac{f(x)}{x-{\rm i}}\right\|_{(\OL)'(\R)}\le\left\|\dfrac{f(x)-f(0)}{x-{\rm i}}\right\|_{(\OL)'(\R)}
+|f(0)|\cdot\|(x-{\rm i})^{-1}\|_{(\OL)'(\R)}\le2+|f(0)|.\quad\bl
$$

\begin{thm}
\label{rinv}
Пусть $h\in\OL'(\R)$. Тогда $h\circ\f\in\OL'(\R)$ для всех дробно-линейных
преобразований $\f\in\dlR$ и
$$
\frac19\|h\|_{\OL'(\R)}\le\|h\circ\f\|_{\OL'(\R)}\le9\|h\|_{\OL'(\R)}.
$$
\end{thm}

\Pf Теорема очевидна, если $\f\in\lR$. В этом случае
$\|h\|_{\OL'(\R)}=\|h\circ\f\|_{\OL'(\R)}=\|h\|_{\OL'(\R)}$.  Таким образом, всё сводится к случаю
$\f(t)=\phi(t)\df t^{-1}$.
Пусть $h=f'$ для некоторой функции $f\in\OL(\R)$ такой, что $f(0)=0$ и $\|f\|_{\OL(\R)}=\|h\|_{\OL'(\R)}$.
Из теоремы \ref{zamf} следует, что $\|x^2f(x^{-1})\|_{\OL(\R)}\le3\|h\|_{\OL'(\R)}$. Следовательно,
$$
\|(x^2f(x^{-1}))'\|_{\OL'(\R)}=\|2xf(x^{-1})-h(x^{-1})\|_{\OL'(\R)}\le3\|h\|_{\OL'(\R)}.
$$
Из теоремы \ref{dxa} вытекает следующее неравенство:
$$
\|xf(x^{-1})\|_{\OL'(\R)}\le\|x^2f(x^{-1})\|_{\OL(\R)}\le3\|h\|_{\OL'(\R)}.
$$
Следовательно,
$$
\|h(x^{-1})\|_{\OL'(\R)}\le\|(x^2f(x^{-1}))'\|_{\OL'(\R)}+2\|xf(x^{-1})\|_{\OL'(\R)}\le9\|h\|_{\OL'(\R)}.
$$
Применяя теперь это неравенство к функции $h(x^{-1})$, получаем:
$$
\frac19\|h(x^{-1})\|_{\OL'(\R)}\le\|h\|_{\OL'(\R)}.\quad\bl
$$

\

\section{\bf Пространства  $\bs{\OL(\R)}$ и $\bs{\OL(\T)}$}
\setcounter{equation}{0}
\label{olriort}

\

Основная цель этого и следующего параграфа состоит в том, чтобы
``пересадить'' теорему \ref{rinv} с прямой на окружность.

Легко видеть, что если  $f\in\OL(\T)$, то $f(e^{{\rm i}t})\in\OL(\R)$ и
$\|f(e^{{\rm i}t})\|_{\OL(\R)}\le\|f\|_{\OL(\T)}$. 
Мы покажем здесь, что верно и обратное, т.е.
любая $(2\pi)$-периодическая
функция $F$ из $\OL(\R)$ представима в виде $F=f(e^{{\rm i}t})$, где
$f\in\OL(\T)$ и $\|f\|_{\OL(\T)}\le\const\|F\|_{\OL(\R)}$. Этот факт 
нетрудно вывести из леммы 9.8 статьи \cite{AP5}, см. также лемму 5.7 в \cite{A1}.

\begin{lem}
\label{111}
Пусть $h(x,y)=\dfrac{x-y}{e^{{\rm i}x}- e^{{\rm i}y}}$. Тогда
$$
\|h\|_{\fM(J_1\times J_2)}\le\frac{3\sqrt2\pi}4
$$
для любых промежутков  $J_1$ и $J_2$ таких, что 
$J_1-J_2\subset[-\frac32\pi,\frac32\pi]$.
\end{lem}


\Pf Рассмотрим $3\pi$-периодическую функцию $\xi$ такую, что 
$\xi(t)=t(2\sin(t/2))^{-1}$
при $t\in[-\frac32\pi,\frac32\pi]$. Тогда
$$
\|h\|_{\fM(J_1\times J_2)}=\big\|e^{\frac{{\rm i}x}2}h(x,y)e^{\frac{{\rm i}y}2}\big\|_{{\fM}(J_1\times J_2)}
=\|\xi(x-y)\|_{{\fM}(J_1\times J_2)},
$$
поскольку $x-y\in[-\frac32\pi,\frac32\pi]$ при $x\in J_1$ и $y\in J_2$. Разложим функцию $\xi$ в ряд Фурье:
$$
\xi(t)=\sum_{n\in\Z}a_ne^{\frac23n{\rm i}t}=a_0+2\sum_{n=1}^\be a_n\cos\frac23nt,
$$
поскольку $a_n=a_{-n}$ при всех $n\in\Z$. Ясно, что $a_0>0$. Заметим, что функция $\xi$
выпукла на промежутке $[-\frac32\pi,\frac32\pi]$. Отсюда следует, что $(-1)^na_n\ge0$
при всех натуральных $n$,  см. теорему 35 в монографии \cite{HR}.
Остаётся заметить, что
\begin{align*}
\|\xi(x-y)\|_{{\fM}(J_1\times J_2)}&\le\|\xi(x-y)\|_{{\fM}(\R\times\R)}
\le\sum_{n\in\Z}|a_n|\cdot\big\|e^{\frac23n{\rm i}x}e^{-\frac23n{\rm i}y}\big\|_{{\fM}(\R\times\R)}\\
&=\sum_{n\in\Z}|a_n|=\xi\left(\frac{3\pi}2\right)=\frac{3\sqrt2\,\pi}4.\quad\bl
\end{align*}

\begin{thm}
\label{2pper}
Пусть $f$ --  непрерывная функция на  $\T$.
Тогда
$$
\|f(e^{{\rm i}x})\|_{\OL(\R)}\le\|f\|_{\OL(\T)}\le\frac{3\sqrt2\pi}2\|f(e^{{\rm i}x})\|_{\OL(\R)}.
$$
\end{thm}

\Pf Как было отмечено выше, первое неравенство очевидно. Докажем второе неравенство.
Положим $g(x)\df f(e^{{\rm i}x})$. Можно считать, что $\|g\|_{\OL(\R)}<\be$.
Тогда функция
$g$ дифференцируема всюду на $\R$. Отсюда следует, что
функция $f$ имеет производную всюду на  $\T$. Заметим,
что из теорем \ref{su}  и \ref{muld} следует, что
$$
\|g\|_{\OL(\R)}=\|\dg g\|_{\fM(\R\times\R)}\;~{\text и}
~\;\|f\|_{\OL(\T)}=\|\dg f\|_{\fM(\T\times\T)}=\|(\dg f)(e^{{\rm i}x},e^{{\rm i}y})\|_{\fM([0,2\pi)\times[-\frac\pi2,\frac{3\pi}2))}
$$
Таким образом, нам нужно доказать, что
$$
\Big\|(\dg f)(e^{{\rm i}x},e^{{\rm i}y})\Big\|_{\fM([0,2\pi)\times[-\frac\pi2,\frac{3\pi}2))}
\le\frac{3\sqrt2\pi}2\|\dg g\|_{\fM(\R\times\R)}.
$$

Обозначим здесь символом $\chi_{jk}$ характеристическую функцию множества
$\bs{J}_{j,k}\df[j\pi,(j+1)\pi)\times[k\pi-\frac\pi2,k\pi+\frac{\pi}2)$,
где $j,k\in\Z$.
Заметим, что
$$
\chi_{jk}(x,y)(\dg f)(e^{{\rm i}x},e^{{\rm i}y})
=\chi_{jk}(x,y)h(x,y)(\dg g)(x,y),
$$
где $h$ обозначает то же, что в лемме \ref{111}.
Отсюда и из леммы \ref{111} следует, что
\bay
\label{n111}
\Big\|(\dg f)(e^{{\rm i}x},e^{{\rm i}y})\Big\|_{\fM(\bs{J}_{j,k})}
\le\frac{3\sqrt2\pi}4\|\dg g\|_{\fM(\R\times\R)}
\ey
при $(j,k)\in\{0,1\}$, $(j,k)\ne(1,0)$.

Случай, когда $j=1$  и $k=0$, следует рассмотреть отдельно, поскольку в этом случае
$J_1-J_2\not\subset[-\frac{3\pi}2,\frac{3\pi}2]$, и  мы не можем применить лемму \ref{111}
напрямую.

Пусть теперь $j=1$ и  $k=0$. Имеем:
$$
\chi_{10}(x+2\pi,y)(\dg f)(e^{{\rm i}x},e^{{\rm i}y})
=\chi_{10}(x+2\pi,y)h(x,y)(\dg g)(x,y).
$$
Теперь, применяя лемму \ref{111},  получаем:
$$
\|(\dg f)(e^{{\rm i}x},e^{{\rm i}y})\|_{\fM(\bs{J}_{1,0})}
=\|(\dg f)(e^{{\rm i}x},e^{{\rm i}y})\|_{\fM(\bs{J}_{-1,0})}\le\frac{3\sqrt2\pi}4\|\dg g\|_{\fM(\R\times\R)}.
$$
Положим также $\bs{J}\df[0,2\pi)\times[-\frac\pi2,\frac{3\pi}2)$. Тогда
\begin{align*}
\big\|(\dg f&)(e^{{\rm i}x},e^{{\rm i}y})\big\|_{\fM(\bs{J})}\le
\big\|(\chi_{00}(x,y)+\chi_{11}(x,y))(\dg f)(e^{{\rm i}x},e^{{\rm i}y})\big\|_{\fM(\bs{J})}\\
&+
\big\|(\chi_{01}(x,y)+\chi_{10}(x,y))(\dg f)(e^{{\rm i}x},e^{{\rm i}y})\big\|_{\fM(\bs{J})}\\
&\le\max\big\{\big\|(\dg f)(e^{{\rm i}x},e^{{\rm i}y})\big\|_{\fM(\bs{J}_{0,0})},
\big\|(\dg f)(e^{{\rm i}x},e^{{\rm i}y})\big\|_{\fM(\bs{J}_{1,1})}\big\}\\
&+
\max\big\{\big\|(\dg f)(e^{{\rm i}x},e^{{\rm i}y})\big\|_{\fM(\bs{J}_{0,1})},
\big\|(\dg f)(e^{{\rm i}x},e^{{\rm i}y})\big\|_{\fM(\bs{J}_{1,0})}\big\}
\le\frac{3\sqrt2\pi}2\|\dg g\|_{\fM(\R\times\R)}.~\;\bl
\end{align*}

\medskip

{\bf Замечание.} Из доказательства теоремы видно, что 
$$
\|f(e^{{\rm i}x})\|_{\OL(\R)}\le\|f\|_{\OL(\T)}\le\frac{3\sqrt2\pi}2\|f(e^{{\rm i}x})\|_{\OL([-\pi,2\pi])}
$$
для любой  функции $f$ из $C(\T)$. Ясно, что промежуток $[-\pi,2\pi]$ может быть
заменён любым промежутком  длины $3\pi$. Можно несколько модифицировать доказательство теоремы \ref{2pper} и доказать, что
$$
\|f(e^{{\rm i}x})\|_{\OL(\R)}\le\|f\|_{\OL(\T)}\le c(|J|)\|f(e^{{\rm i}x})\|_{\OL(J)}
$$
для любой  функции $f$ из $C(\T)$  и для любого промежутка $J$ длины больше, чем $2\pi$.

\

\section{\bf Пространства $\bs{(\OL)'(\R)}$ и $\bs{(\OL)_{\rm loc}'(\T)}$}
\setcounter{equation}{0}
\label{olriolt+}

\

Пространство $(\OL)'(\R)$ было определено в \S\:\ref{difdlp}.  Определим теперь пространство $(\OL)_{\rm loc}'(\T)$ следующим образом:
$$
(\OL)_{\rm loc}'(\T)\df\left\{f: f(e^{{\rm i}t})\in(\OL)'(\R)\right\}\quad\text{и}\quad
\|f\|_{(\OL)_{\rm loc}'(\T)}\df\left\|f(e^{{\rm i}t})\right\|_{(\OL)'(\R)}.
$$
Заметим, что
$$
\|f\|_{L^\be(\T)}=\|f(e^{{\rm i}t})\|_{L^\be(\R)}\le\|f(e^{{\rm i}t})\|_{(\OL)'(\R)}
=\|f\|_{(\OL)_{\rm loc}'(\T)}.
$$

Нам понадобится следующая элементарная лемма. 

\begin{lem} 
\label{l132}
Пусть $f\in\Li(\T)$. Тогда $f\in(\OL)_{\rm loc}'(\T)$
и
$$
\|f\|_{(\OL)_{\rm loc}'(\T)}\le|\widehat f(0)|+\frac\pi{\sqrt3}\|f\|_{\Li(\T)}.
$$
\end{lem}

\Pf Заметим, что $\|f'\|_{L^2(\T)}\le\|f'\|_{L^\be(\T)}\le\|f\|_{\Li(\T)}$ и $\|z^n\|_{(\OL)'(\T)}=1$
при всех $n\in\Z$. Следовательно,
\begin{align*}
\|f\|_{(\OL)_{\rm loc}'(\T)}&\le\sum_{n\in\Z}|\widehat f(n)|\le|\widehat f(0)|+\Big(\sum_{n\ne0}n^2|\widehat f(n)|^2\Big)^{\frac12}
\Big(\sum_{n\ne0}\frac1{n^2}\Big)^{\frac12}\\
&=|\widehat f(0)|\!+\!\frac\pi{\sqrt3}\|f'\|_{L^2(\T)}\le|\widehat f(0)|\!+\!
\frac\pi{\sqrt3}\|f'\|_{L^\be(\T)}
\le|\widehat f(0)|\!+\!\frac\pi{\sqrt3}\|f\|_{\Li(\T)}. ~\bl
\end{align*}

\begin{cor}
\label{132}
Пространство $\OL(\T)$ содержится в  $(\OL)_{\rm loc}'(\T)$ и
$$
\|f\|_{(\OL)_{\rm loc}'(\T)}\le|\widehat f(0)|+\frac\pi{\sqrt3}\|f\|_{\OL(\T)}.
$$
\end{cor}

{\bf Замечание.}  Из доказательство леммы \ref{l132} видно, что имеет место следующее
неравенство:
$$
\|f\|_{(\OL)_{\rm loc}'(\T)}\le|\widehat f(0)|+\frac\pi{\sqrt3}\|f'\|_{L^2(\T)}.
$$

\begin{thm}
\label{rt}
Если  $f\in\OL(\T)$, то $zf'(z)\in(\OL)_{\rm loc}'(\T)$  и
$\|zf'(z)\|_{(\OL)_{\rm loc}'(\T)}\le\|f\|_{\OL(\T)}$.
Если $f\in(\OL)_{\rm loc}'(\T) $ и $\int_\T f(z)\,d\m(z)=0$, то существует функция
$F$ из $\OL(\T)$ такая, что $zF'(z)=f$ и $\|F\|_{\OL(\T)}\le\const\|f\|_{(\OL)_{\rm loc}'(\T)}$.
\end{thm}

\Pf Первое утверждение очевидно, поскольку если $f\in\OL(\T)$, то
$$
\int_0^x e^{{\rm i}t}f'(e^{{\rm i}t})\,dt={\rm i}f(1)-{\rm i}f(e^{{\rm i}x})
$$
и $\|f'\|_{(\OL)_{\rm loc}'(\T)}=\|f(e^{{\rm i}x})\|_{\OL(\R)}\le\|f\|_{\OL(\T)}$.

 Докажем теперь второе утверждение. Положим
 $F(e^{{\rm i}x})\df{\rm i}\int_0^x f(e^{{\rm i}t})\,dt$. Это определение функции
 $F$  корректно, поскольку $\int_0^{2\pi} f(e^{{\rm i}t})\,dt=2\pi\int_\T f(z)\,d\m(z)=0$.
 Ясно, что $zF'(z)=f(z)$.  Остаётся заметить, что $\|f\|_{(\OL)_{\rm loc}'(\T)}=\|F(e^{{\rm i}x})\|_{\OL(\R)}$
 и воспользоваться теоремой \ref{2pper}.  $\bl$

 \begin{cor}
 \label{rt2}
 Функция $f$  на  $\T$  принадлежит пространству
 $(\OL)_{\rm loc}'(\T)$ в том и только в том случае, когда она представима в виде
 $f=\widehat f(0)+zF'(z)$, где $F\in\OL(\T)$. При этом
 $$
 \|f\|_{(\OL)_{\rm loc}'(\T)}
 \le|\widehat f(0)|+\|F\|_{\OL(\T)}\le\const\|f\|_{(\OL)_{\rm loc}'(\T)}.
 $$
 \end{cor}

 \Pf Легко видеть, что $\|1\|_{(\OL)_{\rm loc}'(\T)}=1$.  Отсюда и из теоремы \ref{rt} следует,
 что если  $f=\widehat f(0)+zF'(z)$ для некоторой функции $F$ из $\OL(\T)$,
 то \lb$f\in(\OL)_{\rm loc}'(\T)$ и
 \begin{align*}
\|f\|_{(\OL)_{\rm loc}'(\T)}&\le\|\widehat f(0)+zF'(z)\|_{(\OL)_{\rm loc}'(\T)}\le|\widehat f(0)|+\|zF'(z)\|_{(\OL)_{\rm loc}'(\T)}\\
&\le\|f\|_{(\OL)_{\rm loc}'(\T)}+
\|F\|_{(\OL)(\T)}\le c\|f\|_{(\OL)_{\rm loc}'(\T)}.
 \end{align*}
 Пусть $f\in(\OL)_{\rm loc}'(\T)$. 
 Тогда в силу теоремы \ref{rt} функция $f-\widehat f(0)$
 представима в виде $f-\widehat f(0)=zF'(z)$, где $F\in\OL(\T)$. $\bl$

\begin{cor}
\label{rt3}
Пусть $f\in(\OL)_{\rm loc}'(\T)$. Тогда $z^nf(z)\in(\OL)_{\rm loc}'(\T)$ для всех $n$ из $\Z$.
\end{cor}

\Pf Ясно, что
достаточно рассмотреть случай, когда $f=zF'(z)$,  где $F\in\OL(\T)$.
Тогда
$$
z^nf(z)=z^{n+1}F'(z)=z(z^nF(z))'-nz^nF(z)\in(\OL)_{\rm loc}'(\T),
$$
поскольку $z^nF(z)\in\OL(\T)$ и $\OL(\T)\subset(\OL)_{\rm loc}'(\T)$ в силу следствия \ref{132}. $\bl$

 \begin{cor}
 \label{rt4}
 Функция $f$ на $\T$ принадлежит пространству
 $(\OL)_{\rm loc}'(\T)$ в том и только в том случае, когда она представима в виде
 $f=\widehat f(-1)z^{-1}+F'(z)$, где $F\in\OL(\T)$, при этом
 $$
 \|f\|_{(\OL)_{\rm loc}'(\T)}
 \le|\widehat f(-1)|+\|F\|_{\OL(\T)}\le\const\|f\|_{(\OL)_{\rm loc}'(\T)}.
 $$
 \end{cor}

 \Pf Положим $g(z)\df zf(z)$.  Тогда $\widehat f(-1)=\widehat g(0)$ и
 $g(z)=\widehat g(0)+zF'(z)$. Остаётся сослаться на следствия \ref{rt2} и  \ref{rt3}. $\bl$

\begin{lem}
\label{mol}
Пусть $f,g\in\OL(J)$ где $J$ -- замкнутый ограниченный промежуток
вещественной прямой $\R$. Тогда $fg\in\OL(J)$ и 
$$
\|fg\|_{\OL(J)}\le\Big(\m(J)\|g\|_{\OL(J)}+\max_{J}|g|\Big)\|f\|_{\OL(J)}.
$$
\end{lem}

\Pf Результат вытекает из очевидного неравенства
$$
\|fg\|_{\OL(J)}\le\Big(\m(J)\|g\|_{\OL(J)}+\max_{J}|g|\Big)\|f\|_{\OL(J)}.
\quad\bl
$$

\begin{lem}
\label{zetdzeta}
Пусть $f\in\OL(\T)$ и $\z\in\T$. Тогда 
$\dfrac{f(z)-f(\z)}{z-\z}\in(\OL)_{\rm loc}'(\T)$
и
$$
\left\|\dfrac{f(z)-f(\z)}{z-\z}\right\|_{(\OL)_{\rm loc}'(\T)}
\le\const\|f\|_{\OL(\T)}.
$$
\end{lem}

\Pf Достаточно рассмотреть случай $\z=1$.  Можно считать, что $f(1)=0$.
Нам нужно оценить $\OL(\R)$-полунорму функции $\Phi$,
$$
\Phi(x)\df\int_0^x\frac{f(e^{{\rm i}t})}{e^{{\rm i}t}-1}\,dt.
$$
Ясно, что функция $\Phi$ представима в виде $\Phi(x)=\l x+\Phi_0(x)$,
где $\Phi_0$ -- функция с периодом  $2\pi$. 
Заметим, что
$$
|\l|=\left|\frac1{2\pi}\int_0^{2\pi}\frac{f(e^{{\rm i}t})}{e^{{\rm i}t}-1}\,dt\right|\le
\frac1{2\pi}\int_0^{2\pi}\frac{|f(e^{{\rm i}t})-f(1)|}{|e^{{\rm i}t}-1|}\,dt\le\|f\|_{\Li(\T)}\le\|f\|_{\OL(\T)}.
$$
Таким образом, остаётся оценить $\OL(\R)$-полунорму функции $\Phi_0$.

Оценим сначала $\|\Phi\|_{\OL([-\frac{3\pi}2,\frac{3\pi}2])}$. В силу замечания к теореме
\ref{dxa}  и леммы \ref{mol}
$$
\|\Phi\|_{\OL([-\frac{3\pi}2,\frac{3\pi}2])}\le\left\|\frac{tf(e^{{\rm i}t})}{e^{{\rm i}t}-1}\right\|_{\OL([-\frac{3\pi}2,\frac{3\pi}2])}
\le\const\|f(e^{{\rm i}t})\|_{\OL([-\frac{3\pi}2,\frac{3\pi}2])}
\le\const\|f\|_{\OL(\T)},
$$
поскольку функция $t\mapsto\frac{t}{e^{{\rm i}t}-1}$ бесконечно дифференцируема на $[-\frac{3\pi}2,\frac{3\pi}2]$.
Отсюда 
$$
\|\Phi_0\|_{\OL([-\frac{3\pi}2,\frac{3\pi}2])}\le\const\|f\|_{\OL(\T)}.
$$
Используя теперь замечание к теореме \ref{2pper}, получаем
$$
\|\Phi_0\|_{\OL(\R)}\le\frac{3\sqrt2\pi}2\|\Phi_0\|_{\OL([-\frac{3\pi}2,\frac{3\pi}2])}\le\const\|f\|_{\OL(\T)}.\quad\bl
$$

\begin{thm}
\label{pere}
Пусть $f$ -- функция, заданная на  $\T$, а $\psi$ -- дробно-линейное преобразование
такое, что $\psi(\widehat\R)=\T$.
Тогда $f\in(\OL)_{\rm loc}'(\T)$
в том и только в том случае, когда $f\circ\psi\in(\OL)'(\R)$,
при этом
\bay
\label{0310}
c_1\|f\|_{(\OL)_{\rm loc}'(\T)}
\le\|f\circ\psi\|_{(\OL)'(\R)}\le c_2\|f\|_{(\OL)_{\rm loc}'(\T)},
\ey
где $c_1$ и $c_2$ -- абсолютные положительные константы.
\end{thm}

\Pf  Положим $a=\psi^{-1}(0)$. Тогда легко видеть, что $a\in\C\setminus\R$ и $\psi(z)=\z\dfrac{z-a}{z-\ov a}$
при всех $z\in\widehat\C$, где $|\z|=1$. Ясно, что не умаляя общности, можно считать, что
$\z=1$. Докажем сначала второе неравенство. Пусть $f\in(\OL)_{\rm loc}'(\T)$.    В силу следствия \ref{rt4}
функция $f$ представима в виде $f(z)=\widehat f(-1)z^{-1}+F'(z)$, где
$F\in\OL(\T)$, при этом 
$|\widehat f(-1)|+\|F\|_{\OL(\T)}\le c\|f\|_{(\OL)_{\rm loc}'(\T)}$.
Имеем
\begin{align*}
\|f\circ\psi\|_{(\OL)'(\R)}&=
\left\|\widehat f(-1)\frac1{\psi}+F'\circ\psi\right\|_{(\OL)'(\R)}\\
&\le c\left\|\frac1{\psi}\right\|_{(\OL)'(\R)}\|f\|_{(\OL)_{\rm loc}'(\T)}+\|F'\circ\psi\big\|_{(\OL)'(\R)}.
\end{align*}
Заметим, что
$$
\Big\|\frac1{\psi}\Big\|_{(\OL)'(\R)}=\Big\|\frac{t-\ov a}{t-a}\Big\|_{(\OL)'(\R)}\le1+2|\im a|\cdot\|(t-a)^{-1}\|_{(\OL)'(\R)}
\le3,
$$
что легко следует из примера 2 в \S\:\ref{prim}.
Оценим теперь $\|F'\circ\psi\big\|_{(\OL)_{\rm loc}'(\R)}$.
Функцию $F$  выберем так, чтобы
$F(1)=F(\psi(\be))=0$. Из теоремы \ref{zamrt} следует, что
$$
\left\|\frac{F\circ\psi}{\psi'}\right\|_{\OL(\R)}=\|{\mathcal Q}_\psi F\|_{\OL(\R)}
\le3\|F\|_{\OL(\T)}\le\const\|f\|_{(\OL)_{\rm loc}'(\T)}.
$$
Следовательно,
$$
\left\|(F'\circ\psi)-\frac{\psi''}{(\psi')^2}(F\circ\psi)\right\|_{(\OL)'(\R)}=
\left\|\left(\frac{F\circ\psi}{\psi'}\right)'\right\|_{(\OL)'(\R)}\le\const\|f\|_{(\OL)_{\rm loc}'(\T)}.
$$
Остаётся оценить
$$
\left\|\frac{\psi''}{(\psi')^2}(F\circ\psi)\right\|_{(\OL)'(\R)}=\left\|\frac{\psi''}{\psi'}\,{\mathcal Q}_\psi F\right\|_{(\OL)'(\R)}
=2\left\|\frac{1}{z-\ov a}\,{\mathcal Q}_\psi F\right\|_{(\OL)'(\R)}.
$$
Используя следствие \ref{delab}, получаем:
\begin{align*}
\left\|\frac{1}{z-\ov a}\,{\mathcal Q}_\psi F\right\|_{(\OL)'(\R)}
&\le\left\|\frac{{\mathcal Q}_\psi F-({\mathcal Q}_\psi F)(\re a)}{z-\ov a}\right\|_{(\OL)_{\rm loc}'(\R)}\\
&+|({\mathcal Q}_\psi F)(\re a)|\cdot\|(z-\ov a)^{-1}\|_{(\OL)(\R)}\\
&\le2\|{\mathcal Q}_\psi F\|_{\OL(\R)}+2|F(-1)|\cdot|\im a|\cdot\|(z-\ov a)^{-1}\|_{(\OL)'(\R)}\\
&\le6\|F\|_{\OL(\T)}+2|F(-1)-F(1)|\le10\|F\|_{\OL(\T)}\le10c\|f\|_{(\OL)_{\rm loc}'(\T)}.
\end{align*}

Докажем теперь первое неравенство.
Возьмём функцию $g\in\OL(\R)$ такую, что
$g'(t)\df f(\psi(t))\in(\OL)'(\R)$ и $g(\re a)=0$.
Пусть $\vk$ обозначает дробно-линейной преобразование, обратное преобразованию $\psi$,
т. е. $\vk(z)=\dfrac{a-\ov a z}{1-z}$.
Из теоремы \ref{zamfcom} вытекает, что 
\bay
\label{1-z}
\big\|(2\im a)^{-1}(1-z)^2g(\kappa(z))\big\|_{\OL(\T)}\le 5\|g\|_{\OL(\R)}.
\ey
Следовательно,  
$$
\big\|(\im a)^{-1}(z-1)g(\vk(z))+f(z)\big\|_{(\OL)_{\rm loc}'(\T)}\le5\|g\|_{\OL(\R)}=5\|f\circ\psi\|_{(\OL)'(\R)}
$$
в силу следствия  \ref{rt4}.
Остаётся доказать, что 
$$
\big\|(\im a)^{-1}(z-1)g(\vk(z))\big\|_{(\OL)_{\rm loc}'(\T)}\le\const\|f\circ\psi\|_{(\OL)'(\R)}.
$$
Это мгновенно
вытекает из \rf{1-z} и леммы \ref{zetdzeta}. $\bl$

\begin{thm}
\label{kp}
Пусть $f\in(\OL)_{\rm loc}'(\T)$, а $\f$ -- дробно-линейное преобразование
такое, что $\f(\T)=\T$. Тогда $f\circ\f\in(\OL)_{\rm loc}'(\T)$
и $c^{-1}\|f\|_{(\OL)_{\rm loc}'(\T)}\le\|f\circ\f\|_{(\OL)_{\rm loc}'(\T)}\le c\|f\|_{(\OL)_{\rm loc}'(\T)}$ для некоторого положительного числа $c$.
\end{thm}

\Pf  Эта теорема по существу является аналогом для окружности $\T$ теоремы \ref{rinv},  в которой
речь шла о прямой $\R$. Теорема  \ref{pere} позволяет ``пересадить'' теорему \ref{rinv}
с прямой $\R$ на окружность $\T$.
$\bl$

\

\section{\bf Вокруг достаточного условия Арази--Бартона--Фридмана}
\setcounter{equation}{0}
\label{dostol+}

\

Мы рассмотрим в этом параграфе достаточное условие для операторной липшицевости
на окружности $\T$, найденное в работе Арази--Бартона--Фридмана \cite{ABF}, а также его аналог для прямой $\R$. Следуя \cite{A2}, мы покажем, как вывести эти достаточные условия из теоремы \ref{coint7}. Далее, мы введём понятия карлесоновой меры в сильном смысле и переформулируем эти достаточные условия в терминах карлесоновых мер в сильном смысле.
Мы также покажем, как можно вывести из этих достаточных условий достаточные условия в терминах классов Бесова, см. \S\:\ref{Dost}. Начнём со случая прямой.

Положим $(\CL)'(\C_+)\df\{g':g\in\CL(\C_+)\}$ и
$\|g'\|_{(\CL)'(\C_+)}=\|g\|_{\CL(\C_+)}$.
Ясно, что
$(\CL)'(\C_+)$ -- банахово пространство.
Заметим, что
функции класса $(\CL)'(\C_+)$ заданы всюду на $\clos\C_+\cup\{\be\}$.
Легко видеть, что для любой функции
$g$ из $\CL(\C_+)$ интеграл Пуассона функции $g'\big|\R$ совпадает с функцией
$g'\big|\C_+$. Чтобы убедиться в этом, достаточно заметить, что при всех $t>0$
интеграл Пуассона функции $t^{-1}(g_t-g)\big|\R$ совпадает с функцией $t^{-1}(g_t-g)\big|\C_+$,
где $g_t(z)\df g(z+t)$, и перейти к пределу при $t\to 0^+$.
Обозначим через $(\OL)_+(\R)$ пространство всех функций
$f\in\OL(\R)$, допускающих аналитическое продолжение в верхнюю полуплоскость
$\C_+$, непрерывное вплоть до границы. Положим
$(\OL)'_+(\R)\df\{g':g\in\OL_+(\R)\}$.

Из теоремы \ref{olclr} следует, что пространство $(\OL)'_+(\R)$
отождествляется естественным образом с пространством 
$(\CL)'(\C_+)$. При этом 
$$
\|f'\|_{(\OL)'(\R)}=\|f\|_{\OL(\R)}=\|f\|_{\CL(\C_+)}
=\|f'\|_{(\CL)'(\C_+)}
$$
для всех $f\in(\CL)(\C_+)$.

%

Аналог теоремы Арази--Бартона--Фридмана для полуплоскости можно сформулировать
следующим образом:

\begin{thm}
\label{abfap}
Пусть $f$ -- функция,
аналитическая в верхней полуплоскости. Предположим, что
$$
\sup_{t\in\R}\int_{\C_+}\frac{(\im w)|f'(w)|\,d\m_2(w)}{|t-w|^2}<+\be.
$$
Тогда функция $f$ имеет конечные угловые граничные значения\footnote{Под $f(\be)$ понимается $\lim f(z)$, когда $|z|\to\be$  и $z$ остаётся в любом замкнутом угле с вершиной в $\R$, все точки которого кроме вершины лежат в $\C_+$.} 
всюду на 
$\widehat\R$, которые мы будем обозначать той же буквой $f$, $f\in(\CL)'(\C_+)$,
и
$$
\|f-f(\be)\|_{(\CL)'(\C_+)}\le\frac2\pi\sup_{t\in\R}\int_{\C_+}
\frac{(\im w)|f'(w)|\,d\m_2(w)}{|t-w|^2}.
$$
\end{thm}

\begin{lem}
\label{berg}
Пусть $f$ -- аналитическая в верхней полуплоскости $\C_+$ функция. Предположим,
что 
$$
\int_{\C_+}(\im w)(1+|w|^2)^{-1}|f'(w)|\,d\m_2(w)<+\be.
$$
Тогда функция $f$ имеет конечное угловое значение
$f(\be)$ в бесконечности и
$$
f(z)-f(\be)=\frac{2\rm i}{\pi}\int_{\C_+}\frac{(\im w)f'(w)\,d\m_2(w)}{(\ov w-z)^2}
$$
при всех $z\in\C_+$.
\end{lem}

\Pf Положим
$$
g(z)\df\frac{2{\rm i}}{\pi}\int_{\C_+}\frac{(\im w)f'(w)\,d\m_2(w)}{(\ov w-z)^2}
$$
для $z\in\C_+$.
Ясно, что функция $g$ аналитична в $\C_+$ и
$$
g'(z)=\frac{4{\rm i}}{\pi}\int_{\C_+}\frac{f'(w)\,d\m_2(w)}{(\ov w-z)^3}=f'(z)
$$
для всех $z\in\C_+$, где последнее равенство следует из того, что $4{\rm i}(\pi)^{-1}(\ov w-z)^{-3}$ --
воспроизводящее ядро для пространства Бергмана, состоящего из аналитических в $\C_+$ функций
класса $L^2(\C_+,y\,d\m_2(x+{\rm i}y))$. Это -- хорошо известный и легко проверяемый факт. Остаётся доказать, что некасательный предел функции $g$ в бесконечности равен нулю.
Из равенства
$$
g(z)=\frac{2{\rm i}}{\pi}\int_{\C_+}
\left(\frac{\ov w-{\rm i}}{\ov w-z}\right)^2
\frac{f'(w)\,d\m_2(w)}{(\ov w-{\rm i})^2}
$$
и теоремы Лебега о предельном переходе под знаком интеграла следует,
что сужение функции $g(z)$ на любую полуплоскость
$\e{\rm i}+\C_+$, где $\e>0$, стремится к нулю при $|z|\to\be$. $\bl$

\medskip

{\bf Доказательство теоремы \ref{abfap}.}
Положим
\bey
F(z)\df
\frac{{2\rm i}}\pi\int_{\C_+}(\im w)f'(w)\left(\frac1{\ov w-z}-\frac1{\ov w}\right)\,d\m_2(w)
=\frac{2{\rm i}z}\pi\int_{\C_+}\frac{(\im w)f'(w)\,d\m_2(w)}{\ov w(\ov w-z)}
\eey
для всех $z\in\C$ таких, что $\im z\ge0$.
Сходимость интегралов при вещественных $z$ вытекает из неравенства Коши--Буняковского,
если принять во внимание
следующее очевидное неравенство:
$$
\int_{\C_+}\frac{(\im w)|f'(w)|\,d\m_2(w)}{|z-\ov w|^2}
\le\int_{\C_+}\frac{(\im w)|f'(w)|\,d\m_2(w)}{|\re z-w|^2}.
$$
Заметим, что
\bay
\label{FprimG}
F'(z)
=\frac{2{\rm i}}{\pi}\int_{\C_+}\frac{(\im w)f'(w)\,d\m_2(w)}{(\ov w-z)^2}
=f(z)-f(\be)
\ey
в силу леммы \ref{berg}. Рассмотрим меру Радона $\mu$ в нижней полуплоскости $\C_-$
$$
d\mu(w)\df\frac{2{\rm i}}{\pi}(\im \ov w)f'(\ov w)\,d\m_2(w).
$$
Тогда $F(z)=\widehat\mu_{0}(z)$ при всех  $z$ таких, что $\im z\ge0$, и
$$
\|\mu\|_{{\M}(\C_-)}=\frac2\pi\sup_{z\in\C_+}\int_{\C_+}
\frac{(\im w)|f'(w)|\,d\m_2(w)}{|\ov z-w|^2}=\frac2\pi\sup_{t\in\R}\int_{\C_+}
\frac{(\im w)|f'(w)|\,d\m_2(w)}{|t-w|^2}.
$$

Теперь из теоремы \ref{coint7} следует, что
$$
\|f-f(\be)\|_{(\CL)'(\C_+)}=\|F\|_{\CL(\C_+)}\le\frac2\pi\sup_{t\in\R}\int_{\C_+}
\frac{(\im w)|f'(w)|\,d\m_2(w)}{|t-w|^2}.\quad\bl
$$

\medskip

Обозначим символом  $\mathcal PM(\C_+)$ пространство всех комплексных гармонических
функций $u$, заданных в верхней полуплоскости $\C_+$ таких, что
$$
\|u\|_{\mathcal PM(\C_+)}\df\sup_{y>0}\int_\R|u(x+{\rm i}y)|\,dx<+\be.
$$
Хорошо известно, см., например, \cite{SW}, теоремы 2.3 и 2.5 главы II,
что пространство $\mathcal PM(\C_+)$ совпадает с множеством всех
функций $u$, представимых в виде
$$
u(z)=(P\nu)(z)=\frac1\pi\int_\R\frac{\im z\,d\nu(t)}{|z-t|^2}, \quad z\in\C_+,
$$
где $\nu$ -- комплексная борелевская мера на $\R$, при этом
$\|u\|_{\mathcal PM(\C_+)}=\|\nu\|_{M(\R)}\df|\mu|(\R)$.

Обозначим символом $\mathcal PL(\C_+)$ подпространство пространства
$\mathcal PM(\C_+)$, состоящее из всех функций $u\in\mathcal PL(\C_+)$,
которые соответствуют абсолютно непрерывным мерам $\nu$.

Положительную меру $\mu$ на $\C_+$ будем называть {\it мерой Карлесона в сильном
смысле}, если $\int_{\C_+}|u(z)|\,d\mu(z)<+\be$ для любой функции $u\in\mathcal PM(\C_+)$.
Заметим, что пространство $\mathcal PM(\C_+)$ содержит класс Харди $H^1$ в верхней
полуплоскости $\C_+$. Отсюда следует, что любая мера Карлесона в сильном смысле
является мерой Карлесона в обычном смысле. {\it Обозначим символом 
${\rm CM}_{\rm s}(\C_+)$ пространство всех мер Радона $\mu$ в $\C_+$,
таких, что $|\mu|$ -- мера Карлесона в сильном смысле}.
Положим
$$
\|\mu\|_{{\rm CM}_{\rm s}(\C_+)}\df\sup\left\{\int_{\C_+}|u(z)|\,d\mu(z):u\in\mathcal PM(\C_+), \|u\|_{\mathcal PM(\C_+)}
\le1\right\}.
$$
Легко видеть, что
$$
\|\mu\|_{{\rm CM}_{\rm s}(\C_+)}=\sup\left\{\int_{\C_+}|u(z)|\,d\mu(z):u\in\mathcal PL(\C_+), \|u\|_{\mathcal PM(\C_+)}
\le1\right\}
$$
и
$$
\|\mu\|_{{\rm CM}_{\rm s}(\C_+)}=\frac1\pi\sup_{t\in\R}\int_{\C_+}\frac{(\im w)\, d\mu(w)}{|t-w|^2}
=\frac1\pi\sup_{z\in\C_+}\int_{\C_+}\frac{(\im w)\, d\mu(w)}{|\ov z-w|^2}.
$$

Теперь аналог достаточного условия Арази--Бартона--Фридмана для
полуплоскости можно переформулировать  следующим образом.

\begin{thm} 
Пусть $f$ -- функция,
аналитическая в верхней полуплоскости $\C_+$. Предположим, что 
$|f'|\,d\m_2\in{\rm CM}_{\rm s}(\C_+)$.
Тогда функция $f$ имеет конечные граничные значения всюду на расширенной
числовой прямой, которые мы обозначаем той же буквой $f$, $f\in(\CL)'(\C_+)$,
и
$$
\|f-f(\be)\|_{(\CL)'(\C_+)}\le 2\|f'\,d\m_2\|_{{\rm CM}_{\rm s}(\C_+)}.
$$
\end{thm}


Аналогично можно получить ещё один вариант теоремы Арази--Бартона--Фрид\-мана.
В следующей теореме, как и во всём этом параграфе, $\zh(\n u)(a)\zh$ будет обозначать операторную норму дифференциала 
$d_a u$ функции $u$ в точке $a$.
\begin{thm}
\label{abfhp}
Пусть $u$ -- (комплексная) гармоническая функция
в верхней полуплоскости $\C_+$. Предположим, что $\zh\n u\zh\,d\m_2\in{\rm CM}_{\rm s}(\C_+)$.
Тогда функция $u$ имеет некасательные граничные значения
всюду на $\widehat\R$, $u|\R\in(\OL)'(\R)$ и
$$
\|u-u(\be)\|_{(\OL)'(\R)}\le2\|\zh\n u\zh\,d\m_2\|_{{\rm CM}_{\rm s}(\C_+)}.
$$
\end{thm}

%
%

\Pf Рассмотрим аналитические в $\C_+$ функции $f$ и $g$ такие, что $f+\ov g=u$.
Тогда
$$
f^\prime=\frac{\partial u}{\partial z}=\frac12\left(\frac{\partial u}{\partial x}-{\rm i}\frac{\partial u}{\partial y}\right)
\quad\text{и}\quad \ov{g^\prime}=\frac{\partial u}{\partial \ov{z}}=
\frac12\left(\frac{\partial u}{\partial x}+{\rm i}\frac{\partial u}{\partial y}\right).
$$
Положим
\bey
F(z)\df
\frac{{2\rm i}}\pi\int_{\C_+}(\im w)f'(w)\left(\frac1{\ov w-z}-\frac1{\ov w}\right)\,d\m_2(w)
=\frac{2{\rm i}z}\pi\int_{\C_+}\frac{(\im w)f'(w)\,d\m_2(w)}{\ov w(\ov w-z)}
\eey
и
\bey
G(z)\df
\frac{{2\rm i}}\pi\int_{\C_+}(\im w)g'(w)\left(\frac1{\ov w-z}-\frac1{\ov w}\right)\,d\m_2(w)
=\frac{2{\rm i}z}\pi\int_{\C_+}\frac{(\im w)g'(w)\,d\m_2(w)}{\ov w(\ov w-z)}
\eey
для всех $z\in\C$ таких, что $\im z\ge0$.
Применяя тождество \rf{FprimG} и аналогичное тождество для функции $G$, получаем
\begin{align*}
u(x)-u(\infty)&=F'(x)+\ov G'(x)\\
&=\frac{2{\rm i}}{\pi}\int_{\C_+}\frac{(\im w)f'(w)\,d\m_2(w)}{(\ov w-x)^2}
+\frac{2{\rm i}}{\pi}\int_{\C_+}\frac{(\im w)\ov{g'(w)}\,d\m_2(w)}{(w-x)^2}.
\end{align*}
Применяя теорему \ref{coint7}, получаем:
\begin{multline*}
\|u-u(\be)\|_{(\OL)'(\R)}\le\\
\frac2\pi\sup_{x\in\R}\left(\int_{\C_+}
\frac{(\im w)|f'(w)|\,d\m_2(w)}{|x-\ov w|^2}+\int_{\C_+}
\frac{(\im w)|g'(w)|\,d\m_2(w)}{|x-w|^2}\right)\\
=\sup_{x\in\R}\int_{\C_+}\frac{(\im w)(|f'(w)|+|g'(w)|)\,d\m_2(w)}{|x-w|^2}.
\end{multline*}
Остаётся заметить, что $|f'(w)|+|g'(w)|=\zh (\n u)(w)\zh$ для всех $w$ из $\C_+$,
поскольку операторная норма линейного отображения $h\mapsto\a h+\b\ov h$
равна $|\a|+|\b|$. $\bl$

\begin{cor} 
Пусть  $f\in\Li(\R)$.
Предположим, что $\zh{\rm Hess}\, \mP f\,\zh\,d\m_2\in{\rm CM}_{\rm s}(\C_+)$.
Тогда $f\in\OL(\R)$.
\end{cor}
%

%

Теперь мы покажем, что из достаточного условия Арази--Бартона--Фридмана вытекает достаточное условие для операторной липшицевости, полученное в \cite{Pe1} и \cite{Pe3}, см. теорему \ref{Besdost} в этом обзоре.

Чтобы получить это достаточное условие, нам понадобится следующее элементарное неравенство:
\bay
\label{51}
\|\f\,d\m_2\|_{{\rm CM}_{\rm s}(\C_+)}\le\int_0^\be\ess\sup\{\f(x+{\rm i}y):x\in\R\}\,dy
\ey
для любой измеримой неотрицательной функции $\f$ на $\C_+$.

Переходим теперь к альтернативному доказательству достаточного условия, полученного в \cite{Pe1}, \cite{Pe3}.

\begin{thm} 
Пусть $f\in B^1_{\be 1}(\R)$. Тогда $f\in\OL(\R)$.
\end{thm}

\Pf Ясно, что $f'\in L^\be(\R)$.
Пусть $u$ -- интеграл  Пуассона функции $f'$. Как известно, см. \S\:\ref{Prel},
включение $f\in B^1_{\be 1}(\R)$ равносильно следующему условию:
$$
\int_0^\be \sup_{x\in\R}\zh\n u(x+{\rm i}y)\zh\,dy<+\be.
$$
Остаётся воспользоваться неравенством \rf{51}  и
сослаться на теорему \ref{abfhp}. $\bl$

\medskip

Рассмотрим теперь случай круга. Положим $(\CL)'(\dd)\df\{g':~g\in\CL(\C_+)\}$ и
$\|g'\|_{(\CL)'(\dd)}=\|g\|_{\CL(\dd)}$.

Обозначим через  $\mathcal PM(\dd)$ пространство всех комплексных гармонических
функций $u$, заданных в единичном круге $\dd$, таких, что
$$
\|u\|_{\mathcal PM(\dd)}\df\sup_{0\le r<1}\int_\T|u(r\z)|\,|d\z|<+\be.
$$
Хорошо известно, см., например, \cite{H}, глава 3,
что пространство $\mathcal PM(\dd)$ совпадает с множеством всех
функций $u$, представимых в виде
$$
u(z)=(P\nu)(z)=\frac1{2\pi}\int_\T\frac{(1-|z|^2)\,d\nu(\z)}{|z-\z|^2}, \quad z\in\dd,
$$
где $\nu$ -- комплексная борелевская мера на $\T$, при этом
$\|u\|_{\mathcal PM(\dd)}=\|\nu\|_{M(\T)}\df|\mu|(\T)$.

Обозначим через через $\mathcal PL(\dd)$ подпространство пространства
$\mathcal PM(\dd)$, состоящее из функций $u\in\mathcal PL(\dd)$,
которым соответствует абсолютно непрерывная мера $\nu$.

Положительную меру $\mu$ на $\dd$ будем называть {\it мерой Карлесона в сильном
смысле}, если $\int_{\dd}|u(z)|\,d\mu(z)<+\be$ для любой функции $u\in\mathcal PM(\dd)$.
Заметим, что пространство $\mathcal PM(\dd)$ содержит класс Харди $H^1$ в единичном
круге $\dd$. Отсюда следует, что любая мера Карлесона в сильном смысле
является мерой Карлесона в обычном смысле. {\it Обозначим символом ${\rm CM}_{\rm s}(\dd)$ пространство всех мер Радона $\mu$ в $\dd$,
таких, что $|\mu|$ -- мера Карлесона в сильном смысле}.
Положим
$$
\|\mu\|_{{\rm CM}_{\rm s}(\dd)}\df\sup\left\{\int_{\dd}|u(z)|\,d\mu(z):u\in\mathcal PM(\dd), \|u\|_{\mathcal PM(\dd)}
\le1\right\}.
$$
Легко видеть, что
$$
\|\mu\|_{{\rm CM}_{\rm s}(\dd)}=\sup\left\{\int_{\dd}|u(z)|\,d\mu(z):u\in\mathcal PL(\dd), \|u\|_{\mathcal PM(\dd)}
\le1\right\}.
$$
и
$$
\|\mu\|_{{\rm CM}_{\rm s}(\dd)}=\frac1{2\pi}\sup_{\z\in\T}\int_{\dd}\frac{(1-|w|^2)\, d\mu(w)}{|\z-w|^2}.
$$
Заметим, что
\bay
\label{pm}
\sup_{\z\in\T}\int_{\dd}\frac{(1-|w|^2)\, d\mu(w)}{|\z-w|^2}
=\sup_{z\in\dd}\int_{\dd}\frac{(1-|w|^2)\, d\mu(w)}{|1-z\ov w|^2}.
\ey
Это вытекает из принципа максимума для $L^2$-значных аналитических функций
в круге $\dd$.


%
%
%

Сформулируем теперь в наших обозначениях достаточное условие
Арази--Бартона--Фридмана в случае окружности, см. \cite{ABF}.

\begin{thm}
\label{abfad}
Пусть $f$ -- функция,
аналитическая в единичном круге $\dd$. Предположим, что
$\z^{-1}f'(\z)\,d\m_2(\z)\in{\rm CM}_{\rm s}(\dd)$.
Тогда функция $f$ имеет конечные граничные значения всюду на окружности $\T$,
которые мы обозначаем той же буквой $f$, $f\in(\CL)'(\dd)$
и
$$
\|f-f(0)\|_{(\CL)'(\dd)}\le 2\Big\|\z^{-1}f'(\z)\,d\m_2(\z)\Big\|_{{\rm CM}_{\rm s}(\dd)}.
$$
\end{thm}

Нам понадобится аналог леммы \ref{berg}.

\begin{lem}
\label{bergd}
Пусть $f$ -- аналитическая в единичном круге $\dd$ функция. Предположим,
что $\int_{\dd}(1-|w|^2)|f'(w)|\,d\m_2(w)<+\be$.
Тогда 
$$
f(z)-f(0)=\frac{1}{\pi}\int_{\dd}\frac{(1-|w|^2)f'(w)\,d\m_2(w)}{(1-z\ov w)^2\ov w}
=\frac{1}{\pi}\int_{\C\setminus\dd}\frac{(|w|^2-1)f'(\ov w^{-1})\,d\m_2(w)}{(w-z)^2\ov w^3}
$$
при всех $z\in\C_+$.
\end{lem}

\Pf Мы докажем только первое равенство, поскольку второе равенство
получается из первого при помощи замены переменной $w\mapsto\ov w^{-1}$.
При $z=0$ доказываемое равенство вытекает из теоремы о среднем.
Остаётся заметить, что
$$
f'(z)=\frac{2}{\pi}\int_{\dd}\frac{(1-|w|^2)f'(w)\,d\m_2(w)}{(1-z\ov w)^3}
$$
для всех $z\in\dd$, см., например, следствие 1.5 монографии \cite{HKZ}.  $\bl$.

Отметим ещё, что первое равенство этой леммы может быть получено
дифференцированием по $z$ равенства (4.3) статьи \cite{ABF}.

\medskip

{\bf Доказательство теорема \ref{abfad}.}
Положим
\begin{align*}
F(z)&\df
\frac{1}\pi\int_{\dd}\frac{(1-|w|^2)f'(w)}{\ov w^2}\left(\frac1{1-z\ov w}-1\right)\,d\m_2(w)\\
&=\frac{z}\pi\int_{\dd}\frac{(1-|w|^2)f'(w)\,d\m_2(w)}{\ov w(1-z\ov w)}
\end{align*}
для всех $z$ из $\C$ таких, что $|z|\le1$.
Сходимость интегралов при $z\in\T$ вытекает из неравенства Коши--Буняковского, если принять во внимание
тождество \rf{pm}.
Заметим, что
$$
F'(z)
=\frac{1}{\pi}\int_{\dd}\frac{(1-|w|^2)f'(w)\,d\m_2(w)}{(1-z\ov w)^2\ov w}
=f(z)-f(0)
$$
в силу леммы \ref{berg}. Рассмотрим меру Радону $\mu$ в $\C\setminus\ov\dd$,
$$
d\mu(w)\df\frac{1}{\pi}\ov w^{-3}(|w|^2-1)f'(\ov w^{-1})\,d\m_2(w).
$$
Тогда
\begin{align*}
\widehat\mu_{0}(z)&=\frac1\pi\int_{\C\setminus\dd}\frac{(|w|^2-1)f'(\ov w^{-1})}{\ov w^3}
\left(\frac1{w-z}-\frac1w\right)\,d\m_2(w)\\
&=\frac{1}\pi\int_{\dd}\frac{(1-|w|^2)f'(w)}{\ov w^2}\left(\frac1{1-z\ov w}-1\right)\,d\m_2(w)=F(z).
\end{align*}
Заметим, что
\begin{align*}
\|\mu\|_{{\M}(\C\setminus\clos\dd)}&=\sup_{z\in\clos\dd}\int_{\C\setminus\clos\dd}\frac{d|\mu|(w)}{|w-z|^2}
=\frac1\pi\sup_{z\in\ov\dd}\int_{\C\setminus\clos\dd}\frac{(|w|^2-1)|f'(\ov w^{-1})|}{|w-z|^2|w|^3}\,d\m_2(w)\\
&=\frac1\pi\sup_{z\in\clos\dd}\int_{\dd}\frac{(1-|w|^2)|f'(w)|}{|1-z\ov w|^2|w|}\,d\m_2(w)\\
&=\frac1\pi\sup_{\z\in\T}\int_{\dd}\frac{(1-|w|^2)|f'(w)|}{|\z-w|^2|w|}\,d\m_2(w)
=2\big\|\z^{-1}f'(\z)\,d\m_2(\z)\big\|_{{\rm CM}_{\rm s}(\dd)}.
\end{align*}

Теперь из теоремы \ref{coint7} следует, что
$$
\|f-f(0)\|_{(\CL)'(\dd)}=\|F\|_{\CL(\dd)}\le\|\mu\|_{{\M}(\C\setminus\clos\dd)}
=2\big\|\z^{-1}f'(\z)\,d\m_2(\z)\big\|_{{\rm CM}_{\rm s}(\dd)}.\quad\bl
$$

\begin{cor}
\label{abfadcor}
Пусть $f$ -- функция,
аналитическая в единичном круге $\dd$. Предположим, что
$f'\,d\m_2\in{\rm CM}_{\rm s}(\dd)$.
Тогда функция $f$ имеет конечные граничные значения всюду на окружности $\T$,
которые мы обозначаем той же буквой $f$, $f\in(\CL)'(\dd)$
и
$$
\|f-f(0)\|_{(\CL)'(\dd)}\le\const\|f'\,d\m_2\|_{{\rm CM}_{\rm s}(\dd)}.
$$
\end{cor}

\Pf Достаточно заметить, что для любой непрерывной в $\dd$ функции $h$
условие $h\,d\m_2\in{\rm CM}_{\rm s}(\dd)$ влечёт, что $\z^{-1}h(\z)\,d\m_2(\z)\in{\rm CM}_{\rm s}(\dd)$.
Остаётся сослаться на теорему о замкнутом графике. $\bl$

\medskip

{\bf Замечание.} Можно отказаться от использования теоремы о замкнутом графике
и получить следующую явную оценку:
$$
C_s\Big(|\z|^{-1}h(\z)\,d\m_2(\z)\Big)\le\frac83C_s\big(h\,d\m_2\big)
$$
для любой субгармонической в $\dd$ функции $h$, но нам это не понадобится.

\medskip

\begin{thm}  
Предположим, что  $\zh \n u\zh\,d\m_2\in{\rm CM}_{\rm s}(\dd)$
для некоторой гармонической в $\dd$  функции $u$. Тогда $u$ имеет граничные значения 
всюду на $\T$ и \lb$u\in(\OL)'_{{\rm loc}}(\T)$.
\end{thm}

\Pf 
Функцию $u$  можно представить в виде $u=f+\ov g$, где $f$  и $g$ -- аналитические 
в $\dd$  функции. Из следствия \ref{abfadcor} вытекает, что $f,g\in(\OL)'(\T)$.
Остаётся заметить, что из определения пространства $(\OL)'_{{\rm loc}}(\T)$
мгновенно вытекает его инвариантность относительно комплексного
сопряжения. $\bl$

%

\begin{cor}  Предположим, что  $\zh {\rm Hess}\, u\zh\,d\m_2\in{\rm CM}_{\rm s}(\dd)$
для некоторой гармонической в $\dd$  функции $u$. Тогда $u$ продолжается до
непрерывной функции на $\dd\cup\T$ и
$u\in\OL(\T)$.
\end{cor}

Отсюда легко  вытекает следующий результат работы \cite{Pe1}, доказательство которого приведено в
в \S\:\ref{Dost} этого обзора, см. теорему \ref{Besokr}.

\begin{thm}
\label{612}
Пусть $f\in B^1_{\be 1}(\T)$. Тогда $f\in\OL(\T)$.
\end{thm}


Мы вывели достаточное условие Арази--Бартона--Фридмана  из теоремы \ref{coint7}.
Можно показать, что теорема  \ref{coint7} даёт примеры операторно липшицевых
функций, которые не удовлетворяют аналогу достаточного условия  Арази--Бартона--Фридмана для $\C_+$.
В \cite{A2} построен пример функции $f\in\widehat\M(\clos\C_+)$ 
такой, что  $f''\,d\m_2\not\in{\rm CM}_{\rm s}(\C_+)$. Из теоремы  \ref{coint7} следует,
что такая функция $f$ принадлежит пространству $(\CL)'(\C_+)$,  хотя  достаточное 
условие Арази--Бартона--Фридмана к этой
функции не применимо. Аналогичное утверждение справедливо и для функций в $\dd$.

\medskip

{\bf Замечание.} В работе Арази--Бартона--Фридмана \cite{ABF} было отмечено, что их достаточное
условие операторной липшицевости функции, заданной на единичной ок\-ружности, влечёт достаточное условие,
полученное в работе \cite{Pe1}, см. также \S \ref{Dost} этого обзора. Из результатов статьи \cite{A2} следует, что достаточное условие операторной липшицевости
Арази--Бартона--Фридмана может работать и в том случае, когда функция $f'$ не является
непрерывной. С другой стороны, легко видеть, что если $f\in B^1_{\be 1}(\T)$, то $f'\in C(\T)$.
То же самое можно сказать и по поводу функций класса $B^1_{\be 1}(\R)$ (см.  теорему \ref{Besdiffer}). Действительно,
нетрудно проверить, что функция $f(z)=\exp(-{\rm i}z^{-1})$ удовлетворяет условиям теоремы
\ref{abfap}, хотя её сужение на вещественную прямую терпит разрыв в нуле. Отметим ещё, что
в работе \cite{A2} по существу доказано, что подмножество вещественной прямой может 
быть множеством всех точек разрыва функции $f\big|\R$, где
$f$ удовлетворяет условиям теоремы \ref{abfap}, в том и только в том случае, когда оно 
типа $F_\s$ и не имеет внутренних точек. Можно сделать аналогичное замечание и для   функций на $\T$ и $\dd$.

\medskip

Интересно сравнить достаточное условия операторной липшицевости, приведённое в этом параграфе, с необходимым условием, приведённым в \S\:\ref{Neob}. Комбинация этих условий дана в следующей теореме.

\begin{thm}
Если $f\in\Li(\R)$ и $\zh{\rm Hess}\mP f\zh\,d\m_2$ -- мера Карлесона в сильном смысле, то $f\in\OL(\R)$.
Если же $f\in\OL(\R)$, то $\zh{\rm Hess}\mP f\zh\,d\m_2$ -- мера Карлесона.
\end{thm}

Аналогичное утверждение справедливо и для функций на окружности $\T$.

\

\section{\bf В каких случаях имеет место равенство $\bs{\OL(\fF)=\Li(\fF)}$?}
\setcounter{equation}{0}
\label{orav}

\

\begin{thm}
\label{dal}
Предположим, что $\OL(\fF)=\Li(\fF)$ для некоторого замкнутого
подмножества $\fF$  комплексной плоскости $\C$. Тогда множество $\fF$ конечно.
\end{thm}

\Pf  Предположим, что множество $\fF$ бесконечно. Тогда множество
$\fF$ содержит предельную точку
$a\in\widehat\C\df\C\cup\{\be\}$. Если $a\in\C$, то можно считать, что
$a=0$. Случай $a=\be$
 может быть рассмотрен аналогично. Кроме того, случай $a=\be$ сводится к случаю
 $a=0$ при помощи дробно-линейных преобразований.

 Предположим сначала, что $\fF\subset\R$. Тогда легко построить функцию
 $f\in\Li(\fF)$, у которой нет производной в нуле. Ясно, что $f\not\in\OL(\fF)$.

 Чтобы отказаться от предположения $\fF\subset\R$, нам понадобится следующая
 лемма.

 \begin{lem}
 \label{ldal}
 Пусть $0<q<1$ и пусть
 $\{a_n\}_{n\ge1}$ -- последовательность положительных чисел
 таких, что $a_{n+1}\le q a_n$  при всех $n\ge1$.
 Тогда для любой числовой последовательности $b_n$, удовлетворяющей условию
 $\sum_{n\ge1}|b_n|a^{-1}_n<+\be$, найдётся функция $v\in\OL(\R)$
 такая, что $v(a_n)=b_n$ при всех $n\ge1$.
 \end{lem}

 \Pf Зафиксируем функцию $\f$ класса $C^\be(\R)$ такую, что $\f(0)=1$ и
 $\supp\f\subset[-\delta,\delta]$, где число $\delta$ будет выбрано в конце
 доказательства. Положим $v(t)\df\sum\limits_{n\ge1}b_n\f(a_n^{-1}(t-a_n))$.
 Тогда
 $$
 \|v\|_{\OL(\R)}\le\sum_{n\ge1}|b_n|\cdot\|\f(a_n^{-1}(t-a_n))\|_{\OL(\R)}=\|\f\|_{\OL(\R)}
 \sum_{n\ge1}|b_n|a^{-1}_n<+\be
 $$
и $v(a_n)=b_n$ при всех $n\ge1$, если только число $\delta$ достаточно мало. $\bl$

\medskip

Продолжим доказательство теоремы \ref{dal}. 
Хорошо известно, что любая липшицева функция, заданная на подмножестве
комплексной плоскости $\C$, продолжается до липшицевой функции, заданной
на всей  комплексной плоскости $\C$, см., например, \cite{St}, глава VI, \S\:2, теорема  3. 
Таким образом, достаточно ограничиться случаем,
когда множество $\fF\setminus\{0\}$ состоит из членов последовательности
$\{\l_n\}_{n\ge1}$, как угодно быстро стремящейся к нулю. Пусть $\l_n=a_n+{\rm i}b_n$. Можно считать,
что $\lim\limits_{n\to\be}\frac{\l_n}{|\l_n|}=1$ и
вещественные последовательности $\{a_n\}_{n\ge1}$ и $\{b_n\}_{n\ge1}$
удовлетворяют условиям леммы \ref{ldal}. Положим $h(t)\df t+{\rm i}v(t)$, где
$v$ обозначает то же, что и в лемме \ref{ldal}.
Теперь случай множества $\fF$ сводится к уже разобранному случаю множества
$\re\fF$, поскольку 
$$
\|A-B\|\le\|h(A)-h(B)\|\le(1+\|v\|_{\OL(\R)})\|A-B\|
$$
для любых самосопряжённых операторов $A$ и $B$ таких, что $\s(A),\:\s(B)\subset\re\fF$. $\bl$

\

\begin{center}
\bf\large Заключительные замечания
\end{center}
\label{zaklzam}

\addtocontents{toc}{\vspace*{.03cm}\hspace*{-.35cm}\textbf{Заключительные замечания}\hfill\pageref{zaklzam}}

\

В этом разделе мы вкратце коснёмся некоторых результатов, которые не вошли в основную часть статьи.

\medskip

{\bf 1. Операторные модули непрерывности.} Для непрерывной функции $f$ на $\R$ операторный модуль непрерывности $\O_f$ определяется равенством
$$
\O_f(\d)\df\sup\{\|f(A)-f(B)\|:~A~\mbox{и}~B\;~
\mbox{--~самосопряжённые операторы},~\:\|A-B\|<\d\}.
$$
Операторные модули непрерывности были введены в \cite{AP2}
и подробно изучались в \cite{AP5}. Теорема \ref{prmone}, сформулированная в этой статье, означает, что если $f\in\L_\o(\R)$, где $\o$ -- модуль непрерывности, то
$$
\O_f(\d)\le\const\o_*(\d),\quad\mbox{где}\quad
\o_*(\d)\df\d\int_\d^\be\frac{\o(t)}{t^2}dt
$$
(см. \rf{omzv}).

В работе \cite{AP5} обсуждается точность таких оценок и получены значительно более точные оценки для непрерывных ``кусочно выпукло-вогнутых'' функций $f$. В частности, получена следующая неулучшаемая оценка:
$$
\big\|\,|A|-|B|\,\big\|\le C\|A-B\|\log\left(2+\log\frac{\|A\|+\|B\|}{\|A-B\|}\right)
$$
для ограниченных самосопряжённых операторов $A$ и $B$. Это неравенство значительно улучшает оценку Като, полученную в работе \cite{Kat}.

\medskip

{\bf 2. Коммутаторные оценки для функций нормальных операторов.} Лемма \ref{predosnnaformKotSh} позволяет получить следующую квазикоммутаторную оценку:
$$
\|f(N_1)R-Rf(N_2)\|\le\const\o_*\big(\max\{\|N_1R-RN_2\|,\|N_1^*R-RN_2^*\|\}\big)
$$
для любого модуля непрерывности $\o$, для любой функции $f$ класса $\L_\o(\R)$,
для любого линейного оператора $R$ нормы 1 и для любых нормальных операторов $N_1$ и $N_2$ (см. \cite{APPS}). В работе \cite{AP6} норма квазикоммутатора в левой части неравенства оценивается только через норму $\|N_1R-RN_2\|$. Правда, при этом нужно заменить в правой части неравенства $\o_*$ на $\o_{**}\df(\o_*)_*$:
$$
\|f(N_1)R-Rf(N_2)\|\le\const\o_{**}\big(\{\|N_1R-RN_2\|\}\big).
$$
Отметим, что в случае классов Гёльдера, т.е. $\o(t)=t^\a$, $0<\a<1$, имеет место неравенство $\o_{**}(t)\le\const(1-\a)^{-2}t^\a$. Иными словами, мы получаем коммутаторно гёльдерову оценку.

\medskip

{\bf3. Подход Никольской--Фарфоровской к операторно гёльдеровым функциям.}
В работе \cite{FN} предлагается альтернативный подход к операторно гёльдеровым функциям. Он основан на следующем утверждении:

\medskip

{\it Пусть  $0<\a<1$. Тогда $\L_\a(\Z)\subset\OL(\Z)$.
При этом существует число $c_\a$ такое, что 
$\|f\|_{\OL(\Z)}\le c_\a\|f\|_{\L_\a(\Z)}$.}

\medskip

Из этого утверждения нетрудно вывести теорему \ref{44}, если воспользоваться легко проверяемым неравенством $\O_f(\d)\le2\o_f(\d/2)+2\|f(\d x)\|_{\OL(\Z)}$. Само это утверждение можно доказать посредством интерполирования функции класса $\L_\a(\Z)$ функцией класса $B_{\be,1}^1(\R)$ и применения теоремы \ref{Besdost}, хотя в работе 
\cite{FN} это было доказано совсем другим способом.

\medskip

{\bf4. Функции наборов коммутирующих самосопряжённых операторов.} Изучение функций от нормальных операторов эквивалентно изучению функций от пар коммутирующих самосопряжённых операторов. В работе \cite{NP} удаётся распространить результаты работы \cite{APPS} (см. \S\:\ref{OLnaplBes} этого обзора)
на случай функций от произвольного числа коммутирующих самосопряжённых операторов. При этом используются совсем другие методы.

В работе \cite{A3} получены обобщения на операторно липшицевы функции от $n$ переменных некоторых результатов работы \cite{A1} о дробно-линейных подстановках
(см. \S\:\ref{dlp} этого обзора). В многомерной ситуации роль дробно-линейных преобразований
играют преобразования Мёбиуса, т. е. суперпозиции конечного числа инверсий.

\medskip

{\bf5. Липшицевы функции наборов коммутирующих самосопряжённых операторов.} В работе \cite{KPSS} результаты работы \cite{PS} обобщаются на случай функций наборов $n$ коммутирующих самосопряжённых операторов и получается липшицева оценка в норме $\bS_p$, $1<p<\be$, для липшицевых функций в $\R^n$.

\medskip

{\bf6. Функции пар некоммутирующих самосопряжённых операторов.} Для пары $(A,B)$ не обязательно коммутирующих самосопряжённых операторов в работе \cite{ANP} рассматриваются функции $f(A,B)$, которые определяются с помощью двойных операторных интегралов, и изучается поведение таких функций при возмущении пары.
Оказалось, что в отличие от функций коммутирующих операторов, липшицевы оценки в операторной норме и в ядерной норме сильно отличаются друг от друга. В частности, в работе \cite{ANP} показано, что при $f\in B_{\be,1}^1(\R^2)$ имеет место неравенство
$$
\|f(A_1,B_1)-f(A_2,B_2)\|_{\bS_p}\le\const\|f\|_{B_{\be,1}^1}
\max\{\|A_1-A_2\|_{\bS_p},\|B_1-B_2\|_{\bS_p}\}
$$
при $p\in[1,2]$.
Такое же неравенство было получено ранее в \cite{APPS} для функций коммутирующих операторов. Однако, в случае функций некоммутирующих операторов, {\it такое же неравенство для $p\ge2$ и для операторной нормы неверно}, см. \cite{ANP}.

В качестве инструмента в работе \cite{ANP} используются тройные операторные интегралы и вводятся {\it модифицированные} хогеруповы тензорные произведения пространств $L^\be$.

\medskip

{\bf7. Операторно липшицевы функции и формула следов Лифшица--Крей\-на.} Пусть $A$ и $B$ -- самосопряжённые операторы с ядерной разностью $A-B$. Такой паре соответствует единственная вещественная функция $\xi$ класса $L^1(\R)$, называемая {\it функцией спектрального сдвига}, такая, что для достаточно хороших функций $f$ на $\R$ имеет место формула следов Лифшица--Крейна
$$
\trace(f(A)-f(B))=\int_\R f'(t)\xi(t)\,dt
$$
(см. \cite{Li} и \cite{Kr}). В работе М.Г. Крейна \cite{Kr} было показано, что эта формула справедлива для функций $f$, производная которых является преобразованием Фурье комплексной меры. В работе \cite{Pe3} формула следов распространена на функции $f$ класса Бесова $B_{\be,1}^1(\R)$. В теореме \ref{tsentrez} этого обзора показано, что для того, чтобы оператор $f(A)-f(B)$ был ядерным при условии ядерности оператора $A-B$, необходимо и достаточно, чтобы функция $f$ была операторно липшицевой. Наконец, в недавней работе \cite{Pe7} показано, что для операторно липшицевых функций левая часть формулы следов Лифшица--Крейна не только имеет смысл, но и совпадает с её правой частью. Другими словами, {\it формула Лифшица--Крейна справедлива для любых самосопряжённых операторов $A$ и $B$ с ядерной разностью в том и только в том случае, когда функция $f$ операторно липшицева}.

\medskip

В завершение упомянем недавний обзор \cite{Pe10}, в котором рассматриваются приложения кратных операторных интегралов в различных задачах теории возмущений.

\

\

\noindent
\begin{tabular}{p{9cm}p{15cm}}
А.Б. Александров & В.В. Пеллер \\
Санкт-Петербургское отделение & Department of Mathematics \\
Математический институт Стеклова РАН  & Michigan State University \\
Фонтанка 27, 191023 Санкт-Петербург & East Lansing, Michigan 48824\\
Россия&USA
\end{tabular}

\end{document}